\documentclass[12pt,a4paper]{amsart}
 
\usepackage[T1]{fontenc}
\input xypic\usepackage{amsmath}
\usepackage{amssymb}

\usepackage[usenames]{color}
\usepackage{amssymb}
\usepackage{latexsym}
\usepackage{graphicx}
\usepackage{psfrag}

\newcommand{\observation}[1]{}
\newcommand{\generalizations}[1]{}

\newcommand{\hattuM}{M}
\newcommand{\Sclo}{\mathcal S^{cl}}

\newcommand{\xxi}{\xi}
\newcommand{\zzeta}{\zeta}

\renewcommand{\div}{\hbox{div}}

\newcommand{\Ric}{\hbox{Ric}}
\newcommand{\Ein}{\hbox{Ein}}
\def \Box {$\square$}

\newcommand{\mmmbf}{}
\newcommand{\mattibf}{} 
 
\newcommand{\mltext}{} 

\newcommand{\newtext}{}
\newcommand{\tobecheckedtext}{}

\newcommand{\HOX}[1]{}

\newcommand{\hiddenfootnote}[1]{}

\newcommand{\rhoepsilon}{{\rho}}
\newcommand{\bsequence}{{\bf b}}
\def\X{{\mathcal X}}
\def\Y{{\mathcal Y}}
\def\U{{\mathcal U}}

\def\S{{\mathcal S}}
\def\T{{\mathcal T}}
\def\P{{\mathcal P}}
\def\qP {{\mathcal R}}
 
\def\ear{\hbox{\it E}}
\def\pointear{\hbox{\it e}}

\renewcommand{\H}{{\mathbb H}} 
\newcommand{\R}{{\mathbb R}} 
\newcommand{\D}{{\cal D}} 
\newcommand{\I}{{\cal I}} 
\newcommand{\be}{{\bf e}} 
\newcommand{\E}{{\cal E}} 
\newcommand{\A}{{\cal A}}
\newcommand{\B}{{\cal B}}  
\newcommand{\cC}{{\cal C}}  
\newcommand{\K}{{\cal K}}  

\newcommand{\C}{{\mathbb C}} 
\newcommand{\N}{{\mathbb N}} 
\newcommand{\cal}{\mathcal } 
\renewcommand{\O}{{\mathcal O}} 
 
\renewcommand{\L}{{\mathcal L}} 
\newcommand{\V}{{\cal V}} 
\newcommand{\W}{{\cal W}}

\def\hat{\widehat}
\def\tilde{\widetilde}
\def \bfo {\begin {eqnarray*} }
\def \efo {\end {eqnarray*} }
\def \ba {\begin {eqnarray*} }
\def \ea {\end {eqnarray*} }
\def \beq {\begin {eqnarray}}
\def \eeq {\end {eqnarray}}
\def \supp {\hbox{supp}\,}
\def \dim{\hbox{dim}\,}

\def \re {\hbox{Re}\,}
\def \im {\hbox{Im}\,}

\def \WF {\hbox{WF}\,}

\def \dist {\hbox{dist}}
\def\diag{\hbox{diag }}
\def \det {\hbox{det}}

\def\bra{\langle}
\def\cet{\rangle}

\def \e {\varepsilon}
\def \p {\partial}
\def \a {\alpha}
\def \b {\beta}

\def\M{{\mathcal M}}
\def\F{{\mathcal F}}

\def\Z{{\mathbb Z}}

\newtheorem{definition}{Definition}[section] 
\newtheorem{theorem}[definition]{Theorem} 
\newtheorem{lemma}[definition]{Lemma} 
 
\newtheorem{proposition}[definition]{Proposition} 
\newtheorem{corollary}[definition]{Corollary}

\begin{document}
\title[Determination of space-time]
{Determination of structures in the space-time from local
\\ measurements:  a detailed exposition
}
\date{May 29, 2013}
\author{Yaroslav Kurylev, Matti Lassas,  Gunther Uhlmann}
\address{Yaroslav Kurylev, UCL; Matti Lassas,  University of Helsinki;
Gunther Uhlmann, University of Washington,  and 
 University of Helsinki.  {\rm y.kurylev@ucl.ac.uk, Matti.Lassas@helsinki.fi,  
  gunther@math.washington.edu}}

\email{}

\maketitle

\HOX{ Make also a version of the paper where "observation" and "generalization" macros
are activated so that the  hidden texts are shown. \observation{Observation text is
in quite good shape and is in blue,} \generalizations{Generalization text is messy
and is in red.}}
{\bf Abstract:}
{\it We consider inverse problems for the Einstein equation 
with a time-depending metric on a 4-dimensional  globally hyperbolic
Lorentzian manifold $(M,g)$. We formulate
the concept of active measurements for relativistic models.
We do this by coupling
the Einstein equation with equations for scalar fields
and study the system
$\Ein(g)=T$, $T=T(g,\phi)+F_1$, and $\square_g \phi=F_2+S(g,\phi,F_1,F_2)$.
Here $F=(F_1,F_2)$ correspond to the perturbations of 
the physical fields which we control
and $S$ is a secondary source corresponding to
the adaptation of the system to the perturbation so that the
conservation law $\div_g(T)=0$ will be satisfied.

The inverse problem we study
is the question, do the observation of the solutions $(g,\phi)$
in an open  subset $U\subset M$ of the space-time 
corresponding to sources  $F$ supported 
in  $U$ determine the properties 
of the metric in a larger domain 
$W\subset M$ containing $U$.
To study this problem we define the concept of light observation sets and
show that these sets determine the conformal class of the metric.
This corresponds to
passive observations from a distant area of the space which is  
filled by light sources (e.g. we see light from stars varying in time).
One can apply the obtained 
result to solve inverse problems encountered in general relativity
and in various practical imaging problems.
}


\smallskip
\tableofcontents

\section{Introduction and main results}
\HOX{Paper should be divided to two parts.
We should consider the alternative that the IP for Einstein equations is  the paper number 1
and the geometric result is paper number 2.}
We consider inverse problems for non-linear 2nd order hyperbolic equations 
with time-depending metric on a globally hyperbolic Lorent\-zian manifold $(M,g)$ of dimension
$n\geq 2+1$
and in particular,  the Einstein equation.
Roughly speaking,  we study the problem, can an observer on a Lorentzian manifold, satisfying
certain natural causality conditions, determine
the structure of the surrounding space-time by doing measurements near its
world line. 
Thus the problem we are interested in  this paper is the inverse problem with respect 
to the evolution problem in the general relativity. Let us note that this
last problem has recently attracted much interest in the mathematical
community with many important results been obtained, see e.g.\
\cite{Bernal,ChBook,Chr,Ch-K,Dafermos,Dafermos2,K1,K3,K4,Li1}.

We consider two different kind of inverse problems. 
In the first one we have passive observations:
We detect on an open set $U\subset M$ the wave fronts of the 
waves produced by the point sources  located at points $q$ 
in a relatively compact subset $V\subset M$.
We call such observations the light observation sets $\O_U(q)\subset TU$,
where $TU$ is the tangent bundle on $U$.
We ask, can the conformal class of the metric in $V$ be determined
from these observations. In the second class of problems we consider 
active measurements: We consider non-linear hyperbolic partial differential equations on $M$
assuming that we can control sources supported in $U$ which
produce waves that  we can measure in the same set. We ask,
can the properties of the metric (the metric itself or its conformal class)
or the coefficients of the equation be determined in a suitable larger 
set $W$ containing the set $U$. For instance in the context of 
relativity theory, these correspond to the following examples: 
In the first case we consider
passive observations from a distant area of the space which is  
filled by light sources (e.g. we see light from stars varying in time).
In the second case we assume that one can cause local perturbations in
the stress-energy tensor and measures 
locally the  caused perturbations of the gravitational field. 

\subsection{Notations}\label{sec:notations 1}

Let $(M,g)$ be a $C^\infty$-smooth $n$-dim\-ens\-ion\-al manifold with $C^\infty$-smooth Lorentzian
metric $g$ with a causal structure (For this and other definitions
for Lorentzian manifolds, see the next section.) The tangent
bundle of $M$ is denoted by $TM$ and the projection to the base
 is denoted by $\pi:TM\to M$.

Let us introduce  next some notations needed below.
For $x,y\in M$ we say that $x$ is in the chronological past of $y$ and denote $x\ll y$ if $x\not =y$ and there is a time-like path from $x$ to $y$.
If $x\not=y$ and there is a causal path from $x$ to $y$,
we say that $x$ is  in the causal past of $y$ and  denote $x<y$.
If 
$x< y$ or $x=y$ we denote $x\leq y$.
The chronological future $I^+(p)$ of $p\in M$
consist of all points $x\in M$ such that $p\ll x$,
and the causal future $J^+(p)$ of $p\in M$
consist of all points $x\in M$ such that $p\leq x$.
Similarly chronological past $I^-(p)$ of $p\in M$
consist of all points $x\in M$ such that $x\ll p$ and
the causal past $J^-(p)$ of $p\in M$
consist of all points $x\in M$ such that $x\leq p$.
For a set $A$ we denote $J^\pm(A)=\cup_{p\in A}J^\pm(p)$.
We also denote $J(p,q):=J^+(p)\cap J^-(q)$ and $I(p,q):=I^+(p)\cap I^-(q)$.
If we need to emphasize the metric $g$ which is used to define the causality, 
we denote by $J^\pm_g(q)$ or  $J^\pm_{M,g}(q)$  the past and the future sets
of $q\in M$ with respect to a Lorentzian metric $g$ etc.

Let $\gamma_{x,\xi}(t)=\exp_x(t\xi)$ denote the geodesics in $(M,g)$. 
Also, let $TU=\{(x,\xi)\in TM;\ x\in U\}$.
Let $L_xM$ denote the light-like directions of $T_xM$, 
and $L_x^+M$ and $L_x^-M$  denote the future and past pointing light-like vectors,
correspondingly, and $L_x^{*,+}M$ and $L_x^{*,-}M$ be the 
future and past pointing light-like co-vectors.
Sometimes, to emphasize the metric, we denote $L_x^+M=L_x^+(M,g)$, etc.
We also denote $\L^+_{g}(x)=\exp_x(L^+_xM)\cup\{x\}$ the union of the image of the future light-cone
in the exponential map of $(M,g)$ and the point $x$.

By  \cite{Bernal}, an open Lorenzian manifold $(M,g)$ is globally hyperbolic
if and only if it has a causal structure where there are no closed causal paths in $M$ and
for all $q^-,q^+\in M$ such that $q^-<q^+$ the set 
$J(q^-,q^+)\subset M$ is compact.
We assume throughout the paper that $(M,g)$ is globally hyperbolic. Roughly
speaking, this means that we have no naked singularities which
we could reach by moving along a time-like path starting from a point $q^-$ and
ending to a point $q^+$.


When $g$ is a Lorentzian metric, having eigenvalues $\lambda_j(x)$ and eigenvectors 
$v_j(x)$ in some local coordinates $(U,X)$, we will use also the corresponding
Riemannian metric, denoted $g^+=\hbox{Riem}(g)$ which has the  
eigenvalues $|\lambda_j(x)|$ and the eigenvectors 
$v_j(x)$ in the same local coordinates $(U,X)$.
Let $B_{g^+}(x,r)=\{y\in M;\ d_{g^+}(x,y)<r\}$.

All functions $u(x)$ defined on manifold $M$ etc.\ are real-valued unless otherwise mentioned.
Finally, when $X$ is a set, let  $P(X)=2^X=\{Z;\ Z\subset X\}$ denote the power set of $X$.

Let $\mu_g:[-1,1]\to M$ be a freely falling observer, that is, a time-like geodesic on $(M,g)$.
Let $-1<s_{-2}<s_{-1}<s_{+1}<s_{+2}<s_{+3}<s_{+4}<1$ and denote $p^-=\mu_g(s_{-1})$
and  $p^+=\mu_g(s_{+1})$. Below, we also denote $s_\pm=s_{\pm 1}$.

When  $z_0=\mu_g(s_{-2})$ and $\eta_0=\p_s\mu_g(s)|_{s=s_{-2}}$,  let  $\U_{z_0,\eta_0}(h)$  
be the $h$-neighborhood
of $(z_0,\eta_0)$ in the Sasaki metric of $(TM,g^+)$.  
We use below a small parameter $\hat h>0$.
For $(z,\eta)\in \U_{z_0,\eta_0}(2\hat h)$
we define on $(M,g)$  a freely falling observer 
$\mu_{g,z,\eta}:[-1,1]\to M$,
such that   $\mu_g(s_{-2})=z$, and $\p_s\mu_g(s_{-2})=\eta$.
We assume that $\hat h$ is so small that
for all $(z,\eta)\in \U_{z_0,\eta_0}(2\hat h)$ the geodesic $\mu_{g,z,\eta}([-1,s_{+2}])\subset M$
is well defined and time-like and satisfies 
\beq\label{kaava D}
& &\mu_{g,z,\eta}(s_{-2})\in I_{g}^-(\mu_{g,z_0,\eta_0}(s_{-1})),\quad
\mu_{g,z,\eta}(s_{+2})\in I_{g}^+(\mu_{g,z_0,\eta_0}(s_{+1}),\hspace{-1cm}\\
& &\nonumber \mu_{g,z,\eta}(s_{+4})\in I_{g}^+(\mu_{g,z_0,\eta_0}(s_{+3})).
\eeq

We denote, see Fig.\ 1. 
\beq\label{eq: Def Wg with hat}
U_g=\bigcup_{(z,\eta)\in \U_{z_0,\eta_0}(\hat h)} \mu_{g,z,\eta}((-1,1)).
\eeq 

%

\begin{figure}[htbp]
\begin{center}\label{Fig-1}

\psfrag{1}{$\pi(\mathcal U_{z_0,\eta_0})$}
\psfrag{2}{$U_{\hat g}$}
\psfrag{3}{$J_{\hat g}(\hat p^-,\hat p^+)$}
\psfrag{4}{$W_{\hat g}$}
\psfrag{5}{$\mu_{z,\eta}$}
\includegraphics[width=7.5cm]{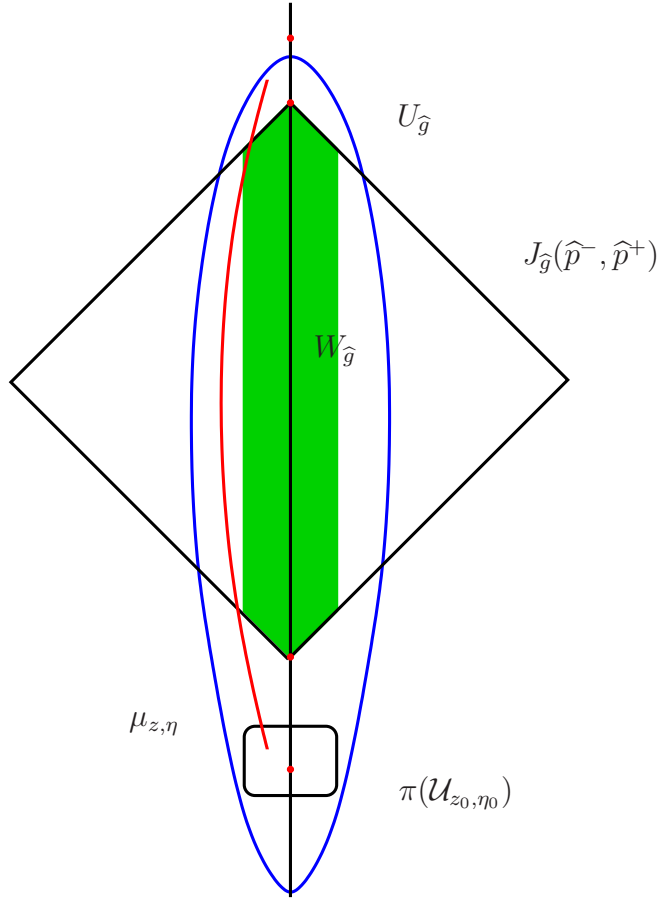}
\end{center}
\caption{Setting throughout the paper: A schematic figure where  the space-time is represented as the  2-dimensional set  $\R^{1+1}$. 
 In the figure the black vertical line is the freely falling
observer $\hat \mu([-1,1])$ and the four red points on it
are $z_0=\hat \mu(s_{-2})$,  $\hat p^-=\hat \mu(s_{-})$, $\hat p^+=\hat \mu(s_{+})$,
and $\hat \mu(s_{+2})$. The rounded black square
is $\pi(\mathcal U_{z_0,\eta_0})$ that is is a neighborhood of $z_0$,
and the red curve starting from $z\in\pi(\mathcal U_{z_0,\eta_0})$ is the
time-like geodesic $\mu_{\hat g,z,\eta}([s_{-2},1))$. 
The boundary of the set $U_{\hat g}$ is shown on blue.
The green area is the set $W_{\hat g}\subset U_{\hat g}$ where the Fermi-type
coordinates are defined, and the black "diamond" is the set
$J_{\hat g}(\hat p^-,\hat p^+)=J_{\hat g}^+(\hat p^-)\cap J^-_{\hat g}(\hat p^+)$.}
 \end{figure}

\subsection{Inverse problem for the light observation sets}

Let  us first consider $M$ with a fixed metric $g$ and denote below $U=U_g$.
Next we define the
 light-observation set
   for point $q$ corresponding
to observations from a light source at the point $q$, see Fig.\ 2. 

\begin{definition}\label{def. O_U}
The light-observation set corresponding to the point $q\in M$ 
and the observation set $U\subset M$ is
\ba
\O_U(q)&=&\{(\gamma_{q,\eta}(r),\dot \gamma_{q,\eta}(r))\in TU;\ r\geq 0,\ \eta\in L_q^{+}M\}.
\ea 
The set of the light-observation points corresponding to $q\in M$ is
\ba
\P_U(q):=\{\gamma_{q,\eta}(r)\in U;\ r\geq 0,\ \eta\in L_q^{+}M\}=\pi(\O_U(q)).
\ea 
\end{definition}

\begin{figure}[htbp] \label{Fig-2}
\begin{center}

\psfrag{1}{$\pi(\mathcal U_{z_0,\eta_0})$}
\psfrag{2}{$U_{\hat g}$}
\psfrag{3}{$J_{\hat g}(\hat p^-,\hat p^+)$}
\psfrag{4}{$W_{\hat g}$}
\psfrag{5}{$\mu_{z,\eta}$}
\includegraphics[width=7.5cm]{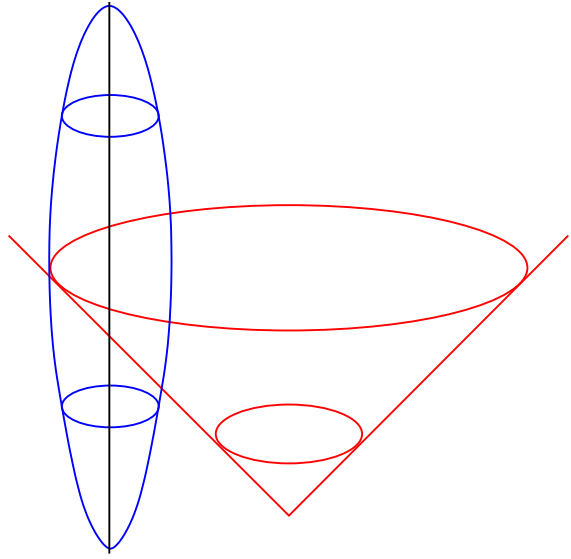}
\end{center}
\caption{A schematic figure where the space-time is represented as the  3-dimensional set $\R^{1+2}$. The future light cone $\mathcal L^+_{\hat g}(q)$ corresponding to the point $q$ is  shown
as a red cone. The point $q$ is the tip of the cone. The boundary of the set $U_{\hat g}$ is shown on blue.
The set of the light observation points $\mathcal P_U(q)$, see Definition 1.1,  is the intersection
of the set $\mathcal L^+_{\hat g}(q)$ and the set $U_{\hat g}$.}
 \end{figure}

In the following, we extend any map $F:A\to A^{\prime}$ to the map
$F:P(A)\to P(A^{\prime})$ by setting $F(B)=\{F(b);\ b\in B\}$ for $B\subset A$.
We call $F:P(A)\to P(A^{\prime})$ the power set extension of $F:A\to A^{\prime}$.

In the following, we will consider 
the collection $\O_U(V):=\{\O_U(q);\ q\in V\}\subset P(TU)$ of all
light observation sets corresponding to the points in an 
open relatively compact set $V\subset M$. Note 
that the collection
 $\O_U(V)$ is considered just as a subset of $P(TU)$
 and for a given element $A\in \O_U(V)$ we do not know what is the point $q$
 for which $A=\O_U(q)$.

Later we will introduce a topology on   
a suitable subset $\mathbb S\subset P(TU)$ that contains $\O_U(V)$. Then we can consider  the image 
of the manifold $V$ in the embedding 
(see Lemma \ref{lemma homeo} below) $\O_U:q\mapsto \O_U(q)$ as an
embedded manifold  in  $\mathbb S$, or in $P(TU)$.
If we assume that $U$ is given, the set $\O_U(V)$ 
will be a homeomorphic representation of the unknown manifold $V$ in
a known topological space $\mathbb S$. For a compact Riemannian
manifold $N$, an analogous representation  
is the metric space $K(N)\subset \hbox{Lip}(N)$ 
is obtained using the Kuratowski-Wojdyslawski embedding, $K:x\mapsto \dist(x,\cdotp)$,
from the metric space  $N$  to space of Lipschitz functions $\hbox{Lip}(N)$ on $N$.
In several inverse problems for Riemannian manifolds with boundary,
a homeomorphic image of the compact manifold $N$ has be obtained by
using 
the 
embedding $R:x\mapsto  \dist(x,\cdotp)$, $R:N\to \hbox{Lip}(\p N)$, see
\cite{AKKLT,Katsuda,KKL,KLU}.

 Our first goal is prove that the collection of the light observation
 sets, $\O_U(V)$,  determine the conformal type of  the Lorentzian manifold $(V,g)$.
 As $ \pi(\O_U(q))= \P_U(q)$, it is enough to show that $\P_U(V)=\{
  \P_U(q);\ q\in V\}$ determines the conformal type of $(V,g)$, see Fig.\ 3.

 \HOX{Gunther's suggestion: Main geometric theorem should be formulated
 to earliest observations.}

\begin{figure}[htbp]\label{Fig-3}
\begin{center}

\psfrag{1}{$\pi(\mathcal U_{z_0,\eta_0})$}
\psfrag{2}{$U_{\hat g}$}
\psfrag{3}{$J_{\hat g}(\hat p^-,\hat p^+)$}
\psfrag{4}{$W_{\hat g}$}
\psfrag{5}{$\mu_{z,\eta}$}
\includegraphics[width=7.5cm]{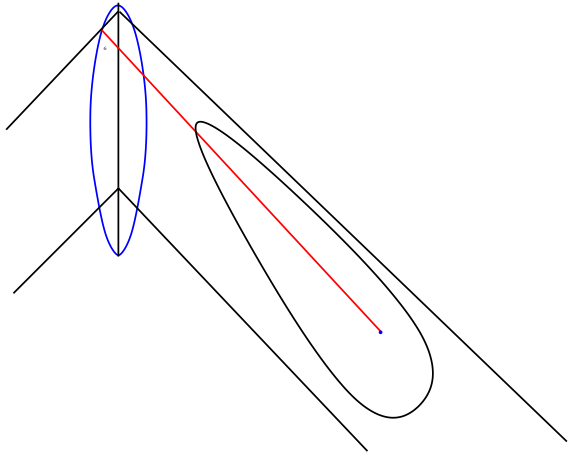}
\end{center}
\caption{A schematic figure where the space-time is represented as the  2-dimensional set $\R^{1+1}$. 
 In Theorem \ref{main thm}, we consider a relatively compact set $V\subset I^-(p^+)\setminus J^-(p^-)$.
The boundary of the set $V$ is shown in the figure with a black curve. The red curve
is a light ray from a point $x\in V$ that intersect the blue set $U_{\hat g}$.
}
 \end{figure}

\begin{theorem}\label{main thm}
Let $(M_j,g_j)$, $j=1,2$ be two open, smooth, globally hyperbolic 
 $(1,n-1)$ Lorentzian manifolds of dimension $n$, $n\geq 3$ and let $p^+_j,p^-_j\in M_j$ be the endpoints of a 
  time-like geodesic $\mu_{g_j}([s_-,s_+])\subset M_j$, that is,
$p_j^ \pm=\mu_{g_j}(s_\pm)$. Let  $U_j\subset M_j$ be open  relatively compact neighborhood
of $\mu_j([s_-,s_+])$ given by (\ref{eq: Def Wg with hat}).  
Moreover, let $V_j$ be open, relatively compact subsets of $I^-(p^+_j)\setminus J^-(p^-_j) \subset M_j$,
$j=1,2$.

Let us denote by
 \ba
\P_{U_j}(V_j)=\{\P_{U_j}(q);\ q\in V_j\}\subset P(U_j)
\ea
the collections of  the observation point
sets on manifold $(M_j,g_j)$ corresponding to the points in the set $V_j$.  
Assume that there is a diffeomorphism
$\Phi:U_1\to U_2$ such that $\Phi(\mu_1(s))=\mu_2(s)$, $s\in [s_-,s_+]$ 
and the the power set extension of $\Phi$ defines a bijection
\ba
\Phi:\P_{U_1}(V_1)\to \P_{U_2}(V_2).
\ea
Then there is a diffeomorphism $\Psi:V_1\to V_2$
and the metric $\Psi^*g_2$ is conformal to $g_1$.
\end{theorem}

Theorem \ref{main thm} can be stated by saying that if 
an observer moves along a path $\mu$ then
the diffeomorphism-type of the neighborhood $U\subset M$ of $\mu$ and the  collection of 
the   observation  point
 sets
 $\{\P_U(q);\ q\in V\}$ determine uniquely the Lorentzian manifold
$(V,g)$ up to a conformal deformation.
We note that by the strong hyperbolicity any  set $V\subset\subset  I(p^-,p^+)$ satisfies
also the condition $V\subset  I^-(p^+)\setminus J^-(p^-)$.

Note that we do not assume that $(M,g)$ is complete, which is crucial in relativity
due to the presence of singularities. 


In the case when metric $g$ is known in $U_g$ and $g$ is Ricci
flat (i.e.\ corresponds to vacuum) in a set $W$ that intersects $U_g$,
after constructing the conformal structure we can find the whole
metric tensor by constructing the conformal factor along light
like geodesics that start at $U_g$ and are subset of $W$ (see Fig.\ 4). This fact is formulated
more precisely  in the following corollary.

\begin{corollary}\label{coro of main thm original Einstein pre}
Assume that $(M _j,g_j)$ and $U_j$, $V_j$,  $j=1,2$ and the conformal map
$\Psi:V_1\to V_2$ are as in  Theorem \ref{main thm}. Moreover, assume that 
 $\Phi=\Psi|_{U_1}:(U_1,g_1)\to (U_2,g_2)$ is an isometry and
 for $j=1,2$  there are sets $W_j\subset M_j$ such that 
the Ricci curvature of $g_j$ is zero 
in $W_j$, and  $V_j\subset W_j \cup U_j$, $\Psi(U_1\setminus W_1)=U_2\setminus W_2$.  
Moreover, assume that any point $x\in V_1\cap W_1$ can be connected to some point $y\in U_1\cap
W_1$ with a light-like geodesic $\gamma_{x,\xi}([0,l])\subset V_1\cap W_1$.
Then
the metric $\Psi^* g_2$ is isometric to $g_1$ in $V_1$.
\end{corollary}

\medskip

\noindent {\bf Example 1.} Theorem \ref {main thm} concerns passive observations
which are idealizations of the
measurements used in astronomy and the engineering sciences. For example,
consider an event when $J$ 
astronomical observatories, located at the points $x_1,x_2,\dots, x_J$ in a neighborhood
$U$ of the path of the Earth in the space
time $M$, observe light from a supernova at directions $\xi_j\in T_{x_j}M$, 
$j=1,2,\dots,J$.
By measuring the properties of the observed frequencies (i.e. spectrum),
 the observatories can draw the conclusion that they are looking at  the same supernova
 and that  there is some point $q$ (the unknown location of the supernova)
 such that $(x_j,\xi_j)\in \O_U(q)$ for all $j=1,2,\dots,J$.
Examples of other point source type events which can possibly be
observed by astronomic measurements, are  
novae, quasars, and pulsars as well as
eclipses of double stars, pulsating stars 
or variable stars with extreme sunspots or flares  \cite{KazS}.
When one considers the Universe in large 
scale, the non-flatness of the space time is clearly visible in the observations. 
For instance,
gravitational lenses produced by massive objects \cite{Henry,Nasa}
cause multiple images of the distant objects and the images
have significant time delays. Mathematically, this
corresponds to the focal points of geodesics. For instance,
 in the observations made on 
the gravitationally lensed quasar Q0957+561, due to the bending of
space, light from the quasar arrive to Earth along two different paths
and the time delay between these two paths has been measured  to be
approximately 
$417$ days \cite{Kundic,Rhee}. In such time scales, many astronomic
events can be considered as point sources and thus the light observation
sets can  be reasonable, although highly idealized, model for the observations.
We note that the above described "point sources"  in astronomy
have very different magnitudes and
with the present telescopes and other astronomic
instruments only some of those are 
detectable e.g.\ from another side of a gravitational lens.

The light observation sets model also observations or theoretical measurements
(or thought experiments) which one could do near astronomic systems
with black holes. For instance, the light observation sets correspond to observations when
a large number of point sources (i.e.\ matter) emitting time-varying
radiation  fall in a black hole(s).

The inverse problem for the light observation sets and for the  Einstein equation
considered in   the next section is related also to practical measurements: The detection of small but rapid perturbations
of the gravitational field, that is, gravitational waves. 
The detection
of gravitational waves is a very rapidly growing field of physics where  new laboratories have be founded
in many countries in the recent years, e.g the  LIGO interferometer  at Washington, US, the 
GEO-600 detector at Germany, and  the VIRGO detector at Italy.
We note that the Einstein equation implies that gravitational waves
exists if strong enough sources exists. Thus the detection of gravitational waves actually
is the question wether  strong sources exists.
Although gravitational radiation has not been directly detected using present measurement devices, there is indirect evidence for its existence. For example, the 1993 Nobel Prize in Physics was awarded for measurements of the Hulse-Taylor binary system which suggests gravitational waves are more than mathematical anomalies.
The detection of gravitational waves can be considered as
a far field measurement. The relation of far-field and near-field measurements is a well studied question
of inverse problems \cite{Ber} and in understanding the far-field observations the near field inverse problem
needs to be studied. This is just the inverse problem for the coupled Einstein-matter field equations we consider below.

The inverse problems analogous to Theorem \ref {main thm} are encountered also 
in several mundane applications,
for instance one can study if the mobile phone signals can be used to determine the 
refractive index of the surrounding urban neighborhood.

\subsection{The inverse problem for Einstein equations}

\subsubsection{How the active measurements could be done}

Let us make a Gedankenexperiment. Assume that we are close to a huge gravitational 
object and want to measure the distortion of the space time.
Let us assume that we can do extremely precise measurements.
Let us use several high-power laser sources (like laser pointers). 
We vary the direction of the laser rays so that some of the rays cross at time $t=T^0$. When the rays cross,
we have an increased density of light at the crossing point, just
 at the crossing time. Because of the non-linearity
of the Einstein-Maxwell equations this creates
an artificial point source of gravitational waves
having a very small amplitude. Using point masses that are
located near the sources of the laser rays and observing
their movements
 we can
in principle detect the waves. Making lots of such
experiments which create artificial sources of gravitational waves at
a large set of points, we get
data analogous to "boundary distance functions", see e.g.\ \cite{KKL}, on the Lorentzian
manifold. Then we ask, is it possible to use this data
to determine the
metric up to a conformal deformation in the portion of the
space time bordered by the possible event horizons.

For sake of mathematical simplicity, we consider next
the Einstein-scalar field equations instead of Maxwell-Einstein equations
and replace the laser rays by gravitational and scalar field waves, in fact
four waves
having only $C^k$-regularity with some finite $k$.

\subsubsection{Inverse problems for non-linear wave equations}

Many physical models lead to non-linear differential equations.  In small perturbations, 
these equations can be approximated by linear equations,
and most of the previous results on hyperbolic inverse problem  concern with
these linear models. Moreover, the existing uniqueness results are limited to
the time-independent or real-analytic coefficients \cite{AKKLT,Bel1,BelKur,Eskin,KKL} as 
these results are based on Tataru's unique 
continuation principle \cite{Tataru1,Tataru2}. Such unique continuation results 
have   been shown to fail for general metric
tensors with are not analytic the time variable \cite{Alinhac}. 
Even some linear inverse problem are not uniquely
solvable.  In fact, the counterexamples for these problems have
been used in the so-called transformation optics. This has led to  
models for  fixed frequency invisibility cloaks, see e.g.\ 
\cite{GKLU1,GKLU2,GLU} and references therein.
These applications give one more motivation 
to study inverse problems for non-linear equations.

\subsubsection{Einstein equations}
In the following, the Einstein tensor of a Lorentzian metric $g=g_{jk}(x)$
 of type $(-,+,+,+)$ on a 4-dimensional manifold $M$ is
\ba
\Ein_{jk}(g)=\Ric_{jk}(g)-\frac 12 (\hbox{tr}\,\Ric (g))
g_{jk}=\Ric_{jk}(g)-\frac 12 (g^{pq}\,\Ric_{pq} (g))
g_{jk}.
\ea
Here, $\Ric_{pq} (g)$ is the Ricci curvature of the metric $g$,
$\hbox{tr}\,\Ric (g)=g^{pq}\,\Ric_{pq} (g)$ is equal to the scalar curvature of $g$.
We define the divergence of a 2-covariant tensor $T_{jk}$ to be in  local coordinates 
$(\div_g T)_k=\nabla_n(g^{nj}T_{jk})$.

Let us consider the
Einstein equation in presence of matter,
\beq\label{Einmat1}
& &\Ein_{jk}(g)=T_{jk},\\
& &\div_{g}T=0,\label{Einmat2}
\eeq
for a Lorentzian metric $g$
and a stress-energy tensor $T$
related to the distribution of mass and energy.
We will consider $g$ which is a small perturbation of
 some unknown background metric $\hat g$
 and the stress-energy tensor $T$ which is 
a perturbation of an unknown background stress-energy tensor $\hat T$.
We recall the essential fact that when (\ref{Einmat1}) has solutions,
then due to the Bianchi's identity $\div_{g}(\Ein(g))=0$ and thus
the equation (\ref{Einmat2}) follows automatically. Equation (\ref{Einmat2})
is called the conservation law for the stress-energy tensor.
 
\subsubsection{Global hyperbolicity and the domain of influence $\K$} 
\label{subsec: Gloabal hyperbolicity} Recall that a Lorentzian
metric $g_1$ dominates the metric $g_2$, if all vectors $\xi$
that are light-like or time-like with respect to the metric $g_2$ are
time-like with respect to the metric $g_1$, and in this case we denote $g_2<g_1$.
Let $(\hattuM ,\hat g)$ be a $C^\infty$-smooth globally hyperbolic Lorentzian
manifold. 
By \cite{Geroch} it holds
that there is a Lorentzian metric $\tilde g$ such that  $(\hattuM ,\tilde g)$
is   globally hyperbolic and $\hat g<\tilde g$. Note that one can
assume 
that the metric $\tilde g$ is smooth (see Appendix C on details). We use the positive definite Riemannian metric $\hat g^+$ to define
norms in  the spaces $C^k_b(\hattuM )$ of functions with bounded $k$ derivatives, sometimes denoted also  by $C^k_b(\hattuM ;\hat g^+)$
and the  Sobolev spaces $H^s(\hattuM )$.

By \cite{Bernal2}, the globally hyperbolic manifold $(\hattuM ,\tilde g)$ has
an isometry $\Phi$  
 to the smooth product manifold $(\R\times N,\tilde  h)$,
where $N$
is a 3-dimensional manifold and the metric $\tilde  h$ can be written as
$\tilde  h=-\beta(t,y) dt ^2+\overline h(t,y)$ where $\beta:\R\times  N\to (0,\infty)$ is a smooth function and 
$\overline h(t,\cdotp)$ is a Riemannian metric on $ N$ depending smoothly
on $t\in \R$, and the submanifolds $\{t^{\prime}\}\times N$ are $C^\infty$-smooth
Cauchy surfaces for all $t^{\prime}\in \R$. We define smooth time function
${\bf t}:\hattuM \to \R$ by setting
 ${\bf t}(x)=t$ if $\Phi(x)\in \{t\}\times N$.
Let us next identify these isometric manifolds,
that is, we denote $\hattuM =\R\times N$.
%
Then 
$S_t=\{x\in \hattuM ;\ {\bf t}(x)=t\}$ are space-like
Cauchy surfaces also for $\tilde g$ and for all metrics $g$ for which $g<\tilde g$.

For $t\in \R$, let $\hattuM (t)=(-\infty,t)\times N$. Let $t_1>t_0>0$
and denote $\hattuM _1=\hattuM (t_1)$ and $\hattuM _0=\hattuM (t_0)$.
When  $\hat p^-\in \hattuM _0$ and $j\in \{0,1\}$, it follows from  by \cite[Cor.\ A.5.4]{BGP}
that $J^+_{\tilde g}(\hat p^-)\cap M(t_j)$ is compact.
We denote below 
\beq\label{Kappa sets}
\K_j:=J^+_{\tilde g}(\hat p^-)\cap  M(t_j).
\eeq
As $\hat  g<\tilde g$,
we see that there exists $\e_0>0$ such that if $\|g-\hat g\|_{C^0_b(\hattuM _1;\hat g^+)}<\e_0$, then
$g|_{\K_1}<\tilde g|_{\K_1}$,
and in particular, we have  $J^+_{g}(p)\cap \hattuM _1\subset \K_1$
for all $p\in \K_1$. Next we consider in particular the manifold $M_0$ and
denote $\K_0=\K$.

\subsubsection{Reduced Einstein tensor}


%
%
%

Let $t_1>t_0>t_{-1}=0$ and  $g^\prime$ be a metric  on $\hattuM (t_1)$
that coincide with $\hat g$ in $\hattuM (t_{-1})$ and  assume 
that $g^\prime$ satisfies the Einstein equation $\Ein(g^\prime)=T^\prime$
on $\hattuM (t_1)$.  
If $g^\prime$  is a small
perturbation of $\hat g$ in a suitable sense (see  Appendix A), there is
a diffeomorphism $f:\hattuM (t_1)\to f(\hattuM (t_1))\subset \hattuM $ that is  a 
$(g^{\prime},\hat g)$-wave map $f:\hattuM (t_1)\to \hattuM $, see (\ref{C-problem 1})-(\ref
{C-problem 1}) in Appendix A,
and satisfies $\hattuM (t_0)\subset f(\hattuM (t_1))$. 
The wave map has the property that
  $\Ein(f_*g^{\prime})= \Ein_{\hat g}(f_*g^{\prime})$,
  where 
$\Ein_{\hat g}(g)$ is the the $\hat g$-reduced Einstein tensor, 
given by formula (\ref{Reduced Einstein tensor}) 
 in Appendix A, and has the form 
\ba
(\Ein_{\hat g} g)_{pq}=-\frac 12 g^{jk}\hat \nabla_j\hat \nabla_k  g_{pq}
+\frac 14( g^{nm}g^{jk}\hat \nabla_j\hat \nabla_k  g_{nm})g_{pq}
+P_{pq}(g,\hat \nabla g),
\ea
and
$\hat \nabla_j$ is the covariant differentation with respect to the metric $\hat g$
and $P_{pq}$ is polynomial of 
$ g_{nm}$, $ g^{nm}$, $ \hat \nabla_jg_{nm}$,  $\hat \nabla_j g^{nm}$,
$\hat g_{nm}$,  $\hat g^{nm}$, and the curvature $\hat R_{jklm}$ of  $\hat g$.
Using
this  map $f$, the tensors $g=f_*g^{\prime}$ and $T=f_*T^{\prime}$ 
satisfy the  Einstein equations $\Ein(g)=T$ on $\hattuM (t_0)$. Moreover, 
as $f$ is a wave map, the Einstein equation for $g$ can be written 
as the $\hat g$-reduced Einstein equations, that is,
\beq\label{eq: EE in correct coordinates}
\Ein_{\hat g}(g)= T\quad\hbox{on }\hattuM (t_0).
\eeq
In the literature, the above is often
stated by saying that 
the reduced
Einstein equation (\ref{eq: EE in correct coordinates})  is
the Einstein equations  written
with the wave-gauge corresponding to the metric $\hat g$. The equation (\ref{eq: EE in correct coordinates}) is a quasi-linear hyperbolic system of equations.
We emphasize that
a solution of the reduced Einstein equation can be a solution
of the original Einstein equation only if the stress energy
tensor satisfies the conservation law $\nabla^{g}_jT^{jk}=0$.
Summarizing the above, when $g$ is close to $\hat g$ and
satisfies the Einstein equation, by changing the coordinates
appropriately (or physically stated, choosing the correct gauge), 
the metric tensor $g$ satisfies  
the $\hat  g$-reduced Einstein equation.
Next we formulate a direct problem for  the $\hat  g$-reduced Einstein equations.

\subsubsection{Motivation for the direct problem} 

We consider metric and physical fields informally on a Lorentzian manifold
$(M,g)$.
We aim to study an inverse problem with active measurements.
As measurements can not be implemented in Vacuum,
we have to add matter fields in the model. Since matter fields and
the metric tensor have to satisfy a standard conservation law, we 
add in the model two types of matter fields. In fact,
we will start with the standard coupled system of the Einstein equation
and scalar filed equations with some sources $\F_1$ and $\F_2$, namely
\beq\label{eq: adaptive model with no source}
& &\Ein_{\hat g}(g) ={\bf T}_{jk}(g,\phi)+\F_1,\quad
\hbox{in }M,
\\ \nonumber
& &{\bf T}_{jk}(g,\phi)=\sum_{\ell=1}^L(\p_j\phi_\ell \,\p_k\phi_\ell 
-\frac 12 g_{jk}g^{pq}\p_p\phi_\ell \,\p_q\phi_\ell-
\frac 12 m\phi_\ell^2g_{jk}),
\\ \nonumber
& &\square_g\phi_\ell - m\phi_\ell=\F_2,\quad \ell=1,2,3,\dots,L,
\eeq
and make the source term $\F_1$ correspond to fluids (or fluid fields) consisting of particles
which energy and moment vectors are controlled and make $\F_2$ to adapt
the changes of $g,\phi,$ and $\F_1$ so that the standard conservation law
is satisfied. Note that  above we could assume that masses $m_\ell$ of the fields $\phi_\ell$
depend on $\ell$, but for simplicity we have assumed that $m_\ell=m$.

\HOX{Maybe we should write first the system of matter-Einstein equations and then
motivate those.}

By an active measurement
we mean a model where we can control some of these fields  and
the other fields adapt to the changes of all fields so that the  conservation law holds. Roughly speaking,
we can consider measurement devices as a complicated
process that changes one energy form to other forms of energy, like a system
of explosives  that transform some potential energy to kinetic energy.
This process creates a perturbation of the metric and the matter fields
that we observe in a subset of the space time.
In this paper,  our aim is not to consider a  physically realistic example but 
a mathematical model that can be rigorously analyzed. 
\HOX{Here was: Roughly speaking,
we consider a theoretical "device" which transforms energy in scalar fields
to matter particles (fluid).
}

To motivate such a model, we start with a non-rigorous discussion.
Following  \cite[Ch. III, Sect. 6.4, 7.1, 7.2, 7.3]{ChBook} 
and  \cite[p. 36]{AnderssonComer}
we start by considering  the 
Lagrangians, associated to gravity, scalar fields $\phi=(\phi_\ell)_{\ell=1}^L$ and non-interacting fluid fields, that is, the number density four-currents ${\bf n}=({\bf n}^\kappa)_{\kappa=1}^J$
(where each ${\bf n}^\kappa$ is a co-vector,  see \cite[p. 33]{AnderssonComer}).
We also add in to the model a Lagrangian  associated with
some scalar valued source fields $S=(S_\ell)_{\ell=1}^L$  and
$Q=(Q_k)_{k=1}^K$. We consider action corresponding to
the coupled Lagrangian\hiddenfootnote{
An alternative Lagrangian would be the one where the term
$P_{pq}g^{pq}$ is replaced by $\frac 12 P_{pq}P_{jk}g^{pj}g^{qk}$.
Then in the stress energy tensor $P_{jk}$ is replaced by
${\mathcal T}^{(1)}_{jk}= F_{jp}F_{kq}g^{qp}$, i.e. ${\mathcal T}^{(1)}=\hat F\hat G\hat F.$
By [Duistermaat, Prop. 1.3.2], 
its derivative with respect to $P_{jk}$, that is,
$((D_F{\mathcal T}^{(1)}_{jk}|_{\hat F})^{ab})$ (COMPUTE THIS)
is a surjective pointwise map is surjective at points $x\in M$ if
$\hat F_{qp}$ is positive definite symmetric (and thus invertible) matrix (CHECK THIS!!).
}
\HOX{Change of L to $\mathcal A$}
\ba
& &\mathcal A=\int_{M}\bigg(L_{grav}(x)+L_{fields}(x)+L_{source}(x)\bigg)\,dV_g(x),\\
& &L_{grav}=R(g),\\
& &L_{fields}=\sum_{\ell=1}^L\bigg(g^{jk}\p_j\phi_\ell \,\p_k\phi_\ell 
-{\mathcal V}( \phi_\ell;S_\ell)\bigg)
+  g^{jk}(\sum_{\kappa=1}^J {\bf n}_j^\kappa {\bf n}_k^\kappa),
\\
& &L_{source}=\e \mathcal H_\e(g,S,Q, {\bf n},\phi),
\ea
where $R(g)$ is the scalar curvature, $dV_g=(-\det g)^{1/2}dx$ is the volume form on $(M,g)$, 
\beq\label{eq: Slava-a}
{\mathcal V}( \phi_\ell;S_\ell)=\frac 12 (m-1)\phi_\ell^2+\frac 12(\phi_\ell -S_\ell)^2
\eeq are energy potentials
of the scalar fields $\phi_\ell$ that depend on $S_\ell$,
 and  $\mathcal H_\e(g,S,Q, {\bf n},\phi)$ is a function modeling the measurement
 device we use.
 We assume that $\mathcal H_\e$ is bounded and   its derivatives with respect
 to $S,Q,{\bf n}$  are very large (like of order $O((\e)^{-2})$) and its derivatives
with respect of $g$ and $\phi$ are bounded when $\e>0$ is small.
We note that the  above Lagrangian for the fluid fields is the sum of the single fluid Lagrangians
where for all fluids the master function  $\Lambda$ is  the identity function, that is, the energy density
of each fluid is given by
$\rho=-\Lambda(n^2)=-n^2$, $n^2=g^{jk}{\bf n}_j{\bf n}_k.$ Note that here ${\bf n}$
is a time-like vector or zero and thus $\rho$ is non-negative.
On fluid Lagrangians, see the discussions in
 \cite[p. 33-37]{AnderssonComer},  \cite[Ch. III, Sect. 8]{ChBook} and
 \cite{Taub} and \cite[p.\ 196]{Felice}.

%
When we compute the critical points of the Lagrangian $L$ and 
neglect the  $O(\e)$-terms, the equation  $\frac {\delta  \mathcal A}{\delta g}=0$ 
gives the Einstein equation with a stress-energy tensor $T_{jk}$
defined below, see (\ref{eq: T tensor1}),  the equation  $\frac {\delta  \mathcal A}{\delta \phi}=0$
gives the wave equations with sources $S_\ell$. We assume that $O(\e^{-1})$ order equations
obtained from 
the equation $(\frac {\delta \mathcal A}{\delta S},\frac {\delta  \mathcal A}{\delta Q},
\frac {\delta  \mathcal A}{\delta {\bf n}})=0$ 
fix the values of the scalar functions $Q$, the tensor 
\ba
\P=\sum_{\kappa=1}^J  {\bf n}^\kappa_j
{\bf n}^\kappa_k dx^j\otimes dx^k
\ea and yields for the sources $S=(S_\ell)_{\ell=1}^L$ equations of the form
$S_\ell=\mathcal S_\ell(\phi, \nabla^g
 \phi,Q,\nabla^g Q,\P,\nabla^g \P,g)$.
 The function
 $\mathcal H_\e$ models the way the measurement device works. 
Due to this we will assume that $\mathcal H_\e$ and thus functions ${\mathcal S}_\ell$ may
be quite complicated.   The interpretation of the above is
that in  each measurement event we use a device that
fixes the values of the scalar functions $Q$ and the tensor $\P$ and gives the equations
$S=\mathcal S(\phi, \nabla^g
 \phi,Q,\nabla^g Q,\P,\nabla^g \P,g)$ that 
 tell how the sources of the $\phi$-fields adapt to these changes
 so that 
the physical conservation laws are satisfied.

 Now we stop the non-rigorous discussion (where the $O(\e)$ terms were
neglected). 

\subsubsection{Formulation of the direct problem} 
To start the rigorous analysis, let us define some physical fields  and introduce a
 model as a system of partial differential
 equations (that hold at a critical point of the above Lagrangian $L$).

We assume that there are $C^\infty$-background fields $\hat g$,  $\hat \phi$,
 $\hat Q$, and  $\hat P$ on $M$.
 We also fix 
 a smooth metric  $\tilde g$ that is a
globally hyperbolic metric on $\hattuM $  such that $\hat g<\tilde g$
and make the the identification $\hattuM =\R\times N$ where
$\{t\}\times N$ are Cauchy surfaces for $\tilde g$. Moreover, we fix $t_0>0$ and 
a point $\hat p^-\in (0,t_0)\times N$ and denote $\hattuM _0=M(t_0)=(-\infty,t_0)\times N$.

{Let $\P=\P_{jk}(x)dx^jdx^k$ be a symmetric tensor on $\hattuM _0$, corresponding below
to a direct perturbation to the stress energy tensor, and 
$Q=(Q_{\ell}(x))_{{\ell}=1}^K$ where  $Q_{\ell}(x)$ are real-valued functions on $\hattuM _0$.
Also, we consider 
a Lorentzian metric  $g$ on $\hattuM _0$ and $\phi=(\phi_\ell)_{\ell=1}^L$ where $\phi_\ell$
are scalar fields on $\hattuM _0$, $L\leq K-1$. The potentials of the fields $\phi_\ell$
are ${\mathcal V}( \phi_\ell;S_\ell)$ given in (\ref{eq: Slava-a}). The way how  $S_\ell$,
called below the adaptive source functions, depend on other 
fields is explained later.

Using the $\phi$ and $\P$ fields, we define the stress-energy tensor
 \beq\label{eq: T tensor1}
\hspace{-5mm} T_{jk}=\sum_{\ell=1}^L(\p_j\phi_\ell \,\p_k\phi_\ell 
-\frac 12 g_{jk}g^{pq}\p_p\phi_\ell \,\p_q\phi_\ell-\mathcal V( \phi_\ell;S_\ell)g_{jk})
+\P_{jk}.\hspace{-14mm}
\eeq
 Below, we introduce
the notation 
\ba
P_{jk}=
\P_{jk}-\sum_{\ell=1}^L \frac 12 S_\ell ^2g_{jk}.
\ea

We assume that $P-\hat P$ and $Q-\hat Q$ are supported on
$\K=J^+_{\tilde g}(\hat p^-)\cap \hattuM _0$. As we will see later, when $P-\hat P$ and $Q-\hat Q$
are small enough, in a suitable sense, and the intersection
of their support and $M_0$ is contained in $\K$, then $g<\tilde g$
and $g-\hat g$, considered as a function on $M_0$, is also supported in $\K$.

Using (\ref{eq: Slava-a}) we can
 write the stress energy tensor (\ref{eq: T tensor1})
 in the form
 \ba
T_{jk}
&=&P_{jk}+Zg_{jk}+{\bf T}_{jk}(g,\phi),\quad Z=\sum_{\ell=1}^L S_\ell\phi_\ell,
\\
{\bf T}_{jk}(g,\phi)&=&\sum_{\ell=1}^L(\p_j\phi_\ell \,\p_k\phi_\ell 
-\frac 12 g_{jk}g^{pq}\p_p\phi_\ell \,\p_q\phi_\ell-
\frac 12 m\phi_\ell^2g_{jk}),
 \ea
where we call $Z$ the stress energy density caused by sources $S_\ell$.

\generalizations{We will also add in our considerations linear, first order differential operators
\beq\label{B-interaction}
B_\ell(\phi)=\sum_{\kappa=1}^L a^\kappa_\ell\phi_\kappa(x)
+ \sum_{\kappa,\alpha,\beta =1}^Lb^{\kappa,\a,\beta}_\ell \phi_\kappa (x)
\phi_\alpha (x)\phi_\beta (x)\hspace{-1cm}
\eeq modeling interaction of matter fields,
 where $a^{\kappa}_\ell,b^{\kappa,\a,\beta}_\ell\in \R$ are constants. 
\HOX{ We have  added here some first order terms
to the wave equation that connect $\phi$-variables together. WE NEED TO ADD 
TERMS IN STRESS ENERGY TENSOR AND IN THE APPENDIX B AND THE MAIN THEOREM}}

Now we are ready to formulate the direct problem for 
the Einstein-scalar field equations. Let $g$ and $\phi$ satisfy
\beq\label{eq: adaptive model}
& &\Ein_{\hat g}(g) =P_{jk}+Zg_{jk}+{\bf T}_{jk}(g,\phi),\quad Z=\sum_{\ell=1}^L S_\ell\phi_\ell,\quad
\hbox{in }\hattuM _0,
\\ \nonumber
& &\square_g\phi_\ell +\mathcal V^{\prime}( \phi_\ell;S_\ell) =0
\generalizations{+B_\ell(\phi)}
\quad
\hbox{in }\hattuM _0,\quad \ell=1,2,3,\dots,L,
\\ \nonumber
& &S_\ell={\mathcal S}_\ell(\phi, \nabla^g
 \phi,Q,\nabla^g Q,P,\nabla^g P,g),\quad
\hbox{in }\hattuM _0,
 \\ \nonumber
& & g=\hat g,\quad \phi_\ell=\hat \phi_\ell,\quad
\hbox{in }\hattuM _0\setminus \K.
\eeq
Above, $\mathcal V^{\prime}(\phi ;s)=\p_\phi \mathcal V(\phi;s)$ so that
$
\mathcal V^{\prime}( \phi_\ell;S_\ell)=m\phi_\ell-S_\ell.
$
We assume that  the background fields $\hat g$,  $\hat \phi$,
 $\hat Q$, and  $\hat P$ satisfy these equations.

We consider  here $P=(P_{jk})_{j,k=1}^4$
and $Q=(Q_\ell)_{\ell=1}^K$   as   fields that we can control. 
As mathematical idealization we will assume that the fields $P-\hat P$ and $Q-\hat Q$ are compactly
supported.
%
%
To obtain a physically meaningful model,
we need to consider how  the adaptive source functions $\mathcal S_\ell$ should be chosen 
so that the physical conservation law in relativity,
\beq\label{conservation law0}
\nabla_k(g^{kp} T_{pq})=0
\eeq
 is satisfied. Here $\nabla=\nabla^g$ is the connection corresponding to the metric $g$. 
  We
note that the conservation law 
is a necessary condition for the equation (\ref{eq: adaptive model})
to have solutions for which
$\Ein_{\hat g}(g)=\Ein(g)$, i.e., that the solutions of
(\ref{eq: adaptive model}) are solutions of  the Einstein field equations.

 The functions ${\mathcal S}_\ell(\phi, \nabla^g
 \phi,Q,\nabla^g Q,P,\nabla^g P,g)$ model the devices that we use to perform active
 measurements. Thus, even though the Assumption S below may appear quite technical,
 this assumption can be viewed as the instructions on  how to build a device that can be used 
 to measure the structure of the space time far away. Outside the support
 of the measurement device (i.e. the union of the supports of $Q$ and $P$) we have just assumed that the standard coupled Einstein-scalar field
 equations hold, c.f. (\ref{S-vanish condition}).

Throughout the paper we assume that the following assumption holds.
\medskip 

{\bf Assumption S}.   \HOX{In the final paper Assumption S should be written in a clearer way.}
Throughout the paper we assume 
that 
the adaptive source functions ${\mathcal S}_\ell(\phi, \nabla^g
 \phi,Q,\nabla^g Q,P,\nabla^g P,g)$ have the following properties: 
\smallskip

{(i) Denoting  $c=\nabla^g
 \phi$, $C=\nabla^g P$, and $H=\nabla^g Q$,
 we assume 
  that ${\mathcal S}_\ell(\phi, c,Q,H,P,C,g)$
 are linear functions of $(Q,H,P,C)$ and in
 particular satisfy
\beq\label{S-vanish condition}
 {\mathcal S}_\ell(\phi,c,0,0,0,0,g)=0.
\eeq
We also assume that 
when $(Q_\ell)_{\ell=1}^K$
and  $(P_{jk})_{j,k=1^4}$ are sufficiently close to
$(\hat Q_\ell)_{\ell=1}^K$
and  $(\hat P_{jk})_{j,k=1^4}$, respectively,  and $\phi$ and $g$ are 
sufficiently close in the $C^1$-topology to
the background fields $\hat g$ and $\hat \phi$ then
the adaptive
source function ${\mathcal S}_\ell(\phi, \nabla^g
 \phi,Q,\nabla^g Q,P,\nabla^g P,g)$ at $x\in \hattuM _0$ is  a smooth function of
 the pointwise values $\phi(x), \nabla^g
 \phi(x),$ $Q(x),\nabla^g Q(x),P(x),\nabla^g P(x)$ and $g_{jk}(x)$.
 
\smallskip

 \HOX{Smoothness $C^{{5}}$ is not optimal. Definitions should be changed
 so that the meaning of the "mysterious" Z-field becomes clearer}

 (ii)  $Q_K=Z$, that is, one of the physical fields we directly control is the 
 density $Z$ of the stress energy tensor caused by the source fields.
 
 \smallskip

(iii)  \HOX{Change wording in (iii). When changing the text to 2 papers,
check that $ U_{\hat g}$ is defined before this.}
We assume that  ${\mathcal S}_\ell$ is independent of $P(x)$ and the dependency of ${\mathcal S}$ on  $\nabla^g P$ 
and $\nabla^g Q$ is only due to the dependency
in the term $g^{pk}\nabla^g_p( P_{jk}+Zg_{jk})=
g^{pk}\nabla^g_p P_{jk}+\nabla^g_jQ_K$, associated to
the divergence of the perturbation of $T$,
that is,
there exist functions $\tilde {\mathcal S}_\ell$ so that
\ba
{\mathcal S}_\ell(\phi, c,Q,H,P,C,g)=\tilde {\mathcal S}_\ell(\phi,c,Q,
R,g),\quad R=(g^{pk}\nabla^g_p ( P_{jk}+Q_Kg_{jk}))_{j=1}^4.\ea
Let $\hat R=\hat g^{pk}\hat\nabla_p \hat P_{jk}+\hat\nabla_j\hat Q_K$.
Moreover,   we assume that  for all
$x\in   U_{\hat g}$ 
 the derivative of $\tilde {\mathcal S}(\hat \phi, \hat \nabla
 \hat \phi,Q,R,\hat g)=(\tilde {\mathcal S}_\ell(\hat \phi, \hat \nabla
 \hat \phi,Q,R,\hat g)_{\ell=1}^L$ with respect to $Q$ and $R$, that is, the map
  \beq\label{surjective 1}
D_{Q,R}\tilde {\mathcal S}(\hat \phi, \hat \nabla
 \hat \phi,Q,R,\hat g)|_{Q=\hat Q, R=\hat R}:\R^{K+4}\to \R^L
 \eeq
 is surjective.
\smallskip
 
(iv)  We assume that 
the adaptive source functions ${\mathcal S}_\ell$ are such 
 that  if $g,\phi$  satisfy (\ref{eq: adaptive model}) with any
 $(Q_\ell)_{\ell=1}^K$
and  $(P_{jk})$ that are sufficiently close to $\hat Q$ and $\hat P$ in $C^{{5}}$-topology then
the conservation
 law (\ref{conservation law0}) is valid.
}
 

 \bigskip

Notice that as sources are small, we need to consider
only local existence of solutions, see Appendix C.
Above, the assumptions on the smoothness of the sources and solutions
are far from optimal. For the local existence results, see \cite{K1,K3,K4,R1}.
The global existence results are considered e.g.\ in \cite{C1,Ch-K,Li1,Li2}.

Below, expect in Corollary \ref{coro of main thm original Einstein},
 we will consider the case when $\hat Q=0$ and $\hat P=0$. This implies
 that for the background fields that adaptive source functions $\mathcal S_\ell$ vanish.

Below we will denote $Q=(Q^{\prime},Q_K)$, $Q^{\prime}=(Q_\ell)_{\ell=1}^{K-1}$.
There are examples when the background fields $(\hat g,\hat \phi)$ and 
the adaptive source functions ${\mathcal S}_\ell(\phi, \nabla^g
 \phi,Q,\nabla^g Q,P,\nabla^g P,g)$ exists and satisfy the Assumption S. 
  This is shown in Appendix B in the case the following condition is valid for the \HOX{We need to check Appendix B carefully.}
   background fields:
 \medskip

 {\bf Condition A}:
Assume that at any $x\in \hat U$ there is
a permutation $\sigma:\{1,2,\dots,L\}\to \{1,2,\dots,L\}$, denoted $\sigma_x$, such that the 
$5\times 5$ matrix $[ B_{jk}^\sigma(\hat \phi(x),\nabla \hat \phi(x))]_{j,k\leq 5}$
is invertible, where
  \ba
[ B_{jk}^\sigma(\phi(x),\nabla \phi(x))]_{k,j\leq 5}=\left[\begin{array}{c}
(\,\p_j  \phi_{\sigma(\ell)}(x))_{\ell\leq 5,\ j\leq 4}\\
(\phi_{\sigma(\ell)}(x))_{\ell\leq 5}\end{array}\right].
 \ea
  \medskip
 
%
%

\subsubsection{Inverse problem for the reduced Einstein equations}

Let us next define the measurements precisely. 

%

 Let $\hat \mu:[-1,1]\to \hattuM _0$
be a freely falling observer of $(\hattuM _0,\hat g)$
and  $z_0=\hat \mu(s_{-2})$
and $\eta_0=\p_s\hat \mu(s_{-2})$. We assume that 
$z_0\in (-\infty,0)\times N$ and $\hat g(\eta_0,\eta_0)=-1$.


We will consider Lorentzian metrics $g$ on $\hattuM _0= (-\infty,t_0)\times N$, $t_0>0$
that is sufficiently close in $C^2_b(\hattuM _0)$ to $\hat g$ and 
 coincides with $\hat g$ in $ (-\infty,0)\times N$. 
For the metric $g$ we will use the notations of an open set $U_g$, 
freely falling observers $\mu_g$ and $\mu_{g,z,\eta}$ with 
$(z,\eta)\in \U_{z_0,\eta_0}$, and numbers $-1<s_{-2}<s_{-1}=s_-<s_{+1}=s_+<s_{+2}<1$,
  that are
defined  near the formula 
(\ref{eq: Def Wg with hat}), with $\hat h$ being small enough and 
independent of $g$ and   so that the set $\U_{z_0,\eta_0}(\hat h)\subset (-\infty,0)\times N$ 
 is the same for all metric tensors $g$ that we consider.
 We assume also that $\hat p^-=\hat\mu(s_-)\in (0,t_0)\times N$.

Let $\mu_{g,z,\eta}:[-1,1]\to M_0$ be geodesics  such that $\mu_{g,z,\eta}(s_{-2})=z$
and $\p_s\mu_{g,z,\eta}(s_{-2})=\eta$. Then for 
$\mu_g:=\mu_{g,z_0,\eta_0}$ we have
$\mu_g(s)=\mu_{\hat g}(s)$ for all $s\leq s_{-1}$ and
 $\mu_{\hat g}=\hat \mu$. 
Moreover, we
denote $\hat p^-=\hat \mu(s_-)$ and $\hat p^+=\hat \mu(s_+)$. 
We assume that $\hat h$ used to define 
$\U_{z_0,\eta_0}(\hat h)$ is independent of the metric and that it is so small that $\pi(\U_{z_0,\eta_0}(\hat h))\subset  
 (-\infty,t_0)\times N$.
  
We note that
when $g$ is close enough the $\hat g$ in the space $C^2_{loc}(\hattuM _0)$,
for all $(z,\eta)\in \U_{z_0,\eta_0}(\hat h)$ we have $\mu_{g,z,\eta}([s_{-2},s_{+2}])\subset U_g$,
$\mu_{g,z,\eta}(s_{-2})\in I_{g}^-(\mu_g(s_-))$,
and 
$\mu_{g,z,\eta}(s_{+2})\in I_{g}^+(\mu_g(s_+))$.

Moreover, for $r>0$ let
\beq\label{observer neighborhood with hat}
& &W_g(r)=\bigcup_{s_-<s<s_+-r} I_{M,g}( \mu_g(s),\mu_g(s+r))
\eeq
and let $r_0\in (0,1)$ be so small that  $W_{\hat g}(2r_0)\subset U_{\hat g}.$ 
%
We denote next $W_{g}=W_{g}(r_0)$.

%
%

%

Let us  use causal Fermi-type coordinates:
Let $Z_j(s)$, $j=1,2,3,4$
be a parallel frame of linearly independent time-like vectors on $\mu_g(s)$ such
that $Z_1(s)=\dot\mu_g(s)$. Let $\Phi_g:(t_j)_{j=1}^4\mapsto
\exp_{\mu_g(t_1)}(\sum_{j=2}^4 t_jZ_j(t_1))$.  
We assume that $r_0>0$ used above is so small that 
$\Psi_{\hat g}=\Phi_{\hat g}^{-1}$ defines coordinates in $W_{\hat g}(2r_0)$.
Then, when $g$ is sufficiently close to $\hat g$ in the $C^2_b$-topology,
$\Psi_g=\Phi_g^{-1}$ defines coordinates in $W_g(2r_0)$ that we call
the Fermi-type coordinates.
We define 
the norm-like functions
\ba
& &\mathcal N_{\hat g}^{(k)}(g)=
\|(\Psi_{g})_*g-(\Psi_{\hat g})_*\hat g)\|_{C^{k}(\overline {\Psi_{\hat g}(W_{\hat g})})},\\
& &\mathcal N^{(k)}(F)=
\|(\Psi_{g})_*F\|_{C^{k}(\overline {\Psi_{\hat g}(W_{\hat g})})},
\ea
where $k\in \N$,
that measures the $C^{k}$ distance of $g$ from $\hat g$ and $F$ from zero in
a Fermi-type coordinates. As we have assumed 
that the background metric $\hat g$ and the field $\hat \phi$
are $C^\infty$-smooth, we can  consider as smooth sources
as we wish. Thus, for clarity, we use below  smoothness assumptions on sources 
that are far from the optimal ones.
\medskip

}

Let us  define  (recall that here $\hat Q=0$ and $\hat P=0$)
the source-to-observation 4-tuples\hiddenfootnote{We need to study, how  source-to-observation pairs determine data analogous to the 
second derivative of "measurement map". 
} 
  corresponding to me\-as\-ur\-em\-ents in  $U_{g}$ with sources $F=(P,Q)$ 
 supported in $U_{g}$. 
We define 
\beq\nonumber
{\cal D}(\hat g,\hat \phi,\e)=\{[(U_g,g|_{U_g},\phi|_{U_g},F|_{U_g})]&;&(g,\phi,F)\hbox{ are smooth 
solutions}\\ 
& &  \label{eq: main data}
\hspace{-4cm}\hbox{ of  
(\ref{eq: adaptive model}) with $F=(P,Q)$, }F\in C^\infty_0(W_{g};\B^K),
\\  \nonumber
& &\hspace{-6cm}\hbox{$J_g^+(\supp(F))\cap J_g^-(\supp(F))\subset 
W_g, $ $\mathcal N^{({ {16}})}(F)<\e$, $\mathcal N^{({{16}})}_{\hat g}(g)<\e$}\}.
\eeq
Above, the sources $F$ are considered as sections of the bundle $\B^K$,
where $\B^K$ is a vector bundle on $\hattuM $ that is the product bundle of
the bundle  symmetric $(2,0)$-tensors and the trivial vector bundle with the fiber $\R^K$.
Above, $[(U_g,g,\phi,F)]$ denotes  equivalence class of
all
Lorenztian manifolds $(U^{\prime},g^{\prime})$ and  functions $\phi^{\prime}=(\phi^{\prime}_\ell)_{\ell=1}^L$ and the  tensors $F^{\prime}$ defined on 
$C^\infty$-smooth manifold $U^{\prime}$,
such that there is $C^\infty$-smooth 
diffeomorphism $\Psi:U^{\prime}\to U_g$ satisfying $\Psi_*g^{\prime}=g$, $\Psi_*\phi_\ell^{\prime}=\phi_\ell$  and $\Psi_*F=F$.

Note that as $[(U_{\hat g},\hat g,\hat \phi,0)]$ is the intersection of all collections ${\cal D}(\hat g,\hat \phi,\e)$, $\e>0$, 
we see that
the collection $ {\cal D}(\hat g,\hat \phi,\e)$ determines the isometry type of $(U_{\hat g},\hat g)$.

Our goal is to prove the following result:

\begin{theorem}\label{main thm Einstein}
Let 
$(\hattuM _j ,\hat g_j)$, $j=1,2$ be two 
open,  $C^\infty$-smooth, globally hyperbolic   Lorentzian manifolds
and  $\hat \phi^{(j)}$, $j=1,2$ be background values of the scalar fields on these manifolds,
and let $\mathcal S_\ell$ be the adaptive source functions satisfying  (\ref{S-vanish condition})
and  Assumption S, and
assume that background source fields vanish,  $\hat P^{(j)}=0$ and $\hat Q^{(j)}=0$.

Let $\mu_{\hat g_j}:[-1,1]\to \hattuM _j $ be freely falling observers on  $(\hattuM _j  ,\hat g_j)$
and $U_{\hat g_j}$ be the open sets
defined by formula (\ref{eq: Def Wg with hat}), and $\hat p_j^\pm=\mu_{\hat g_j}(s_{\pm})$.

Assume that there is a $C^{17}$-isometry $\Psi_0:(U_{\hat g_1},\hat g_1)\to (U_{\hat g_2},\hat g_2)$.
We identify these isometric sets and denote
\ba
\hat U=U_{\hat g_1}=U_{\hat g_2},\quad 
\hat g|_{\hat U}=\hat g_1|_{U_{\hat g_1}}=\hat g_2|_{U_{\hat g_2}}.
\ea
Let $\e>0$ and assume that source-to-observation 4-tuples ${\cal D}(\hat g_j,\hat \phi^{(j)},\e)$ 
for the manifolds $(\hattuM _j ,\hat g_j)$ and fields $\hat \phi^{(j)} $ 
satisfy
\beq\label{eq: data are the same}
{\cal D}(\hat g_1,\hat \phi^{(1)},\e)= {\cal D}(\hat g_2,\hat \phi^{(2)},\e).
\eeq

Then  there is a diffeomorphism $\Psi:I_{\hattuM  _1 ,\hat g_1}(\hat p^-_1,\hat p^+_1)\to 
I_{\hattuM _2 ,\hat g_2}(\hat p^-_2,\hat p^+_2)$,
and the metric $\Psi^*\hat g_2$ is conformal to $\hat g_1$ in $I_{\hattuM _1 ,\hat g_1}(\hat p^-_1,\hat p^+_1)$.
\end{theorem}

Note that above we have assumed that the adaptive source functions $\mathcal S_\ell$ are
the same on $(\hattuM _1 ,\hat g_1)$ and $(\hattuM _2 ,\hat g_2)$.
%

Recall that above $I_{\hattuM _1 ,\hat g_1}(\hat p^-_1,\hat p^+_1)=I_{\hattuM _1 ,\hat g_1}^+(\hat p^-_1)\cap I_{\hattuM _1 ,\hat g_1}^-(\hat p^+_1)$.

   \observation{
\medskip

{\bf Remark 1.1.} The result of \HOX{Check the Remark 3.3 very carefully
or remove this and the other remarks.}
Theorem \ref{main thm Einstein} can be improved in
the case when the adaptive source functions $\mathcal S_\ell$ are
the ones constructed in Appendix B: 
As the background source fields vanish, $\hat Q=0$ and $\hat P=0$
we can consider the sources where $Q=0$ and we observe
only the $g$-components of the wave. 
This means that instead of assuming
that we are given the data set ${\cal D}(\hat g,\hat \phi,\e)$, we can assume that
we are given 
\beq\nonumber
{\cal D}_0(\hat g,\hat \phi,\e)&=&\{
[(U_g,g|_{U_g},F|_{U_g})]
;\ [(U_g,g|_{U_g},\phi|_{U_g},F|_{U_g})]
\in {\cal D}(\hat g,\hat \phi,\e),\\ 
& &  \label{eq: main data 2}
\hbox{  where $F=(P,Q)$ is such that $Q=0$} \}.
\eeq
Then a result  similar to Theorem \ref{main thm Einstein}
can be obtained. We explain the  tools needed for this result, see e.g. 
Remark 5.1 below, but give the detailed proof elsewhere.\bigskip }

The measurements in  Theorem \ref{main thm Einstein} provides a 
subset of all possible  sources $(\F_1,\F_2)$ for equations
(\ref{eq: adaptive model with no source}).  Thus, as we see later,  using
Theorem \ref{main thm Einstein} 
we can prove the following
result where we assume  that we have information on
 measurements with a larger class of sources than was used in Theorem \ref{main thm Einstein}: 
\HOX{T3: We should consider if the "alternative" formulation on main
theorem, Thm. \ref{alternative main thm Einstein}, is better than Thm. 1.2.}

\medskip

 \begin{theorem}\label{alternative main thm Einstein} Assume that $(\hattuM ,\hat g)$ is a globally hyperbolic
manifold, and that in the open set $U_{\hat g}$ the  Condition A is valid. 
Assume that we given
the set ${\cal D}^{alt}(\hat g,\hat \phi,\e)$ of the 
 the equivalence classes
$[(U_g,g|_{U_g},\phi|_{U_g},\F|_{U_g})]$  where $g$ and $\phi$ and $\F=(\F_1,\F_2)$  satisfy the equations 
\HOX{T4: 
The formulation and the proof of Thm. \ref{alternative main thm Einstein} need to be checked
 and indexes of the norms of sources have to be checked in the claim.}
(\ref{eq: adaptive model with no source}), the conservation law
$\nabla^g_j({\bf T}^{jk}(g,\phi)+\F_1^{jk})=0$, 
the sources $\F_1$ and $\F_2$ are supported in
$U_g$, and satisfy $\mathcal N^{({ {15}})}(\F)<\e$ and $\mathcal N^{({{15}})}_{\hat g}(g)<\e$. Then these data 
determine
the conformal type of $\hat g$ in $I_{\hattuM  ,\hat g}(\hat p^-,\hat p^+)$.
 \end{theorem}

\medskip

The above result means that if the manifold $(\hattuM _0,\hat g)$ is unknown, then
the source-to-observation pairs corresponding to
freely falling sources which are near the 
freely falling observer $\mu_{\hat g}$ and the measurements of 
the metric tensor and the scalar fields in a neighborhood $U_{\hat g}$ of  $\mu_{\hat g}$,
 determine the  metric tensor up to conformal transformation
in the set $I_{\hattuM _0,\hat g}(\hat p^-,\hat p^+)$.


We want to point out that by the main theorem, if we have
two non-isometric  space times, a generic
measurement gives different results
on these manifolds. In particular, this implies that the perfect 
space-time cloaking, see \cite{Fridman,McCall}, with a smooth metric  in a globally hyperbolic universe  is not possible.

Also, one can ask if 
one can make an approximative image of the space-time knowing
 only one measurement.  In general, in many
inverse problems several measurements can be packed together to 
one measurement.
For instance, for the wave equation with a  time-independent simple metric
this is done in \cite{HLO}. Similarly, Theorem \ref{main thm Einstein} and its proof
make it possible to do approximate reconstructions in a suitable class of
manifolds with only one measurement.  We will discuss this
in detail in forthcoming articles.

\medskip

Using Theorem \ref{main thm Einstein}  we  will show that the metric tensor can be determined in the domain which can
be connected to a measurement set with light-like geodesics through vacuum.

Below, we consider the case when  $\hat Q$ and $\hat P$ are non-zero, and we define
\ba
{\cal D}^{mod}(\hat g,\hat \phi,\e)=\{[(U_g,g|_{U_g},\phi|_{U_g},F|_{U_g})]&;&(g,\phi,F)\hbox{ are smooth 
solutions,}\\
& & 
\hspace{-6cm}\hbox{ of  
(\ref{eq: adaptive model}) with $F=(P-\hat P,Q-\hat Q)$, }F\in C^\infty_0(W_{g};\B^K),
\\
& &\hspace{-6cm}\hbox{$J_g^+(\supp(F))\cap J_g^-(\supp(F))\subset 
W_g, $ $\mathcal N^{({ {16}})}(F)<\e$, $\mathcal N^{({{16}})}_{\hat g}(g)<\e$}\}.
\ea
Next we consider the case when $\hat Q^{(j)}$ and  $\hat P^{(j)}$ are not assumed to be zero,
see Fig.\ 4. 

\begin{corollary}\label{coro of main thm original Einstein}
Assume that $\hattuM _j$, $\hat g_j$, $\hat \phi^{(j)}$, $j=1,2$ and $\Psi_0$ are as in  Theorem \ref{main thm Einstein} and
 (\ref{eq: data are the same}) is not assumed to be 
valid. Assume that 
also that for $j=1,2$  there are sets $W_j\subset M_j$ such that 
$\hat \phi^{(j)}$, $\hat Q^{(j)}$ and $\hat P^{(j)}$ are zero (and thus
the metric tensors $\hat g_j$ have vanishing Ricci curvature)
in $W_j$ and that  $I_{\hat g_j}(\hat p^-_j,\hat p^+_j)\subset W_j \cup U_{\hat g}$. 
If
\ba
{\cal D}^{mod}(\hat g^{(1)},\hat \phi^{(1)},\e)={\cal D}^{mod}(\hat g^{(2)},\hat \phi^{(2)},\e)
\ea
then the metric $\Psi^*\hat g_2$ is isometric to $\hat g_1$ in $I_{\hat g_1}(\hat p^-_1,\hat p^+_1)$.
\end{corollary}
\medskip

\begin{figure}[htbp]
\begin{center}\label{Fig-4}

\psfrag{1}{}
\psfrag{2}{}
\psfrag{3}{}
\psfrag{4}{$W_{\hat g}$}
\psfrag{5}{$\mu_{z,\eta}$}
\includegraphics[width=4.5cm]{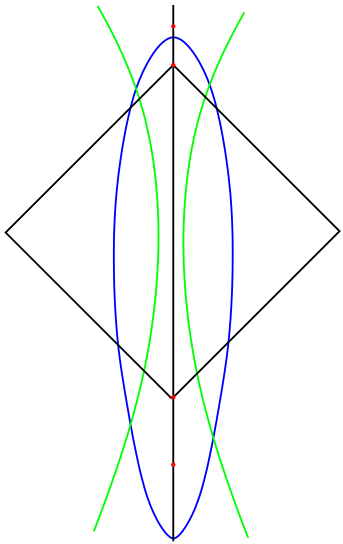}
\end{center}
\caption{A schematic figure where the space-time is represented as the  2-dimensional set $\R^{1+1}$ on the setting in Corollary \ref{coro of main thm original Einstein}.
The set $J_{\hat g}(\hat p^-,\hat p^+)$, i.e., the diamond with the black boundary, is contained in the union of the blue set $U_{\hat g}$
and the set $W$. The set $W$ is in the figure the area outside of the green curves. 
The sources are controlled in the set $U_{\hat g}\setminus W$ and  the set $W$ consists of vacuum.
}
 \end{figure}

In the setting of Corollary \ref{coro of main thm original Einstein} the set
 $W$ is such $I(\hat p^-,\hat p^+)\cap (M\setminus W)\subset U$. This means that if we restrict to the domain $I(\hat p^-,\hat p^+)$
then we have Vacuum Einstein equations in the unknown domain $I(\hat p^-,\hat p^+)\setminus U$
and have matter only in the domain $U$ where we implement our measurement (c.f.\ a space
ship going around in a system of black holes). This could be considered as an "Inverse problem
for the vacuum Einstein equations". 
%

\medskip

\noindent
{\bf Example 3:} 
Consider a black hole modeled by a Schwarzschild metric in 
$(\R^3\setminus \{0\})\times \R$
with the event horizon $\p B(R)\times\R$, $B(R)=\{y\in \R^3:\ |y|=R\}$, and add to the stress
energy tensor the (small) effect caused by the "space ship" doing the measurements.
By Theorem \ref{main thm Einstein} the measurements in
$U=(\R^3\setminus B(r))\times \R$ with any $r>R$ 
determines uniquely the metric outside the event horizon,
that is, in $W=(\R^3\setminus B(R))\times \R$. Physically a more interesting
example could be obtained  by considering several black holes
or by adding matter to the system.
\medskip

%

\medskip
We note the techniques used to prove Theorem \ref{main thm Einstein} are suitable for
studying other equations, e.g.\ Einstein-Maxwell system, where strong electromagnetic
waves and the gravitational field have an interaction. However, we do not study
these equations in this paper.

\medskip
\section{Proofs for inverse problem for light observation sets}
\subsection{Notations and definitions}


%

A smooth $n$-dim\-ens\-ion\-al manifold with smooth Lorentzian, type  $(1,n-1)$ metric $g$ has a causal structure, if there is a
 globally defined smooth vector field $X_c$
such that $g(X_c(y),X_c(y))<0$ for all $y\in M$. We say that
a piecewise smooth path $\alpha(t)$
is time-like if $g( \dot\alpha(t),\dot\alpha(t))< 0$
for almost every $t$. Also,
a piecewise smooth path $\alpha(t)$
is causal (or non-space-like) if $\dot\alpha(t)\not =0$ and $g( \dot\alpha(t),\dot\alpha(t))\leq 0$
for almost every $t$.

By \cite{Bernal}, a globally hyperbolic manifold, as defined in the
introduction, satisfies the  strong causality condition:
\beq\label{strong claim}
& &\hbox{For every  $z\in M$ and every neighborhood $V\subset M$ of $z$ there is 
}\\
\nonumber & &\hbox{a neighborhood $V^{\prime}\subset M$ of $z$ that if $x,y\in V^{\prime}$ and  $\alpha\subset M$
is }
\\
\nonumber & &\hbox{a causal path connecting $x$ to $y$ then $\alpha\subset V$.}
\eeq


Recall that $g^+$ is a Riemannian metric determined by a Lorentzian metric $g$ on $M$. In the following,
for open $V\subset M$, denote $TV=\{(x,\xi)\in TM:\ x\in V\} $
and $S^{g^+}V=\{(x,\xi)\in TM;\ x\in V,\ \|\xi\|_{g^+}=1\}$.

The metric $g^+$ on $M$ induces a Sasaki metric on $TM$
which we denote by $g^+$, too.

In the next sections, $g$ will be fixed and we denote $U=U_g$, $J^\pm_g(x)=J^\pm(x)$, etc.
Let us next define  the earliest observation functions on a geodesic $\gamma_{g,z,\eta}$,
$(z,\eta)\in \U_{z_0,\eta_0}$.
If $V\subset TU$ and $W=\pi(V)\subset U$ are such that
 $\mu_{g,z,\eta}\cap W\not=\emptyset$,
we define ${\bf s}_{z,\eta}(W)=
\inf\{s\in [-1,1];\ 
\mu_{g,z,\eta}(s)\cap W\}$ and
\beq\label{Earliest element sets}
\ear_{z,\eta}(V)&=&\{(x,\xi)\in \hbox{cl}\,(V);\ x=\mu_{g,z,\eta}({\bf s}_{z,\eta}(\pi(V)))\},\\
\nonumber
\pointear_{z,\eta}(V)&=&\pi(\ear_{z,\eta}(V)),\\
\nonumber
\pointear_{z,\eta}(W)&=&\{x\in \hbox{cl}\,(W);\  x=\mu_{g,z,\eta}({\bf s}_{z,\eta}(W))\},
\eeq
where $\hbox{cl}\,(V)$ denotes the closure of $V$.
Moreover, we denote
$\E_{z,\eta}(q)=\ear_{z,\eta}(\O_U(q))$, $\P_{z,\eta}(q)=\pointear_{z,\eta}(\O_U(q))$.

Finally,
 ${}^{\flat }:TM\to T^*M$ is the lowering index operator by the metric tensor
$(a^j\frac \p{\p x^j})^{\flat }=g_{jk}a^j dx^k$ and ${}^{\sharp }:T^*M\to TM$ its inverse.


%

\subsection{Determination of the conformal class of the metric}

\begin{lemma} \label{A: lem 1} Let $z\in M$.
Then
there  is a neighborhood $V$ of $z$ so that
\smallskip 

\noindent
(i) If the geodesics $\gamma_{y,\eta}([0,s])\subset V$ and $\gamma_{y,\eta^{\prime}}([0,s^{\prime}])\subset V$,
$s,s^{\prime}>0$
satisfy $\gamma_{y,\eta}(s)=\gamma_{y,\eta^{\prime}}(s^{\prime})$,
then $\eta= c\eta^{\prime}$  and $s^\prime=cs$ with some  $c>0$. 
\smallskip 

\noindent
(ii) For any $y\in V$, $\eta\in T_yM\setminus 0$ there is $s>0$ such that $\gamma_{x,\eta}(s)\not \in V$.
\end{lemma}

\noindent
{\bf Proof.} The property (i) follows from 
\cite[Prop.\ 5.7]{ONeill}.  Making  $s>0$ so small that $\overline V\subset B_{g^+}(z,\rho)$ with 
a sufficiently small $\rho$, the claim (ii) from \cite[Lem.\ 14.13]{ONeill}. 
%
\hfill \Box \medskip


Let $q^-,q^+\in M$. As $M$ is globally hyperbolic,   the set $J(q^-,q^+)=J^-(q^+)\cap J^+(q^-)$ is compact.

Let us consider points $x, y\in M$. If  $x< y$,
we define  the time separation function $\tau(x,y)\in [0,\infty)$ from $x$ to $y$
to be the supremum of the lengths 
\ba
L(\alpha)=\int_0^1 \sqrt{-g(\dot\alpha(s),\dot\alpha(s))}\,ds
\ea
of the piecewise smooth
causal paths $\alpha:[0,1]\to M$ from $x$ to $y$. If  the condition $x< y$ does not
hold, we define $\tau(x,y)=0$. 
We note that $\tau(x,y)$ satisfies the reverse triangle inequality
 \beq\label{eq: reverse}
 \tau(x,y)+\tau(y,z)\leq \tau(x,z)\quad\hbox{for }x\leq y\leq z.
 \eeq

As $M$ is globally hyperbolic,
the time separation function $(x,y)\mapsto \tau(x,y)$  is continuous in  $M\times M$
by \cite[Lem.\ 14.21]{ONeill}.

 Let us consider points $x, y\in M$, $x< y$. 
  As the set $J(x,y)$
 it compact,  by \cite[Prop.\ 14.19]{ONeill}
 the points $x$ and $y$ can be connected by a causal geodesic whose
 length is $\tau(x,y)$.   
 In particular, this implies that if 
 $y$ can be connected to
$x$ with a causal path which is not a light-like geodesic then $\tau(x,y)>0$,
see  \cite[Prop.\ 10.46]{ONeill}. 

The above facts can be combined as follows:
consider a path which is the union of  the future going light-like geodesics
$\gamma_{x_1,\theta_1}([0,t_1])\subset M$  and $\gamma_{x_2,\theta_2}([0,t_2])\subset M$,
where $x_2=\gamma_{x_1,\theta_1}(t_1)$ and $t_1,t_2>0$. Let
$\zeta=\dot \gamma_{x_1,\theta_1}(t_1)$ and $x_3=\gamma_{x_2,\theta_2}(t_2)$. 
Then, if there are no $c>0$ such that $\zeta= c\theta_2$, the union of
these geodesic is not a light-like geodesics and
thus  $\tau(x_1,x_3)>0$. In particular, then there exists a time-like geodesic 
from $x_1$ to $x_3$. 
In the following we call this kind of argument for a union of light-like geodesics a short-cut
argument.

When $(x,\xi)$ is a  light-like vector, we define  $\T(x,\xi)$ to
be the length of
the maximal interval on which $\gamma_{x,\xi}:[0,\T(x,\xi))\to M$ is defined.

When $(x,\xi_+)$ is a  future pointing light-like vector,
and $(x,\xi_-)$ is a  past pointing light-like vector,
we define the modified cut locus functions, c.f.\ \cite[Def.\ 9.32]{Beem},
$\rho(x,\xi_\pm)=\rho_g(x,\xi_\pm)$,
\beq\label{eq: max time}
& &\rho(x,\xi_+)=\sup\{s\in [0,\T(x,\xi)):\ \tau(x,\gamma_{x,\xi_+}(s))=0\},\\
\nonumber
& &\rho(x,\xi_-)=\sup\{s\in [0,\T(x,\xi)):\ \tau(\gamma_{x,\xi_-}(s),x)=0\}.
\eeq
The point $\gamma_{x,\xi}(s)|_{s=\rho(x,\xi)}$ is called the cut point
on the geodesic  $\gamma_{x,\xi}$.
Below, we  say that  $s=\rho(x,\xi)$ is a cut value
on $\gamma_{x,\xi}([0,a])$.

By \cite[Thm. 9.33]{Beem}, the function $\rho(x,\xi)$ is lower semi-continuous on a globally hyperbolic
Lorentzian manifold $(M,g)$. We note that  by \cite[Thm. 9.15]{Beem},
a cut point $\gamma_{x,\xi}(s)|_{s=\rho(x,\xi)}$ is either a conjugate point or a null cut point, that is, an intersection point of two light geodesics
starting from the point $x$.
Note that by \cite[Cor.\ 10.73]{Beem} the infimum of
the  null cut points 
is smaller or equal to the first null conjugate point.

%

Let $p^\pm=\mu_g(s_\pm)$.
By \cite[Lem.\ 14.3]{ONeill},
the set $I(p^-,p^+)=I^-(p^+)\cap  I^+(p^-)$ is open.
We need the following first observation time function $f^+_\mu$.

Recall that by  formula (\ref{kaava D}),
$p^\pm=\mu_g(s_\pm)$ satisfy $p^\pm\in I^\mp(\mu_{g,z,\eta}(s_{\pm 2}))$
for all ${z,\eta}\in \U_{z_0,\eta_0}$.

\begin{definition}\label{def: f functions}
Let $\mu=\mu_{g,z,\eta}$, ${z,\eta}\in \U_{z_0,\eta_0}$. 
For $x\in J^-(p^+)\setminus I^-(p^-)$ we define $f_\mu^+(x)\in [-1,1]$
by setting 
\ba
& &
h_+(s)=\tau(x,\mu(s)), \quad
A_\mu^+(x) =\{s\in (-1,1):\ h_+(s)>0\}\cup\{1\},
\\
& & f_\mu^+(x) =\inf A_\mu^+(x).
\ea
Similarly, 
for $x\in J^+(p^-)\setminus I^+(p^+)$ we define $f_\mu^-(x)\in [-1,1]$
by setting 
\ba
& &
h_-(s)=\tau(\mu(s),x),\quad  A_\mu^-(x)=\{s\in (-1,1):\ h_-(s)>0\}\cup\{-1\},
\\
& &f_\mu^-(x) =\sup A_\mu^-(x).
\ea 
\end{definition}

We need the following simple properties of these functions.

\begin{lemma} \label{B: lemma} 
Let $\mu=\mu_{z,\eta}$, ${z,\eta}\in \U_{z_0,\eta_0}$,
and   $x\in  J^-(p^+)\setminus I^-(p^-)$.
Then

\smallskip 

\noindent
(i) The function $s\mapsto \tau(\mu(s),x)$ is non-decreasing on the interval $s\in [-1,1]$
and strictly increasing on $s\in [f_\mu^+(x),1]$.
\smallskip 

\noindent
(ii) It holds that
 $s_{-2}< f_\mu^+(x)<s_{+2}$.
\smallskip 

\noindent (iii)  Let $y=\mu(f_\mu^+(x))$. Then $\tau(x,y)=0$.
Also, if $x\not\in \mu$, there is a light-like geodesic $\gamma([0,s])$ in $M$ from $x$ to  $y$ with no conjugate points
on $\gamma([0,s))$.

\smallskip 

\noindent (iv)
The map $f^+_\mu:J^-(p^+)\setminus I^-(p^-)\to (-1,1)$ is continuous. 

\smallskip 

\noindent (v)  For $q\in J^-(p^+)\setminus I^-(p^-)$ the map $F:\U_{z_0,\eta_0}\to \R;$
$F(z,\eta)=f^+_{\mu(z,\eta)}(q)$ is continuous. 

\end{lemma}

Note that above we consider also a single point path as a light-like geodesic.

\noindent
{\bf Proof.} (i) \HOX{The proof should be shortened in the final paper version.}
As $\mu$ is a time like-path, we have $\tau(\mu(s),\mu(s^{\prime}))>0$ for $s<s^{\prime}$. 
When $x\leq \mu(s)$ and $s<s^{\prime}$, the 
 reverse triangle inequality  (\ref{eq: reverse}) 
yields that $ \tau(x,\mu(s))<\tau(x,\mu(s^{\prime}))$. As $\tau(x,\mu(t))=0$ for $t\leq f^+_\mu(x)$, (i) follows. 

(ii) Recall that  $p^\pm\in I^\mp(\mu(s_{\pm 2}))$ by (\ref{kaava D}). As $x\in  I^-(p^+)$,
we see that $\tau(x,\mu(s_{+2}))>0$ .
Due to the global hyperbolicity, $\tau$ is continuous in $M\times M$,
see \cite[Lem.\ 14.21]{ONeill}. Therefore,
$h_+(s)=\tau(x,\mu(s))>0$ for some $s<s_{+2}$ and we see that $s_+=f_\mu^+(x)<s_{+2}$. 

Assume next that $f^+_\mu(x)\leq s_{-2}$. Since then $h_+$ is
strictly increasing, we have
$h_+(s)=\tau(x,\mu(s))>0$ for all  $s\in (s_{+2},1]$.
Let  $s_j>s_{-2}$ be such that $s_j\to s_{-2}$ as $j\to \infty$ and 
$h_+(s_j)>0$, so that $x<\mu(s_j)$. As $J^+(x)\cap J^-(p^+)$ 
is closed, then  
$x\leq q_-=\mu(s_{-2})<p^-$. This is not possible since 
$x\in  J^-(p^+)\setminus I^-(p^-)$.
Thus (ii) is proven.

(iii) 
 Let $s_x^+=f_\mu^+(x)$ and
$s_j^+< s_x^+$ be such that $s_j^+\to s_x^+$ as $j\to \infty$.
Then, $\tau(x,\mu(s_j^+))=0$ and by continuity of $\tau$, $\tau(x,y)=0$.
 On the other hand, let $s_j^+\in A_+(x)$ be such that $s_j^+\to s_x^+$ as $j\to \infty$.
Then 
$x\leq \mu(s_j^+)$, and by closedness of $J^+(x)$ (see \cite[Lem.\ 14.22]{ONeill}),
$x\leq y=\mu(s_x^+)$. As $x\not \in \mu$,  $x<y$
and by  \cite[Lem.\ 10.51]{ONeill}
there is a light-like geodesic  from $x$ to  $y$ with no conjugate points before $y$. This proves (iii).

(iv)  Assume that $x_j\to x$ in $J^-(p^+)\setminus I^-(p^-)$ as $j\to\infty$.
Let
$s_j=f_\mu^+(x_j)$ and $s=f_\mu^+(x)$.  As $\tau$ is continuous,
 for any $\e>0$ we have 
$\lim_{j\to\infty} \tau(x_j,\mu(s+\e))=\tau(x,\mu(s+\e))>0$ and thus
for $j$ large enough $x_j<\mu(s+\e).$ Thus $\lim_{j\to \infty}s_j\leq s$.
Assume next that $\liminf_{j\to \infty}s_j=\tilde s< s$ and denote $\e=\tau(\mu(\tilde s),\mu(s))>0$. Then
$\liminf_{j\to\infty}\tau(x_j,\mu(s))\geq \e$, and as $\tau$ is continuous in $M\times M$,
we obtain $\tau(x,\mu(s))\geq \e$, which is not possible as $s=f_\mu^+(x)$.
Hence $s_j\to s$ as $j\to \infty$. This proves (iv).

(v) Observe that as $J^+(q)$ is a closed set,  $F(z,\eta)$ is equal to the smallest value $s\in [-1,1]$
such that $\mu_{z,\eta}(s)\in J^+(q)$.
Let $(z_j,\eta_j)\to (z,\eta)$ in $(TM,g^+)$ 
as $j\to \infty$ and $s_j=F(z_j,\eta_j)$
and $\underline s=\liminf_{j\to \infty}s_j$.
 As the map $(z,\eta,s)\mapsto \mu_{z,\eta}(s)$
is continuous, we see that  for a suitable subsequence
 $\mu_{z,\eta}(\underline s)=\lim_{k\to \infty} \mu_{z_{j_k},\eta_{j_k}}(s_{j_k})\in J^+(q)$
 and hence $F(z,\eta)\leq \underline s=\liminf_{j\to \infty}F(z_j,\eta_j)$.
 This shows that $F$ is lower-semicontinuous.
 
 On the other hand, let $\overline s=F(z,\eta)$. As
 $\mu_{z,\eta}$ is a time-like geodesic, we see that  for any $\e\in (0,1-\overline s)$
 we have $\tau(q,\mu_{z,\eta}(\overline s+\e))>0$. As $\tau$
 and the map $(z,\eta,s)\mapsto \mu_{z,\eta}(s)$
are continuous, we see that there is $j_0$ such that 
 if $j>j_0$ then $\tau(q,\mu_{z_j,\eta_j}(\overline s+\e))>0$.
 Hence, $F(z_j,\eta_j)\leq \overline s+\e$. Thus
 $\limsup_{j\to \infty}F(z_j,\eta_j)\leq \overline s+\e$,
 and as $\e>0$ can be chosen to be arbitrarily small, we have  $\limsup_{j\to \infty}F(z_j,\eta_j)\leq \overline s=
 F(z,\eta)$. Thus $F$ is also upper-semicontinuous that proves (v).
\hfill \Box \medskip

Similarly, under the assumptions of Lemma \ref {B: lemma},
we see that if $x^{\prime}\in  I^+(p^-)\setminus J^+(p^+)$ then
 the function $s\mapsto \tau(\mu(s),x^{\prime})$ is non-increasing on the interval $s\in [-1,1]$
 and strictly decreasing on $[-1,f_\mu^-(x)]$. In addition, for
$s_-=f_\mu^-(x)$ we have
 $\tau(\mu(s_-),x^{\prime})=0$   and if $x\not \in \mu$,
there is a light-like geodesic from  $\mu(s^-)$ to $x^{\prime}$.
Moreover, $ f_\mu^-:J^+(p^-)\setminus I^+(p^+)\to \R$ 
is continuous.

In the following, let us consider the light observation set $\O_U:q\mapsto \O_U(q) $
as a map $\O_U:I^-(p^+)\setminus  J^-(p^-)\to P(L^+U)$.

\HOX{Improve style of text here}
Let $q\in I^-(p^+)\setminus  J^-(p^-)$ and $\mu=\mu_{z,\eta}$, $(z,\eta)\in \U_{z_0,\eta_0}$.
By  Lemma \ref{B: lemma} (iii),
$Z=\O_U(q)\subset P(TU)$ is non-empty.

First, let us consider the case when
$q\in M\setminus \mu$.
Then 
 there is a sequence  $s_j\in A_-(q)$, $s_j\to f_\mu^+(q) $ as $j\to \infty$,
and  $\zeta_j\in  S^{g^+}_{z_j}M$, $z_j=\mu(s_j)$, such that 
$(z_j,\zeta_j)\in O_U(q)$. Then  $z_j$ converge in the metric of $(M,g^+)$
to $z_0=\mu (f_\mu^+(q) )$ as $j\to \infty$.
As  the set  $J^+(q)\cap J^-(p^+)$ is closed, 
we see that $q\leq z_0$. Moreover,  we see as above that 
$\tau(q,z_0)=\lim_{j\to \infty}\tau(q,\mu(s_j)) =0$. 
Hence there is a light-like geodesic
$\gamma_{q,\theta}([0,l])$ from $q$ to $z_0$. 
Second,
 in the case when $q\in \mu$ we see that $L^+_qM\subset O_U(q)$ and $z_0=\mu (f_\mu^+(q) )=q$.
These cases show that for any  $q\in M$  
there are $z_0=\mu (f_\mu^+(q) )$ and  $\zeta_0\in L^+_{z_0}M$ 
such that $(z_0,\zeta_0)\in O_U(q)$.
Also, $z_0\in \P_U(q)$.

\begin{definition}\label{def: earliest element}
Let $\mu=\mu_{z,\eta}$, $(z,\eta)\in \U_{z_0,\eta_0}$ and
$q\in I^-(p^+)\setminus  J^-(p^-)$.
%
{\tobecheckedtext We  define the set of the earliest observations,  c.f. (\ref{Earliest element sets}),
\ba
\mathcal E_U(q)&=&\bigcup_{(z,\eta)\in \U_{z_0,\eta_0}}\E_{z,\eta}(q),
\ea 
and 
the map 
 $\be_{U}:J(p^-,p^+)\to \mathbb P(U)$ given by $\be_{U}(q)=\bigcup
 _{(z,\eta)\in \U_{z_0,\eta_0}}\be_{z,\eta}(q)$,  where
 $\be_{z,\eta}(q)= \pointear_{z,\eta}(\P_U(q))$.

}
\end{definition}
Above we have seen that  $ z_0=\mu (f_\mu^+(q) )\in \P_U(q)\subset \pi(O_U(q))$
satisfies  $\tau(q,z_0)=0$ and on the other hand, $\tau(q,\gamma_{q,\xi}(s))>0$
for $s< \rho(q,\xi)$. Using these, we see that
\beq\label{eq: two definitions of E}
\mathcal E(q)=\{\dot\gamma_{q,\xi}(s)\in U_g;\ 0\leq s\leq \rho(q,\xi),\ \xi\in L^+_qM\}.
\eeq
Below,  denote $T(I^-(x))=\{(x,\xi)\in TM_0;\ x\in I^-(x_0)\}$.

{\tobecheckedtext

\HOX{Check later if the claim can simplified
by replacing $x_0$ by $p^+$.}
\begin{lemma}
\label{lemma: from P(q) to E(q)}  
Assume that we are given $(U,g)$, $U=U_g$, $x_0\in U$, and  a set $F\subset U$
satisfying $\be_U(q_0)\cap I^-(x_0)\subset F\subset \P_U(q_0)$ for
some $q_0\in I^-(p^+)\setminus J^-(p^-)$. These data determine the set $\E_U(q_0)
\cap T(I^-(x_0))$. If we are given the above data for all $x_0\in U$,
we can determine   $\E_U(q_0)$.
\end{lemma}

\noindent
{\bf Proof.} Given $F$, we can find $\be_U(q_0)\cap I^-(x_0)$ by taking union
of points that can be represented as $\pointear_{z,\eta}(F)\in I^-(x_0)$ with some $(z,\eta)\in \U_{z_0,\eta_0}$.
Thus we may assume that we are given $\be_U(q_0)\cap I^-(x_0)$.

Below, let $y=\be_U(q_0)\cap I^-(x_0)$ be such that $y\in I^-(x_0)$.
  
  Consider a vector $\theta\in L^-_yM$  
such that $(y,\theta)\in \E_U(q_0)$. If $r>0$ is so small that
 $\tilde y=\gamma_{y,\theta}(-r)\in U_g$,
 a short cut argument shows that $\gamma_{y,\theta}([-r,0])$ 
 is the only causal geodesic in $U$
connecting $\tilde y$ and $y$. Moreover,
there is $(\tilde  z, \tilde \eta)\in \U_{z_0,\eta_0}$
such that for $\tilde \mu=\mu_{\tilde  z, \tilde \eta}$
we have $\tilde y\in \tilde \mu$. Moreover, as $\tilde y\in U$, using the
reverse triangle 
inequality we see that then
$\tilde y=\pi(\E_{\tilde z,\tilde \eta}(q_0))\in \be_U(q_0)$. Clearly, $\tilde y\in I^-(x_0)$.

On the other hand,  if there exists $ \theta_1\in L^-_yM$
such that for some $r>0$ we 
have $\gamma_{y, \theta_1}([-r,0])\subset U$
and $\tilde y_1=\gamma_{y,\theta_1}(-r)$
is such that  $\tilde y_1\in \tilde \mu_1 \cap (\be_U(q_0)\cap I^-(x_0))$ for 
some  $(\tilde  z_1, \tilde \eta_1)\in \U_{z_0,\eta_0}$
and   $\tilde \mu_1=\mu_{\tilde  z_1, \tilde \eta_1}$,
 a short cut argument shows that then  we must have
$(y,\theta_1)\in \O(q_0)$. 
This shows that knowing $\be_U(q_0)\cap I^-(x_0)$ we can determine
$ \E_U(q)\cap T(I^-(x_0))$. 

Finally, we observe that if we can find  $ \E_U(q)\cap T(I^-(x_0))$  for all $x_0\in U$,
by taking union of these  sets
we find $ \E_U(q)$.
\hfill \Box \medskip
}

\begin{lemma}
\label{lemma: global injectivity}  
Let  
$q_1,q_2\in I^-(p^+)\setminus J^-(p^-)$ be such that  $\pointear_{z,\eta}(\O_U(q_1))=
\pointear_{z,\eta}(\O_U(q_2))$ where  $(z,\eta)\in \U_{z_0,\eta_0}$. 



If $\pointear_{z,\eta}(\O_U(q_1))$
has a neighborhood $V\subset M$ such that
\ba
\E_U(q_1)\cap TV=
\E_U(q_2)\cap TV,
\ea 
then  $q_1=q_2$.

In particular, if $q_1,q_2\in I^-(p^+)\setminus J^-(p^-)$ are such that
$\O_U(q_1)=\O_U(q_2)$,
then  $q_1=q_2$.
\end{lemma} 

The above lemma can be also stated as follows: a germ
 of the set $\O_U(q)$
near $\pointear_{z,\eta}(\O_U(q))$ determine $q\in I^-(p^+)\setminus J^-(p^-)$ uniquely.\medskip

\noindent
{\bf Proof.} Let $\mu=\mu_{g,z,\eta}$ and 
$Z:=\E_U(q_1)\cap TV=\E_U(q_2)\cap TV$. 

Let us assume that $q_1\not=q_2$.
Let $(y,-\zeta)\in \ear_{z,\eta}(Z)$. Then
$f^+_\mu(q_1)=f^+_\mu(q_2)$ and
$y=\mu(f_\mu^+(q_1))$. Let $s_1,s_2\geq 0$ be such that
$q_1=\gamma_{y,\zeta}(s_1)$ and $q_2=\gamma_{y,\zeta}(s_2)$.
Without loss of generality, we can assume that $s_1<s_2$.
Let $\theta_j= -\dot \gamma_{y,\zeta}(s_j)\in L^+_{q_j}M$, $j=1,2$,
and let consider a vector $\theta^{\prime}_2\in L^+_{q_2}M$ 
satisfying $\theta_2^{\prime}\not=\theta_2$ and  $\|\theta^{\prime}_2\|_{g^+}=\|\theta_2\|_{g^+}$.
We assume that $\theta_2^{\prime}$  is so close to $\theta_2$ in the Sasaki metric
of $(TM,g^+)$ 
 and that $s_2^{\prime}<s_2$ is so close to $s_2$
that $y^{\prime}=\gamma_{q_2,\theta_2^{\prime}}(s_2^{\prime})\in V.$ 
Moreover, as the function $(x,\xi)\mapsto \rho(x,\xi)$ is lower semicontinuous,
and $\rho(q_2,\theta_2)\geq s_2$ we can assume that $\theta_2^{\prime}$ and $s_2^{\prime}$
are so close to $\theta_2^{\prime}$ and $s_2^{\prime}$, correspondingly, that
that $\rho(q_2,\theta_2^{\prime})>s_2^{\prime}$.

Let us define 
 $ \zeta^{\prime}=-\gamma_{q_2,\theta_2^{\prime}}(s_2^{\prime})\in L^-_{y^{\prime}}M$.
Then $(y^{\prime},-\zeta^{\prime})\in Z=\O_U(q_2)\cap TV$, and $q_2= \gamma_{y^{\prime},\zeta^{\prime}}(s_2^{\prime})$.
As  then $(y^{\prime},-\zeta^{\prime})\in Z=\O_U(q_1)\cap TV$, too, we see that there is $s_1^{\prime}\geq 0$
such that
\ba
q_1=\gamma_{y^{\prime},\zeta^{\prime}}(s_1^{\prime}).
\ea 
Then $s_1^{\prime}<s_2^{\prime}$, as otherwise the union of 
 $\gamma_{q_2,\theta_2}([0,s_2-s_1])$ and  $\gamma_{y^{\prime},\zeta^{\prime}}([s_2^{\prime},s_1^{\prime}])$,
 oriented in the opposite direction,
 would be a closed causal path.

Let us consider a path which is the union of  the  light-like geodesics
$\gamma_{q_2,\theta_2^{\prime}}([0,s_2^{\prime}-s_1^{\prime}])$  and $\gamma_{q_1,\theta_1}([0,s_1])$.
As $\theta^{\prime}_2\not =\theta_2$, this path is not a light-like geodesics and
using a short-cut argument we see that 
$\tau(q_2,y)>0$. Hence $f_\mu^+(q_2)<f_\mu^+(q_1)$ which
is a contradiction with the fact that $\E_U(q_1)\cap TV=\E_U(q_2)\cap TV$,
proving that $q_1=q_2$.

Finally, consider the case when  $q_1 \in \mu$ and $\E_U(q_1)\cap TV=\E_U(q_2)\cap TV$.
Let $(q_1,\xi)\in \E_U(q_1)$. 
Then,
 if $q_2\not=q_1$, we have $(q_1,\xi)\in \E_U(q_2)$ and thus there is $s>0$
 so that $q_2=\gamma_{q_1,\xi}(-s)$. When $\tilde s>0$ is sufficiently
 small, we see that $(\tilde x,\tilde \xi)=(\gamma_{q_1,\xi}(-\tilde s),\dot \gamma_{q_1,\xi}(-\tilde s)) 
 \in \E_U(q_2)\cap TV$ but $\gamma_{q_1,\xi}(-\tilde s)<q_1$ and
 thus $(\tilde x,\tilde \xi)\not \in  \E_U(q_1)\cap TV$ that is is a contradiction.
 Hence $q_1=q_2$.
%
%
%

The last claim follows from the fact that if $\O_U(q_1)=\O_U(q_2)$
then $\E_U(q_1)\cap TV_1=\E_U(q_2)\cap TV_1$ for some neighborhood
$V_1\subset U$ of  $\pointear_{z,\eta}(\O_U(q_1))$.
\hfill \Box \medskip

\HOX{Find  a direct reference for  Lemma \ref{A: lem 1}  (ii) and Lemma \ref{lemma compact geodesics}. The proof is probably standard and should be omitted in
the final paper}

\begin{lemma} \label{lemma compact geodesics} 
Let $K\subset M$ be a compact set. Then there is $R_1>0$ such that
if $\gamma_{y,\theta}([0,l])\subset K$ is a light-like geodesic
with $\|\theta\|_{g^+}=1$, then $l\leq R_1$. In the case when $K=J(p^-,p^+)$,
with $q^-,q^+\in M$ we have $\gamma_{y,\theta}(t)\not \in J(q^-,q^+)$ for $t>R_1$.
\end{lemma}

The proof of this lemma is standard, but we include it for the convenience of the reader.

\noindent
{\bf Proof.}
Assume that there are no such $R_1$.
Then
there are geodesics $\gamma_{y_j,\theta_j}([0,l_j])\subset K$, $j\in \Z_+$
such that $\|\theta_j\|_{g^+}=1$ and $l_j\to \infty$ as $j\to \infty$. 
Let us choose a subsequence $(y_j,\theta_j)$ which
converges to some point $(y,\theta)$ in $(TM,g^+)$. 
As $\theta_j$ are light-like, also $\theta$ is light-like.

Then, we observe that for all $R_0>0$ the functions $t\mapsto \gamma_{y_j,\theta_j}(s),$
converge in $C^1([0,R_0];M)$ to $s\mapsto \gamma_{y,\theta}(t),$ as $j\to \infty$. 
As $\gamma_{y_j,\theta_j}([0,l_j])\subset K$ for all $j$,
we see that $\gamma_{y,\theta}([0,R_0])\subset K$ for all $R_0>0$.
Let $z_n=\gamma_{y,\theta}(n)$, $n\in \Z_+$. As $K$ is compact,
we see that there is a subsequence $z_{n_k}$ which converges
to a point $z$ as $n_k\to \infty$. Let now $V\subset M$ be a small convex neighborhood
of $z$ such that each geodesic starting from $V$ exits the set $V$ (cf.\ Lemma \ref{A: lem 1}). 
Let $V^{\prime}\subset M$ be a neighborhood
of $z$  so that the strong causality condition (\ref{strong claim})
is satisfied for $V$ and $V^{\prime}$. Then we see that there is $k_0$ such that
if $k\geq k_0$ then $z_{n_k}\in V^{\prime}$, implying that  
$\gamma_{y,\theta}([n_{k_0},\infty))\subset V$. 
This is a contradiction and thus the claimed $R_1>0$ exists.

Finally, in the case when $K=J(q^-,q^+)$, $q^-,q^+\in M$ we see that
if $q(s)=\gamma_{y,\theta}(s) \in K$ for some $s>R_1$, then for all $\tilde s\in [0,s]$
we have 
$q(\tilde s)\leq q(s) \leq q^+$ and $q_-\leq q(0)\leq q(\tilde s)$. Thus 
$q(\tilde s)\in K$ for all $\tilde s\in [0,s]$. As $t>R_1$, this is not possible by
the above reasoning, and thus the last assertion follows.
\hfill \Box \medskip

Assume next that $y_j=\gamma_{q,\eta_j}(s_j)$
and $\zeta_j=\dot \gamma_{q,\eta_j}(s_j)$, where
$q\in M$, $\eta_j\in L^+_qM$, $\|\eta_j\|_{g^+}=1$, and $t_j>0$, are such that  $(y_j,\zeta_j)\to (y,\zeta)$ in $TM$ 
as $j\to\infty$. There exists $p\in M$ such that $p\in I^+(y)$. Then  for sufficiently large $j$ 
we have $y_j\in J(q,p)$ and we see  that by Lemma \ref{lemma compact geodesics},
 $s_j$ are uniformly bounded. Thus
there exist  subsequences $s_{j_k}$ and $\eta_{j_k}$  satisfying $s_{j_k}\to t$ and $\eta_{j_k}\to \eta$ as $k\to \infty$.
Then $(y,\zeta)=(\gamma_{q,\eta}(s),\dot \gamma_{q,\eta}(s))$. This shows that the light observation
set $\O_U(q)$ is a closed subset $TU$. 


{\tobecheckedtext
 %
%
 %
%
 \hiddenfootnote{
Alternatively to the presented construction, we consider the collection
  $\mathcal C_{\mu({z,\eta})}$ of 
relatively closed sets $K\subset TU_{\hat g}$ intersecting $\pi^{-1}(\mu({z,\eta}))$. When 
$K\subset TU_{g}$ is relatively closed we say
that the earliest intersection point of $K$ and $\pi^{-1}(\mu({z,\eta}))$,
denote by $\pointear_{\mu({z,\eta})}(K)$, is the point $z\in \pi(K)\cap \mu({z,\eta})$ for which we have
$z\leq z^{\prime}$ for all $z^{\prime}\in \pi(K)\cap\mu({z,\eta})$. 
Next we consider the germs of the closed
sets $K$ intersecting $\mu({z,\eta})$ near the earliest intersection point.
When $K\in \mathcal C_{\mu({z,\eta})}$ we define 
 $H_K:\R_+\to \mathcal C_{\mu({z,\eta})}$ to be
 $H_K(r)=K\cap TB_{g^+}(\pointear_{\mu({z,\eta})}(K),r)$.
We denote the map $H=H_K$ below by
$(H(r))_{r>0}$ and say that two such the maps $H=(H(r))_{r>0}$ and $\tilde H=(\tilde H(r))_{r>0}$
are equivalent and denote $H\sim \tilde H$, if there is $r_0>0$
so that $H(r)=\tilde H(r)$ for all $0<r<r_1$.
We call these equivalence classes
the germs of the sets in    $\mathcal C_{\mu({z,\eta})}$
and denote their set by  $ \mathcal G_{z,\eta}$.
We also consider maps $\mathcal H:\U_{z_0,\eta_0}\to \bigcup_{(z,\eta)\in \U_{z_0,\eta_0}}  \mathcal G_{z,\eta}$
having the property that $\mathcal H(z,\eta)\in \mathcal G_{z,\eta}$
and denote their set by $\mathcal G$.

%


We will consider consider  the map
$\O^\sim: J(p^-,p^+)\to  \mathcal G$, 
 given by
\ba
& &O^\sim(z,\eta,q)=O^\sim_{z,\eta}(q),\\
 & &(\O^{\sim}_{z,\eta}(q))(r)=\O^{r}_{z,\eta}(q):=\O_U(q)\cap T
 B_{g^+}(\pointear_{\mu({z,\eta})} ( \O_U(q)),r).
 \ea
We note that it is difficult to define topology on $ \mathcal G$
that would be such that the map $\O^\sim:\U_{z_0,\eta_0}\times J(p^-,p^+)\to \mathcal G$,
would be continuous and that the topology on $ \mathcal G$
would have the Hausdorff property. Because of this we consider in the main text the earliest observations sets.}
%
%
}

{\tobecheckedtext 
Let $\mathbb S$ be the collection of relatively closed
sets $K$ in $U_{g}$ that intersect all geodesics
$\gamma_{z,\eta},$ $(z,\eta)\in \U_{z_0,\eta_0}$ precisely once. 
We endow $\mathbb S$  with the topology $\tau_e$ 
that is the weakest topology for which all maps $\pointear_{z,\eta}$, $K\mapsto 
\pointear_{z,\eta}(K)$, parametrized
by $(z,\eta)\in \U_{z_0,\eta_0}$, are continuous. 
Note that $ \P_U(q)\in 
\mathbb S$ for all $q\in I^-(p^+)\setminus  J^-(p^-)$ by
 formula (\ref{kaava D}).  We use below the  continuous map 
 $F:\mathbb S\to \prod_ {(z,\eta)\in \U_{z_0,\eta_0}} \R=\R^{\U_{z_0,\eta_0}}$ 
 that is defined
 for $Z\in \mathbb S$ by setting
$F(Z)=(s_{z,\eta})_{(z,\eta)\in \U_{z_0,\eta_0}}$ where $\mu_{z,\eta}(s_{z,\eta})\in \pi (Z)$.
}

\begin{lemma} \label{lemma homeo} 
Let $V\subset I^-(p^+)\setminus J^-(p^-)$ be relatively compact
open set 
Then the map 
\ba
 \be_U:\overline V\to \mathbb S
 \ea
defines a homeomorphism $\be_{U}:\overline V\to \be_{U}(\overline V)$.
\end{lemma}
\noindent
{\bf Proof.} Let $\mu=\mu_{g,z,\eta}$.
{\tobecheckedtext 
First we note that as the map $x\mapsto f^+_\mu(x)$ is continuous
in $I^-(p^+)\setminus J^-(p^-)$ for all $\mu=\mu(z,\eta)$, the map $\be_{U}$ is continuous. 

Let us next show that the relative topology on $\be_{U}(J(p^-,p^+))$ determined
by $(\mathbb S,\tau_e)$ 
is a Hausdorff topology.

If $(\be_{U}(J(p^-,p^+)),\tau_e)$ is not Hausdorff space 
then there are  $q_1,q_2\in J(p^-,p^+)$, $q_1\not =q_2$ such that $f^+_{\mu(z,\eta)}(q_1)=  
 f^+_{\mu(z,\eta)}(q_2)$ for all $(z,\eta)\in \U_{z_0,\eta_0}$
 and there is $(x,\zeta_1)\in 
 \E_{U}(q_1)\setminus  \E_{U}(q_2)$.
 By Lemma \ref{lemma: global injectivity} , 
 we have either that the set
 $ \E_{U}(q_1)\setminus  \E_{U}(q_2)$
 or  $ \E_{U}(q_1)\setminus  \E_{U}(q_2)$ is non-empty.
 Next, consider the case where  $ \E_{U}(q_1)\setminus  \E_{U}(q_2)\not=\emptyset$.

 Let $t>0$ be so small that $\gamma_{x_1,\zeta_1}([-t,0])\subset U_{g}$ and
 $\gamma_{x_1,\zeta_1}([-t,0])$  is the only
 light-like geodesic from $x_1$ to $x_1^{\prime}=\gamma_{x_1,\zeta_1}(-t)$,
 see Lemma \ref{A: lem 1}  (i). Let
 $\zeta^{\prime}_1=\dot \gamma_{x_1,\zeta_1}(-t)$.
 Then $(x^{\prime}_1,\zeta^{\prime}_1)\in 
 \E_U(q_1)$ and $\tau(x_1^{\prime},x_1)=0$.
  Let
 $\mu_k(s)=\mu(z_k,\eta_k;s)$, $k=1,2$ be such  that 
 $x_1=\mu_1(s_1)$ and $x_1^{\prime}=\mu_2(s_2)$.
 Then $f^+_{\mu_k}(q_2)=f^+_{\mu_k}(q_1)=s_k$ for
  $k=1,2$  implies  that $\tau(q_2,x_1)=\tau(q_2,x_1^{\prime})=0$.
  As  $\tau(x_1^{\prime},x_1)=0$, this implies by a short cut argument 
  that the union 
  of a light-like geodesic $\tilde \gamma$ from $q_2$ to $x_1^{\prime}$
  and $\gamma_{x_1,\zeta_1}([-t,0])$ from $x_1^{\prime}$ to $x_1$
is a light-like
  geodesic and there is $s_2$ such that $\gamma_{x_1,\zeta_1}(-s_2)=q_2$.
   As $\tau(q_2,x_1)=0$ we see that $\gamma_{x_1,\zeta_1}([-s_2,0])$
   is a longest possible curve between its end points
   and thus  we have to have $(x_1,\zeta_1)\in \E_U(q_2)$.
   This contradiction shows that $\tau_e$
   induces on $\E_U(J(p^-,p^+))$
  a Hausdorff topology.

%
  The claim follows then  from the  general fact that a continuous bijective map
 from a compact Hausdorff space onto a Hausdorff space is a homeomorphism.}
%
%
%
\hfill \Box \medskip

In the next lemma we consider coordinates associated with light observations, see Fig.\ 5. 

\begin{lemma} \label{lemma coordinates} Let $q_0\in  I^-(p^+)\setminus J^-(p^-) $.
Then there is a neighborhood $W\subset M$ of $q_0$
and  time-like paths $\mu_{z_j,\eta_j}( (-1,1))\subset U_g$,
$({z_j,\eta_j})\in \U_{z_0,\eta_0}$, $j=1,2,\dots,{n}$, such that if $Y^j(q)\in \R$
are determined by the equations
\ba
  \mu_j(Y^j(q))=\pointear_{z_j,\eta_j}(\O_U(q)),\quad
  \hbox{or equivalently,}\quad
  Y^j(q)  =f^+_{\mu_j}(q),
\ea  
then $Y(q)=(Y^j(q))_{j=1}^{ n}$ define coordinates $Y:W\to \R^{ n}$
which are compatible with the differentiable structure of $M$.
Moreover, if  $(\tilde z,\tilde \eta)\in \U_{z_0,\eta_0}$ is given,
the points $({z_j,\eta_j})\in \U_{z_0,\eta_0}$, $j=1,2,\dots,{n}$
can be chosen in an arbitrary open neighborhood of $(\tilde z,\tilde \eta)$.

\end{lemma}
%

\begin{figure}[htbp] \label{Fig-5}
\begin{center}

\psfrag{1}{$q$}
\psfrag{2}{$x$}
\includegraphics[width=5.5cm]{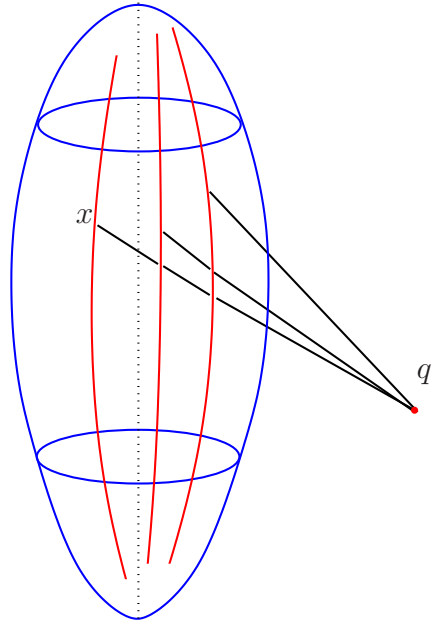}
\end{center}
\caption{A schematic figure where the space-time is represented as the  3-dimensional set $\R^{2+1}$.
The red curves are paths $\mu_{z_j,\eta_j}=\mu(z_j,\eta_j)$. In a neighborhood of a point $q_0\in
I^-(\hat p^+)\setminus I^-(\hat p^-)$ we can choose $(z_j,\eta_j)\in \mathcal U_{z_0,\eta_0}$,
$j=1,2,3,4$ such that $q\mapsto ( f^+_{\mu(z_j,\eta_j)}(q))_{j=1}^4$ defines
coordinates near $q_0$.  In the figure, light lines are light-like geodesics
and $x=\mu_{z_j,\eta_j}(  f^+_{\mu(z_j,\eta_j)}(q))$.}
 \end{figure}

\noindent
{\bf Proof.}
Let $q_0\in  I^-(p^+)\setminus J^-(p^-)$,  $(y_0,\xi_0)\in \ear_{\mu}(\O_U(q_0))$, 
and
$\gamma_{y_0,-\xi_0}([0,t_0])$ be a light-like geodesic
from $y_0$ to $q_0$. Let  $\vartheta_0=-\dot\gamma_{y_0,-\xi_0}(t_0)$.
As $\gamma_{q_0,\vartheta_0}([0,t_0])$ is a longest causal curve between
its end points (note that it needs  to be unique),
 for any $t_1\in (0,t_0)$ (see  \cite[Thm.\ 10.51]{ONeill} or 
\cite[Prop 4.5.12]{HE}) there are no conjugate points on the geodesic
$\gamma_{q_0,\vartheta_0}([0,t_1])$. 
Let $\V_r= \V_{r}(q_0,t_1\vartheta_0)\subset TM$ be the $r$-neighborhood of 
$(q_0,t_1\vartheta_0)$ in the Sasaki metric of $(TM,g^+)$.
We see by using \cite[Prop.\ 10.10]{ONeill} that, for all $0<t_1<t_0$, 
there is $r_1=r_1(t_1)>0$ such that 
the exponential map $\Phi:(x,\zeta)\mapsto (x,\gamma_{x,\zeta}(1))$
is a diffeomorphism $\Phi:\V_{r_1}\to \Phi(\V_{r_1})\subset M\times M$. 

In the following, let  $t_1<t_0$ be so close to $t_0$ and $r_2\in (0,r_1(t_1))$ be so small that for 
all $(x,\xi)\in  \V_{r_2}$ we have $\gamma_{x,\xi}(1)\in U_g$.

Next we show that there  is $r_3\in (0,r_2)$ such that, if
$(x,\xi)\in L^+M\cap \V_{r_3}$, then 
the geodesic $\gamma_{x,\xi}([0,1])$
is the unique causal curve between its endpoints $x$ and $y=\gamma_{x,\xi}(1)\in U$.
To show this, assume that there are no such $r_3$.
Then 
there is a sequence $r_j\to 0$, points $(x_j,\xi_j)\in \V_{r_j}\cap L^+M$, and
$\zeta_j\in T_{x_j}M$, $\zeta_j\not=\xi_j$
 so that  $\gamma_{x_j,\zeta_j}([0,1])$ is some other causal geodesic
between the points $\gamma_{x_j,\xi_j}(0)$ and $\gamma_{x_j,\xi_j}(1)$
having at least the same length as  $\gamma_{x_j,\xi_j}([0,1])$.
Note that $\gamma_{x_j,\xi_j}(0)\to q_0$ and $\gamma_{x_j,\xi_j}(1)\to
\gamma_{q_0,\vartheta_0}(t_1)$ as $j\to \infty$.
Then, the sequence $(x_j,\zeta_j)$ has
a subsequence which converges to some point $(q_0,\zeta)$ in the Sasaki metric,
see Lemma \ref{lemma compact geodesics}. 
If $\zeta=t_1\vartheta_0$, we see that there the map $\Phi$ is
not a local diffeomorphism  near the point $(q_0,t_1\vartheta_0)$
which is not possible. On the other
hand, if  $\zeta\not=t_1\vartheta_0$, we see using a short cut
argument for the union of the geodesics $\gamma_{q_0,\zeta}([0,1])$ and
$\gamma_{q_0,\vartheta_0}([t_1,t_0])$, that the geodesic
$\gamma_{q_0,\vartheta_0}([0,t_0])$ is not a longest causal curve
between its end points $q_0$ and $y_0$. This is in contradiction with the 
earlier assumptions, and thus we see that the claimed $r_3>0$
exists. 

For $(x,y)\in \Phi(\V_{r_4})$, $r_4=r_3/2$  we denote $(x,\zeta_{x,y})=\Phi^{-1}(x,y)$.
If $\zeta=\zeta_{x,y}\in L^+M$, then the light-like geodesic $\gamma_{x,\zeta}([0,1])$
 is the unique causal curve 
between $x$ and $y$. Note that for $\e>0$ small enough, e.g., $\e=r_4$,  also the geodesic
$\gamma_{x,\zeta}([0,1+\e])$ with endpoints in $\Phi(\V_{r_3})$ is the unique causal curve between
its end points and then by  \cite[Thm.\ 10.51]{ONeill} or 
\cite[Prop 4.5.12]{HE}) there are no conjugate points on the geodesic
$\gamma_{x,\zeta}([0,1])$.

Let us next fix  $t_1<t_0$ and $r_3$  as above,
and consider $t_1\vartheta_1\in L^+_{q_0}M\cap \V_{r_3}$, $\|\vartheta_1\|_{g^+}=\|\vartheta_0\|_{g^+}$,
$\vartheta_1\not=\vartheta_0$, and
$y_1=\gamma_{q_0,\vartheta_1}(t_1)\in U_g$. Let $(x,\eta)\in \U_{z_0,\eta_0}$
be such the $y_1\in \mu_{z,\eta}$. Then  $\dot \mu_{z,\eta}$ and $\theta_1=\dot \gamma_{q_0,\vartheta_1}(t_1)$ are both 
light-like vectors in $T_{y_1}M$ and hence $g(\dot\mu_{z,\eta}(0),\theta_1)\not=0$.

Let us recall that the energy of a piecewise smooth curve $\alpha:[0,l]\to M$
is defined by
\ba
E(\a)=\frac 12\int_0^l g(\dot \alpha(t),\dot \alpha(t))\,dt.
\ea
Then the energy of the geodesic $\gamma_{x,\zeta}([0,1])$,
$\zeta=\zeta_{x,y}$ connecting the
points $x$ and $y$ is equal to $g(\zeta_{x,y},\zeta_{x,y})$,
and thus is a $C^\infty$ smooth function in
$(x,y)\in \Phi(\V_r)$.

Next, let $X\in T_{q_0}M$, and $z(s)=\gamma_{q_0,X}(s)$, $s\in [-s_0,s_0]$, $s_0>0$. Let $T:\R\to \R$ be a smooth function
so that $T(0)=f^+_{\mu_{z,\eta}}(q_0)$.
Let $\gamma_s(t)=\gamma_{z(s),\zeta(s)}(t)$,
where $(z(s),\zeta(s))=\Phi^{-1}(z(s),\mu_{z,\eta}(T(s)))$. 
Denote also $\theta(s)=\dot \gamma_{z(s),\zeta(s)}(1)\in T_{\mu_{z,\eta}(T(s))}M$.

Using \cite[Prop.\ 10.39]{ONeill} for the  
first variation of the energy of the geodesics $\gamma_s(t)$,
we see that
\ba
\frac d{ds}{ E}(\gamma_s)
=g(\theta(s),\dot \mu_{z,\eta}(T(s)))T^{\prime}(s)-g(\dot \gamma_{q_0,X}(s),\zeta(s)).
\ea
Thus, $ {\mathcal E}(\gamma_s)=0$ for all $s\in [-s_0,s_0]$ if and only if
\beq\label{ODY1}
T^{\prime}(s)&=&\frac{g(\dot\gamma_{q_0,X}(s),\zeta(s))}{g(\theta(s),\dot \mu_{z,\eta}(T(s)))}\\
&=&\left.\frac{g(\dot\gamma_{q_0,X}(s),\zeta(s))}{g(\dot \gamma_{z(s),\zeta(s)}(1),\dot \mu_{z,\eta}(T(s)))}
\right|_{ (z(s),\zeta(s))=\Phi^{-1}(z(s),\mu_{z,\eta}(T(s)))}\nonumber
\eeq
for  all $s\in [-s_0,s_0]$. Now the  equation (\ref{ODY1})
can be considered as
an ordinary differential
equation for $T(s)$, which has a unique solution with  initial data $T(0)=f^+_{\mu_{z,\eta}}(q_0)$
on the interval  $s\in [-s_0,s_0]$ when $s_0>0$ is sufficiently small.
Let us denote the solution by $T_1(s)$.
Using in the above $T(s)=T_1(s)$, we see that $ {E}(\gamma_s)=0$, implying 
$\zeta(s)\in L^+_{z(s)}M$. Thus, when $|s|$ is small
enough $(z(s),\mu_{z,\eta}(T(s)))\in  \Phi(\V_{r_3})$ and the light-like  geodesic
 $\gamma_{z(s),\zeta(s)}([0,1])$ is a longest geodesic
between   $z(s)$ and $ \mu_{z,\eta}(T(s))$.
Hence, for $s_0>0$ small enough  $T=T_1(s)$, $s\in [-s_0,s_0]$ is
the solution of the equation
\ba
  \mu_{z,\eta}(T)=\pointear_{{z,\eta}}(\O_U(z(s))),
\ea  
that is, $T_1(s)=f^+_{\mu_{z,\eta}}(z(s))$.
In particular, we see that for $s=0$,
\ba
0=
\frac d{ds}{ E}(\gamma_s)|_{s=0}=g(\theta_1,\dot \mu_{z,\eta}(0))T_1^{\prime}(0)-g(X,\vartheta_1),
\ \ \hbox{i.e.}\
T_1^{\prime}(0)=\frac{g(X,\vartheta_1)}{g(\theta_1,\dot \mu_{z,\eta}(0))}.
\ea
Let us now denote $Y^1(z(s)):=T_1(s)$. Above $T^\prime_1(0)=
g(\hbox{grad}_g\, Y^1|_{q_0},X)$ where $X=\dot z(0)$.
Thus, by varying the vector $X\in T_{q_0}M$ in
the above construction, we see that $q_0$ has a neighborhood $W$ in $M$, 
$Y^1:W\to \R$ is well defined function and
\ba
\hbox{grad}_g\, Y^1|_{q_0}&:=&g^{jk}\frac {\p Y^1(x)}{\p x^j} \frac \p{\p x^k}=c_1\vartheta_1,\\
c_1&=&\frac 1{g(\theta_1,\dot \mu_{z,\eta}(0))},
\quad \vartheta_1=\vartheta_1^j\frac \p{\p x^j}.
\ea
Now, choose such $\vartheta_j\in L^+_{q_0}M$, 
 $\|\vartheta_j\|_{g^+}=\|\vartheta_0\|_{g^+}$ and $t_j<t_0$  so close to $t_0$,
  $j=2,3,\dots,{ n}$,  such that $(q_0,t_j\vartheta_j)\in \V_{r_3}$ and 
  $\vartheta_1,\vartheta_2,\dots,\vartheta_{ n}$ are linearly independent.
Then    $y_j=\gamma_{q_0,\vartheta_j}(t_j)\in U_g$
and there are no conjugate points on the geodesic
$\gamma_{q_0,\vartheta_j}([0,t_j+\e])$ with some $\e>0$.
Let $(z_j,\eta_j)\in \U_{z_0,\eta_0}$, $j=2,3,\dots,n$ and $\mu_j=\mu_{g,z_j,\eta_j}$
be a smooth timelike curve going through the point $y_j=\mu_j(0)$.
Again, $g(\dot\mu_j(f^+_{\mu_j}(q_0)),\theta_j)\not=0,$
$\theta_j=\dot \gamma_{q_0,\vartheta_j}(t_j)$.
Finally, let
 $Y^j(z)$ be such that
\ba
  \mu_j(Y^j(z))=\pointear_{z_j,\eta_j}(\O_U(z)),\quad j=2,3,\dots { n}.
\ea  
Since the vectors $\hbox{grad}_gY_j|_{q_0}=c_j\vartheta_j$, $j=1,2,\dots,{ n}$,
$c_j\not=0$ are
linearly independent, the function $Y:W\to \R^{ n},$
$Y(z)=(Y^1(z), Y^2(z),\dots,Y^{ n}(z))$ gives coordinates near $q_0$.

\HOX{
We need to check if the last observation is not needed later.}
Finally, we note that as 
$(y_0,\xi_0)$ was in above an arbitrary point of $\ear_{\mu}(\O_U(q_0))$,
we see easily that in the above construction 
$({z_j,\eta_j})\in \U_{z_0,\eta_0}$, $j=1,2,\dots,{n}$,
can be chosen in an arbitrary open neighborhood of $(\tilde z,\tilde \eta)$.
\hfill \Box \medskip

Let $W\subset V$ be a neighborhood of $q_0$ and $Y:W\to \R^n$ be local coordinates 
 considered in Lemma \ref{lemma coordinates}.
Let us 
consider the metric ${\bf g}=(\E_U)_*g$ which makes $\E_U:(W,g)\to (\E_U(W),{\bf g})$ an  isometry.
Next we show that we can determine
the conformal class of ${\bf g}$.

\HOX{We could mention the double foliations and the incidence relation}Let $(z,\zeta)\in TM$ and define
\ba
\Gamma(z,\zeta)= \{\E_U(q^{\prime})\in \E_U(W);\ (z,\zeta)\in \E_U(q^{\prime})\}.
\ea
Then
$
\Gamma(z,\zeta)=\E_U(\gamma_{z,\zeta}((-\infty,0])\cap W),
$
that is, $\Gamma(z,\zeta)$ is the image of a geodesic 
$\gamma_{z,\zeta}((-\infty,0])$ on $\E_U(W)$. 
Thus  when $\E_U(W)$ is given, we can find
 in the
$Y$-coordinates
 all light-like geodesics in $\E_U(W)$ that go through a
specified point $\E_U(q)$. Consider the tangent space $T_Q(\E_U(W))$ of 
the manifold $\E_U(W)$ at $Q=\E_U(q)$.
Then  $\E_U(W)$  determines
an open subset of the light cone (with respect to the metric ${\bf g}$) in $T_Q(\E_U(W))$.
 As the light cone is a quadratic surface in $T_Q(\E_U(W))$, we can determine
the whole light cone in $T_Q(\E_U(W))$ in  local coordinates.
 Thus we can determine all light-like vectors
in the tangent space at the point $\E_U(q)$, where  $q\in W$ is arbitrary.
By \cite[Thm.\ 2.3]{Beem}, these collections of  light-like vectors determine uniquely the conformal
class of the tensor ${\bf g}$ in $\E_U(W)$.
This proves Theorem \ref{main thm}.\hfill  $\square$

\subsection{Geometric preparations for analytic results}\label{subsec; analytic}

Next we prove some auxiliary geometrical results needed later in the analysis
Einstein equations.

As the modified
cut locus function $\rho(x,\xi)$ on $(M,{g})$ is lower semi-continuous,
there is $\rho_0>0$, depending on $(M,g)$ and $p^+,p^-$, such that $\rho(x,\xi)\geq \rho_0$ for all $x\in J_{g}(p^-,p^+)$
and $\xi\in L^+_x(M,{ g})$ with $\|\xi\|_{g^+}=1$.
Below, let $R_1>0$ be  given in Lemma \ref {lemma compact geodesics}.  On the setting
of the lemma, see Fig. 6. 


 \begin{lemma} \label{lem: detect conjugate 0} 
 Denote  $\mu_{g,z_0,\eta_0}=\mu$. 
  There are 
$\vartheta_0,\kappa_0,\kappa_1,\kappa_2,\kappa_3\in (0,\rho_0)$ such that 
for  all $y\in \mu([s_-,s_+])$, $y=\mu(r_1)$ with $r_1\in [s_-,s_+]$,
$\zeta\in L^+_yM$, $\|\zeta\|_{g^+}=1$  
and 
$(x,\xi)\in L^+M$ satisfying $d_{g^+}((y,\zeta),(x,\xi))\leq \vartheta_0$
we have  $x\in J(\mu(-1),\mu(+1))$ and the following holds:

(i)
If $t>R_1$ and $t<\T(x,\xi)$, then  $\gamma_{x,\xi}(t)\not \in J^-(\mu(s_{+2}))$.
Moreover, if $0<t\leq R_1$ and  $\gamma_{y,\zeta}(t) \in J^-(\mu(s_{+2}))$,
then $\gamma_{x,\xi}([0,t+\kappa_0])\subset I^-(\mu_g(s_{+3}))$.
Finally,  if $0<t\leq 10\kappa_1$, then $\gamma_{x,\xi}(t)\in U_g$.

(ii) Let $t_0\in [\kappa_1,6\kappa_1]$ and $t_2(x,\xi,t_0)
:=\rho(\gamma_{x,\xi}(t_0),\dot \gamma_{x,\xi}(t_0))\in (0,\infty ]$.
If $\rho(y,\zeta)\geq R_1+\kappa_0$, then $t_2(x,\xi,t_0)+t_0>R_1$.
If $\rho(y,\zeta)\leq R_1+\kappa_0$, 
then $t_2(x,\xi,t_0)+t_0\geq \rho(y,\zeta)+3\kappa_2$.

(iii) 
Let $t_1=\rho(y,\zeta).$  Assume that $t_2\in [0,\T(y,\zeta))$ is such that 
$t_2-t_1\geq \kappa_2$ and $p_2=\gamma_{y,\zeta}(t_2)\in J^-(\mu(s_{+2}))$. 
Then     
 $r_2=f^-_\mu(p_2)$ satisfies $r_2-r_1>3\kappa_3$.

Moreover, let $p_3=\gamma_{x,\xi}(t_2)$. Then either 
$p_3\not \in J^-(\mu(s_{+1}))$ or $r_3=f^-_\mu(p_3)$ satisfies $r_3-r_1>2\kappa_3$.


%

(iv) Above, $\kappa_1>0$ can
be chosen so that  \beq\label{t0-s0 assumption}
 & &\gamma_{ \mu(s^\prime),\xi}([0,4\kappa_1])\cap J^+( \mu(s^{\prime\prime}))=\emptyset\\
 & &\nonumber \quad
 \hbox{when }s_-\leq s^\prime<s^{\prime\prime}\leq s_+,
\ \ \xi\in L^+_{\mu(s^\prime)}M_0.
 \eeq
 
\end{lemma}

\begin{figure}[htbp] \label{Fig-6}
\begin{center}

\psfrag{1}{$y$}
\psfrag{2}{$x$}
\psfrag{3}{$z$}
\psfrag{4}{$q_1$}
\psfrag{5}{$p_1$}
\psfrag{6}{$p_2$}
\includegraphics[width=7.5cm]{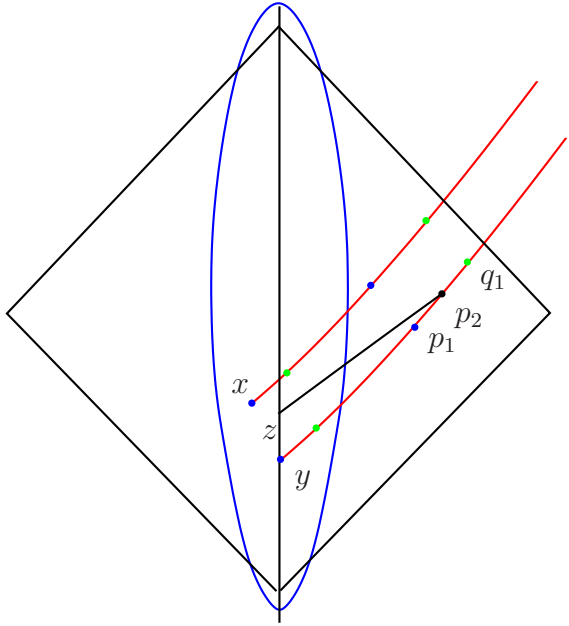}
\end{center}
\caption{A schematic figure where the space-time is represented as the  2-dimensional set $\R^{1+1}$.
The  figure shows the situation in Lemma \ref{lem: detect conjugate 0}. The point $y=\hat \mu(r_1)$ is on the time-like
geodesic $ \hat \mu$ shown as a black line. 
The black diamond is the set $J_{\hat g}(\hat p^-,\hat p^+)$.
The point $y$
is the starting point
of the light-like geodesic $\gamma_{y,\xi}$ and $(x,\zeta)$ is a light-like direction close to $(y,\xi)$. The geodesics  $\gamma_{y,\xi}$ and  $\gamma_{x,\zeta}$ are shown as
red curves and the blue points on them are the first cut-point on
 $\gamma_{y,\xi}([0,\infty))$, that is,
 $p_1=\gamma_{y,\xi}(t_1)$, $t_1=\rho(y,\xi)$ and the first cut point on $\gamma_{x,\zeta}([0,\infty))$.
The green points are $\gamma_{y,\xi}(t_0)$ and  $\gamma_{x,\zeta}(t_0)$
and the first cut point $q_1=\gamma_{y,\xi}(t_c)$  on
  $\gamma_{y,\xi}([t_0,\infty))$ and  the first cut point on $\gamma_{x,\zeta}([t_0,\infty))$.
In Lemma 2.10, $t_0+t_c\geq t_1+2\kappa_2$.
   The black   point  $p_2=\gamma_{y,\xi}(t_2)$ is such that  $t_2\geq t_1+\kappa_2$. Note that if
 $z=\hat \mu(r_2)$ is such that $r_2=f^-_{\hat \mu}(p_2)$, then $r_2-r_1>3\kappa_3$.
}
 \end{figure}

{\tobecheckedtext 

{\bf Proof.} Below we denote $\mu=\mu_{g,z_0,\eta_0}$. We start with the claim (iv).

(iv) We observe that if $\kappa_1>0$ is small enough (v) follows 
from the definition of the Fermi type coordinates. Below we assume
that $\kappa_1$ is so small that (iv) is valid.


Next we prove (i)-(iii).

 (i) $\vartheta_0,$ $\kappa_0$, and $\kappa_1$ can be chosen using
 the compactness of the sets $J(p^-,\mu(s_{+j}))$, $j\leq 3$ and $S=\{(y,\zeta)\in L^+M;\ y\in\mu_g([s_-,s_+]),\ \|\zeta\|_{  g^+}=1\}$.

(ii) 
Let us fix $\kappa_0>0$. As
$\rho(x,\xi)$ is lower semi-continuous, we see  that if 
 $\vartheta_0>0$ is chosen first enough,
 the first claim in (ii) holds. 
 
 To show the second claim in (ii),
assume that the opposite holds. Then  there are $(y^n,\zeta^n)\in L^+(\mu( [s_-,s_+]),$ $\|\zeta^n\|_{g^+}=1$, and $t_0^n
\geq \kappa_1$, $n\in \Z_+$ such that 
$\lim_{n\to\infty}(t_2(y^n,\zeta^n,t_0^n)+t_0^n)=\lim_{n\to\infty}\rho(y^n,\zeta^n)$.
As $\rho(x,\zeta)$ is lower semi-continuous,
we see that for some subsequences of
$(y^n,\zeta^n)$ and $t_0^n$  have limits
$(\overline y,\overline \zeta)$ and $\overline t_0$ 
such that $t_2(\overline y,\overline\zeta,\overline t_0)+\overline t_0=\rho(\overline y,\overline \zeta)$.
This is not possible and hence (ii) is proven.

(iii) Denote $p_{+2}=\mu(s_{+2})$. Let $T_{+}(x,\zeta)=\sup \{t\geq 0\ ;\ \gamma_{x,\zeta}(t)\in J^-(p_{+2})\}$
and
\ba
& &
K=\{(r,\xi);\ r\in [s_-,s_+],\ \xi\in L^+_{\mu (r)}\hattuM _0,\ \|\xi\|_{ g^+}=1\},\\ 
& &K_0=\{(r,\xi)\in K; \ \rho(\mu(r),\xi)+\kappa_2\leq T_{+}(\mu(r),\xi)\},
\quad K_1=K\setminus K_0.
\ea 
Then  the map $L: G_0=\{(r,\xi,t)\in K\times \R_+;\ 
t\leq T_{+}(\mu(r),\xi)\}\to \R,$
\ba
L(r,\xi,t)= f^-_\mu(\gamma_{\mu(r),\xi}(t))-r
\ea
is continuous.
\HOX{We could choose better notations}
%
%
%
%
%
%
%
%
Also, we define 
 a map  $H:K\to \R$,
\ba
H(r,\xi) =\begin{cases}
L(r,\xi,\rho(\mu(r),\xi)+\kappa_2),&\hbox{for $(r,\xi)\in K_0$},\\
 \quad \quad \quad 3,
& \hbox{for $(r,\xi)\in K_1$}.
\end{cases}
\ea
%
%
%
%
Note that if  $(r,\xi)\in K_0$ then 
\ba
L(r,\xi,\rho(\mu(r),\xi)+\kappa_2)=f^-_\mu(P(r,\xi))-r\leq 2,
\ea
where $P(r,\xi)=\gamma_{\mu(r),\xi}(\rho(\mu(r),\xi)+\kappa_2)\in J^-(p_{+2})$. 
Above, $\rho(x,\xi)$ is lower semi-continuous
and $T_{+}(x,\xi)$ is upper semi-continuous on the set $L^+(J(p^-,p_{+2}))$
and thus sets $K_0$  and $G_0$ are closed. Moreover, as
$\rho(x,\xi)$ is lower semi-continuous and
the function $t \mapsto L(r,\xi,t)$
is increasing, we see that $H:K\to \R$ is a lower semi-continuous
function.

%
%

For $(r,\xi)\in K_0$, we see that  
$\tau(q,\gamma_{q,\xi}(\rho(q,\xi)+\kappa_2))>0$, where
$q=\mu(r)$. Hence
$H$ is strictly positive.  
If $K_0=\emptyset$,
it is obvious that the claim is valid as the condition $p_2\in J^-(p_{+2})$ never
holds. Assume next  that $K_0\not=\emptyset$.
Then $H$ obtains it minimum $\e_1:= H(r_0,\xi_0)\in (0,3)$
at some point $(r_0,\xi_0)\in K_0$.
As $t \mapsto L(r,\xi,t)$
is increasing, the first claim follows by choosing $\kappa_3<\e_1/3$.

As $f^-_\mu:J(\mu(-1),\mu(1))\to \R$ is 
 continuous and 
$J(\mu(-1),\mu(1))$ is compact, this function is uniformly continuous.
Thus there exists  $\delta_0>0$ such that if $d_{g^+}(x_1,x_2)<\delta$,
then $|f^-_\mu(x_1)-f^-_\mu(x_2)|<\kappa_3$.
%

Let 
\ba
B_0&=&\{(y^\prime ,\zeta^\prime ,t_1^\prime ,x^\prime ,\xi^\prime,t_2^\prime )    ;\ (r^\prime ,\zeta^\prime ,t_1^\prime )\in G_0,\ y^\prime =\mu(r^\prime ),\\
& &\quad  
(x^\prime ,\xi^\prime )\in L^+M,\ d_{g^+}((y^\prime ,\zeta^\prime ),(x^\prime ,\xi^\prime ))\leq \vartheta_0,\ |t_1^\prime-t_2^\prime|\leq \kappa_0\},
\ea
where $G_0$ is the closed set defined above. The set $B_0$ is 
compact and by claim (i) for $(y^\prime ,\zeta^\prime ,t_1^\prime ,x^\prime ,\xi^\prime,t_2^\prime )  \in B_0$ we have $\gamma_{x^\prime ,\xi^\prime}(t_2^\prime)\in  J(\mu(-1),\mu(+1))$. Thus by using Lipschitz continuity of exponential map in
$ T(J(\mu(-1),\mu(+1)))$, we see that there is $L_0>0$ such that
for $(y^\prime ,\zeta^\prime ,t_1^\prime ,x^\prime ,\xi^\prime,t_2^\prime )  \in B_0$ we have 
\beq\label{Lip estimate kiire}
& &|d_{g^+}(\gamma_{y^\prime ,\zeta^\prime}(t_1^\prime),\gamma_{x^\prime ,\xi^\prime}(t_2^\prime))|\leq L_0(d_{g^+}((y^\prime ,\zeta^\prime ),(x^\prime ,\xi^\prime ))+
|t_1^\prime-t_2^\prime|),
\eeq
and moreover, there are $L_1,L_2>0$ such that
\beq \nonumber
& &|{\bf t}(\gamma_{y^\prime ,\zeta^\prime}(t_1^\prime))-{\bf t}(\gamma_{x^\prime ,\xi^\prime}(t_2^\prime))|\leq L_1(d_{g^+}((y^\prime ,\zeta^\prime ),(x^\prime ,\xi^\prime ))+
|t_1^\prime-t_2^\prime|),\\
\label{Lip L1L2 est} & &\frac \p{\p s}{\bf t}(\gamma_{x^\prime ,\xi^\prime}(s))|
_{s=t_1^\prime}\geq L_2
\eeq
where ${\bf t}:M\to \R$ is the smooth time function on $M$ used to introduce 
the identification $M=\R\times N$.

Let $\delta_1=\dist_{g^+}(J(\mu(s_{-1}),\mu(s_{+1}),M\setminus J(\mu(s_{-2}),\mu(s_{+2}))$.
Let us now assume that $\vartheta_0<\min(\delta_0,\delta_1)/L_0$.

Then we see that if $\gamma_{y,\zeta}(t_2)\not \in J^-(\mu(s_{-2}))$,
then $\gamma_{x,\xi}(t_2)\not \in J^-(\mu(s_{-1}))$.
On the other hand, if $p_2=\gamma_{y,\zeta}(t_2) \in J^-(\mu(s_{-2}))$,
then $f^-_\mu(p_2)>r_1+3\kappa_3$ and we see that 
$p_3=\gamma_{x,\xi}(t_2)$ satisfies $f^-_\mu(p_3)>r_1+3\kappa_3-\kappa_3=
r_1+2\kappa_3.$
%
%
 \hfill \Box \medskip}

%
%
%
%

Note that for proving the unique solvability of the inverse problem we 
need to consider  two manifolds, $(M^{(1)},\hat g^{(1)})$ and $(M^{(2)},\hat g^{(2)})$
 with same data. For these manifolds, we can choose $R_1,\vartheta_0,\kappa_j$ so that they are
same for the both manifolds.
\medskip


{\bf Proof.}\  (of Corollary \ref{coro of main thm original Einstein pre})
By our assumptions, $\Psi:(V_1,g^{(1)})\to (V_1,g^{(2)})$
is conformal, the
Ricci curvature of $g^{(j)})$ is zero 
in $W_j$,  and $V_j\subset W_j \cup U_j$, $\Psi(U_1\setminus W_1)=U_2\setminus W_2$.
Moreover, 
  any point $x\in V_1\cap W_1$ can be connected to some point $y\in U_1\cap
W_1$ with a light-like geodesic $\gamma_{x,\xi}([0,l])\subset V_1\cap W_1$.
As $\Psi:V_1\to V_2$  and  $\Psi:U_1\to U_2$  are bijections and
$\Psi(U_1\setminus W_1)=U_2\setminus W_2$, we have
  $\Psi(V_1\cap W_1)=V_2\cap W_2$.
 This implies that any point $\Psi(x)\in V_2\cap W_2$ 
can be connected to the point $\Psi(y)\in U_2\cap
W_2$ with a light-like geodesic $\Psi(\gamma_{x,\xi}([0,l]))\subset V_2\cap W_2$.

By the above there is $f:V_1\to \R$ such that $ g^{(1)}=e^{2f} \Psi^* g^{(2)}$ on $V_1$. To simplify notations we  denote $\hat g= g^{(1)}$ and 
$g= \Psi^*g^{(2)}$. Next consider how function $f:V_1\to \R$ can be constructed
when $\hat g$ and $g$ are given in $U_1$ and   $\hat g_{jk}=e^{2f}g_{jk}$ 
corresponds to the vacuum, i.e., its Ricci curvature vanishes
in a set $W_1$.
 By \cite[formula (2.73)]{Rendall}, the
 Ricci tensors $ \hbox{Ric}_{jk}(g)$ of $g$ and ${\hbox{Ric}}_{jk}(\hat g)$ of $ \hat g$ satisfy
 in  $W_1$ 
 \ba
& &0={\hbox{Ric}}_{jk}(\hat g)={\hbox{Ric}}_{jk}(g)_{jk}-2\nabla^{}_{j}\nabla_{k}^{}f+2(\nabla_{j}^{}f)(\nabla_{k}^{}f)\\
& &\hspace{2cm}- (g^{pq}\nabla^{}_{p}\nabla_{q}^{}f+2g^{pq}(\nabla_{p}^{}f)(\nabla_{q}^{}f))g_{jk}
\ea
where $\nabla=\nabla^g$. For scalar curvature this yield
 \ba
0=e^{2f}\hat g^{jk}{\hbox{Ric}}_{jk}(\hat g)=
g^{jk}{\hbox{Ric}}_{jk}(g)-3g^{jk}\nabla^{}_{j}\nabla_{k}^{}f.
\ea
Combining these, we obtain
\ba
2\nabla^{}_{j}\nabla_{k}^{}f-2(\nabla_{j}^{}f)(\nabla_{k}^{}f)
+2g^{pq}(\nabla_{p}^{}f)(\nabla_{q}^{}f)g_{jk}
={\hbox{Ric}}_{jk}(g)- \frac 13 g^{pq}{\hbox{Ric}}_{pq}(g) g_{jk}.
\ea
%
Let us next use local coordinate neighborhood that cover the geodesic $\gamma(t)=(\gamma^j(t))_{j=1}^4$.
The above  equation gives a system of second order ordinary differential equations
which
 can be solved along the light-like geodesics $\gamma(t)$ in $ W_1$ starting from $U_{\hat g}$.
Indeed, we can take a contraction of the above equations with  $\dot \gamma^j(t)$
and obtain, on each piece of the geodesic that belongs in one coordinate neighborhood, an  ordinary differential equation for $[(\nabla_{k}^{}f)(\gamma(t))]_{k=1}^n$. 
Since  by our assumption any point $x\in V_1\cap W_1$ can be connected to $y\in U_1\cap W_1$
with a light-like geodesic that is a subset of $V_1\cap W_1$, we can
use the above system of the second order ordinary differential equations 
and the fact that we know $f$ in $U_1$ to determine $f(x)$.
This proves the claim.\hfill \Box \medskip



%
%
%
%
%
%

\hiddenfootnote{ 
SLAVA ASKED THAT WE INCLUDE THE FOLLOWING LEMMA TO SIMPLIFY THE PROOF. NOTE THAT THIS LEMMA HAS BEEN INCLUDED IN A PROOF IN THE "COMBINING RESULTS" SECTION.\begin{lemma}\label{lemma:Slava1}
Let  $x_1\in I^-(p^+)\setminus J^-(p^-)$.  Let $\gamma_1([0,t_1])$ and $\gamma_2([0,t_2])$  be two future going light-like geodesic with 
$\gamma_1(0)=\gamma_2(0)=x_1$ and that there are no $c>0$ such that
 $\dot \gamma_1(0)=c\dot \gamma_2(0)$. Assume that 
$x_2=\gamma_1(s_1)=\gamma_2(s_2)\in I^-(p^+)$ with $s_1>0$. 
Let $z_j=\pointear_\mu(x_j)\in \mu$, $j=1,2$. Then either $z_1\ll z_2$ or $x_2=z_1$.
\end {lemma}

\noindent
{\bf Proof.}  As $s_1>0$, the fact that there are no closed causal curves implies that $s_2>0$ and
that $z_1\leq z_2$. As $z_1,z_2\in \mu$ and $\mu$  is time-like curve, we then have either $z_1\ll z_2$
or $z_1=z_2$. If the former one is true, the claim holds. In the case when $z_1=z_2$,
and let us assume, opposite to our claim, that $z_1\not =x_2$. 

Then $\tau(x_1,z_1)=0$ and $\tau(x_2,z_2)=0$. As $z_1=z_2\not =x_2$, this 
implies that there is 
a light-like geodesics, denoted $\gamma_{x_2,\xi}([0,a])$, from $x_2$ to $z_1$.
Moreover, let
$\eta_j$, $j=1,2$ be a curve from $x_1$ to $z_1$ which 
is the union of  $\gamma_j([0,s_j])$
from $x_1$ to $x_2$ and 
the light-like geodesics $\gamma_{x_2,\xi}([0,a])$, from $x_2$ to $z_1$. 
Then either $\eta_1$ or $\eta_2$ is a causal curve  
which is not a light-light geodesic from $x_1$ to $z_1$. Hence by 
\cite[Prop.\ 10.46, 10.51]{ONeill} there is a time-like curve connecting 
 $x_1$ to $z_1$ and thus $\tau(x_1,z_1)>0$. However, this contradicts  the equation $\tau(x_1,z_1)=0$,
 and hence $z_1=z_2\not =x_2$ is not possible. This proves the claim.
\hfill \Box}

\bigskip

\section{Analysis of Einstein equations in wave coordinates}

\subsection{Notations}\label{introduction 1}

Let $X$ be a manifold of dimension $n$
and $\Lambda\subset T^*X\setminus \{0\}$ be
a Lagrangian submanifold.
Let  
$\phi(x,\theta)$, $\theta\in \R^N$ be a non-degenerate
phase function that locally parametrizes $\Lambda$.
We say that a distribution $u\in {\cal D}^{\prime}(X)$ 
is a Lagrangian distribution associated with $\Lambda$ and denote
$u\in \I^m(X;\Lambda)$, if it can locally be represented 
 as
\ba
u(x)=\int_{\R^N}e^{i\phi(x,\theta)}a(x,\theta)\,d\theta,
\ea
where $a(x,\theta)\in S^{m+n/4-N/2}(X; \R^N\setminus 0)$,
see  \cite{GU1,H4,MU1}.

In particular, when $S\subset X$ be a submanifold,  its conormal bundle
$N^*S=\{(x,\xi)\in T^*X\setminus \{0\};\ x\in S,\ \xi\perp T_xS\}$ is a Lagrangian submanifold.

Let us next consider the case when $X=\R^n$ and
let
$(x^1,x^2,\dots,x^n)=(x^{\prime},x^{\prime\prime},x^{\prime\prime\prime})$ 
be the Euclidean coordinates with $x^{\prime}=(x_1,\dots,x_{d_1})$,
  $x^{\prime\prime}=(x_{d_1+1},\dots,x_{d_1+d_2})$,
 $x^{\prime\prime\prime}=(x_{d_1+d_2+1},\dots,x_{n})$. 
 If $S_1=\{x^{\prime}=0\}\subset \R^n$, $\Lambda_1=N^*S_1$
and $u\in \I^m(X;\Lambda_1)$,
then
\ba
u(x)=\int_{\R^{d_1}}e^{i x^{\prime}\cdotp \theta^{\prime}}a(x,\theta^{\prime})\,d\theta^{\prime},
\quad
a(x,\theta^{\prime})\in S^{\mu}(X;\R^{d_1}\setminus 0)
\ea
where $\mu=m-d_1/2+n/4$. 
 
Next we recall the definition of $\I^{p,l}(X;\Lambda_1,\Lambda_2)$, the space of
the distributions $u$ in ${\cal D}^{\prime}(X)$ associated to two
cleanly intersecting Lagrangian manifolds $\Lambda_1,\Lambda_2\subset T^*X
\setminus \{0\}$, see \cite {GU1,MU1}.
Let us start on the case when $X=\R^n$.

Let $S_1,S_2\subset \R^n$ be the linear subspaces of codimensions
$d_1$ and $d_1+d_2$, respectively, and
$S_2\subset S_1$, given by
$S_1=\{x^{\prime}=0\}$, $S_2=\{x^{\prime}=x^{\prime\prime}=0\}$.
Let us denote $\Lambda_1=N^*S_1,$  $\Lambda_2=N^*S_2$. Then
$u\in \I^{p,l}(\R^n;N^*S_1,N^*S_2)$ if and only if
\ba
u(x)=\int_{\R^{d_1+d_2}}e^{i (x^{\prime}\cdotp \theta^{\prime}+
x^{\prime\prime}\cdotp \theta^{\prime\prime})}a(x,\theta^{\prime},\theta^{\prime\prime})\,d\theta^{\prime}d\theta^{\prime\prime},
\ea
where $
a(x,\theta^{\prime},\theta^{\prime\prime})$ belongs to the product type symbol class
$S^{\mu^{\prime},\mu^{\prime\prime}}(\R^n; (\R^{d_1}\setminus 0)\times \R^{d_2})$
that is the space of function  $a\in C^{\infty}(\R^n\times \R^{d_1}\times \R^{d_2})$ that satisfy
\beq\label{product symbols}
|\p_x^\gamma\p_{\theta^{\prime}}^\alpha \p_{\theta^{\prime\prime}}^\beta a(x,\theta^{\prime},\theta^{\prime\prime})|
\leq C_{\alpha\beta\gamma K}(1+|\theta^{\prime}|+|\theta^{\prime\prime}|)^{\mu-|\alpha|}
(1+|\theta^{\prime\prime}|)^{\mu^{\prime}-|\beta|}\hspace{-1.7cm}
\eeq
for all $x\in K$, multi-indexes  $\alpha,\beta,\gamma$, and
compact sets $K \subset \R^n$. Above,
$\mu=p+l-d_1/2+n/4$ and
$\mu^{\prime}=-l-d_2/2$.

When $X$ is a manifold of dimension $n$ and 
$\Lambda_1,\Lambda_2\subset T^*X
\setminus \{0\}$ are two
cleanly intersecting Lagrangian manifolds, 
we define the class $\I^{p,l}(M;\Lambda_1,\Lambda_2)\subset \mathcal D^\prime(X)$ to consist
of locally finite sums of functions $u=Au_0$, where 
$u_0\in \I^{p,l}(\R^n;N^*S_1,N^*S_2)$ and $S_1,S_2\subset \R^n$ are the linear subspace of codimensions
$d_1$ and $d_1+d_2$, respectively, such that
$S_2\subset S_1$, 
and $A$ is a Fourier integral
operator of order zero with a canonical relation $\Sigma$ for which
$ \Sigma\circ (N^*S_1)^\prime\subset \Lambda_1^\prime$ and 
$\Sigma\circ (N^*S_2)^\prime\subset \Lambda_2^\prime$. 
Here, $ \Lambda^\prime=\{(x,-\xi)\in T^*X;\ (x,\xi)\in \Lambda\}$.

In most cases, below $X=M$. We denote then
  $\I^{p}(M;\Lambda_1)=\I^{p}(\Lambda_1)$
and 
 $\I^{p,l}(M;\Lambda_1,\Lambda_2)=\I^{p,l}(\Lambda_1,\Lambda_2)$, etc.
Also, $\I(\Lambda_1)=\cup_{p\in \R}\I^{p}(\Lambda_1)$ etc.

By \cite{GU1,MU1}, microlocally away from $\Lambda_2$ and $\Lambda_1$,
\beq\label{microlocally away}
& &\I^{p,l}(\Lambda_1,\Lambda_2)\subset \I^{p+l}(\Lambda_1\setminus \Lambda_2)\quad
\hbox{and}\quad \I^{p,l}(\Lambda_1,\Lambda_2)\subset \I^{p}(\Lambda_2\setminus \Lambda_1),
\eeq
correspondingly.
Thus the principal symbol of $u\in \I^{p,l}(\Lambda_0,\Lambda_1)$
is well defined on $\Lambda_0\setminus \Lambda_1$ and  $\Lambda_1\setminus \Lambda_0$.

Below, when $\Lambda_j=N^*S_j,$ $j=1,2$ are conormal bundles of smooth cleanly
intersecting submanifolds
$S_j\subset M$ of codimension $d_j$, where $\dim(M)=n$,
 we
use the traditional notations,
\beq\label{eq: traditional}\quad\quad
\I^\mu(S_1)=\I^{\mu+1/2-n/4}(N^*S_1),\quad \I^{\mu_1,\mu_2}(S_1,S_2)
=\I^{p,l}(N^*S_1,N^*S_2),\hspace{-1cm}
\eeq
where $p=\mu_1+\mu_2+d_1/2-n/4$ and $l=-\mu_2-d_2/2$, \HOX{Slava asked about notations: Is it that
$\mu_2=\mu^\prime$, $\mu_1=\mu^{\prime\prime}$, and $\mu=\mu_1+\mu_2$. Discuss more about this!}
and call such distributions the conormal distributions  associated to $S_1$ or 
product type  conormal distributions  associated to
$S_1$ and $S_2$.  {We note that
$\I^\mu(X;S_1)\subset L^p_{loc}(X)$ for $\mu< -d_1(p-1)/p$, $1\leq p<\infty$, see \cite{GU1}.}

By \cite{MU1}, see also \cite{K2}, a classical pseudodifferential operator $P$ of real
principal type and order $m$  on $M$ has a parametrix $Q\in \I^{p,l}(\Delta^\prime_{T^*M},\Lambda_P)$, 
$p=\frac 12-m$, $l=-\frac 12$ ,
where $\Delta_{T^*M}=N^*(\{(x,x);\ x\in M\})$ and $\Lambda_P\subset T^*M\times T^*M$
is the Lagrangian manifold associated to the canonical relation of  the operator $P$,
that is,
\ba
\Lambda_P=\{(x,\xi,y,-\eta);\ (x,\xi)\in\hbox{Char}\,(P),\ (y,\eta)\in \Theta_{x,\xi}\},
\ea
where $\Theta_{x,\xi}\subset T^*M$ is the bicharacteristic of $P$ containing $(x,\xi)$.  

For the
wave operator $\square_{g}$ on the globally hyperbolic manifold $(M,{g})$
$\hbox{Char}\,(\square_{g})$ is the set of light-like vectors with respect to $g$,
and $(y,\eta)\in \Theta_{x,\xi}$ if and only if there is $t\in\R$ such
that $(y,a)=(\gamma^{g}_{x,b}(t),\dot \gamma^{g}_{x,b}(t))$
where $\gamma^{g}_{x,b}$ is a light-like geodesic with resect to the metric $g$
with the initial data $(x,b)\in TM$, $a=\eta^\flat$, $b=\xi^\flat$.
For $P=\square_{g}+B^0+B^j\p_j$, where $B^0$ and $B_j$
are tensors, we denote  $\Lambda_{P}=\Lambda_g$. When
 $(M,g)$ is globally hyperbolic manifold, the operator 
$P$ has a causal inverse operator, see e.g.\  \cite[Thm.\ 3.2.11]{BGP}.
We denote it by  $P^{-1}$ and by \cite{MU1}, we have
$P^{-1} \in \I^{-3/2,-1/2}(\Delta^\prime_{T^*M},\Lambda_g)$.
We will repeatedly use the fact (see \cite[Prop. 2.1]{GU1}) that 
if $F\in \I^{p}(\Lambda_0)$ and $\Lambda_0$ intersects Char$(P)$ transversally
so that all bicharacterestics of $P$ intersect $\Lambda_0$ only finitely many
times, then $(\square_{g}+B^0+B^j\p_j )^{-1}F\in \I^{p-3/2,-1/2}(\Lambda_0,\Lambda_1)$
where $\Lambda_1^{\prime}=\Lambda_g\circ \Lambda_0^{\prime}$ is called the flowout from $\Lambda_0$ on Char$(P)$,
that is,
\ba
\Lambda_1=\{(x,-\xi);\ 
(x,\xi,y,-\eta)\in \Lambda_g,\ (y,\eta)\in \Lambda_0 \}.
\ea


\subsubsection{Notations used to consider Einstein equations}

For Einstein equation, we will consider a smooth background metric $\hat g$ on $\hattuM $
and the smooth metric $\tilde g$ for which $\hat g<\tilde g$ and $(\hattuM ,\tilde g)$
is globally hyperbolic. We also use the notations defined in Section \ref{subsec: Gloabal hyperbolicity}. 
In particular, we identify $M=\R\times N$ and consider the metric tensor $g$ on $\hattuM _0=
(-\infty,t_0)\times N$, $t_0>0$ that
coincide with $\hat g$ in $(-\infty,0)\times N$. Recall also that we consider
a freely falling observer $\hat \mu=\mu_{\hat g}:[-1,1]\to \hattuM _0$
for which $\hat \mu(s_-)=\hat p^-\in  (0,t_0)\times N$. 
We denote the cut locus function on $(M_0,\hat g)$ by $\rho(x,\xi)=\rho_{\hat g}(x,\xi)$,
 denote $L_x^+M_0=L^+_x(M_0,\hat g)$ and $L^+M_0=L^+(M_0,\hat g)$, 
and denote by $\hat U=U_{\hat g}$ the neighborhood of the geodesic $\hat \mu=\mu_{\hat g}$,
and denote by $\hat \gamma_{x,\xi}(t)$ the geodesics of $(M_0,\hat g)$.
To simplify notations,
we  study below the case when $\hat P=0$ and $\hat Q=0$,
that corresponds to the setting of Theorem \ref{main thm Einstein}. 


%

%

\subsection{Direct problem}\label{subsec: Direct problem}

{Let  we consider the solutions $(g,\phi)$ of the equations  (\ref{eq: adaptive model})
with  source $F=(P,Q)$.
To consider their existence, let us write them in the form
\beq\label{u and solutions}
& &(g,\phi)=(\hat g,\hat \phi)+u,\quad\hbox{where }
u=((G_{jk})_{j,k=1}^4,(\Phi_\ell)_{\ell=1}^L).
\eeq
 Next we formulate an appropriate equation for $u$ and review the
 existence and uniqueness results for it.

Let us assume that 
 $F$ is small enough in the norm  $C^5_b(M_0)$ and that it is  supported in a compact set  $\K=J_{\tilde g}(\hat p^-)\cap  [0,t_0]\times N\subset M_1$.
Then using the fact that $\p_jg^{nm}=-g^{na}(\p_j g_{ab})g^{bm}$ we can write the equations (\ref{eq: adaptive model}) for $u$ appearing in (\ref{u and solutions}) in the form
\beq\label{eq: notation for hyperbolic system 1}
& &P_{g(u)}(u)={R}_{\hat g}(x,u(x),\p u(x)) F,\quad x\in \hattuM _0,\\
& & \nonumber \supp(u)\subset \K,
\eeq
where
\beq\label{eq: notation for hyperbolic system 1b B}
& &P_{g(u)}(u):=g^{jk}(x;u)\p_j\p_k  u(x)+H(x,u(x),\p u(x)),\\
& &(g^{jk}(x;u))_{j,k=1}^4=((\hat g_{jk}(x)+G_{jk}(x))_{j,k=1}^4)^{-1},\nonumber
\eeq
{and  $(x,v,w)\mapsto H(x,v,w)$ is a smooth function which is a second order polynomial in $w$ 
with coefficients being smooth functions of $v$ and derivatives of $\hat g$
and ${R}_{\hat g}(x,u(x),\p u(x))F$, when represented in local coordinates, is a first order  linear differential operator in $F$, i.e., at the point $x$ 
it depends linearly on $F(x)$ and $\p F(x)$,
which coefficients  are smooth functions depending on  $u(x)$, $\p u(x)$,
and derivatives of $\hat g$.}

\observation{
\medskip

{\bf Remark 3.1.}
We note that  the Laplace-Beltrami
operator can be written  for 
as $\square_g \psi_\ell =g^{jk}\p_j\p_k \phi_\ell-g^{pq} \Gamma^n \p_n \phi_\ell=
g^{jk}\p_j\p_k \phi_\ell- \Gamma^n \p_n \phi_\ell$
and thus  in the $\hat g$-wave map coordinates we have
\beq\label{eq: wave operator in wave gauge}
\square_g \phi_\ell=g^{jk}\p_j\p_k \phi_\ell-g^{pq}\hat \Gamma^n_{pq} \p_n \phi_\ell.
\eeq
Thus, the scalar field equation $\square_g \phi_\ell+m\phi_\ell=0$
does not involve derivatives of $g$.  As this can be the case in
 (\ref{eq: notation for hyperbolic system 1b B}),
the system (\ref{eq: notation for hyperbolic system 1b B}) is
a slight generalization of (\ref{eq: adaptive model}).
\medskip
}

 Below we endow $N$
with the Riemannian metric $h:=\iota^*\tilde g$ inherited from the embedding
$\iota: N\to \{0\}\times  N\subset M$ and
use it to define the Sobolev spaces $H^s(N)$.

Let $s_0\geq 5$ be an odd integer. 
Below we will consider the solutions $u=(g,\phi)$ as sections
of the bundle $\B^L$ on $M_0$ and the sources $F=(Q,P)$ as sections
of the bundle $\B^K$ on $M_0$.
We will consider these functions as elements
of section valued Sobolev spaces $H^s(M_0;\B^L)$ and  $H^s(M_0;\B^K)$   etc.
Below, we
omit the bundle $\B^K$ in these notations and denote
just $H^s(M_0;\B^L)=H^s(M_0)$ etc. We use the same convention
for the spaces
\ba
E^s=\bigcap_{j=0}^s C^j([0,t_0];H^{s-j}(N)),\quad s\in \N.
\ea
Note that $E^s\subset C^p([0,t_0]\times N))$ when $0\leq p<s-2$.
 Using standard techniques
for quasi-linear equations developed e.g.\ in 
  \cite{Kato1975} or \cite{HKM},
or \cite[Section 9
]{Ringstrom} (for details, see  Appendix C),
we see that 
when $F$ is supported in the compact set $\K$ and $\|F\|_{E^{s_0}}<c_0$, then 
there exists a unique function $u$ satisfying equation (\ref{eq: notation for hyperbolic system 1}) 
on $M_0$  with the source $F$ and 
\beq\label{eq: Lip estim}
\|u\|_{{E^{s_0-1}}}\leq
C_1\|F\|_{E^{s_0}}.
\eeq

%

For convenience of the reader
we give the proofs of these facts in Appendixes A and C.


\subsection{Asymptotic expansion for non-linear wave equation}


Let us  consider a small parameter $\e>0$ and the sources
 $P=\e {\bf p}$ and  $Q^{\prime}=(Q_\ell)_{\ell=1}^{K-1}=\e{\bf q^{\prime}}$, $Q_K=\e{\bf z}$ (cf.\
Assumption S, (ii)), and denote ${\bf q}=({\bf q^{\prime}},{\bf z})$.
This corresponds to the source $F=\e f=(\e {\bf p},\e{\bf q^{\prime}},\e{\bf z})$ in 
 (\ref{eq: notation for hyperbolic system 1}).
Below, we always assume that ${\bf p}$ and ${\bf q}$ are supported in $\K$ and
${\bf p},{\bf q}\in E^s$, where $s\geq 13$ is an odd integer.
We consider the solution $u=u_\e$ of (\ref{eq: notation for hyperbolic system 1})\HOX{Some details  should
 be removed from final version}
 and write it in the form
\beq\label{eq: epsilon expansion}
 u_\e(x)=\sum_{j=1}^4\e^j w^j(x)+w^{res}(x,\e).
 \eeq
%
To obtain equations for $ w^j$, we use the
 representation   for the $\hat g$-reduced Ricci tensor
 given
 in the Appendix A, see (\ref{Reduced Ric tensor}) 
to write an analogous representation for the  $\hat g$-reduced Einstein tensor,
see (\ref{Reduced Einstein tensor})  or (\ref{hat g reduced einstein and reduced einstein}),
 and
substitute the expansion (\ref {eq: epsilon expansion})
in to the equation (\ref{eq: notation for hyperbolic system 1}).
This, Assumption S and 
equation (\ref{S-vanish condition}) 
imply
 that $ w^j$, $j=1,2,3,4$ are
 given in term of sources
 of $ \mathcal H^1= \mathcal H^1(\hat g,f),$
  $ \mathcal H^2= \mathcal H^2(\hat g,f,w^1),$ 
   $ \mathcal H^3= \mathcal H^3(\hat g,f,w^1,w^2),$ and 
   $ \mathcal H^4= \mathcal H^4(\hat g,f,w^1,w^2,w^3)$ that
   satisfy 
   \beq\nonumber 
 w^j&=&(g^j,\phi^j)={\bf Q}_{\hat g}\mathcal H^j,\quad j=1,2,3,4,\hbox{ where }\\
\nonumber
 \mathcal H^1&=&R_{\hat g}(x,0,0)f,\\ \nonumber
  \mathcal H^2&=&(\mathcal G_2,0)
+\A_2^{(2)}( w^1,\p  w^1;f,\p f),\\
\nonumber
\mathcal G_2&=&2\hat g^{jp}w^1_{pq}\hat g^{qk}\p_j\p_k   w^1,\\
  \mathcal H^3&=& (\mathcal G_3,0)+\A_2^{(3)}( w^1,\p  w^1, w^2,\p  w^2;f,\p f)+
 \label{w1-w3}
\\ \nonumber
& &\ \ +
\A_3^{(3)}( w^1,\p  w^1;f,\p f),
\\
\nonumber
\mathcal G_3&=&
-6\hat g^{jl}w^1_{li}\hat g^{ip}w^1_{pq}\hat g^{qk}\p_j\p_k   w^1+\\
& &\ \ +
3\hat g^{jp}( w^2)_{pq}\hat g^{qk}\p_j\p_k   w^1+3\hat g^{jp}w^1_{pq}\hat g^{qk}\p_k\p_j   w^2
\nonumber
\eeq
where $f\mapsto R_{\hat g}(x,0,0)f$ is a first order linear differential
 operator appearing in (\ref{eq: notation for hyperbolic system 1}),
and \HOX{We could simplify notations by writing 
$\A_2^{(4)}+\A_3^{(4)}+\A_4^{(4)}=\A^{(4)}$ etc.}
\beq  \nonumber
\hspace{-10mm}  \mathcal H^4&=& (\mathcal G_4,0)+\A_2^{(4)}( w^1,\p  w^1, w^2,\p  w^2, w^3,\p  w^3;f,\p f)\\
\nonumber
& &\ \ +\A_3^{(4)}( w^1,\p  w^1, w^2,\p  w^2;f,\p f) +\A_4^{(4)}( w^1,\p  w^1;f,\p f),\\
\nonumber
\hspace{-10mm} \mathcal G_4 &=&24\hat g^{js}w^1_{sr}\hat g^{rl}w^1_{li}\hat g^{ip} w^1_{pq}\hat g^{qk}\p_k\p_j   w^1
+\\
  \nonumber
& &\ \ 
-18\hat g^{jl}w^1_{li}\hat g^{ip} w^2_{pq}\hat g^{qk}\p_k\p_j   w^1
-12\hat g^{jl}w^1_{li}\hat g^{ip} w^1_{pq}\hat g^{qk}\p_k\p_j   w^2+
\hspace{-15mm}\\
\label{w4}
& &\ \ +
3\hat g^{jp}w^3_{pq}\hat g^{qk}\p_k\p_j   w^1+
3\hat g^{jp}w^1_{pq}\hat g^{qk}\p_k\p_j   w^3+\\
  \nonumber
& &\ \ +
6\hat g^{jp}w^2_{pq}\hat g^{qk}\p_k\p_j   w^2.
  \nonumber
\eeq
Here, we have used the notation $ w^j=(( w^j)_{pq})_{p,q=1}^4,(( w^j)_\ell)_{\ell=1}^L)$
with $(( w^j)_{pq})_{p,q=1}^4$ being the $g$-component of $ w^j$ and 
$(( w^j)_\ell)_{\ell=1}^L$
being the $\phi$-component of $ w^j$. Moreover,
${\bf Q}_{\hat g}=(\square_{\hat g}+V(x,D))^{-1}$
is the causal inverse of the operator $\square_{\hat g}+V(x,D)$
where $V(x,D)$ is a first order differential operator with coefficients
depending on $\hat g$ and its derivatives,
$R(x,D)$ is  a first order linear differential operator
 with coefficients
depending on $\hat g$ and its derivatives,
 and
 $\A_m^{(\a)}$, $\a=2,3,4$ denotes a multilinear operators  of order $m$ 
 having at a point $x$ the 
representation
\beq\label{eq: mixed source terms}
& &\hspace{1cm}(\A_m^{(\a)}(v^1,\p v^1,v^2,\p v^2,v^3,\p v^3;f,\p f))(x)
\\ &=& 
\nonumber
\sum_{} \bigg(a^{(\a)}_{abcijkP_1P_2P_3pqn\ell_1\ell_2 }(x)\,
 (v^1_a(x))^i(v^2_b(x))^j (v^3_c(x))^k \cdotp \\
 \nonumber
 & &\quad \quad\quad \cdotp P_1 ( \p v^1(x))P_2 ( \p  v^2(x)) P_3 ( \p v^3(x))  ( f_p(x))^{\ell_1}(\p_n f_{q} (x))^{\ell_2}
 \bigg)\hspace{-1cm} 
\eeq
where 
$(v^1_a(x))^i$ denotes the $i$-th power of $a$-th component of $v^1(x)$, etc, and
 the sum is taken over the indexes $a,b,c,p,q,n,$ integers
$ i,j,k,\ell_1,\ell_2$,
and the homogeneous monomials $P_d(y)=y^{\beta_d}$,
$\beta_d=(b_1,b_2,\dots,b_{4(10+L)})\in \N^{4(10+L)}$, $d=1,2,3$
 having orders $|\beta_d|$, correspondingly,  that satisfy
\beq
\label{Eq: F condition1}
& &\ell_1+\ell_2\leq 1,\\
\label{Eq: F condition1a}
& &i+2j+3k+|\beta_1|+2|\beta_2|+3|\beta_3|+\ell_1+\ell_2= \alpha,\\
\label{Eq: F condition2}
& &i+j+k+|\beta_1|+|\beta_2|+|\beta_3|+\ell_1+\ell_2= m,\quad\hbox{and}\\
& &\hbox{if $\ell_1=0$ and $\ell_2=0$ then $|\beta_1|+|\beta_2|+|\beta_3|\leq 2$.}
\label{Eq: F condition3}\eeq
Here, condition (\ref{Eq: F condition1}) means
that $\A_m^{(\a)}$
is an affine function of $f$ and its first derivative, 
condition (\ref{Eq: F condition1a}) implies that the term $\A_m^{(\a)}$ produces
a term of order $O(\e^\alpha)$ when $v^j=w^j$, 
condition (\ref{Eq: F condition2}) that $\A_m^{(\a)}$ is multilinear of 
order $m$, and
condition (\ref{Eq: F condition3}) means
that for $x\not\in \supp(f)$ the non-vanishing 
terms in $\A_m^{(\a)}(v^1,\p v^1,v^2,\p v^2,v^3,\p v^3;f,\p f)$
contain only terms where the total power of derivatives
of $v^1,v^2$, and $v^3$ is at most two. 
We note that the inequalities (\ref{Eq: F condition1})-(\ref{Eq: F condition3}) follow from
 Assumption S and 
equation (\ref{S-vanish condition}).

 

 By \cite[App.\ III, Thm.\ 3.7]{ChBook}, or alternatively, the proof of \cite [Lemma 2.6]{HKM} adapted for manifolds, for details, see 
 (\ref{U bounds}) and   (\ref{U bounds2}) in the Appendix C,
the estimate  $\| {\bf Q}_{\hat g}H\|_{E^{s_1}}\leq C_{s_1}\|H\|_{E^{s_1}}$
 holds  for all   $H$ that are supported in $\K=\K_0$, having the form (\ref{Kappa sets}) and satisfy 
  $H\in E^{s_1}$, $s_1\in \Z_+$.
Note that we are interested only on the local solvability of the Einstein
equations and that singularities can appear on 
sufficiently large time intervals, see \cite{Dafermos}.
  
Recall that $s\geq 13$ is an odd integer and $f\in E^{s+1}$.
 By defining $ w^j$ via the above equations with $f\in E^{s+1}$ we obtain 
  $ w^j\in E^{s+2-2j}$,
for  $j=1,2,3,4$. Thus, by using Taylor's expansion of the coefficients
  in the equation (\ref{eq: notation for hyperbolic system 1})
 we see that the approximate 4th order expansion
 $u^{app}_\e=\e  w^1+\e^2  w^2+\e^3  w^3+\e^4  w^4$
 satisfies an equation of the form
 \beq\label{eq: notation for hyperbolic system 1b}
& &P_{g(u^{app}_\e)}(u^{app}_\e)
\hspace{-1mm}=\hspace{-1mm}{R}(x,u^{app}_\e(x),\p u^{app}_\e(x)) 
(\e f) \hspace{-1mm}+\hspace{-1mm}H^{res}(\,\cdotp,\e),\ \ x\in \hattuM _0,
\hspace{-5mm}\\
& & \nonumber \supp(u^{app}_\e)\subset \K,
\eeq
%
%
%
such that 
 $$\|H^{res}(\,\cdotp,\e)\|_{E^{s-8}}
 \leq c_1\e^5.  $$ 
 Using the Lipschitz continuity of the solution
 of the equation (\ref{eq: notation for hyperbolic system 1b}), see Appendix C,
 we see that
 there are $c_1,c_2>0$ such that for all $0<\e<c_1$ the function
 $w^{res}(x,\e)=u_\e(x)-u^{app}_\e(x)$ satisfies 
  \ba
 \|u_\e-u^{app}_\e\|_{E^{s-8}}
 \leq c_2\e^5.
 \ea
Thus $ w^j=\p_\e^j u_\e|_{\e=0}\in E^{s-8}$, $j=1,2,3,4.$}
%
%
%
 \subsection{Linearized conservation law and divergence condition}
 
 \subsubsection{Linearized Einstein equation}

We will below consider  sources $Q^{\prime}=(Q_\ell)_{\ell=1}^{K-1}=\e{\bf q^{\prime}}$, $Q_K=\e{\bf z}$ and 
 $P=\e {\bf p}$.
We denote 
${\bf q}=({\bf q^{\prime}},{\bf z})$.

To analyze the linearized waves, we denote $ w^1=u^{(1)}$.
We see that 
$u^{(1)}$ satisfies the linear equation
\beq\label{lin wave eq}
& &\square_{\hat g}  u^{(1)}+V(x,\p_x) u^{(1)}={\bf h},
\eeq
where $v\mapsto V(x,\p_x)v$ is a linear first order partial differential operator with
coefficients depending on the derivatives of $\hat g$ 
{and, 
${\bf h}=H(x;{\bf p},{\bf q})$, where 
%
%
\beq
 \nonumber
H(x;{\bf p},{\bf q})
&=&\left(\begin{array}{c}
{\bf p}\\
 0\end{array}\right)+  \left(\begin{array}{c}
M_{(1)}(x) \,{\bf q^{\prime}}(x)
\\  
 M_{(2)}(x) \, {\bf q^{\prime}}(x)
  \end{array}\right)+
  \\
  & &\label{lin wave eq source B}\hspace{-3cm}
 + \left(\begin{array}{c}
L_{(1)}(x)\,{\bf z}(x)
\\  
 L_{(2)}(x)\,{\bf z}(x)
  \end{array}\right)
+
   \left(\begin{array}{c}
N^j_{(1)} (x)\, \hat g^{lk}\hat \nabla_l  ({\bf p}_{jk}+{\bf z}\hat g_{jk})
\\ 
N^j_{(2)}(x) \, \hat g^{lk}\hat \nabla_l ( {\bf p}_{jk}+{\bf z}\hat g_{jk})
\end{array}\right),\hspace{-2cm}
%
  \eeq
  where $M_{(k)}=M_{(k)}(\hat \phi(x),\hat \nabla  \hat \phi(x),\hat g(x))$,
  $  L_{(k)}=L_{(k)}(\hat \phi(x),\hat \nabla  \hat \phi(x),\hat g(x))$,
  and $N^j_{(k)}=N^j_{(k)}(\hat \phi(x),\hat \nabla  \hat \phi(x),\hat g(x))$
  are, in local coordinates, matrices whose elements  are smooth functions 
of  $\hat \phi(x),\hat \nabla  \hat \phi(x),$ and $\hat g(x)$.
Later we use the fact that  by the form of the equation
(\ref{eq: adaptive model}), $M_{(1)}=0$, $N^j_{(1)}=0$, and $L_{(1)}=\hat g$.}
We see using Assumption S (iii) that  the union of the image spaces of  the matrices 
$M_{(2)}(x)$ and
$L_{(2)}(x) $,  
and $N^j_{(2)}(x) $,  $j=1,2,3,4$,
span the space $\R^L$  for all $x\in \hat U$.
%
 %

\subsubsection{The linearized conservation law for sources} 
Assume that $Y\subset \hattuM _0$
is a 2-dimensional space-like submanifold and consider local coordinates defined 
in  $V\subset \hattuM _0$. Moreover, assume that 
in these local coordinates $Y\cap V\subset \{x\in \R^4;\ x^jb_j=0,\  x^jb^\prime_j=0\}$,
where $b^\prime_j$ are some constants and let ${\bf p}\in \I^{n-3/2}(Y)$, $n\leq n_0=-15$,  be defined by 
\beq\label{eq: b b-prime}
{\bf p}_{jk}(x^1,x^2,x^3,x^4)=
\re \int_{\R^2}e^{i(\theta_1b_m+\theta_2b^\prime_m)x^m}c_{jk}(x,\theta_1,\theta_2)\,
d\theta_1d\theta_2.\hspace{-2cm}
\eeq
Here, we assume that 
 $c_{jk}(x,\theta)$, $\theta=(\theta_1,\theta_2)$ 
 are  classical symbols and we denote their principal symbols by
 $\sigma_p({\bf p}_{jk})(x,\theta)$. When 
$x\in Y$ and $\xi=(\theta_1b_m+\theta_2b^\prime_m)dx^m$ so that $(x,\xi)\in N^*Y$,
we denote the value of the principal symbol $\sigma_p({\bf p})$ at $(x,\theta_1,\theta_2)$ by 
 $\tilde c^{(a)}_{jk}(x,\xi)$, that is,  $ \tilde c^{(a)}_{jk} (x,\xi)=\sigma_p({\bf p}_{jk})(x,\theta_1,\theta_2)$,
 and say that it is the principal symbol of ${\bf p}_{jk}$ at $(x,\xi)$, associated to the local $X$-coordinates
 and the phase function $\phi(x,\theta_1,\theta_2)=(\theta_1b_m+\theta_2b^\prime_m)x^m$.
The above defined principal symbols can be defined invariantly, see \cite{GuU1}.

We assume that also ${\bf q}^\prime,{\bf z}\in \I^{n-3/2}(Y)$ have 

representations (\ref{eq: b b-prime}) with classical symbols.
{Let
 us denote the principal symbols of ${\bf p},{\bf q^{\prime}},{\bf z}\in \I^{n-3/2}(Y)$ 
 by  $\tilde c^{(a)}(x,\xi),$ $\tilde d_{1}(x,\xi)$,
 $\tilde d_2^{(a)}(x,\xi)$, correspondingly
  and let $\tilde c^{(b)}$ and $\tilde d_2^{(b)}(x,\xi)$ denote the sub-principal
symbols of ${\bf p}$ and ${\bf z}$, correspondingly, at $(x,\xi)\in N^*Y$.


We will below consider what happens when  $({\bf p}_{jk}+{\bf z}\hat g_{jk})\in \I^{n-3/2}(Y)$ satisfies 
\beq\label{eq: lineariz. conse. law}
\hat g^{lk} \nabla_l^{\hat g} ({\bf p}^{(a)}_{jk}+{\bf z}\hat g_{jk})\in \I^{n-3/2}(Y),\quad j=1,2,3,4.
\eeq
When (\ref{eq: lineariz. conse. law}) is valid, we say
that the leading order of singularity of the wave satisfies
the {\it linearized conservation
law}. This  corresponds to the assumption that the principal
symbol of the sum of divergence of the first two terms appearing in the stress energy tensor
on the right hand side of (\ref{eq: adaptive model}) vanishes.
}


When  (\ref{eq: lineariz. conse. law})  is valid, we have
\beq\label{eq: lineariz. conse. law symbols}
& &\hat g^{lk}\xi_l(\tilde c^{(a)}_{kj} (x,\xi)
+\hat g_{kj}(x)\tilde d_2^{(a)}(x,\xi))=0,\quad\hbox{for 
$j\leq 4$ and $\xi\in N^*_x Y$}.
\eeq
We say that this is the {\it linearized conservation law for principal symbols}.

\subsubsection{A divergence condition for linearized solutions}

Above, we have considered the $\hat g$-reduced Einstein equations.
As discussed in the Appendix A, if we pull back 
a solution of the Einstein equation via a $(\hat g,g)$-wave map,
the metric tensor will satisfy the $\hat g$-reduced Einstein equations.
Roughly speaking, this means that we are in fact considering
the metric tensor $g$ in special coordinates associated to the 
metric $\hat g$. Next we recall some well known consequences of this. 

Assume now that $(g,\phi)$ satisfy equations (\ref{eq: adaptive model}).
By Assumption S (iv), the conservation law (\ref{conservation law0})
is valid. 
As discussed in Appendix A.5,  
the conservation law (\ref{conservation law0}) and 
the $\hat g$-reduced Einstein equations (\ref{eq: adaptive model}) imply 
 that the harmonicity functions $\Gamma^j=g^{nm} \Gamma^j_{nm}$
 satisfy 
 \beq\label{harmonic condition}
g^{nm} \Gamma^j_{nm}= g^{nm}\hat \Gamma^j_{nm}.
\eeq
Next we discuss the implications  of this  for the metric
component of the solution of
the linearized Einstein equation.

 \HOX{We have to  look literature
 if this subsubsection can be removed from the final paper. Footnote should
 be shortened.}

Let us next do calculations  in local coordinates of $M_0$ and 
denote $\p_k=\frac\p{\p x^k}$.
Direct calculations show that 
$
h^{jk}=g^{jk}\sqrt{-\det(g)}
$
satisfies $\p_k h^{kq}= - \Gamma^q_{kn}h^{nk}$.
Then (\ref{harmonic condition})
implies that 
\beq\label{harmonic condition - alternative2}
\p_k h^{kq}=
-\hat \Gamma^q_{kn}h^{nk}.
\eeq
We call (\ref{harmonic condition - alternative2}) the {\it  harmonicity condition} for  metric $g$.

Assume now that  $g_\e$ and $\phi_\e$ satisfy (\ref{eq: adaptive model})
with sources
 $P=\e {\bf p}$ and $Q=\e {\bf q}$ where $\e>0$
 is a small parameter.
 We define $h_\e^{jk}=g_\e^{jk}\sqrt{-\det(g_\e)}$
 and denote $\dot g_{jk}=\p_\e (g_\e)_{jk}|_{\e=0}$,
$\dot g^{jk}=\p_\e (g_\e)^{jk}|_{\e=0}$,
and $\dot h^{jk}=\p_\e h_\e^{jk}|_{\e=0}$.

%
%
%

The equation (\ref{harmonic condition - alternative2})
yields then\footnote{The treatment on this de Donder-type gauge condition is known in
folklore of the field. For a similar gauge condition to (\ref{harmonic condition for dots1}) 
in harmonic coordinates, see \cite[pages 6 and 250]{Maggiore},  or  
\cite[formulas 107.5, 108.7, 108.8]{Landau}, or \cite[p.\ 229-230]{HE}.
}
\beq\label{harmonic condition for dots1}
\p_k\dot h^{kq}=-\hat \Gamma^q_{kn}\dot h^{nk}.
\eeq
A direct computation  shows that 
\ba
\dot h^{ab}
&=&(-\det(\hat g))^{1/2}\kappa^{ab},
\ea
where $\kappa^{ab}=\dot g^{ab}-\frac 12 \hat g^{ab}\hat g_{qp}\dot g^{pq}$.
Thus (\ref{harmonic condition for dots1}) gives
\beq\label{extra 2a}
\p_a((-\det(\hat g))^{1/2}\kappa^{ab})
=
-\hat \Gamma^b_{ac}(-\det(\hat g))^{1/2}\kappa^{ac}
\eeq
that implies 
$
\p_a\kappa^{ab}
+\kappa^{nb}\hat\Gamma^a_{an}+\kappa^{an}\hat \Gamma^b_{an}=0,
$
or equivalently, 
\beq\label{extra 2b}
\hat \nabla_a\kappa^{ab}=0.
\eeq
We call (\ref{extra 2b}) the {\it  linearized divergence condition} for $g$. 
Writing this for $\dot g$, we obtain
\beq\label{divergence condition}
& &\hspace{-1cm}-\hat g^{an}\p_a \dot g_{nj}+\frac 12 \hat g^{pq}  \p_j\dot g_{pq}=m^{pq}_j\dot g_{pq}
\eeq
where $m_j$  depend on $\hat g_{pq}$ and its derivatives.
On similar conditions for the polarization tensor, see \cite[form.\ (9.58) and example  9.5.a, p. 416]{Padmanabhan}.

\subsubsection{Properties of the principal symbols of the waves}

Let $K\subset \hattuM _0$ be a light-like submanifold of dimension 3.
We use below local coordinates  $X:V\to \R^4$ defined
in a neighborhood $V$ of a point $x_0\in M_0$.
We use the Euclidian coordinates
 $x=(x^k)_{k=1}^4$ on $X(V)$. We assume these coordinates are such that 
in $V$ the submanifold $K$ is given by $K\cap V\subset \{x\in \R^4;\ b_kx^k=0\}$, where $b_k\in \R$ are constants.
Assume that the solution $u^{(1)}=(\dot g,\dot \phi)$ of the 
linear wave equation (\ref{lin wave eq}) with the right hand
side vanishing  in $V$ is such that  $\dot g_{jk}\in \I^\mu (K)$ with a suitable $\mu\in \R$.
Let us write    $\dot g_{jk}$  as an oscillatory integral using 
 a phase function $\varphi(x,\theta)=b_kx^k\theta$, and
a symbol $a_{jk}(x,\theta)\in S^\mu (\R^4,\R)$, 
\beq\label{eq: oscil}
\dot g_{pq}(x^1,x^2,x^3,x^4)=
\re \int_{\R}e^{i(\theta b_mx^m)}a_{pq}(x,\theta)\,
d\theta.
\eeq
We assume that $a_{jk}(x,\theta)$
is a classical symbol and denote
its  (positively homogeneous) principal symbol by 
$\sigma_p(\dot g_{pq})(x,\theta)$. 
When $x\in K$ and $\xi=\theta b_kdx^k$ so that $(x,\xi)\in N^*K$,
we denote the value of $\sigma_p(\dot g_{pq})$ at $(x,\theta)$ by 
 $\tilde a_{jk}(x,\xi)$, that is,  $\tilde a_{jk}(x,\xi)=\sigma_p(\dot g_{pq})(x,\theta)$
 and say that it is the principal symbol of $\dot g_{pq}$ at $(x,\xi)$, associated to the local $X$-coordinates
 and the phase function $\phi(x,\theta)=\theta b_mx^m$.

%
%
Then, 
if $\dot  g_{jk}$ satisfies the divergence condition (\ref{divergence condition}),
its principal symbol $\tilde a_{jk}(x,\theta)$ satisfies 
\beq\label{divergence condition for symbol}
& &\hspace{-1cm}-\hat g^{mn}(x)\xi_m v_{nj}+\frac 12  \xi_j(\hat g^{pq}(x) v_{pq})=0,
\quad v_{pq}=\tilde a_{pq}(x,\xi),
\eeq
where $j=1,2,3,4$ and  $\xi=\theta b_kdx^k\in N_x^*K$.
 If (\ref{divergence condition for symbol}) holds, we say that the {\it  divergence condition for the symbol} is satisfied for $ \tilde a(x,\xi)$
at $(x,\xi)\in N^*K$.

\subsection{Pieces of spherical waves satisfying linear wave equation}

\subsubsection{Solutions which are singular on hypersurfaces}

 Next we consider a 
piece of  spherical wave whose singular support is
concentrated near a  geodesic. 

\begin{figure}[htbp] \label{Fig-7}
\begin{center}

\psfrag{1}{$y_0$}
\psfrag{2}{$\Sigma$}
\psfrag{3}{$\Sigma_1$}
\psfrag{4}{$y^\prime$}
\includegraphics[width=9.5cm]{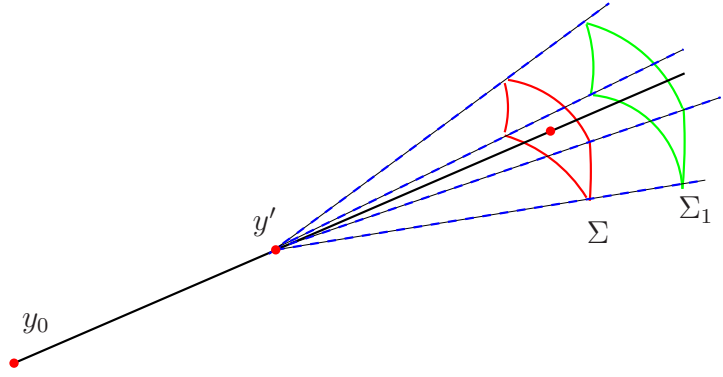}
\end{center}
\caption{A schematic figure where in the   3-dimensional Euclidean space  $\R^{3}$
we describe the route of the piece of the spherical wave that
propagates near the geodesic $\gamma_{x_0,\zeta_0}((0,\infty))$.
The geodesic is the black line in the figure. The spherical wave is the solution
$u_1$ that is singular on the surface $K(x_0,\zeta_0;t_0,s_0)\subset \R^{1+3}$ that is
 a subset of a light cone centered at $x^\prime=\gamma_{x_0,\zeta_0}(t_0)$. 
 When $P:\R^{1+3}\to \R^{3}$ is the projection to the space component, $P(t,y)=y$,
 the figure shows $P(K(x_1,\xi_1;t_0,s_0))$ and the points $y_0=P(x_0)$
 and $y^\prime=P(x^\prime)$. 
 The piece of spherical wave is sent from the surface $\Sigma=P(Y(x_0,\zeta_0;t_0,s_0))$,
 (the surface with red boundary) at time $T={\bf t}(\gamma_{x_0,\zeta_0}(2t_0))$.
 It starts to propagate and at a later time $T_1>T$ its singular support is the surface $\Sigma_1$
 shown in the figure in green.
%
}
 \end{figure}

We define the 3-submanifold
$K(x_0,\zeta_0;t_0,s_0)\subset M_0$  associated to  $(x_0,\zeta_0)\in L^{+}(M_0,\hat g)$, $x_0\in U_{\hat g}$ and  parameters $t_0,s_0\in \R_+ $
as 
\beq\label{associated submanifold}
& &K(x_0,\zeta_0;t_0,s_0)=\{\gamma_{x^{\prime},\eta}(t)\in M_0;\ \eta\in \W,\ t\in (0,\infty)\}.\hspace{-1cm}
\eeq
where  $(x^{\prime},\zeta^{\prime})=(\gamma_{x_0,\zeta_0}(t_0),\dot \gamma_{x_0,\zeta_0}(t_0))$
and $\W\subset L^{+}_{x^{\prime}}(M_0,\hat g)$ is a  neighborhood of $\zeta^{\prime}$
consisting of vectors $\eta\in L^+_{x^{\prime}}(M_0,\hat g)$ that satisfy
\ba
\left \|\eta -\zeta^{\prime}\right \|_{\hat g^+}<s_0,
\ea
where $\hat g^+$ is the Riemannian metric corresponding to the Lorentzian metric $\hat g$.
Note that $K(x_0,\zeta_0;t_0,s_0)\subset \L^+_{\hat g}(x^{\prime})\cup  \L^-_{\hat g}(x^{\prime})\cup \{x^\prime\}$ is a subset
of the light cone starting from $x^{\prime}=\gamma_{x_0,\xi_0}(t_0)$ and that it 
is singular
near the point $x^{\prime}$. Moreover, 
$\bigcap_{s_0>0}K(x_0,\zeta_0;t_0,s_0)=\gamma_{x_0,\zeta_0}((-t_0,\infty))\cap M_0$.
Let $S=\{x\in \hattuM _0;{\bf t}(x)={\bf t}(\gamma_{x_0,\zeta_0}(2t_0))\}$ be a Cauchy surface which intersects 
$\gamma_{x_0,\zeta_0}(\R)$ transversally at the point $\gamma_{x_0,\xi_0}(2t_0)$.
When $t_0>0$  is  small enough,
$Y(x_0,\zeta_0;t_0,s_0)=S\cap K(x_0,\zeta_0;t_0,s_0)$ is a smooth 2-surface that is
a subset of $U_{\hat g}$.

%
Let $\Lambda(x_0,\zeta_0;t_0,s_0)$ be the lagranginan manifold that is the  flowout from  
$N^*Y(x_0,\zeta_0;t_0,s_0)\cap  N^*K(x_0,\zeta_0;t_0,s_0)$ on Char$(\square_{\hat g})$ in the future direction.
%

Note that if $K_{s}$ is the smooth 3-dimensional manifold that is subset of $K(x_0,\zeta_0;t_0,s_0)$,
then $N^*K_s\subset \Lambda(x_0,\zeta_0;t_0,s_0)$.

\begin{lemma}\label{lem: lagrangian 1} 
Let $n\leq n_0=-15$ be an integer,  $t_0,s_0>0$, $Y=Y(x_0,\zeta_0;t_0,s_0)$,
$K=K(x_0,\zeta_0;t_0,s_0)$, and $\Lambda_1=\Lambda(x_0,\zeta_0;t_0,s_0)$.
Assume that ${\bf h}=({\bf h}_1,{\bf h}_2)\in \I^{n-3/2}(Y)$,
  is a $\B^L$-valued conormal
distribution  that is supported in  a neighborhood 
$V\subset \hattuM _0$ of  $\gamma_{x_0,\zeta_0}\cap Y=\{\gamma_{x_0,\zeta_0}(2t_0)\}$ and has  a $\R^{10+L}$-valued 
classical symbol.
Denote 
the principal symbol of ${\bf h}$ by 
$\tilde h(x,\xi)=(\tilde h_k(x,\xi))_{k=1}^{10+L}
$.
 Assume that the symbol of ${\bf h}$ vanishes near
the  light-like directions
in $N^*Y\setminus N^*K$.

Let $u^{(1)}=(\dot g,\dot\phi)$ be a solution of the 
linear wave equation (\ref{lin wave eq})
with the source ${\bf h}$. Then $u^{(1)}$,
considered as a vector valued lagrangian distribution on the set $M_0\setminus Y$,  satisfies
$
u^{(1)}\in \I^{n-1/2} ( M_0\setminus Y;\Lambda_1), 
$
and its principal symbol  
$\tilde a(y,\eta)=(\tilde a_{1}(y,\eta),\tilde a_{2}(y,\eta))$ at $(y,\eta)\in  \Lambda_1$ is given by 
\beq\label{eq: R propagation}\tilde a(y,\eta)=\sum_{k=1}^{10+L} R_j^k(y,\eta,x,\xi)\tilde h_{k}(x,\xi)
,\hspace{-1cm}
\eeq
where the pairs $(x,\xi)$ and $(y,\eta)$ are on the same bicharacteristics of $\square_{\hat g}$, \HOX{The arxiv submission v1 contained here a misprint: x and y were changed.}
and $x<y$, that is, $((x,\xi),(y,\eta))\in \Lambda_{\hat g}^\prime$, and in addition,
$(x,\xi)\in N^*Y$.
%
 Moreover, the matrix $(R_j^k(y,\eta,x,\xi))_{j,k=1}^{10+L}$
is invertible. 
\end{lemma}
We call the solution $u^{(1)}$   a piece of spherical wave that is associated
to the submanifold $K(x_0,\zeta_0;t_0,s_0)$. 
Below, we will consider interaction of spherical waves $u^{(1)}\in \I^{n-1/2} ( M_0\setminus Y;\Lambda_1)$.
The interaction terms involve several derivatives, and thus these solutions need to be
smooth enough. This is why we chose above $ n_0=-15$.

{\bf Proof.}\footnote{We note that Nils Dencker's results for the polarization sets, see
\cite{Dencker} are closely related to this result.} 
By \cite{MU1}, the parametrix  of the scalar wave equation  satisfies $(\square_{\hat g}+v(x,D))^{-1}\in
 \I^{-3/2,-1/2}(\Delta^\prime_{T^*M_0},\Lambda_{\hat g})$,
where $v(x,D)$ is a 1st order differential operator, $\Delta_{T^*M_0}$ is the conormal bundle of the diagonal of $M_0\times M_0$
and $\Lambda_{\hat g}$ is the flow-out of the canonical relation of $\square_{\hat g}$.

Let $x=(x^1,x^{\prime})\in \R^4$ denote local coordinates of $\R^4$
and $\p_1=\frac \p{\p x^1}$.
By 
\cite[Prop.\ 26.1.3]{H4} 
there are elliptic Fourier integral operators $\Phi_1$ and $\Phi_0$, of order $1$ and $0$,
having the same  canonical relation  so that
$\square_{\hat g}=\Phi_1\p_1\Phi_0^{-1}$.  
Thus an operator matrix
$A=(A_{jk})_{j,k=1}^{10+L}$ with
$A_{jk}=\square_{\hat g}\delta_{jk}+B_{jk}^p\nabla_p+C_{jk}$ can
be written as $A=\Phi_1\tilde A\Phi_0^{-1},$ 
$\tilde A=(\tilde A_{jk})_{j,k=1}^{10+L}$ with
$\tilde A_{jk}=\p_1\delta_{jk}+\tilde R_{jk},$ $\tilde R_{jk}=\Phi_0(B_{jk}^p\nabla_p+C_{jk})\Phi_1^{-1}$.

Furthermore, $\Phi_1$ and $\Phi_0$ have the same canonical
relation and we conclude that 
$\tilde R$ is zeroth order pseudodifferential operator.
Thus the parametrix for $A$ is
\ba
A^{-1}=\Phi_0(\p_1I+\tilde R)^{-1}\Phi_1^{-1},
\ea
where $(\p_1I+\tilde R)^{-1}\sim \sum_{j=0}^\infty
(\p_1I)^{-1}(\tilde R(\p_1I)^{-1})^j.$ 
Here $\sim$ denotes an asy\-mptotic sum.
This implies for also  the matrix valued wave operator, $\square_{\hat g}I+V(x,D)$,
when  $V(x,D)$ is the 1st order differential operator, that
 $(\square_{\hat g}I+V(x,D))^{-1}\in \I^{-3/2,-1/2}(\Delta^\prime_{T^*M_0},\Lambda_{\hat g})$.
By \cite[Prop.\ 2.1]{GU1}, this yields that $u^{(1)}\in \I^{n-1/2}(\Lambda_1)$\hiddenfootnote{
{Let $m=n+3/2$. Then we have ${\bf h}=({\bf h}_1,{\bf h}_2)\in \I^{m-1/2}(\Lambda(x_0,\zeta_0;t_0,s_0))$,
Thus $u^{(1)}=Q\bf h}\in \I^{m-1/2-3/2,-1/2}(\Lambda(x_0,\zeta_0;t_0,s_0),N^*K_1)^{10+L}$,
so that away from $Y(x_0,\zeta_0;t_0,s_0))$ we have
 $u^{(1)}\in \I^{m-1/2-3/2}(N^*K_1))^{10+L}$, or
$u^{(1)}\in \I^{m-1/2-3/2+1/2}(K_1))^{10+L}$
}
and the formula (\ref{eq: R propagation}) where
$R=(R_j^k(y,\eta,x,\xi))_{j,k=1}^{10+L}$ is related to the symbol
of $(\square_{\hat g}I+V(x,D))^{-1}$ on $\Lambda_{\hat g}$. Making a similar consideration
for the adjoint of the $(\square_{\hat g}I+V(x,D))^{-1}$, i.e., considering
the propagation of singularities using reversed causality, we see that
the matrix $R$ is invertible.

 \hfill \Box \medskip

\observation{{\bf Remark 3.2.}
Let us write the above principal symbol
$\tilde h(y,\eta)\in \R^{10+L}$ as
$\tilde h(y,\eta)=(\tilde h^{\prime}(y,\eta),\tilde h^{\prime\prime}(y,\eta))\in \R^{10}\times
\R^L$, where  $\tilde h^{\prime}(y,\eta)$ corresponds
to the principal symbol of the  $g$-component of the wave
and $\tilde h^{\prime\prime}(y,\eta)$ the  $\phi$-component of the wave.
Recall that due to (\ref{eq: wave operator in wave gauge}),
the equations $\square_g \phi_\ell+m\phi_\ell=0$
does not involve derivatives of $g$. 
This implies that the  
matrix $R(y,\eta,x,\xi)=[R_j^k(y,\eta,x,\xi)]_{j,k=1}^{L+10}$ in (\ref{eq: R propagation}) is such that
if $\tilde h(y,\eta)=(\tilde h^{\prime}(y,\eta),\tilde h^{\prime\prime}(y,\eta))$
is such that $\tilde h^{\prime\prime}(y,\eta)=0$, then
\beq
\label{new eq: R propagation}
\tilde a(y,\eta)=\left(\begin{array}{c}
\tilde a^{\prime}(y,\eta)\\
0
\end{array}\right)
=
R(y,\eta,x,\xi)\left(\begin{array}{c}
{\tilde h}^{\prime}(y,\eta) \\
\tilde h^{\prime\prime}(y,\eta)
\end{array}\right)
,\hspace{-1cm}
\eeq
that is, $\tilde a^{\prime\prime}(y,\eta)=0$. Roughly speaking,
this means that during the propagation along a bicharacteristic, 
the leading order singularities can
not move from the $g$-component of the wave to $\phi$-component.
\medskip
}

Let  $\frak S^n_{y,\eta}$ be the space of 
the elements in $\B^L_y$  satisfying
the divergence condition for the symbols (\ref{divergence condition for symbol})
at $({y,\eta})$.

\begin{lemma}\label{lem: lagrangian} 
Let $n\leq n_0$, $t_0,s_0>0$, $Y=Y(x_0,\zeta_0;t_0,s_0)$,
$K=K(x_0,\zeta_0;t_0,s_0)$, and $\Lambda_1=\Lambda(x_0,\zeta_0;t_0,s_0)$.
Let us consider
 ${\bf p},{\bf q^{\prime}},{\bf z}\in \I^{n-3/2}(Y)$ 
 that  have classical symbols 
 with principal symbols $\tilde c^{(a)}(x,\xi),$ $\tilde d_{1}(x,\xi)$,
 $\tilde d_2^{(a)}(x,\xi)$, correspondingly, at $(x,\xi)\in N^*Y$.
 Moreover, assume that the principal symbols of ${\bf p}$ and ${\bf z}$ 
 satisfy the linearized conservation law 
 for the principal symbols, that is, (\ref{eq: lineariz. conse. law symbols}),  at  a light-like co-vector  $(x,\xi)\in N^*Y$.
 
 Let $(y,\eta)\in T^*M_0$, $y\not \in Y$
be a light-like co-vector such that $(x,\xi)\in\Theta_{y,\eta}\cap  N^*Y$.
Then  the principal symbol $\tilde a(y,\eta)=(\tilde a_1(y,\eta),\tilde a_2(y,\eta))$ 
of the wave $u^{(1)}\in \I^{n-1/2} ( M_0\setminus Y;\Lambda_1)$, defined in Lemma \ref{lem: lagrangian 1},  varies in  $\frak S_{y,\eta}$. 
{
Moreover,  by varying ${\bf p},{\bf q^{\prime}},{\bf z}$  so that the linearized conservation law  (\ref{eq: lineariz. conse. law symbols})  for principal symbols 
is satisfied, the principal symbol  $\tilde a(y,\eta)$ at $(y,\eta)$
 achieves all values in the $(L+6)$ dimensional space $\frak S_{y,\eta}$.}
\end{lemma}

 Below, we denote ${\bf f}\in  {\mathcal I_{S}}^{n-3/2}
(Y(x_0,\zeta_0;t_0,s_0))$ when
\beq\label{eq: Z-sources}
{\bf f}=({\bf p},{\bf q}),\quad {\bf q}=({\bf q^{\prime}},{\bf z})
\eeq
and ${\bf p},{\bf q^{\prime}},{\bf z}\in \I^{n-3/2}
(Y(x_0,\zeta_0;t_0,s_0)) $  and
 the principal symbols of ${\bf p}$ and ${\bf z}$ 
 satisfy the linearized conservation law 
 for principal symbols, that is, equation (\ref{eq: lineariz. conse. law symbols}).

\noindent 
{\bf Proof.} 
{
Let us use local coordinates $X:V\to \R^4$ where $V\subset M_0$ is a neighborhood of $x$.
In these coordinates, let $\tilde c^{(b)}(x,\xi)$ and $\tilde d_{2}^{(b)}(x,\xi)$ denote the sub-principal
symbols of ${\bf p}$ and ${\bf z}$, respectively, at $(x,\xi)\in N^*Y$. 
{
 Moreover, let $\tilde c^{(c)}_j(x,\xi)=\frac \p {\p x^j}\tilde c^{(a)}_j(x,\xi)$ and $\tilde d^{(c)}_j(x,\xi)=\frac \p {\p x^j}\tilde d^{(a)}_2(x,\xi)$, $j=1,2,3,4$ be the $x$-derivatives of the principal symbols
and let us denote $\tilde c^{(c)}(x,\xi)=(\tilde c^{(c)}_j(x,\xi))_{j=1}^4$ and $\tilde d^{(c)}(x,\xi)=(\tilde d^{(c)}_{j}(x,\xi))_{j=1}^4$.}


 Let ${\bf h}=({\bf h}_1,{\bf h}_2)=H(x;{\bf p},{\bf q})$ 
 be defined by
 (\ref{lin wave eq source B}), where we recall that $M_{(1)}=0$, $N^j_{(1)}=0$, and $L_{(1)}=\hat g$.
 Then ${\bf h}\in \I^{n-3/2}(Y)$ has the principal symbol
  $\tilde h(x,\xi)=(\tilde h_{1}(x,\xi),\tilde h_{2}(x,\xi))$ at $(x,\xi)$, given by  
   \beq\nonumber 
\left(\begin{array}{c}
\tilde  h_{1}(x,\xi)\\
 \tilde  h_{2}(x,\xi)\end{array}\right)\hspace{-2mm}&=&\hspace{-2mm}\left(\begin{array}{c}
(\tilde c^{(a)}+\hat g\tilde d_2^{(a)})(x,\xi)\\
  K_{(2)}(x,\xi)(\tilde c^{(a)}+\hat g\tilde d_2^{(a)})(x,\xi)
 \end{array}\right)+ \left(\begin{array}{c}
0 
\\  
 M_{(2)}(x)\tilde d_{1}(x,\xi)  \end{array}\right) \hspace{-5mm}
 \\    \label{lin wave eq source symbols} 
  & &\hspace{-20mm}+ 
 \left(\begin{array}{c}
0\\
  J_{(2)}(x,\xi)(\tilde c^{(c)}+\hat g\tilde d^{(c)})(x,\xi)
 \end{array}\right)
  +\left(\begin{array}{c}
0 
\\  
L_{(2)}(x)\tilde d_2^{(a)}(x,\xi)  \end{array}\right)
 +
\\
\nonumber
  & &\hspace{-20mm}+
     \left(\begin{array}{c}
0
\\ 
N^j_{(2)}(x) \, \hat g^{lk}\xi_l (\tilde c_{1}^{(b)}+\hat g\tilde d_2^{(b)})_{jk}(x,\xi)
\end{array}\right),\hspace{-5mm}
%
  \eeq 
  where $ K_{(2)}(x,\xi)(\tilde c^{(a)}+\hat g\tilde d_2^{(a)})(x,\xi)$ 
  and  $ J_{(2)}(x,\xi)(\tilde c^{(c)}+\hat g\tilde d_2^{(c)})(x,\xi)$  is related to the sub-principal
  symbol of the term $\hat g^{lk} \nabla_l^{\hat g} ({\bf p}_{jk}+{\bf z}\hat g_{jk})$. 
  Observe that  the map $(c^{(b)}_{jk})\mapsto (\hat g^{lk}\xi_lc^{(b)}_{jk})_{j=1}^4$, 
 defined as
   $\hbox{Symm}(\R^{4\times 4})\to \R^4$, is surjective.
  Denote 
  $\tilde m^{(a)}=(\tilde c^{(a)}+\hat g\tilde d_2^{(a)})(x,\xi)$,
  $\tilde m^{(b)}=(\tilde c^{(b)}+\hat g\tilde d_2^{(b)})(x,\xi)$,
  and $\tilde m^{(c)}=(\tilde c^{(c)}+\hat g\tilde d^{(c)})(x,\xi)$.
 As noted above,
by Assumption S (iii),   the union of the image spaces of  the matrices 
$M_{(2)}(x)$ and
$L_{(2)}(x) $,  
and $N^j_{(2)}(x) $,  $j=1,2,3,4$,
span the space $\R^L$  for all $x\in \hat U$.  
%
Hence the map 
  \ba
  {\bf A}:(\tilde m^{(a)},\tilde m^{(b)},\tilde m^{(c)},
  \tilde d_{1}, \tilde d_2^{(a)})|_{(x,\xi)}\mapsto
  (\tilde h_1(x,\xi),\tilde h_2(x,\xi)),
  \ea 
  given by (\ref{lin wave eq source symbols}), considered as a 
  map ${\bf A}:Y=(\hbox{Symm}(\R^{4\times 4}))^{1+1+4}\times \R^{K}\times \R\to
  \hbox{Symm}(\R^{4\times 4})\times \R^L$,
  is surjective.
 Let $\mathcal X$ be  the set of
  elements $(\tilde m^{(a)}(x,\xi),\tilde m^{(b)}(x,\xi),\tilde m^{(c)}(x,\xi),
  \tilde d_{1}(x,\xi),$ $\tilde d_2^{(a)}(x,\xi))\in Y$
  where  
   $\tilde m^{(a)}(x,\xi)=(\tilde c^{(a)}+\hat g\tilde d_2^{(a)})(x,\xi)$ is such that the  pair  $(\tilde c^{(a)}(x,\xi),\tilde d_2^{(a)}(x,\xi))$  satisfies
 the linearized conservation law for principal symbols, see (\ref{eq: lineariz. conse. law symbols}). 
 Then    $\mathcal X$ has codimension 4 in $Y$, we see
 that the image  ${\bf A}(\mathcal X)$ has co-dimension less or equal
 to 4, that is, it has at least dimension $(L+6)$.}
%

%
%
%
 
Let $({\bf h}_1,{\bf h}_2)$ be a source with 
the principal symbol  $\tilde h(x,\xi)$ that corresponds
to  functions $({\bf p},{\bf q})$.
When $(P,Q)=(\e {\bf p},\e{\bf q})$ and $\e>0$
is small enough, it
 follows from Assumption $S$  (iv), that the $\hat g$-reduced Einstein
 equations (\ref{eq: adaptive model}) have a solution $(g_\e,\phi_\e)$
that satisfies the conservation law (\ref{conservation law0}) and thus $g_\e$  satisfies
the  harmonicity condition (\ref{harmonic condition}). Hence $\dot g=\p_\e g_\e|_{\e=0}$ satisfies
the linearized divergence condition (\ref{divergence condition}). 
Observe that the metric component of the solution $u^{(1)}$ of the linearized
Einstein equation corresponding to the functions $({\bf p},{\bf q})$
is equal to  $\dot g$.
Thus the principal symbol $\tilde a_1(y,\eta)$ 
of $u^{(1)}$ at $(y,\eta)$ satisfies 
the divergence condition for the symbols (\ref{divergence condition for symbol}),
that is, four linear conditions, and also satisfies  $\tilde a_1(y,\eta)\in \frak S_{y,\eta}$.

By Lemma \ref{lem: lagrangian 1}  the  $R:\tilde h(x,\xi)\mapsto \tilde a(y,\eta)$
maps a subspace of $\R^{10+L}$ whose has codimension $4$
onto some subspace of $\R^{10+L}$ whose codimension is $4$.

{
The above imply that the map $R\circ {\bf A}$, that maps the principal and sub-principal symbols
$(\tilde m^{(a)}(x,\xi),\tilde m^{(b)}(x,\xi),\tilde m^{(c)}(x,\xi),
  \tilde d_{1}(x,\xi), \tilde d_2^{(a)}(x,\xi))$ 
   to the
 principal symbol $\tilde a(y,\eta)$
of the solution $u^{(1)}$ at $(y,\eta)$, is such
that $R\circ {\bf A}$ maps $\mathcal X$ onto a subspace  of $\R^{10+L}$ which
codimension is at most 4 and is a subspace of $\frak S_{y,\eta}$. 
As $\frak S_{y,\eta}$ has codimension 4, we see that
$R\circ {\bf A}(\mathcal X)$ coincides with  $\frak S_{y,\eta}$  and has thus
dimension $(L+6)$.}
%
\hfill \Box \medskip

 \subsection{Interaction of non-linear waves}
 
 Next we consider interaction of four $C^k$-smooth 
 waves on a $(1+3)$ dimensional manifold having conormal singularities, where $k\in \Z_+$ is sufficiently large. 
 Interaction of such waves produces a "corner point" in the space
 time. On the related microlocal tools to consider scattering by corners, see \cite{Vasy1,Vasy2}.
 Earlier related interaction of three waves 
 has been studied by Melrose and Ritter \cite{MR1,MR2}
 and Rauch and Reed,  \cite{R-R}  for a non-linear hyperbolic equations in $\R^{1+2}$
 where the non-linearity appears in the lower order terms. 
   
\subsubsection{Interaction of non-linear waves on a general manifold}

Next, we introduce a vector of four $\e$ variables
denoted by $\vec \e=(\e_1,\e_2,\e_3,\e_4)\in \R_+^4$. 
Let $s_0,t_0>0$ and consider the solution $u_{\vec \e}=(g_{\vec\e}-\hat g,\phi _{\vec\e}-\hat \phi)$
where $g_{\vec \e}$ and $\phi_{\vec \e}$ solve the equations (\ref{eq: adaptive model})
with $(P,Q)={\bf f}_{\vec \e}$ being 
\beq\label{eq: f vec e sources}
 {\bf f}_{\vec \e}:=\sum_{j=1}^4\e_j {\bf f}_{j},\quad
 {\bf f}_{j}\in {\mathcal I_{S}}^{n-3/2} (Y(x_j,\zeta_j;t_0,s_0)),
 \eeq
 see (\ref{eq: Z-sources}), where  $(x_j,\zeta_j)$ are light-like
vectors with $x_j\in \hat U.$
Moreover,  we assume that for some $0<r_1<r_0$  and $s_-+r_1<s_1<s_+$ the sources satisfy
 \beq\label{eq: source causality condition}
& & \supp({\bf f}_{j})\cap J ^+_{\hat g}(\supp({\bf f}_{k}))=\emptyset,\quad\hbox{for all }j\not=k,\\
\nonumber
& &\supp({\bf f}_{j})\subset  I_{\hat g}( \mu_{\hat g}(s_1-r_1),\mu_{\hat g}(s_1)),\quad \hbox{for all }j=1,2,3,4,
%
 \eeq
 where $r_0$ is the parameter introduced after (\ref{observer neighborhood with hat})
 to define $W_{\hat g}=W_{\hat g}(r_0)$.
  The first condition implies that the supports of the sources are causally independent.

The sources   $ {\bf f}_{j}$ give raise for  $\B^L$-section valued solutions
of the linearized wave equations, denoted
\ba
u_j:=u^{(1)}_{j}=QR \,{\bf f}_{j}\in \I( \Lambda(x_0^{(j)},\zeta_0^{(j)};t_0,s_0)),
\ea
where $R$ is a first order differential operator depending on $\hat g$
and
   ${\bf Q}={\bf Q}_{\hat g}=(\square_{\hat g} +V(x,D))^{-1}$ is  the causal inverse of
the wave equation where $V(x,D)$ is a first order differential operator.

  In the following we use the notations
\ba
& &\p_{\vec \e}^1u_{\vec \e}|_{\vec \e=0}:=\p_{\e_1}u_{\vec \e}|_{\vec \e=0},\\
& &
\p_{\vec \e}^2 u_{\vec \e}|_{\vec \e=0}:=\p_{\e_1}\p_{\e_2}u_{\vec \e}|_{\vec \e=0},\\
& &
\p_{\vec \e}^3 u_{\vec \e}|_{{\vec \e}=0}:=\p_{\e_1}\p_{\e_2}\p_{\e_3} u_{\vec \e}|_{{\vec \e}=0}
,\\
& &
\p_{\vec \e}^4 u_{\vec \e}|_{\vec \e=0}:=\p_{\e_1}\p_{\e_2}\p_{\e_3}\p_{\e_4} u_{\vec \e}|_{{\vec \e}=0}.
\ea

\begin{figure}[htbp] \label{Fig-8}
\begin{center}

\psfrag{1}{\hspace{-3mm}$x_1$}
\psfrag{2}{\hspace{-3mm}$x_2$}
\psfrag{3}{$q$}
\includegraphics[width=5.5cm]{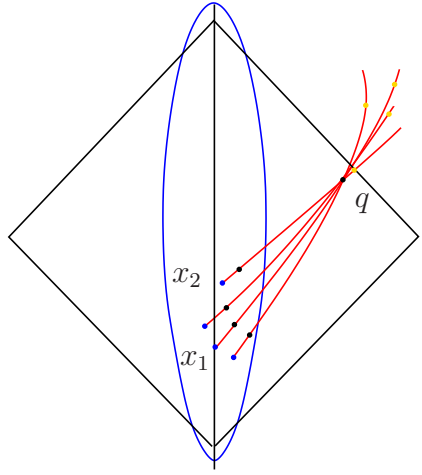}
\end{center}
\caption{A schematic figure where the space-time is represented as the  2-dimensional set $\R^{1+1}$. 
The four light-like geodesics $\gamma_{x_j,\xi_j}([0,\infty))$, $j=1,2,3,4$ starting at
the blue points $x_j$
intersect at $q$ before the first cut points of $\gamma_{x_j,\xi_j}([t_0,\infty))$,
denoted by the golden points. We send pieces
of spherical waves from the black points $\gamma_{x_j,\xi_j}(t_0)$ that propagate near these geodesics
and produce an artificial point source at the point $q$.}
 \end{figure}

Next we denote the waves produced by the $\ell$-th order interaction by
\ba
\M^{(\ell)}:=\p_{\vec \e}^\ell u_{\vec \e}|_{\vec \e=0},\quad \ell\in \{1,2,3,4\}
\ea
in $\hat U$, see Fig.\ 8. 
Below,  
$\B_j^{\beta}$, $j,\beta \in \Z_+$ are differential operators having in local coordinates the form 
\beq\label{extra notation1}
& &\hbox{$\B^\beta _j:(v_{p})_{p=1}^{10+L}\mapsto (b^{r,(j,\beta )}_{p}(x)\p_x^{\vec k(\beta,j)} v_{r}(x))_{p=1}^{10+L}$, and}\\ \nonumber
& &\hbox{$S^\beta _j={\bf Q}$ or $S^\beta _j=I$},
\eeq 
 where the coefficients $b^{r,(j,\beta )}_{p}(x)$
depend on the derivatives of $\hat g_{jk}$.

 Computing $\e_j$ derivatives of the equations (\ref{eq: mixed source terms})
with the sources ${\bf f}_{\vec \e}$, and taking into account
 the condition (\ref{eq: source causality condition}), we obtain 
  \beq\label{w1-3 solutions}
& & \hspace{-2cm}\M^{(1)}=u_1,\\ \nonumber
  & &\hspace{-2cm}\M^{(2)}=\sum_{\sigma\in \Sigma(2)}\sum_{\beta\in J_2 }{\bf Q}(\B^{\beta}_2u_{\sigma(2)}\,\cdotp \B^{\beta}_1u_{\sigma(1)}),\hspace{-2cm}\\ \nonumber
 & &\hspace{-2cm}\M^{(3)}=\sum_{\sigma\in \Sigma(3)}\sum_{\beta\in J_3}{\bf Q}(\B_3^{\beta}u_{\sigma(3)}\,\cdotp \mathcal C_1^{\beta} S^{\beta}_1(\B^{\beta}_2u_{\sigma(2)}\,\cdotp \B^{\beta}_1u_{\sigma(1)})),\hspace{-2cm} \nonumber
 \\ \nonumber
& &\hspace{-2cm}\M^{(4)}=
{\bf Q}\F^{(4)},\quad  \F^{(4)}=
\sum_{\sigma\in \Sigma(4)}\sum_{\beta\in J_4}
(\mathcal G^{(4),\beta}_\sigma+\tilde {\mathcal G}^{(4),\beta}_\sigma),
\eeq 
where \HOX{We have to explain  $n_2,n_3,n_4$ that denote bounds for index $\beta$.}
$\Sigma(\ell)$  is the set of permutations, that is, bijections
 $\sigma:\{1,2,\dots,\ell\}\to \{1,2,\dots,\ell\}$, and $J_2,J_3,J_4\subset \Z_+$ are finite sets.  Moreover,
the operators $\mathcal C_j^{\beta}$ are of the same form 
as the operators $\mathcal B_j^{\beta}$. 
The orders
$k_j^\beta=ord(\B^\beta_j)$ 
of the differential operators $\B^\beta_j$  and
 the orders 
$\ell_j^\beta=ord(\mathcal C^\beta_j)$ 
satisfy  $k_1^\beta+k_2^\beta\leq 2$ and 
 $k_3^\beta+\ell_1^\beta\leq 2$.
Moreover, 
\beq\label{M4 terms}
& & {\mathcal G}^{(4),\beta}_\sigma=\B_4^{\beta}u_{\sigma(4)} \,\cdotp \mathcal C_2^{\beta}S_2^{\beta}(\B_3^{\beta}u_{\sigma(3)}\,\cdotp \mathcal C_1^{\beta} S^{\beta}_1(\B^{\beta}_2u_{\sigma(2)}\,\cdotp \B^{\beta}_1u_{\sigma(1)}))\hspace{-2cm}
\\ & &{\M}^{(4),\beta}_\sigma={\bf Q}{\mathcal G}^{(4),\beta}_\sigma,
\nonumber 
\eeq
where  the orders of the differential operators
satisfy $k_4^\beta+\ell_2^\beta\leq 2$,  $k_3^\beta+\ell_1^\beta\leq 2$, 
 $k_2^\beta+k_1^\beta\leq 2$, 
and finally
\beq\label{tilde M4 terms}
& &\tilde  {\mathcal G}^{(4),\beta}_\sigma=\mathcal C_2^{\beta}S_2^{\beta} (\B_4^{\beta}u_{\sigma(4)} \,\cdotp \B_3^{\beta}u_{\sigma(3)})\,\cdotp \mathcal C_1^{\beta}S^{\beta}_1(\B^{\beta}_2u_{\sigma(2)}\,\cdotp \B^{\beta}_1u_{\sigma(1)}),\hspace{-2cm}
\\ & &\tilde {\M}^{(4),\beta}_\sigma={\bf Q}\tilde {\mathcal G}^{(4),\beta}_\sigma,
\nonumber 
\eeq
where $\ell_1^\beta+\ell_2^\beta\leq 2$,  $k_4^\beta+k_3^\beta\leq 2$, 
 $k_2^\beta+k_1^\beta\leq 2$.
 Note that due to the conditions (\ref{eq: source causality condition}), for  $\ell=2,3,4$,
 $\M^{(\ell),\beta}_\sigma$ and $\tilde \M^{(\ell),\beta}_\sigma$
  do not contain
terms that involve sources $f_j$, for example,  using (\ref{eq: mixed source terms}) we obtain the formula
\ba
\M^{(2)}={\bf Q}\bigg(B_1(u_1)f_2+B_2(u_2)f_1+K(u_1,u_2)\bigg),
\ea
where the terms $B_1(u_1)f_2$ and $B_2(u_2)f_2$ vanish due to  (\ref{eq: source causality condition}) and $K$ is a bilinear operator.  

We make also the observation that when $\vec S_\beta=(  {\bf Q},  {\bf Q})$,
the terms  $  \M^{(4),\beta}_\sigma$ and $\tilde  \M^{(4),\beta}_\sigma$
can be written in the form
\beq\label{eq: tilde M1}
& & \M^{(4),\beta}_\sigma={\bf Q}(A[  u_{\sigma(4)}, {\bf Q}(A[ u_{\sigma(3)},{\bf Q}(A[ u_{\sigma(2)}, u_{\sigma(1)}])]),
\\ \label{eq: tilde M2}
& &\tilde  \M^{(4),\beta}_\sigma={\bf Q}(A[ {\bf Q}(A[ u_{\sigma(4)}, u_{\sigma(3)}]),{\bf Q}(A[ u_{\sigma(2)}, u_{\sigma(1)}])])
\eeq
where 
$A[V,W]$ is a generic notation (i.e., its exact form can vary even inside the formula) for a 2nd  order multilinear operator in $V$ and $W$ having the form
\beq\label{A form}
A[V,W]=\sum_{|\a|+|\gamma|\leq 2}a_{\a\gamma}(x)(\p_x^\a V(x))\,\cdotp(\p_x^\gamma W(x)).
\eeq
We use in particular two such bilinear forms that are given
for $V=(v_{jk},\phi)$, 
and  $W=( w^{jk},\phi^{\prime})$,
by
\beq\label{A alpha decomposition2}
& & A_1[V,W]=-\hat g^{jb} v_{ab}\hat g^{ak}\p_j\p_k  w_{pq},\quad A_2[V,W]=A_1[W,V].
\eeq

By considering the terms that we obtain by substituting formulas (\ref{w1-w3}) in   (\ref{w4}),
%
we see that 
the case when the quadratic forms $A$ appearing
 in the formulas (\ref{eq: tilde M1}) and (\ref{eq: tilde M2}) 
 have second order derivatives, that is,
when the terms  with
 $|\gamma|=2$ or $|\a|= 2$ in (\ref{A form})   are non-zero, can happen only when 
$A$ is either the bilinear form $A_1$ or $A_2$. 

%
%




%


%

\subsubsection{On the singular support of the non-linear interaction  of three waves}
\label{subsection: sing supp interactions}


Below, we use  for a pair $(x,\xi)\in L^+\hattuM _0$  the notation 
\beq\label{eq: x(h) notation}
& &(x(h),\xi(h))=(\gamma_{x,\xi}(h),\dot \gamma_{x,\xi}(h)).
\eeq

Let us next consider four light-like future pointing directions $(x_j,\xi_j)$, $j=1,2,3,4$, and
use below the notation
\ba
(\vec x,\vec\xi)=((x_j,\xi_j))_{j=1}^4.
\ea
We will consider the case when we send pieces of spherical waves
propagating on surfaces $K(x_j,\xi_j;t_0,s_0)$,  $t_0,s_0>0$, cf.\ (\ref{associated submanifold}), and
these waves interact. Let us use for $(x_j,\xi_j)$ the notation (\ref{eq: x(h) notation})
and assume that $x_k(t_0)\not\in K(x_j,\xi_j;t_0,s_0)$ for $k\not =j$, see Fig.\ 9. 

Next we consider the 3-interactions of the waves, see Fig.\ 9 on analogous considerations
in $(1+2)$ dimensional Lorentz space. 
For $1<j_1<j_2<j_3\leq 4$,
let $K_{j_p}=K(x_{j_p},\xi_{j_p};t_0,s_0)$, $p=1,2,3$. 
We define 
 \beq\label{eq: 3-interactions}
\X(j_1,j_2,j_3;t_0,s_0)=\hspace{-4mm}\bigcup_{z\in K_{j_1}\cap K_{j_2}\cap K_{j_3}}\hspace{-4mm}(
N_zK_{j_1}
+N_zK_{j_2}+ N_zK_{j_3})\cap
L^+_z \hattuM _0.\hspace{-1.3cm}
\eeq
Note that $K_{123}=K_1\cap K_2\cap K_3$ is a space like curve and
$N^*_zK_{123}=N^*_zK_1+N^*_zK_2+N^*_z K_3$.  
Moreover, we define ${\mathcal Y}(j_1,j_2,j_3;t_0,s_0)$ to be the set of  all  $y\in M_0$
such that there are $z\in K_{j_1}\cap K_{j_2}\cap K_{j_3}$,
$\zeta\in \X (j_1,j_2,j_3;t_0,s_0)$, and
 $t\geq 0$ such that $\gamma_{z,\zeta}(t)=y$. \HOX{We should change all $\Y$ to $\K$ etc.}
We use below  the 3-interaction sets (See Figs.\ 10 and 11) 
\beq\label{K and X sets}
& &{\mathcal Y}((\vec x,\vec\xi);t_0,s_0)=\bigcup_{1\leq j_1<j_2<j_3\leq 4}{\mathcal Y}(j_1,j_2,j_3;t_0,s_0),\\
& &{\mathcal Y}((\vec x,\vec\xi);t_0)=\bigcap_{s_0>0}{\mathcal Y}((\vec x,\vec\xi);t_0,s_0)\subset M_0,\nonumber \\
& &\X((\vec x,\vec\xi);t_0,s_0)=\bigcup_{1\leq j_1<j_2<j_3\leq 4}\X(j_1,j_2,j_3;t_0,s_0),\nonumber\\
& &\X((\vec x,\vec\xi);t_0)=\bigcap_{s_0>0}\X((\vec x,\vec\xi);t_0,s_0)\subset TM_0.\nonumber
\eeq
For instance in  Minkowski space, when three plane waves (which singular supports 
are hyper-planes) collide,
the intersections of the hyperplanes is a 1-dimensional space-like line $K_{123}$ in the 4-dimensional space-time.
This corresponds to a point  moving continuously in time. Roughly speaking, the point 
seem to move at a higher speed than light (i.e.\ it appears like a tachyonic
point-like object)
and  produces a shock wave type of singularity that moves
on the set ${\mathcal Y}((\vec x,\vec\xi);t_0,s_0)$ in the space-time.

\begin{figure}[htbp] \label{Fig-9}
\begin{center}

\psfrag{1}{$x^0$}
\psfrag{2}{$x^1$}
\psfrag{3}{$x^2$}
\includegraphics[width=7.5cm]{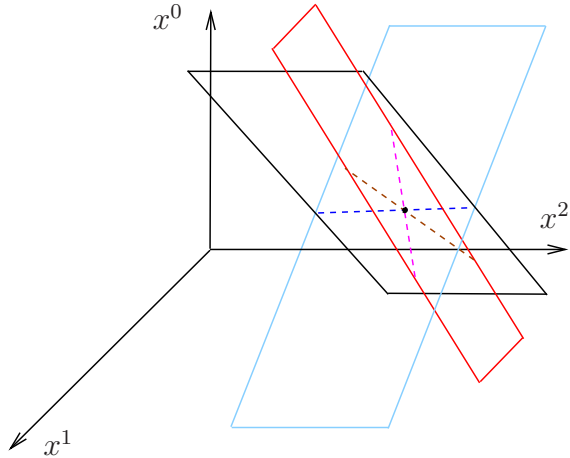}
\end{center}
\caption{A schematic figure where the space-time is represented as the  3-dimensional set $\R^{2+1}$. 
In the figure 3 pieces or plane waves  have singularities 
on strips of hyperplanes (in fact planes) $K_1,K_2,K_3$, colored
by light blue, red, and black. These planes  have intersections, and
in the figure the sets $K_{12}=K_1\cap K_2$, $K_{23}=K_2\cap K_3$,
and $K_{13}=K_1\cap K_3$ are shown as dashed lines with
dark blue, magenta, and brown colors. These dashed lines  intersect
at a point $\{q\}=K_{123}=K_1\cap K_2\cap K_3$.}
 \end{figure}

\begin{figure}[htbp]  \label{Fig-10}
\begin{center}

\psfrag{1}{$\pi(\mathcal U_{z_0,\eta_0})$}
\psfrag{2}{$U_{\hat g}$}
\psfrag{3}{$J_{\hat g}(\hat p^-,\hat p^+)$}
\psfrag{4}{$W_{\hat g}$}
\psfrag{5}{$\mu_{z,\eta}$}
\includegraphics[width=7.5cm]{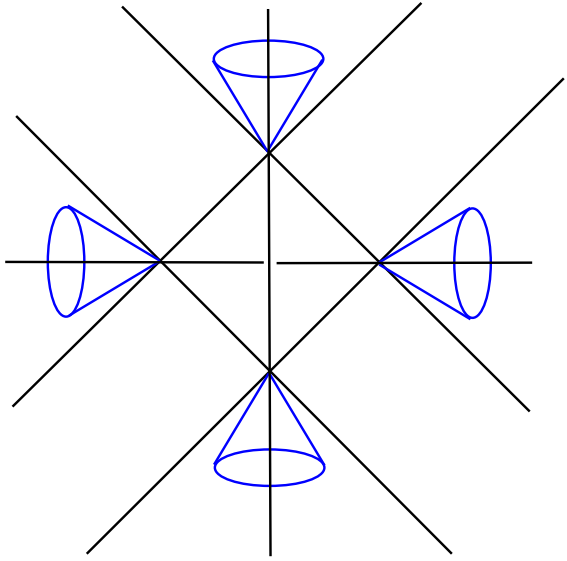}
\end{center}
\caption{In the section \ref{subsection: sing supp interactions} we consider four colliding pieces of spherical waves.
In the figure we consider the Minkowski space $\R^{1+3}$ and the figure corresponds 
to a "time-slice" $\{T_1\}\times \R^3$.
We assume that the spherical waves are sent 
so far away that the waves look like pieces of plane waves. The plane waves $u_j\in \mathcal I(K_j)$, $j=1,2,3,4$,
are conormal distributions that are  solutions of the linear wave equation and their singular supports are the sets $K_j$,
that are pieces of 3-dimensional planes in the space-time. The sets $K_j$ are not shown in the figure.
The 2-wave interaction wave $\mathcal M^{(2)}$ is singular on the set $\cup_{j\not =k}K_j\cap K_j$.
There are 6 intersection sets $K_j\cap K_j$ that are shown as black line segments. Note that these
lines have 4 intersection points, that is, the vertical and the horizontal black
lines do not intersect  in $\{T_1\}\times \R^3$.
The four intersection points of the black lines are the sets $(\{T_1\}\times \R^3)\cap (K_j\cap K_j\cap K_n)$.
These points correspond to  points moving in time (i.e., they are curves in the space-time) that produce singularities of the 3-interaction wave 
$\mathcal M^{(3)}$. The points
 seem to move faster than the speed of the light (similarly, as  a shadow of a
far away  object may seem to  move faster than the speed of the light).  Such
point sources produce "shock waves", and due to this, $\mathcal M^{(3)}$ is
singular on the sets $\mathcal Y((\vec x,\vec \xi),t_0,s_0)$ defined in 
formulas (65)-(66). The set $(\{T_1\}\times \R^3)\cap \mathcal Y((\vec x,\vec \xi),t_0,s_0)$
is the union of the four blue cones shown in the figure. 
}
 \end{figure}

%
%
%
%

\begin{figure}[htbp] \label{Fig-12} 
\begin{center}

\psfrag{1}{$\pi(\mathcal U_{z_0,\eta_0})$}
\psfrag{2}{$U_{\hat g}$}
\psfrag{3}{$J_{\hat g}(\hat p^-,\hat p^+)$}
\psfrag{4}{$W_{\hat g}$}
\psfrag{5}{$\mu_{z,\eta}$}
\includegraphics[width=7.5cm]{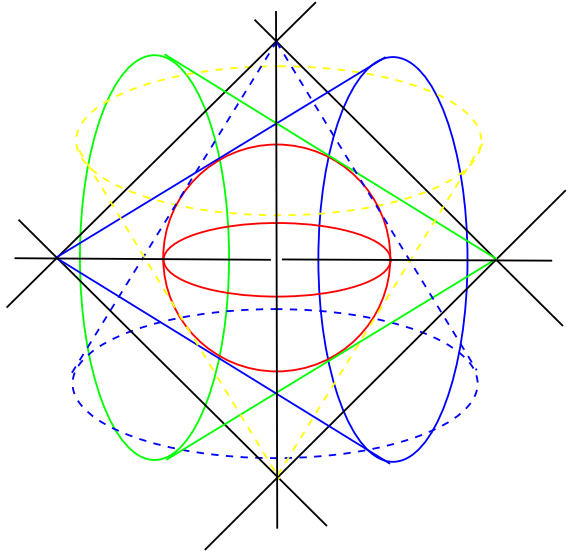}
\end{center}
\caption{The same situation that was described in Fig.\ 9 
is shown
at a later time, that is, the figure shows the time-slice $\{T_2\}\times \R^3$ with
$T_2>T_1$, when the parameter $s_0$ is quite large. The four  pieces of the spherical waves have now collided and
produced a point source in the space-time at a point $q\in K_1\cap K_2\cap K_3\cap K_4$,
that produces the singularities of the 4-interaction wave $\mathcal M^{(4)}$.
In the figure $T_1<t<T_2$, where $q$ has the time coordinate $t$.
The four cones in the figure, shown with solid blue and green curves
and dashed blue and yellow curves are   the intersection of  the time-slice $\{T_2\}\times \R^3$ 
and the set $\mathcal Y((\vec x,\vec \xi),t_0,s_0)$. Inside the cones the red sphere
is the set $\mathcal L^+(q)\cap (\{T_2\}\times \R^3)$ that
corresponds to the spherical wave produced by the point source at $q$.
}
 \end{figure}

\begin{figure}[htbp]\label{Fig-13} 
\begin{center}

\psfrag{1}{$\pi(\mathcal U_{z_0,\eta_0})$}
\psfrag{2}{$U_{\hat g}$}
\psfrag{3}{$J_{\hat g}(\hat p^-,\hat p^+)$}
\psfrag{4}{$W_{\hat g}$}
\psfrag{5}{$\mu_{z,\eta}$}
\includegraphics[width=5.5cm]{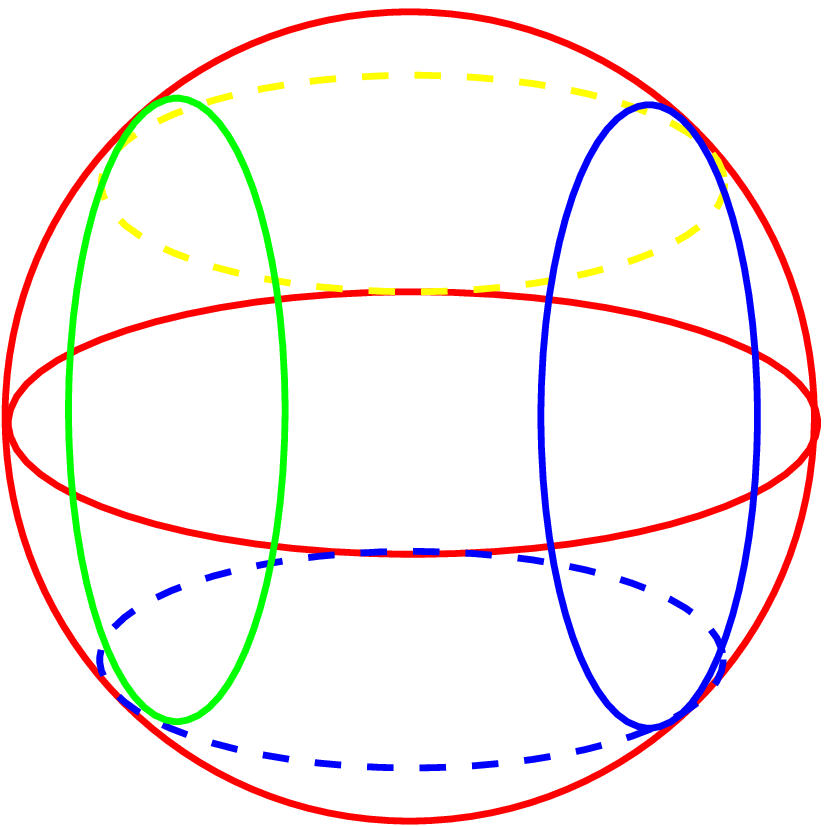}
\end{center}
\caption{The same situation that was described in Fig.\ 11, 
that is,  the figure shows in 
the time-slice $\{T_2\}\times \R^3$ the singularities produced by four colliding
spherical waves, when the parameter $s_0$ is very small. In this case the truncated cones
degenerate to circles.
}
 \end{figure}

\subsubsection{Gaussian beams}

Our aim is to consider interactions of  4 waves to produce  a new source, and
to this end we use test sources that produce gaussian beams.

Let $y\in \hat U$ and  $\eta\in T_{y}M$ be a future pointing light-like vector. 
We choose a complex function $p\in C^\infty(\hattuM _0)$
such that
 $\im  p(x) \geq 0$ and $\im  p(x)$ vanishes
only at $y$, 
$ p(y)=0$, $d(\re p)|_{y}=\eta^\sharp$, 
$d(\im  p)|_y=0$ and
the Hessian of $\im  p$ at $y$ is positive definite.
To simplify notations, we use below also complex sources and waves.
 The physical linearized waves can be obtained by taking the real part of the 
corresponding complex wave. We use a large parameter $\tau$ and define a  test source 
\beq\label{Ftau source}
F_\tau(x)=\tau^{-1}
\exp (  i\tau  p(x))
 h(x) 
\eeq
where $h$ is section on $\hbox{sym}(\Omega^2M)\times \R^L$ supported
in a small neighborhood $W$ of $y$. The construction of $p(x)$ and $F_\tau$ is discussed later.

\HOX{Change $u^\tau$ to $u_\tau$ systematically in the final version of the paper.}

We consider both the usual causal
solutions and the  solutions
for which time is reversed, that is, we use the anti-causal parametrix
${\bf Q}^*={\bf Q}^*_{\hat g}$ instead of the usual causal   parametrix ${\bf Q}=(\square_{\hat g}+V(x,D))^{-1}$. 
The wave   $u_\tau={\bf Q}^*F_\tau$. It satisfies by \cite{Ralston}
\beq\label{gaussian beam 1}
\|u_\tau-u^N_{\tau}\|_{C^k(J(p^-,p^+))}\leq C_N\tau^{-n_{N,k}}
\eeq
where $n_{N,k}\to \infty$ as $N\to \infty$ and $u^N_\tau$ is a formal Gaussian beam \cite{Ralston} of order $N$ having
the form
\beq\label{gaussian beam 2}
u^N_{\tau}(x ) = 
\exp (i\tau \varphi (x) )\left(
\sum _{n=0}^N U_n(x) \tau^{-n}\right),
\eeq
where
$
\varphi (x)=A(x)+iB(x)
$
and $A(x)$ and $B(x)$ are real functions, $B(x)\geq 0$, and $B(x)$ vanishes
only on $\gamma_{y,\eta} (\R)\cap J^-_{\hat g}(W)$, and
for $z= \gamma_{y,\eta}(t)$, and
$\zeta=\dot \gamma_{y,\eta}^\flat (t)$, $t<0$ we have
$dA|_{z}=\zeta$,  $dB|_{z}=0$, and
the Hessian of $B$ at $z$  restricted to the orthocomplement 
of $\zeta$ (with respect to $\hat g^+$) is positive definite.
Above, functions $h$ and $U_n$ can be chosen to be smooth and supported in
any neighborhood $V$ of $y$ and
any neighborhood of  
$\gamma_{y,\eta} ((-\infty,0])$. 

The source function $F_\tau$ can
be constructed in local coordinates using the asymptotic representation of the 
gaussian beam, namely, by considering first 
\ba
F^{(series)}_\tau(x)=c\tau^{-1/2}\int_\R e^{-\tau s^2}\phi(s)f_\tau(x;s)ds,
\ea
where
$f_\tau(x;s)=\square_{\hat g}(H(s-x^0)u_\tau(x))$, $\phi\in C^\infty_0(\R)$ 
has value 1 in some neighborhood of zero, and $H(s)$ is the Heaviside function,
and writing 
\ba
e^{-i\tau p(x)}F^{(series)}_\tau(x)=e^{-i\tau p(x)}F^{(1)}(x)\tau^{-1}+e^{-i\tau p(x)}F^{(2)}(x)\tau^{-2}+O(\tau^{-3}).
\ea Then,
$F_\tau$  can be defined by $F_\tau=F^{(1)}(x)\tau^{-1}.$ 
{
Indeed,\HOX{The details on computation of $F^{(series)}_\tau(x)$ can be omitted in the final paper version}
we see that
\ba
\square_{\hat g}^{-1}F^{(series)}_\tau&=&
\tau^{-1}\Phi_\tau(x)\,u_\tau(x),\\
\Phi_\tau(x)&=&
c\tau^{1/2}\int_\R e^{-\tau s^2}\phi(s)H(s-x^0) ds\\
&=&
\left\{\begin{array}{cl}
O(\tau^{-N}),&\hbox{for }x_0<0,\\
1+O(\tau^{-N}),&\hbox{for }x_0>0.
\end{array}\right.
\ea
Moreover, by writing $\square_{\hat g}=\hat g^{jk}\p_j\p_k+\hat \Gamma^j\p_j$,
$\hat \Gamma^j=\hat  g^{pq}\hat \Gamma^j_{pq}$ we see that
\ba
f_\tau(x;s)&=&
\hat g^{00}u_\tau(x)\delta^\prime(s-x^0)\\
& &-2\sum_{j=1}^3\hat g^{j0}(\p_ju_\tau(x))\,\delta(s-x^0)
-\hat \Gamma^0u_\tau(x)\,\delta(x^0-s).
\ea
Let us use  normal coordinates centered at the point $y$ so that
$\Gamma^j(0)=0$ and $\hat g^{jk}(0)$ is the Minkowski metric.
Then we define $p(x)=(x^0)^2+\varphi(x)$ and see that the leading
order term of  $F^{(series)}_\tau$ is given by
\ba
F_\tau= c\tau^{-1/2}e^{-\tau p(x)}\left(
2x^0\phi(x^0)U_0(x)+2i(\p_0\varphi(x))U_0(x)\phi(x^0)
\right),
\ea
where $\p_0\varphi|_y$ does not vanish. 

}


\subsubsection{Indicator function for singularities produced by interactions}

Let $y\in U$ and $\eta\in T_{x_0}M$ be  a future pointing light-like vector.
We will next make a test to see if  
$(y,\eta^\flat)\in\WF(\M^{(\ell)})$ with $\ell\leq 4$.

Using the functions $\M^{(\ell)}$ defined
in (\ref{w1-3 solutions})    with
the  pieces of plane waves 
$u_j\in \I(\Lambda(x_{j},\xi_{j};t_0,s_0))$, $j\leq 4$,
and the source $F_\tau$ in (\ref{Ftau source}),
 we define indicator functions
\beq\label{test sing}
\Theta_\tau^{(\ell)}=\bra F_{\tau},\M^{(\ell)}\cet_{L^2(U)},\quad \ell=1,2,3,4,
\eeq
We can write $\Theta_\tau^{(\ell)}$ 
is a sum of terms $T_{\tau,\sigma}^{(\ell),\beta}$
and $\tilde T_{\tau,\sigma}^{(\ell),\beta}$, where $\beta\in \Z_+$
are just numbers indexing terms, and $\sigma:\{1,2,\dots,\ell\}\to \{1,2,\dots,\ell\}$
is in the set of permutations of $\ell$ indexes, 
\ba
\Theta_\tau^{(\ell)}=\sum_{\beta\in J_\ell}\sum_{\sigma\in \Sigma(\ell)}(T_{\tau,\sigma}^{(\ell),\beta}+\tilde T_{\tau,\sigma}^{(\ell),\beta}).
\ea
To define the terms
 $T_{\tau,\sigma}^{(\ell),\beta}$ and $\tilde T_{\tau,\sigma}^{(\ell),\beta}$ 
appearing above,
we use  generic notations where we drop
the index $\beta$, that is, we denote $\B_j=\B_j^\beta$ and $S_j=S_j^\beta.$
Then
$T_{\tau,\sigma}^{(2),\beta}$ are terms of the form \HOX{In the main text we do not need
$\M^{(2)}$ or $\M^{(3)}$ and the analysis of these terms can be removed from the final paper.}
\beq\label{T2 term}
T^{(2),\beta}_{\tau,\sigma}&=&
\bra F_\tau,{\bf Q}(\B_2u_{\sigma(2)}\,\cdotp \B_1u_{\sigma(1)})\cet_{L^2(\hattuM _0)}
\\ \nonumber
&=&\bra {\bf Q}^* F_\tau,\B_2u_{\sigma(2)}\,\cdotp \B_1u_{\sigma(1)}\cet_{L^2(\hattuM _0)}
\\ \nonumber
&=&\bra u_\tau,\B_2u_{\sigma(2)}\,\cdotp \B_1u_{\sigma(1)}\cet_{L^2(\hattuM _0)},
\eeq
and $\tilde T_{\tau,\sigma}^{(2),\beta}=0$.
Moreover, $T_\tau^{(3),\beta}$ are of the form 
\beq\label{T3 term}
T^{(3),\beta}_{\tau,\sigma}
&=&\bra F_\tau,{\bf Q}(\B_3u_{\sigma(3)}\,\cdotp \cC_1 S_1(\B_2u_{\sigma(2)}\,\cdotp \B_1u_{\sigma(1)}))\cet_{L^2(\hattuM _0)}\\
\nonumber
&=&\bra u_\tau,\B_3u_{\sigma(3)}\,\cdotp  \cC_1 S_1(\B_2u_{\sigma(2)}\,\cdotp \B_1u_{\sigma(1)})\cet_{L^2(\hattuM _0)},
\eeq
and  $\tilde T_{\tau,\sigma}^{(3),\beta}=0$.
When we consider the 4th order interaction terms
$T^{(4),\beta}_{\tau,\sigma}$ we change our notations by 
commuting $S_j$ and $\cC_j$ and using
 Leibniz rule. For instance, if $\cC_1$ is a first order operator, we write  
\ba
 \cC_1^\beta S_1^\beta (\B_2^\beta u_{2}\,\cdotp \B_1^\beta u_{1})&=&
[\cC_1^\beta  , S_1^\beta ](\B_2^\beta u_{2}\,\cdotp \B_1^\beta u_{1})+\\
& &+
S_1((\cC_1^\beta  \B_2^\beta u_{2})\,\cdotp  \B_1^\beta u_{1})+
S_1( \B_2^\beta u_{2}\,\cdotp ((\cC_1^\beta \B_1^\beta u_{1})).
\ea
Using this we can eliminate the operators $\cC_j^\beta $ and increase the order
of the differential operators $\B_j^\beta $ and allow $S_j^\beta $ also be a commutator of ${\bf Q}$ and $\cC_j^\beta $,
that is, below we allow $\B^\beta _j$ and $S^\beta _j$ to be of the form
\beq\label{extra notation}
& &\hbox{$\B^\beta _j:(v_{p})_{p=1}^{10+L}\mapsto (b^{r,(j,\beta )}_{p}(x)\p_x^{\vec k(\beta,j)} v_{r}(x))_{p=1}^{10+L}$, and}\\ \nonumber
& &\hbox{$S^\beta _j={\bf Q}$ or $S^\beta _j=I$, or $S^\beta _j=[{\bf Q},a(x)D^\a]$},
\eeq
where $|\vec k(\beta)|=\sum_{j=1}^4k(\beta,j)\leq 6$. To simplify notations, we also 
enumerate again the obtained \HOX{In final paper, we need to check that explanation on indexes $\beta$ is clear.}
terms with the indexes   $\beta$, that is, we consider $\beta$ as an index running over
a non-specified finite set.
To simplify notations, we sometimes drop the super-index $\beta$ below. With these notations,
\beq\nonumber
T^{(4),\beta}_{\tau,\sigma}&=&\bra F_\tau,{\bf Q}(\B_4u_{\sigma(4)} \,\cdotp  S_2(\B_3u_{\sigma(3)}\,\cdotp   S_1(\B_2u_{\sigma(2)}\,\cdotp \B_1u_{\sigma(1)})))\cet_{L^2(\hattuM _0)}\\ 
\label{T-type source}
&=&\bra {\bf Q}^* F_\tau,\B_4u_{\sigma(4)} \,\cdotp   S_2(\B_3u_{\sigma(3)}\,\cdotp  S_1(\B_2u_{\sigma(2)}\,\cdotp \B_1u_{\sigma(1)}))\cet_{L^2(\hattuM _0)}\\ \nonumber
&=&\bra (\B_4u_{\sigma(4)}) \,\cdotp u_{\tau}, S_2(\B_3u_{\sigma(3)}\,\cdotp  S_1(\B_2u_{\sigma(2)}\,\cdotp \B_1u_{\sigma(1)}))\cet_{L^2(\hattuM _0)}
\\ \nonumber
&=&\bra (\B_3u_{\sigma(3)})\,\cdotp S_2^*(  (\B_4u_{\sigma(4)}) \,\cdotp u^\tau), S_1(\B_2u_{\sigma(2)}\,\cdotp \B_1u_{\sigma(1)})\cet_{L^2(\hattuM _0)}
\eeq
and
\beq\nonumber
\tilde T^{(4),\beta}_{\tau,\sigma}&=&\bra F_\tau,{\bf Q}( S_2(\B_4u_{\sigma(4)} \,\cdotp \B_3u_{\sigma(3)})\,\cdotp  S_1(\B_2u_{\sigma(2)}\,\cdotp \B_1u_{\sigma(1)}))\cet_{L^2(\hattuM _0)}\\ \label{tilde T-type source}
&=&\bra {\bf Q}^* F_\tau, S_2(\B_4u_{\sigma(4)} \,\cdotp  \B_3u_{\sigma(3)})\,\cdotp  S_1(\B_2u_{\sigma(2)}\,\cdotp \B_1u_{\sigma(1)})\cet_{L^2(\hattuM _0)}\\ \nonumber
&=&\bra  S_2(\B_4u_{\sigma(4)} \,\cdotp  \B_3u_{\sigma(3)}) \,\cdotp u_\tau,  S_1(\B_2u_{\sigma(2)}\,\cdotp \B_1u_{\sigma(1)})\cet_{L^2(\hattuM _0)}.
\eeq



When  $\sigma$ is the identity, we will omit it in our notations
and denote $T^{(4),\beta}_{\tau}=T^{(4),\beta}_{\tau,id}$, etc.

In particular,   the term
$\bra F_\tau,{\bf Q}(A_2[u_4,{\bf Q}(A_2[u_3,{\bf Q}(A_2[u_2,u_1])])])\cet$,
where we recall that $A_2[v,w]=\hat g^{np}\hat g^{mq}v_{nm}\p_p\p_q w_{jk}$,
can be written  as a sum of terms of the type 
 $ T^{(4),\beta}_{\tau,\sigma}$, 
 we obtain one term that will later cause the leading order asymptotics.
 This term
 corresponds to $\sigma=Id$ and the indexes $k_1=6$,
 $k_2=k_3=k_2=0$, and we  enumerate this term to  correspond $\beta=\beta_1:=1$ (we define these 
 indexes used later)
 \beq\label{beta0 index}
\hbox{$\vec S_{\beta_1}=({\bf Q},{\bf Q})$ and $k^{\beta_1}_1=6$, 
 $k^{\beta_1}_2=k^{\beta_1}_3=k^{\beta_1}_2=0$.}
 \eeq


%
%

\subsubsection{Properties of indicator functions on a general manifold}
To consider the properties of 
the indicator functions related to the sources
$f_j\in \I(Y(x_j,\xi_j;t_0,s_0))$, considered as sections of
the bundle $\B^K$,
where $(x_j,\xi_j)\in L^+M_0,$ $x_j\in\hat U$, $j\leq 4$,
and the source $F_\tau$ determined by  
 $(y,\eta)\in L^+M_0$,  $y\in\hat U$
 we denote in the following 
$(x_5,\xi_5)=(y,\eta)$ and continue to use the notation 
$(\vec x,\vec \xi)=((x_j,\xi_j))_{j=1}^4$.
%

\begin{definition}\label{def:  4-intersection of rays}
We say that the geodesics corresponding to $(\vec x,\vec \xi)=((x_j,\xi_j))_{j=1}^4$
intersect  and  the intersection takes place at the point $q$ 
if there is $q\in \hattuM _0$ such that
 for all $j=1,2,3,4$, there are $t_j\in (0,{\bf t}_j)$, ${\bf t}_j=\rho(x_j,\xi_j)$ such that 
 $q=\gamma_{x_j,\xi_j}(t_j)$. In this case  that  such $q$ exists
 and is of the form  $q=\gamma_{x_5,\xi_5}(t_5)$, $t_5<0$, we say that
   \HOX{In the final paper, clarify or remove
 the definition of " $(x_5,\xi_5)$ comes from the 4-intersection of rays".}
 $(x_5,\xi_5)$ comes from the 4-intersection of rays
corresponding to $(\vec x,\vec \xi)=((x_j,\xi_j))_{j=1}^4$.
Also, we say that  $q$ is the intersection point corresponding
 to $(\vec x,\vec \xi)$ and $(x_5,\xi_5)$.
 \end{definition}
 
  Let  $\Lambda_q^+$ be the lagrangian manifold 
  \ba
  \Lambda_q^+=\{(x,\xi)\in T^*M_0;\ x=\gamma_{q,\zeta}(t),\ 
  \xi^\sharp=\dot\gamma_{q,\zeta}(t),\ \zeta\in L^+_qM_0,\ t>0\}
  \ea
  Note that the projection of $ \Lambda_q^+$
 on $\hattuM _0$ is  the light cone $\L^+_{\hat g}(q)\setminus \{q\}$. 

\begin{figure}[htbp]
\begin{center}  \label{Fig-15}

\psfrag{1}{$x_1$}
\psfrag{2}{$x_2$}
\psfrag{3}{$q$}
\psfrag{4}{$p_1$}
\psfrag{5}{$p_2$}
\psfrag{6}{\hspace{-2mm}$x_6$}
\psfrag{7}{$x_5$}
\includegraphics[width=8cm]{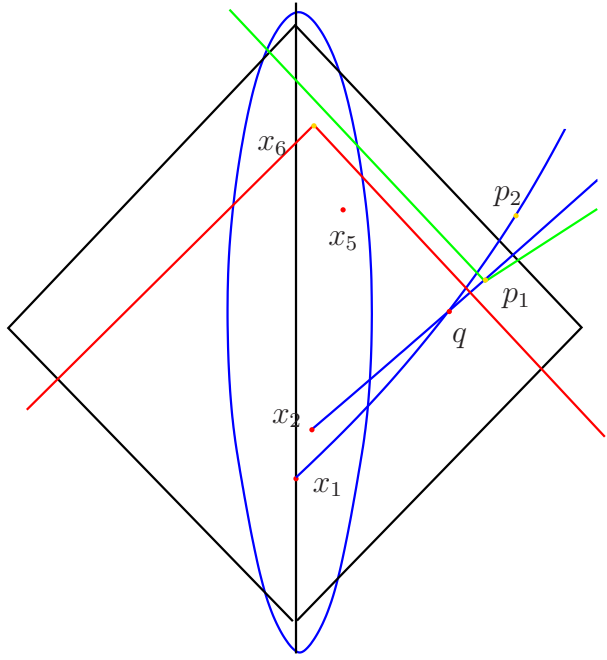}
\end{center}

\caption{ A schematic figure where the space-time is represented as the  2-dimensional set $\R^{1+1}$. 
The figure shows the configuration in formulas (\ref{eq: summary of assumptions 1}) and  (\ref{eq: summary of assumptions 2}).
The  points $x_1$ and $x_2$, marked with red dots,
are the points from where we send light-like geodesics $\gamma_{x_j,\xi_j}([t_0,\infty))$,
$j=1,2.$ The geodesics $\gamma_{x_j,\xi_j}([t_0,\infty))$,
$j=3,4$ are not in the figure.  We assume the these  geodesics intersect
at point $q$  that is shown as a  red point. The geodesics $\gamma_{x_j,\xi_j}([t_0,\infty))$
have cut points and the first cut points
  $p_j$ are shown as golden points. The point $x_6\in U_{\hat g}$
is such that no cut points $p_j$ of geodesics are in the causal past $J^-(x_6)$ of the point $x_6$
shown with red lines.
Observations are done at the point $x_5\in  J^-(x_6)\cap U_{\hat g}$. 
The set  $J^-(x_6)$ has a neighborhood $\mathcal V=\mathcal V((\vec x,\vec \xi),t_0)$ that is the complement
of the set $\bigcup_{j=1}^4 J^+(p_j)$. The boundary of $\mathcal V$ is shown with green lines.
The the black "diamond" is again  the set
$J_{\hat g}(\hat p^-,\hat p^+)$ and the black line is the geodesic $\hat \mu$.
The  domain inside the blue curve is the set $U_{\hat g}$ where the sources are supported and the observations are done.
}
\end{figure}

 Below we say that a function $F(z)$ is a meromorphic function of the variables $z\in \R^m$
 if $F(z)=P(z)/Q(z)$ where $P(z)$ and $Q(z)$ are real-analytic functions and 
 $Q(z)$ does not vanish identically.
Next we consider  $(\vec x,\vec \xi)=((x_j,\xi_j))_{j=1}^4$ and $\vartheta_0,t_0>0$ that satisfy (see Fig.\ 13)
\beq\label{eq: summary of assumptions 1}& &\hspace{-1.0cm}(i)\ x_j\in \hat U,\ \xi_j\in L^+_{x_j}\hattuM _0,\
 x_j\not \in J^+_{\hat g}(x_k),
  \hbox{ for $j,k\leq 4$, $j\not =k$,}  \hspace{-5mm}\\
 & &\nonumber
 \hspace{-1.0cm}(ii)\ 
\hbox{For all $j,k\leq 4$,
 $d_{\hat g^+}((x_j,\xi_j),(x_k,\xi_k))<\vartheta_0$,} \hspace{-5mm}\\
 & &\nonumber
 \hspace{-1.0cm}(iii)\ 
\hbox{There is $y\in  \hat \mu$ such that \ 
for all $j\leq 4$,
 $d_{\hat g^+}(y,x_j)<\vartheta_0$,} \hspace{-5mm}\\
& &\nonumber
\hspace{-1.0cm}(iv) \hbox{ For $j,k\leq 4$, $j\not =k$, 
we have $x_j(t_0)\not \in \gamma_{x_k,\xi_k}(\R_+)$.} \hspace{-5mm}
 \eeq
 Above, 
$(x_j(h),\xi_j(h))$ are defined in (\ref{eq: x(h) notation}).
 {\newtext We also consider a point $x_6\in U_{\hat g}$ that satisfies for  
 $(\vec x,\vec \xi)$ and $t_0$ satisfies the condition \HOX{Change notation $\V$ to $\M$ as Slava suggested.
}
 \beq\label{eq: summary of assumptions 2}
& & \hbox{For all $j\leq 4$, }x_j^{cut}=\gamma_{x_j(t_0),\xi_j(t_0)}({\bf t}_j)\not \in J^-_{\hat g}(x_6), 
 \\
 \nonumber
& & \hbox{where }{\bf t}_j:=\rho(x_j(t_0),\xi_j(t_0)).\hspace{-1cm}\\ \nonumber
& & \hbox{Then, we denote }\V((\vec x,\vec \xi),t_0)= M_0\setminus \bigcup_{j=1}^4
J^+_{\hat g}(\gamma_{x_j(t_0),\xi_j(t_0)}({\bf t}_j)).
 \eeq
 Observe that $ \V((\vec x,\vec \xi),t_0)$ is an open neighborhood of $J^-_{\hat g}(x_6)$.
 \HOX{In the final paper,
 check if $ \V((\vec x,\vec \xi),t_0)$  can be replaced by $I^-(x_6)$ and ${\bf p}(y,...)$ by $x_6$.}

%

 The condition  (\ref{eq: summary of assumptions 2}) implies that when we consider the past of $x_6$,
 then no geodesics $\gamma_{x_j(t_0),\xi_j(t_0)}$ has
 conjugate or cut points.} Note that
 two such geodesics $\gamma_{x_j(t_0),\xi_j(t_0)}([0,\infty))$ can intersect only once in $\V((\vec x,\vec \xi),t_0)$.

In the proposition below we will consider four spherical
waves $u_j\in \I^{n-1/2}(\Lambda(x_{j},\xi_{j};t_0,s_0))$, $j=1,2,3,4$, where the order 
$n$ is low enough, that is $n\leq -n_1$. A suitable value of $n_1$ is $n_1=12$, and the reason
for this is that we need to consider fourth order interaction terms and each interaction involves
two derivativs.

  \begin{proposition}\label{lem:analytic limits A}    Let $(\vec x,\vec \xi)=((x_j,\xi_j))_{j=1}^4$  be future pointing light-like vectors 
   and $x_6\in U_{\hat g}$ 
   satisfing (\ref{eq: summary of assumptions 1})-(\ref{eq: summary of assumptions 2}).
  Let $x_5\in \V((\vec x,\vec \xi),t_0)\cap U_{\hat g}$ and $(x_5,\xi_5)$ be a future pointing light-like vector
such that
 $x_5\not\in \gamma_{x_j,\xi_j}(\R)$ for $j\leq 4$, $t_0>0$,  and 
  $x_5\not\in{\mathcal Y}={\mathcal Y}(((x_{j},\xi_{j}))_{j=1}^4;t_0)$,
  see (\ref{associated submanifold}) and (\ref{eq: 3-interactions}).


There exists $n_1\in \Z_+$ such that the following  holds:
When $s_0>0$ is small enough, 
the function $\Theta^{(\ell)}_\tau$, see (\ref{test sing}),  corresponding to  the linear waves
$u_j\in \I^{n-1/2}(\Lambda(x_{j},\xi_{j};t_0,s_0))$, $j\leq 4$, 
defined in Lemma \ref{lem: lagrangian 1} with $n\leq -n_1$,
and the source $F_\tau$ 
satisfies the following:

%
%
 \medskip

 \noindent
 (i) If 
$(x_5,\xi_5)$ does not come from the 4-intersection
of  rays corresponding to $(\vec x,\vec \xi)$, we have
   $|\Theta^{(4)}_\tau|\leq C_N\tau^{-N}$
 for all $N>0$.

(ii) if
$(x_5,\xi_5)$ comes from the 4-intersection
of  rays corresponding to $(\vec x,\vec \xi)$ and 
$q$ is the corresponding intersection point, $q=\gamma_{x_j,\xi_j}(t_j)$,
then 
\beq\label{indicator}
\Theta^{(4)}_\tau\sim  
\sum_{k=m}^\infty s_{k}\tau^{-k}
\eeq
as $\tau\to \infty$ where   $m=-4n+2$. 
Here we use $\sim$ to denote that the terms
have the same  asymptotics up an error $O(\tau^{-N})$ for all $N>0$.


Moreover,
let $b_j=(\dot\gamma_{x_j,\xi_j}(t_j))^\flat$ and
 $\bsequence=(b_{j})_{j=1}^5\in (T^*_q\hattuM _0)^5$,
 $w_j$ be the principal symbols of the waves $u_j$
 at $(q,b_j)$, and ${\bf w}=(w_j)_{j=1}^5$. 
 Then there is 
  a real-analytic function $\mathcal G(\bsequence,{\bf w})$ such that
 the leading order term in (\ref{indicator})  satisfies
 \beq\label{definition of G}
s_{m}=
\mathcal  G(\bsequence,{\bf w}).
\eeq
\end{proposition}

\medskip

\noindent
{\bf Proof.}  
Below, to simplify notations,
we denote $K_j=K(x_{j},\xi_{j};t_0,s_0)$ and 
$K_{123}=K_1\cap K_2\cap K_3$ and $K_{124}=K_1\cap K_2\cap K_4$, etc.
We will denote $ \Lambda_j=\Lambda(x_j,\xi:j;t_0,s_0)$
to consider also the singularities of $K_j$ related to conjugate points.

Below we will consider separately the case when the
following linear independency condition,
\medskip

(LI) Assume if that  $J\subset \{1,2,3,4\}$ and $y\in J^-(x_6)$
are such that for all $j\in J$  we have $\gamma_{x_j,\xi_j}(t^\prime_j)=y$
with some $t^\prime_j\geq 0$, then  the vectors
$\dot\gamma_{x_j,\xi_j}(t_j^\prime )$, $j\in J$ are linearly independent.
\medskip

\noindent is valid and the case when (LI) is not valid. \HOX{Jan.\ 4, 2013: 
Maybe we do not need to analyze (LI), see
condition (D)  below, and thus (LI) can be moved to claim of the proposition.}

Let us first  consider  the case when (LI) is valid.

By the definition of ${\bf t}_j$,  if
the intersection $\gamma_{x_5,\xi_5}(\R_-) \cap(\cap_{j=1}^4\gamma_{x_j,\xi_j}((0,{\bf t}_j)))
$ is non-empty, it can
contain only one point. In the case that such a
point exists, we denote it by $q$. When $q$ exists,
the intersection of $K_j$ at this point are transversal and
we see that when $s_0$ is small enough, the set
$\cap_{j=1}^4 K_j$ consists only of the point $q$.
Next we consider so small $s_0$ that this is true.

We consider    two local coordinates $Z:W_0\to \R^4$
and  $Y:W_1\to \R^4$ such that $W_0,W_1\subset \V((\vec x,\vec \xi),t_0)$,
 see definition after (\ref{eq: summary of assumptions 2}).
 We assume that  local coordinates are such that
$K_j\cap W_0=\{x\in W_0;\ Z^j(x)=0\}$ and $K_j\cap W_1=\{x\in W_1;\ Y^j(x)=0\}$ for $j=1,2,3,4$.
In the Fig.\ 14, $W_0$ is a neighborhood of $z$ and
$W_1$ is the neighborhood of $y$. 
We note that the origin $0=(0,0,0,0)\in \R^4$ does not
necessarily belong to the set $Z(W_0)$ or the set $Y(W_1)$,
for instance in the case when the four geodesic $\gamma_{x_j,\xi_j}$ do not intersect.
However, this is the case when the  four geodesic intersect at the point $q$, we need
to consider the case when $W_0$ and $W_1$ are neighborhoods of $0$. 
 %
Note that
$W_0$ and $W_1$ do not contain any cut points of the geodesic
 $\gamma_{x_j,\xi_j}([t_0,\infty)$.

We will denote $z^j=Z^j(x)$.
We assume that  $Y:W_1\to \R^4$ are similar coordinates
and denote  $y^j=Y^j(x)$. We also denote
 below $dy^j=dY^j$ and  $dz^j=dZ^j$.
Let $\Phi_0\in C^\infty_0(W_0)$ and 
 $\Phi_1\in C^\infty_0(W_1)$.  
 
 Let us next considering the map ${\bf Q}^*:C^\infty_0(W_1)\to
 C^\infty(W_0).$  By \cite{MU1}, ${\bf Q}^*\in I(W_1\times W_0;
   \Delta_{T\hattuM _0}^{\prime}, \Lambda_{\hat g})$ is an  operator 
   with a classical symbol and its canonical
   relation 
$\Lambda_{{\bf Q}^*}^{\prime}$ is associated to a union of two  intersecting lagrangian manifolds,
$\Lambda_{{\bf Q}^*}^{\prime}= \Lambda_{\hat g}^{\prime}\cup \Delta_{T\hattuM _0}$,
intersecting cleanly \cite{MU1}.
Let $\e_2>\e_1>0$ and $B_{\e_1,\e_2}$ be a pseudodifferential
operator on $\hattuM _0$ which is microlocally smoothing
operator (i.e., the full  symbol vanishes  in local coordinates)
outside in the $\e_2$-neighborhood $\V_2\subset T^*\hattuM _0$ of the set of the light like covectors
$L^*\hattuM _0$ and for which 
$(I-B_{\e_1,\e_2})$ is microlocally smoothing
operator 
 in the $\e_1$-neighborhood $\V_2$ of $L^*\hattuM _0$.
 The neighborhoods here are defined 
 with respect to the Sasaki metric of $(T^*\hattuM _0,\hat g^+)$ and
$\e_2,\e_1$ are chosen later in the proof. 
 Let us decompose the operator ${\bf Q}^*={\bf Q}^*_1+{\bf Q}^*_2$
 where  ${\bf Q}^*_1={\bf Q}^*(I-B_{\e_1,\e_2})$ and   ${\bf Q}^*_2={\bf Q}^*B_{\e_1,\e_2}$.
 As $\Lambda_{{\bf Q}^*}= \Lambda_{\hat g}\cup \Delta_{T\hattuM _0}^{\prime}$,
we see that then there is 
a neighborhood $ \W_2=\W_2(\e_2)$ of $L^*\hattuM _0\times L^*\hattuM _0\subset (T^*\hattuM _0)^2$
such that the Schwartz
kernel ${\bf Q}^*_2(r,y)$ of the operator ${\bf Q}^*_2$ satisfies
\beq\label{eq: W_2 neighborhood}
\hbox{WF}({\bf Q}^*_2)\subset \W_2.
\eeq
Moreover,  $\Lambda_{{\bf Q}^*_1}\subset \Delta_{T\hattuM _0}^{\prime}$ implying
that ${\bf Q}^*_1$ is a pseudodifferential operator with a classical symbol,
 ${\bf Q}^*_1\in I(W_1\times W_0;
   \Delta_{T\hattuM _0}^{\prime})$,
 and
 ${\bf Q}^*_2 
 \in I(W_1\times W_0;
   \Delta_{T\hattuM _0}^{\prime}, \Lambda_{\hat g})$
   is a Fourier integral operator (FIO) associated
   to two cleanly intersecting lagrangian manifolds, similarly to ${\bf Q}^*$. 
In the case when  $p=1$ we can write ${\bf Q}^*_p$ as 
\beq\label{eq: representation of Q* p}
({\bf Q}^*_1v)(z)=\int_{\R^{4+4}}e^{i\psi_1(z,y,\xi)}q_1(z,y,\xi)v(y)\,dyd\xi,
\eeq
where 
\beq\label{eq: psiDO phase function}
\psi_1(z,y,\xi)=(y-z)\,\cdotp \xi,\quad \hbox{for }(z,y,\xi)\in W_1\times W_0\times \R^4,
\eeq
{
and a classical symbol $q_1 (z,y,\xi)\in S^{-2}(W_1\times W_0; \R^4)$, having a real
valued principal symbol
\ba
\tilde q_1 (z,y,\xi)=\frac {\chi(z,\xi)}{g_z(\xi,\xi)}
\ea
where $\chi(z,\xi)\in \C^\infty$ is a cut-off function vanishing
in a neighborhood of  the set where $g_z(\xi,\xi)=0$. Note that
then $Q_1-Q_1^*\in \Psi^{-3}(W_1\times W_0)$.}

Furthermore, let us decompose 
${\bf Q}^*_1={\bf Q}^*_{1,1}+{\bf Q}^*_{1,2}$ corresponding to the decomposition
$q_1(z,y,\xi)=q_{1,1}(z,y,\xi)+q_{1,2}(z,y,\xi)$
 of the symbol,
where
 \beq\label{eq: q_{1,1} and q_{1,2}} \hspace{-.3cm}
 q_{1,1}(z,y,\xi)=q_{1}(z,y,\xi)\psi_R(\xi),\quad q_{1,2}(z,y,\xi)=q_{1}(z,y,\xi)(1-\psi_R(\xi)),
 \hspace{-2cm}
 \eeq
  where $\psi_R\in C^\infty_0(\R^4)$
is a cut-off function that is equal to one in  a ball $B(R)$ of radius $R$ specified below.

Next we start to consider the terms $T_{\tau}^{(4),\beta}$ and $\tilde T_{\tau}^{(4),\beta}$ of the type
 (\ref{T-type source}) 
and
(\ref{tilde T-type source}). In these terms,
we can represent the gaussian beam $u_{\tau}(z)$ in $W_1$ in the form 
\beq\label{eq:representation of the wave}
u_\tau(y)=e^{i\tau \varphi (y)}a_5(y,\tau)
\eeq
where the function $ \varphi$ is a complex phase function having non-negative imaginary
part such that $\im   \varphi$, defined on $W_1$, vanishes exactly on the geodesic
 $\gamma_{x_5,\xi_5}\cap W_1$. Note that 
  $\gamma_{x_5,\xi_5}\cap W_1$ may be empty. Moreover,
 $a_5\in
S^{0}_{clas}(W_1;\R)$ is a classical symbol.

 Also, for $y=\gamma_{x_5,\xi_5}(t)\in W_1$
we have that $d\varphi(y)=c\dot \gamma_{x_5,\xi_5}(t)^\flat$, with $c\in \R\setminus \{0\}$, is light-like.



We consider first the asymptotics of 
 terms $T_\tau^{(4),\beta}$ and $\tilde T^{(4),\beta}_\tau$ of the type
 (\ref{T-type source}) 
and
(\ref{tilde T-type source})
   where $S_1=S_2={\bf Q}$ and, symbols $a_j(z,\theta_j)$, $j\leq 3$ and  
    $a_j(y,\theta_j)$, $j\in \{4,5\}$ are scalar valued
   symbols written in the $Z$ and $Y$ coordinates,
$\B_1,\B_2,\B_3$ are multiplication operators with $\Phi_0$,
and $\B_4$ is a  multiplication operators with $\Phi_1$
 and consider section-valued symbols and general operators later.

Let us consider functions
$U_j\in \I^{p_j}(K_j)$, $j=1,2,3,$ supported in $W_0$ and 
$U_4\in \I(K_4)$,  supported in $W_1$,  having classical symbols, 
 \beq\label{U1term}
& &U_j(x)=\int_{\R}e^{i\theta_jx^j}a_j(x,\theta_j)d\theta_j,\quad a_j\in
S^{p_j}_{clas}(W_{k(j)};\R), 
\eeq
for all $j=1,2,3,4$ (Note that here the phase function is  $\theta_jx^j=\theta_1x^1$ for $j=1$ etc,
that is, there is no summing over index $j$). We may assume that $a_j(x,\theta_j)$ vanish near $\theta_j=0$.

As  $x_5\in \V((\vec x,\vec \xi),t_0)\cap U_{\hat g}$ and 
$W_0,W_1\subset \V((\vec x,\vec \xi),t_0)$, we see that
an example of functions (\ref{U1term}) are $U_j(z)=\Phi_{k(j)}(z)u_j(z)$, $j=1,2,3,4$,
where $u_j(z)$ are the pieces of the spherical waves. Here and below, 
 $k(j)=0$ for $j=1,2,3$  and $k(4)=1$ and we
also denote $k(5)=1$.
Note that $p_j=n$ correspond to the case when $U_j\in \I^n(K_j)=\I^{n-1/2}(N^*K_j)$.

Denote $\Lambda_j=N^*K_j$ and $\Lambda_{jk}=N^*(K_j\cap K_k)$. By \cite[Lem.\ 1.2 and 1.3]{GU1}, the  
pointwise product 
$U_2\,\cdotp U_1\in \I(\Lambda_1,\Lambda_{12})+ \I(\Lambda_2,\Lambda_{12})$  and thus
by \cite[Prop.\ 2.2]{GU1},
${\bf Q}(U_2\,\cdotp U_1)\in \I(\Lambda_1,\Lambda_{12})+ \I(\Lambda_2,\Lambda_{12})$ 
and  it
can be written as 
\beq\label{Q U1U2term}
{\bf Q}(U_2\,\cdotp U_1)=\int_{\R^2}e^{i(\theta_1z^1+\theta_2z^2)}c_1(z,\theta_1,\theta_2)d\theta_1d\theta_2.
\eeq
Note that here $c_1(z,\theta_1,\theta_2)$ is sum of product type symbols, see (\ref{product symbols}).
As $ N^*(K_1\cap K_2)
\setminus (N^*K_1\cup N^*K_2)$ consist of vectors which
are non-characteristic for 
the wave operator, that is, the wave operator is elliptic in
a neighborhood of this subset of the cotangent bundle,   the principal symbol 
$\tilde c_1$ of $c_1$ on  $N^*(K_1\cap K_2)
\setminus (N^*K_1\cup N^*K_2)$ is  given by
\beq\label{product type symbols}
& &\tilde c_1(z,\theta_1,\theta_2)\sim
 s(z,\theta_1,\theta_2)a_1(z,\theta_1)a_2(z,\theta_2),
 \\ \nonumber
 & &\hspace{-1cm}s(z,\theta_1,\theta_2)
=1/{ g(\theta_1b^{(1)}+\theta_2b^{(2)},\theta_1b^{(1)}+\theta_2b^{(2)})}
 =1/({2 g(\theta_1b^{(1)},\theta_2b^{(2)})}). \hspace{-1cm}
 \eeq
Note that $s(z,\theta_1,\theta_2)$ is a smooth function  
on $N^*(K_1\cap K_2)
\setminus (N^*K_1\cup N^*K_2)$ and homogeneous of order $(-2)$ in $\theta=(\theta_1,\theta_2)$.
Here, we use $\sim$ to denote that the symbols have the same principal symbol.
%
%
%
%
Let us next
make computations in the case when 
 $a_j(z,\theta_j)\in C^\infty(\R^4\times \R)$ is positively homogeneous for $|\theta_j|>1$, that is, we have 
 $ a_j(z,s)=a^{\prime}_j(z)s^{p_j}$, where $p_j\in \N$  and $|s|>1$.
 We consider $T_\tau^{(4),\beta}=\sum_{p=1}^2T_{\tau,p}^{(4),\beta}$
 where $T_{\tau,p}^{(4),\beta}$ is defined as 
 $T_{\tau}^{(4),\beta}$ by replacing the term ${\bf Q}^*(U_4\,\cdotp u_{\tau})$
  by ${\bf Q}^*_p(U_4\,\cdotp u_{\tau})$.

 Let us now consider the case $p=2$ and choose the parameters $\e_1$ and $\e_2$
that determine the
decomposition ${\bf Q}^*={\bf Q}^*_1+{\bf Q}^*_2$. 
First, we observe that for $p=2$ we can write using $Z$ and $Y$ coordinates
 \beq\label{eq; T tau asympt, p=2}
T_{\tau,2}^{(4),\beta}
&=&
\tau^{4}
\int_{\R^{12}}e^{i\tau \Psi_2(z,y,\theta)}
c_1(z,\tau \theta_1,\tau \theta_2)\cdotp\\ \nonumber
& &\hspace{-1.3cm}
\cdotp  a_3(z,\tau \theta_3)
{\bf Q}^*_2(z,y)
a_4(y,\tau \theta_4) a_5(y,\tau )\,d\theta_1d\theta_2d\theta_3d\theta_4 dydz,\\
\nonumber \Psi_2(z,y,\theta)
&=&\theta_1z^1+\theta_2z^2+\theta_3z^3+\theta_4y^4+\varphi(y).
\eeq

Denote $\theta=(\theta_1,\theta_2,\theta_3,\theta_4)\in \R^{4}$.
Consider the case when $(z,y,\theta)$ is a critical point of $\Psi_2$
satisfying $\im \varphi(y)=0$. Then we have
 $\theta^{\prime}=(\theta_1,\theta_2,\theta_3)=0$
and $z^{\prime}=(z^1,z^2,z^3)=0$, 
$y^4=0$,
  $d_y\varphi(y)=(0,0,0,-\theta_4)$,
implying that $y\in K_4$ and 
$(y,d_y\varphi(y))\in N^*K_4$.
Since $\im \varphi(y)=0$, we have that 
  $y=\gamma_{x_5,\xi_5}(t_0)$ with some $t_0\in\R_-$.
 As we have $\dot \gamma_{x_5,\xi_5}(t_0)^\flat=d_y\varphi(y)\in N^*_yK_4$,
we obtain $\gamma_{x_5,\xi_5}([t_0,0])\subset K_4$. However, 
this 
is not possible by our assumption 
$x_5\not \in \cup_{j=1}^4 K(x_{j},\xi_{j};s_0)$ when $s_0$ is small enough.
Thus  the phase function 
$ \Psi_2(z,y,\theta)$ has no critical points satisfying $\im \varphi(y)=0$.

%

When the orders $p_j$ of the symbols $a_j$ are small enough, the integrals
in the $\theta$ variable in  (\ref{eq; T tau asympt, p=2}) are convergent in the classical sense.
Next we use  properties of wave front set to 
compute the asymptotics of  an oscillatory integrals
and to this end we introduce the function
\beq
\tilde {\bf Q}^*_2(z,y,\theta)={\bf Q}^*_2(z,y),
\eeq
that is, consider ${\bf Q}^*_2(z,y)$ as a constant function
in $\theta$. Below, denote $\psi_4(y,\theta_4)=\theta_4y^4$ and $r=d\varphi(y)$. 
Note that then $d_{\theta_4}\psi_4=y^4$ and $d_y\psi_4=(0,0,0,\theta_4)$.
Then  in $W_1\times W_0\times \R^4$
\ba
d_{z,y,\theta}\Psi_2&=&(\theta_1,\theta_2,\theta_3,0;
r+d_y\psi_4(y,\theta_4),
z^1,z^2,z^3,d_{ \theta_4}\psi_4(y,\theta_4))
\\
&=&(\theta_1,\theta_2,\theta_3,0;
d\varphi(y)+(0,0,0,\theta_4),
z^1,z^2,z^3,y^4)
\ea
and we see that if $((z,y,\theta),d_{z,y,\theta}\Psi_2)\in \hbox{WF}(\tilde {\bf Q}_2^*)$
and  $\im \varphi(y)=0$, we have $(z^1,z^2,z^3)=0$, $y^4=d_{\theta_4}\psi_4(y,\theta_4)=0$
and $y\in  \gamma_{x_5,\xi_5}$. Thus 
$z\in K_{123}$ and $y\in  \gamma_{x_5,\xi_5}\cap K_4$.

Let  us use the following notations
\beq\label{new notations}
& &z\in K_{123},\quad \omega_\theta:=(\theta_1,\theta_2,
 \theta_2,0)=\sum_{j=1}^3\theta_j dz^j\in T^*_z\hattuM _0,\\
 \nonumber
& &y\in K_4\cap \gamma_{x_5,\xi_5},\quad (y,w):=(y,d_y\psi_4(y,\theta_4))\in N^*K_4,
\\
 \nonumber
& & r=d\varphi(y)=r_jdy^j\in T^*_y\hattuM _0,\quad \kappa:=r+w.
 \eeq
Then, $y$ and $ \theta_4$ satisfy $y^4=d_{ \theta_4}\psi_4(y,\theta_4)=0$ and
$w=(0,0,0,\theta_4)$.
 
 Note that  by definition of the $Y$  coordinates $w$ is a light-like covector.
 By definition of the  $Z$ coordinates,
 $\omega_\theta\in 
 N^*K_1+N^*K_2+N^*K_3=N^*K_{123}$.  
 
 Let us first consider what happens if $\kappa=r+w=d\varphi(y)+(0,0,0,\theta_4)$ is light-like.
In this case, all vectors $\kappa$, $w$, and $r$ are light-like
and satisfy $\kappa =r+w$. This is possible
only if  $r\parallel w$, i.e., $r$ and $w$ are parallel, see \cite[Cor.\ 1.1.5]{SW}.  Thus  $r+w$ is light-like if and only
if $r$ and $w$ are parallel.

Consider next the case when  $(x,y,\theta)\in W_1\times W_0\times \R^4$ 
is such that $((x,y,\theta),d_{z,y,\theta}\Psi_2)\in \hbox{WF}(\tilde {\bf Q}_2^*)$
and   $\im \varphi(y)=0$.
Using the above notations (\ref{new notations}), we obtain
$
d_{z,y,\theta}\Psi_2=(\omega_\theta,
r+w;(0,0,0,d_{ \theta_4}\psi_4(y, \theta_4)))=(\omega_\theta,d\varphi(y)+(0,0,0,\theta_4);(0,0,0,y^4))$,
 where $y^4=d_{ \theta_4}\psi_4(y, \theta_4)=0$,
and thus we have
\ba
((z,\omega_\theta),(y,
r+w))\in \hbox{WF}( {\bf Q}_2^*)=\Lambda_{{\bf Q}^*_2}.
\ea

\begin{figure}[htbp] \label{Fig-2bb}
\begin{center}

\psfrag{1}{$(x_1,\xi_1)$}
\psfrag{2}{\hspace{-5mm}$(x_2,\xi_2)$}
\psfrag{3}{\hspace{-5mm}$(x_3,\xi_3)$}
\psfrag{4}{\hspace{-5mm}$(x_4,\xi_4)$}
\psfrag{5}{$(x_5,\xi_5)$}
\psfrag{6}{$z$}
\psfrag{7}{$y$}
\psfrag{8}{$r$}
\psfrag{9}{$w$}
\psfrag{0}{$\omega_\theta$}
\includegraphics[width=7.5cm]{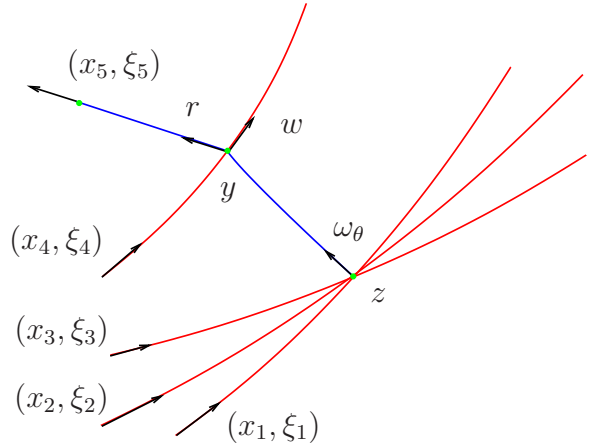}
\end{center}
\caption{A schematic figure where the space-time is represented as the  3-dimensional set $\R^{1+1}$. 
In the figure we consider the case A1 where three geodesics intersect at $z$ and 
the waves propagating near these geodesics interact and create a wave that hits
 the fourth geodesic at the point $y$, the produced singularities are detected by the gaussian beam
source at the point $x_5$. Note that $z$ and $y$ can be conjugate points on the geodesic connecting
them. In the case A2 the points $y$ and $z$ are the same.
}
 \end{figure}

As $\Lambda_{{\bf Q}^*_2}\subset \Lambda_{\hat g}\cup \Delta_{T\hattuM _0}^{\prime}$,
this implies that one of the following conditions are valid:
\ba
& &(A1)\ ((z,\omega_\theta),(y,
r+w))\in \Lambda_{\hat g}, \hspace{5cm} \
\\ & &\hspace{-0.5cm}\hbox{or}\\ 
& &(A2)\ ((z,\omega_\theta),(y,
r+w))\in \Delta_{T\hattuM _0}^{\prime}.
\ea
Let
$\gamma_0$ 
be the geodesic with
$
\gamma_0(0)=z,\ \dot \gamma_0(0) =\omega_\theta^\sharp.
$
Then (A1) and (A2) are equivalent to the following conditions:
\ba
& &(A1)\ \hbox{There is $t_0\in \R$ such
that $(\gamma_0(t_0), \dot \gamma_0(t_0)^\flat) =(y,\kappa)$ and},\\
& &\quad\quad \ \ \hbox{
 the vector $\kappa$ is  light-like,}
\\ & &\hspace{-0.5cm}\hbox{or}\\ 
& &(A2)\ z=y\hbox{ and }
\kappa=-\omega_\theta.
\ea

Consider next the case when (A1)  is valid. 
As $\kappa$ is light-like,  $r$ and $w$ are parallel. 
Then, as $(\gamma_{x_5,\xi_5}(t_1),\dot \gamma_{x_5,\xi_5}(t_1)^\flat)=(y,r)$ we see that $\gamma_0$ is a continuation of the geodesic $ \gamma_{x_5,\xi_5}$,
that is, for some $t_2$ we
have $(\gamma_{x_5,\xi_5}(t_2),\dot \gamma_{x_5,\xi_5}(t_2))=(z,\omega_\theta)
\in N^*K_{123}$. This implies that $x_5\in {\mathcal Y}$  that is not possible
by our assumptions. Hence (A1) is not possible.

Consider next the case when (A2)  is valid. If  we then would also have 
that $r\parallel w$  then $r$ is parallel to $\kappa=-\omega_\theta\in N^*K_{123}$,
 and as $(\gamma_{x_5,\xi_5}(t_1),\dot \gamma_{x_5,\xi_5}(t_1)^\flat)=(y,r)$ we
 would have  $x_5\in {\mathcal Y}$. As this is not possible by our assumptions,
 we see that $r$ and $w$ are not parallel. This implies that
 $\omega_\theta=-\kappa$ is not light-like.

For any given  $(\vec x,\vec \xi)$  and $(x_5,\xi_5)$ there exists
 $\e_2>0$ so that $(\{(y, \omega_\theta)\}\times T^*M)\cap \W_2(\e_2)=\emptyset$,
 see (\ref{eq: W_2 neighborhood}), and thus
 \beq\label{eq: y and omega}
 (\{(y, \omega_\theta)\}\times T^*M)\cap \hbox{WF}({\bf Q}_2^*)^{\prime}=\emptyset.
 \eeq
 Next we assume that $\e_2>0$ and also $\e_1\in (0,\e_2)$ are chosen
 so that  (\ref{eq: y and omega}) is valid.
 Then there are no $(z,y,\theta)$ such that 
  $((z,y,\theta),d\Psi_2(z,y,\theta))\in  \hbox{WF}({\bf Q}_2^*)$ 
 and  $\im \varphi(y)=0$. Thus by 
  Corollary 1.4 in \cite{Delort} or
  \cite[Lem.\ 4.1]{Ralston}
 yields  $T_{\tau,2}^{(4),\beta}=O(\tau^{-N})$ for all $N>0$.
Alternatively, one can use 
 the complex version of \cite[Prop. 1.3.2]{Duistermaat}, 
 obtained using combining the proof of  \cite[Prop. 1.3.2]{Duistermaat} and 
 the method of stationary phase with a complex phase, see \cite[Thm.\ 7.7.17]{H1}.
 
 Thus to analyze the asymptotics of  $T_{\tau}^{(4),\beta}$ 
 we need to consider only $T_{\tau,1}^{(4),\beta}$.
Next, we analyze the case when $U_4$ is a conormal
distribution and has the form (\ref{U1term}).

%
%

 Let us thus consider the case $p=1$. Now
 \beq
U_4(y)\,\cdotp u_{\tau}(y)&=&
\int_{\R^{1}}e^{i\theta_4y^4+i\tau \varphi(y)}
a_4(y,\theta_4)a_5(y,\tau )
\,d\theta_4. \label{QU4Ut term PRE}
\eeq
 Denoting by $\psi_1(z,t,\xi)=(y-z)\,\cdotp \xi$
 the phase function of the pseudodifferential operator ${\bf Q}^*_1$  we obtain
by (\ref{eq: representation of Q* p}) 
\beq\nonumber
({\bf Q}^*_1(U_4\,\cdotp u_{\tau}))(z)&=&
\int_{\R^{9}}e^{i(\psi_1(z,y,\xi)+\theta_4y^4+\tau \varphi(y))}q_1(z,y,\xi)\cdotp\\
& &\quad
\cdotp a_4(y,\theta_4)a_5(y,\tau )
\,d\theta_4\,dy d\xi. \label{QU4Ut term}
\eeq
Then
 $T_{\tau,1}^{(4),\beta}=T_{\tau,1,1}^{(4),\beta}+T_{\tau,1,2}^{(4),\beta}$, 
 cf.\ (\ref{eq: q_{1,1} and q_{1,2}}),
 where \HOX{In the final version
 of the paper we could shorten the discussion and just use $q_1(x,\xi)$
 instead of decomposition  (\ref{eq: q_{1,1} and q_{1,2}}).}
 \beq\label{eq; T tau asympt pre}
& &\hspace{-.4cm}T_{\tau,1,k}^{(4),\beta}
=
\int_{\R^{16}}e^{i(\theta_1z^1+\theta_2z^2+\theta_3z^3+\psi_1(z,y,\xi)+\theta_4y^4+\tau \varphi(y))}
c_1(z,\theta_1,\theta_2)\cdotp
\hspace{-1cm}
\\ \nonumber
& &\cdotp  a_3(z,\theta_3)
q_{1,k}(z,y,\xi)a_4(y,\theta_4) a_5(y,\tau )\,d\theta_1d\theta_2d\theta_3d\theta_4dydz d\xi,
\eeq
or
 \beq\label{eq; T tau asympt}
& &\hspace{-.1cm}T_{\tau,1,k}^{(4),\beta}=\tau^{8}
\int_{\R^{16}}e^{i\tau (\theta_1z^1+\theta_2z^2+\theta_3z^3+\psi_1(z,y,\xi)+\theta_4y^4+\varphi(y))}
c_1(z,\tau \theta_1,\tau \theta_2)\cdotp\hspace{-1cm}\\ \nonumber
& &\hspace{-0.3cm}
\cdotp  a_3(z,\tau \theta_3)
q_{1,k}(z,y,\tau \xi)a_4(y,\tau \theta_4) a_5(y,\tau )\,d\theta_1d\theta_2d\theta_3d\theta_4dydz d\xi.\hspace{-1cm}
\eeq

%
Let $(z,\theta,y,\xi)$ be a critical point of
 the phase function
  \beq\label{Psi3 function}
 \Psi_3(z,\theta,y,\xi)=\theta_1z^1+\theta_2z^2+\theta_3z^3+(y-z)\,\cdotp \xi+\theta_4y^4+\varphi(y).\hspace{-1cm}
 \eeq
Then
 \beq\label{eq: critical points} 
  & &\p_{\theta_j}\Psi_3=0,\ j=1,2,3\quad\hbox{yield}\quad z\in K_{123},\\
  \nonumber
& &\p_{\theta_4}\Psi_3=0\quad\hbox{yields}\quad y\in K_4,\\
  \nonumber
 & &\p_z\Psi_3=0\quad\hbox{yields}\quad\xi=\omega_\theta,\\
   \nonumber
& &\p_{\xi}\Psi_3=0\quad\hbox{yields}\quad y=z,\\
  \nonumber
 & &\p_{y}\Psi_3=0\quad\hbox{yields}\quad \xi=
 -\p_y\varphi(y)-w.
  \eeq
  The critical points we need to consider for the asymptotics satisfy also
\beq\label{eq: imag.condition} \\ \nonumber
\im \varphi(y)=0,\quad\hbox{so that }y\in \gamma_{x_5,\xi_5},\ \im d\varphi(y)=0,\
\re d\varphi(y)\in L^{*,+}_y\hattuM _0.\hspace{-1cm}
\eeq


Next\hiddenfootnote{Alternatively, to analyze  the term $T_{\tau,1,2}^{(4),\beta}$
we could use the fact that that  $e^{\Psi_3}=|\xi|^{-2} (\nabla_z-\omega_\theta)^2 e^{\Psi_3}$
where $\omega_\theta=(\theta_1,\theta_1,\theta_1,0)$. This can be used in the 
integration by parts in the $z$ variable. 
Possibly, we could use this to analyze directly
$T_{\tau,1}^{(4),\beta}$ and omit the decomposition ${\bf Q}^*_1={\bf Q}^*_{1,1}+{\bf Q}^*_{1,2}$.  
} we analyze      the terms $T_{\tau,1,k}^{(4),\beta}$ starting with $k=2$.
 Observe that 
 the 3rd and 5th equations in (\ref{eq: critical points})  imply that at the critical points
 $\xi=( \p_{y_1}\varphi(y), \p_{y_2}\varphi(y),$ $ \p_{y_3}\varphi(y),0)$.
 Thus the critical points are bounded in the $\xi$ variable. 
 Let us now fix the parameter $R$ determining $\psi_R(\xi)$ in
 (\ref{eq: q_{1,1} and q_{1,2}}) so that 
 $\xi$-components of the critical
 points are in a ball $B(R)\subset \R^4$.
  Using the identity  $e^{\Psi_3}=|\xi|^{-2} (\nabla_z-\omega_\theta)^2 e^{\Psi_3}$
where $\omega_\theta=(\theta_1,\theta_1,\theta_1,0)$
we can include the operator $|\xi|^{-2}(\nabla_z-\omega_\theta)^2$
in
the integral (\ref{eq; T tau asympt}) with $k=2$
and integrate by
parts. Doing this  two times we can show
that this oscillatory  integral (\ref{eq; T tau asympt}) with $k=2$ becomes an integral of a Lebesgue-integrable
 function. 
 Then, by using method of stationary phase
 and  
%
  the fact that $\psi_R(\xi)$ vanishes at all critical points of $\Psi_3$
 where $\im \Psi_3$ vanishes, we see that   $T_{\tau,1,2}^{(4),\beta}=O(\tau^{-n})$
 for all $n>0$.  
%
%


Above, we have shown that  the term $T_{\tau,1,1}^{(4),\beta}$  has the same
asymptotics as $T_\tau^{(4),\beta}$. Next we analyze  this term.
Let $(z,\theta,y,\xi)$ be a critical point
of $ \Psi_3(z,\theta,y,\xi)$  such that $y$ satisfies (\ref{eq: imag.condition}).
Let us next use the same notations (\ref{new notations}) which
we used above.
Then (\ref{eq: critical points}) and (\ref{eq: imag.condition}) imply
\beq\label{eq: psi1 consequences}
z=y\in \gamma_{x_5,\xi_5}\cap \bigcap_{j=1}^4 K_j,\ \xi=\omega_\theta=-r-w.
\eeq
Note that in this case all the four geodesics $\gamma_{x_j,\xi_j}$ intersect
at the point $q$ and by our assumptions, $r=d\varphi(y)$ is such a co-vector
that in the $Y$-coordinates $r=(r_j)_{j=1}^4$ with $r_j\not =0$ for all $j=1,2,3,4$.
In particular, this shows that
 the existence of the critical point  of $ \Psi_3(z,\theta,y,\xi)$
implies that there exists 
an intersection point of $\gamma_{x_5,\xi_5}$ and $ \bigcap_{j=1}^4 K_j$.
Equations (\ref{eq: psi1 consequences}) imply also that 
\ba
r=\sum_{j=1}^4r_jdy^j=-\omega_\theta-w=-\sum_{j=1}^3\theta_jdz^j-
\theta_4 dy^4.
\ea

%


To consider the case when $y=z$, let us assume for a while that that $W_0=W_1$
and that the $Y$-coordinates and $Z$-coordinates coincide, that is, $Y(x)=Z(x)$.
Then the covectors $dz^j=dZ^j$ and $dy^j=dY^j$ coincide for $j=1,2,3,4$.
Then we have
\beq\label{eq: r is -theta}
r_j=-\theta_j,\quad\hbox{i.e.,\ } \theta:=\theta_jdz^j=-r=r_jdy^j\in T^*_y\hattuM _0.
\eeq


Let us apply
the method of stationary phase to $T^{(4),\beta}_{\tau,1,1}$ as $\tau\to \infty$. 
Note that
as $c_1(z,\theta_1,\theta_2)$ is a product type symbol, we need to use the fact $\theta_1\not=0$ and $\theta_2\not=0$
for the critical points as we have by (\ref{eq: r is -theta}) 
and the fact that $r=d\varphi(y)|_{y}\not \in N^*K_{234}\cup N^*K_{134}$
as $x_5\not \in  {\mathcal Y}$ and assuming that $s_0$ is small enough.

{
In local $Y$ and $Z$ coordinates where $z=y=(0,0,0,0)$ we
can use the method of  stationary phase,
 similarly to proofs of \cite[Thm.\ 1.11 and 3.4]{GrS},
 to compute 
the asymptotics of  (\ref{eq; T tau asympt}) with $k=1$.
Let us explain this  computation in detail.
To this end, let us start with some preparatory considerations.

 \HOX{Jan. 4, 2013: Remove these details from the final version of
the paper. }
Let $\phi_1(z,y,\theta,\xi)$ be a smooth bounded
 function in $C^\infty(W_0\times W_1\times \R^4\times \R^4)$ that is 
homogeneous of degree zero in the $(\theta,\xi)$ variables in the set  $\{|(\theta,\xi)|>R_0\}$ with some $R_0>0$. 
Assume that $\phi_1(z,y,\theta,\xi)$ is equal to one  in 
a conic neighborhood, with respect to $(\theta,\xi)$, of the points where some of the 
$\theta_j$ or $\xi_k$ variable is zero. Note that in this set 
the positively homogeneous  functions $a_j(z,\theta_j/|\theta_j|)$ may be non-smooth.
Let $\Sigma\subset W_0\times W_1\times \R^4\times \R^4$ 
be a conic 
 neighborhood 
of the critical points of the phase function $\Psi_3(y,z,\theta,\xi)$.
Also, assume that the function  $\phi_1(y,z,\theta,\xi)$ vanishes in 
the intersection of $\Sigma$ 
and the set  $\{|(\theta,\xi)|>R_0\}$.
Let  \HOX{Slava asked why we need $\Psi_{(\tau)}$ and not only $\Psi_3$,
the reason is to do some computations with the formula (\ref{eq; T tau asympt pre}) and
and some with  (\ref{eq; T tau asympt}). It seems that this can be done when
$a_j$ are smooth near origin. How about smoothness of $c$? }
\ba
\Psi_{(\tau)}(z,y,\theta,\xi)= 
\theta_1z^1+\theta_2z^2+\theta_3z^3+\psi_1(z,y,\xi)+\theta_4y^4+\tau \varphi(y)\ea 
 be the phase function appearing in (\ref{eq; T tau asympt pre}).
 Note that $\Sigma$ contains  the critical points of $\Psi_{(\tau)}$ for all $\tau>0$.
 Let
$$L_\tau=\frac {\phi_1(z,y,\tau^{-1}\theta,\tau^{-1} \xi)}{ |d_{z,y,\theta,\xi}
\Psi_{(\tau)}(z,y,\theta,\xi)|^{2}}\, (d_{z,y,\theta,\xi}
\overline{\Psi_{(\tau)}(z,y,\theta,\xi)})\cdotp d_{z,y,\theta,\xi},$$
so that 
\beq\label{kaav A}
L_\tau \exp(\Psi_{(\tau)})= \phi_1(z,y,\tau^{-1} \theta,\tau^{-1}\xi)\,\exp(\Psi_{(\tau)}).
\eeq
Note that if $\im \varphi(y)=0$, then $y\in \gamma_{x_5,\xi_5}$
and hence $d\varphi(y)$ does not vanish and that when $\tau$ is large enough,  the function $\Psi_{(\tau)}$ has 
no critical points in 
the support of $\phi_1$. Using
these we see that $ \gamma_{x_5,\xi_5}\cap W_1$
has a neighborhood $V_1\subset W_1$ where $|d\varphi(y)|>C_0>0$ and
there are $C_1,C_2,C_3>0$ so that if
  $\tau>C_1$, $y\in V_1$, and $(z,y,\theta,\xi)\in \supp(\phi_1)$
then
\beq\label{str est1}
\,|d_{z,y,\theta,\xi}
\Psi_{(\tau)}(z,y,\theta,\xi)|^{-1}\leq \frac {C_2}{\tau-C_3}.
\eeq
After these preparatory steps, we are ready to  compute the asymptotics of
 $T^{(4),\beta}_{\tau,1,1}$. To this end, we first transform the integrals in  (\ref{eq; T tau asympt})
 to an integral of a Lebesgue integrable function by using integration by parts of $|\xi|^{-2}(\nabla_z-\omega_\theta)^2$
 as explained above. Then
 we decompose  $T^{(4),\beta}_{\tau,1,1}$ into three
 terms $T^{(4),\beta}_{\tau,1,1}=I_1+I_2+I_3$. To obtain the first term $I_1$ we include the factor
 $(1-\phi_1(z,y,\theta,\xi))$ in the  integral  
 (\ref{eq; T tau asympt}) with $k=1$. 
  The integral $I_1$ can then be computed
 using the method of stationary phase similarly to the proof of \cite[Thm.\ 3.4]{GrS}. 
 Let $\chi_1\in C^\infty(W_1)$ be a function that 
 is supported in $V_1$ and vanishes on $\gamma_{x_5,\xi_5}$. 
 The terms $I_2$ and $I_3$
 are obtained by 
 including the factor
 $\phi_1(z,y,\tau^{-1}\theta,\tau^{-1}\xi)\chi_1(y)$ and  $\phi_1(z,y,\tau^{-1}\theta,\tau^{-1}\xi)(1-\chi_1(y))$
  in the  integral 
 (\ref{eq; T tau asympt pre}) with $k=1$, respectively.
(Equivalently, 
 the terms $I_2$ and $I_3$
 are obtained by 
 including the factor
 $\phi_1(z,y,\theta,\xi)\chi_1(y)$ and  $\phi_1(z,y,\theta,\xi)(1-\chi_1(y))$
  in the  integral 
 (\ref{eq; T tau asympt}) with $k=1$, respectively.)
 Using integration by parts in integral (\ref{eq; T tau asympt pre}) and inequalities  (\ref{kaav A}) and (\ref{str est1}), we see that that $I_2=O(\tau^{-N})$
 for all $N>0$. Moreover, the fact that $\im \varphi(y)>c_1>0$ in  
$W_1\cap V_1$ implies that $I_3=O(\tau^{-N})$
 for all $N>0$.}


 Combining the above we obtain
the asymptotics
 \beq\label{definition of xi parameter1}
& &\quad \quad T^{(4),\beta}_\tau
\sim \tau^{4+4-16/2-2+\rho-2}\sum_{k=0}^\infty c_k
\tau^{-k} =\tau^{-4+\rho}\sum_{k=0}^\infty c_k\tau^{-k}
, \\ \nonumber
& &\hspace{-8mm}c_0=h(z^{(0)})c_1(0,-r_1,-r_2)\tilde a_3(0,-r_3)
\tilde q_1(0,0,-(r_1,r_2,r_3,0))\tilde a_4(0,-r_4)\tilde a_5(0,1),\hspace{-1cm}
\eeq
where $\rho=\sum_{j=1}^5 p_j$
and $\tilde a_j$ is the principal
symbol of $a_j$ etc.   The factor $h(z^{(0)})$ is non-vanishing and is
determined by the determinant of the Hessian of the phase function
$\varphi$ at $q$.
A direct computation shows that $\det(\hbox{Hess}_{z,y,\theta,\xi} 
\Psi_{3}(z^{(0)},y^{(0)},\theta^{(0)},\xi^{(0)})=1$. Above,
$(z^{(0)},y^{(0)},\theta^{(0)},\xi^{(0)})$ is the critical point satisfying
(\ref{eq: critical points}) and (\ref{eq: imag.condition}), 
where in the local coordinates $(z^{(0)},y^{(0)})=(0,0)$ and
$h(z^{(0)})$ is constant times powers of values of the cut-off functions $\Phi_0$
and $\Phi_1$ at zero. Recall that we considered above the case
when $\mathcal B_j$ are multiplication operators with these cut-off functions.
 The term $c_k$ depends on the 
derivatives of the symbols $a_j$ and $q_1$ of order less
or equal to $2k$ at the critical point.
If  $ \Psi_3(z,\theta,y,\xi)$  has no critical points, that is,
 $q$ is not an intersection point, we
obtain the asymptotics $T^{(4),\beta}_{\tau,1,1}=O(\tau^{-N})$ for all $N>0$.

 For  future reference we note that if we use the method of  stationary
 phase in the last integral of (\ref{eq; T tau asympt}) only in the integrals with respect
 to $z$ and $\xi$, yielding that at the critical point we have  $y=z$ and $\xi=\omega_\beta(\theta)=(\theta_1,\theta_2,\theta_3,0)$, we 
 see that 
 $T_{\tau,1,1}^{(4),\beta}$ can
 be written as
 \beq \label{eq; T tau asympt modified}
& &T_{\tau,1,1}^{(4),\beta}
=c\tau^{4}
\int_{\R^{8}}e^{i(\theta_1y^1+\theta_2y^2+\theta_3y^3+\theta_4y^4)+i\tau \varphi(y)}
c_1(y, \theta_1, \theta_2)\cdotp\\ \nonumber
& &
\cdotp a_3(y, \theta_3)
q_{1,1}(y,y,\omega_\beta(\theta))a_4(y,\theta_4) a_5(y,\tau )\,d\theta_1d\theta_2d\theta_3d\theta_4dy\\ \nonumber
& &=c\tau^{8}
\int_{\R^{8}}e^{i\tau (\theta_1y^1+\theta_2y^2+\theta_3y^3+\theta_4y^4+\varphi(y))}
c_1(y, \tau \theta_1, \tau \theta_2)\cdotp\\ \nonumber
& &
\cdotp a_3(y,\tau  \theta_3)
q_{1,1}(y,y,\tau  \omega_\beta(\theta))a_4(y,\tau \theta_4) a_5(y,\tau )\,d\theta_1d\theta_2d\theta_3d\theta_4dy.
\eeq

Next, we consider the terms $\tilde T^{(4),\beta}_\tau$ of the type
(\ref{tilde T-type source}). Such term
is an integral of the product of $u_\tau$ and two other factors 
${\bf Q}(U_2\,\cdotp U_1)$ and ${\bf Q}(U_4\,\cdotp U_3)$.
As the last two factors can be written in the form (\ref{Q U1U2term}),
one can see using the method of stationary phase that  $\tilde T^{(4),\beta}_\tau$ 
has similar asymptotics to  $T^{(4),\beta}_\tau$ as $\tau\to \infty$,
with the leading order coefficient $\tilde c_0=\tilde h(z^{(0)})\tilde c_1(0,-r_1,-r_2)
\tilde c_2(0,-r_3,-r_4)$ where $\tilde c_2$ is given as in (\ref{product type symbols})
with symbols $\tilde a_3$ and $\tilde a_4$, and moreover,
$\tilde h(z^{(0)})$ is a constant times powers of values of the cut-off functions $\Phi_0$
and $\Phi_1$ at zero.

This proves the claim in the special case where
$u_j$ are conormal distributions supported in the coordinate neighborhoods $W_{k(j)}$,
$a_j$ are positively homogeneous scalar valued symbols, $S_j={\bf Q}$, and $\B_j$ are 
multiplication functions with smooth cut-off functions.
By using a suitable partition of unity and 
summing the results of the  above computations, 
similar results to the above follows when  
$a_j$ are general classical symbols that are $\mathcal B$-valued and 
the waves $u_j$ are supported on
$J^+_{\hat g}(\supp({\bf f}_j))$.
Also, $S_j$ can be replaced by operators of type (\ref{extra notation}) and
$\B_j$ can replaced by differential operator without other essential changes
expect that the highest order power of $\tau$ changes.
Then, in the asymptotics of terms $T^{(4),\beta}_\tau$  the function
$h(z^{(0)})$ in (\ref{definition of xi parameter1}) is a section in  dual bundle $(\B_L)^4$. The coefficients
of $h(z^{(0)})$ in local coordinates are polynomials of $\hat g^{jk}$, $\hat g_{jk}$, 
$\hat \phi_{\ell}$, and 
their derivatives at $z^{(0)}$. Similar representation is obtained for the asymptotics
 of terms $\tilde T^{(4),\beta}_\tau$.

As we integrated by parts two times the operator
$(\nabla_z-\omega_\theta)^2$ and the total order of 
of $\B_j$ is less or equal to 6, we see that it is
enough to assume above that the symbols $a_j(z,\theta_j)$
are of  order $(-12)$ or less.
The leading order asymptotics come from the term \HOX{Clarify explanation on $p_j$.} 
where the sum of orders of $\B_j$ is 6 and $p_j=n$ for $j=1,2,3,4$,
$p_5=0$ so that $m=-4-\rho+6=4n+2$. 
We also see that the terms containing permutation $\sigma=\sigma_\b$ of the indexes
of the spherical waves can be analyzed analogously. This proves (\ref{indicator}).

Making the above computations explicitly, we obtain an explicit formula for
the leading order coefficient $s_{m}$ in (\ref{indicator})  in terms of
$\bsequence$ and ${\bf w}$,
multiplied with the power $(-1/2)$ of the determinant of the Hessian of the phase function
$\Psi_3$ at the critical point $q$ in the $Z$-coordinates determined 
by $(\vec x,\vec \xi)$ and $(x_5,\xi_5)$. This show that $s_{m}$
coincides with some real-analytic function $G(\bsequence,{\bf w})$ multiplied
by a non-vanishing function $\mathcal R(p,(\vec x,\vec \xi),(x_5,\xi_5))$,
corresponding to the power of the Hessian, that depends on the phase function $\varphi$ in the $Z$ coordinates.
%
%
This proves the claim in the case when the linear independency condition (LI) is valid.

Next, consider the case when the linear independency condition (LI) is not valid.
Again, by the definition of ${\bf t}_j$,  if
the intersection $\gamma_{x_5,\xi_5}(\R_-) \cap(\cap_{j=1}^4\gamma_{x_j,\xi_j}((0,{\bf t}_j)))
$ is non-empty, it can
contain only one point. In the case that such a
point exists, we denote it by $q$.

When (LI) is not valid, we have that the linear space $\hbox{span}(b_j;\ j=1,2,3,4)\subset T_q^*M_0$ has dimension  3 or less.
We use the facts that for $w\in I(\Lambda_1,\Lambda_2)$
we have WF$(w)\subset \Lambda_1\cup \Lambda_2$
and the fact, see \cite[Thm.\ 1.3.6]{Duistermaat}  
 \ba
\hbox{WF}(v\,\cdotp w)\subset 
\hbox{WF}(v)\cup \hbox{WF}(w)
\cup\{(x,\xi+\eta);\ (x,\xi)\in \hbox{WF}(v),\ (x,\eta)\in \hbox{WF}(w)\}.
\ea
Let us next consider the terms corresponding to the permutation $\sigma=Id$.
The above facts imply that $\tilde {\mathcal G}^{(4),\beta}$ 
in (\ref{w1-3 solutions}) satisfies
\ba
\hbox{WF}(\tilde {\mathcal G}^{(4),\beta})
\cap T_q^*M_0\subset \mathcal Z_{s_0}:=\X_{s_0}\cup \bigcup _{1\leq j\leq 4} N^*K_j\cup \bigcup_{1\leq j<k\leq 4} N^*K_{jk}, 
\ea 
where $\X_{s_0}=\X((\vec x,\vec\xi);t_0,s_0)$.
Also, for
\ba
w_{123}=\B_3^{\beta}u_{3}\,\cdotp \mathcal C_1^{\beta} S^{\beta}_1(\B^{\beta}_2u_{2}\,\cdotp \B^{\beta}_1u_{1}),
\ea
appearing in (\ref{M4 terms}), we have WF$(w_{123})
\subset \mathcal Z_{s_0}$ and thus using H\"ormander's theorem \cite[Thm.\ 26.1.1]{H4}, we see
that  WF$(S_2^{\beta}(w_{123}))\subset \Lambda^{(3)}$,
where  $\Lambda^{(3)}$ is the flowout of $\mathcal Z_{s_0}$ in the canonical relation 
of ${\bf Q}$. Then
\ba
\pi(\Lambda^{(3)})\subset \Y_{s_0}\cup \bigcup _{1\leq j\leq 4} K_j\cup \bigcup_{1\leq j<k\leq 4} K_{jk}, 
\ea 
where $\Y_{s_0}=\Y((\vec x,\vec\xi);t_0,s_0)$ and $\pi:T^*M_0\to M_0$ is the projection to the
base point.

Observe that $E=\hbox{span}(b_j;\ j=1,2,3,4)\subset T_q^*M_0$
has dimension 3 or less, $\Lambda^{(3)}\cap T_q^*M_0\subset E$
and WF$(u_{4})\cap T_q^*M_0\subset E$. Thus,
${\mathcal G}^{(4),\beta}=\B_4^{\beta}u_{4}\,\cdotp C_2^{\beta}S_2^{\beta}(w_{123})$ 
satisfies WF$({\mathcal G}^{(4),\beta})
\cap T_q^*M_0\subset E$.  Now, $E\subset \mathcal Z_{s_0}.$ By our assumption,
$(q,b_5)\not \in \mathcal X((\vec x,\vec\xi);t_0)$, and thus, cf.\ (\ref{w1-3 solutions}),  we see that 
\ba
\bra u_\tau,{\bf Q}(\sum_{\b\leq n_1} {\mathcal G}^{(4),\beta}+\tilde {\mathcal G}^{(4),\beta})\cet=O(\tau^{-N})
\ea 
for all $N>0$ when $s_0$ is small enough.  The terms where
the permutation $\sigma$ is not the identity can be analyzed similarly.
This proves the claim in the case when the linear independency condition (LI) is not valid.
\hfill \Box \medskip

\generalizations{
 \begin{proposition}\label{lem:analytic limits B}
 Let the assumptions of Proposition \ref{lem:analytic limits A} are valid.
  Moreover, let $\ell\leq 3$. Then
   $|\Theta^{(\ell)}_\tau|\leq C_N\tau^{-N}$
 for all $N>0$.
 \end{proposition}

\noindent {\bf Proof.} 
The cases $\ell=1$ and
$\ell=2$ follow from the fact that  when $s_0$ is small enough, $x_5\not\in K_j$ for $j=1,2,3,4$
and that  the waves $u_j$ restricted in the set 
$\V((\vec x,\vec \xi),t_0)$, see  (\ref{eq: summary of assumptions 2}),
 are 
conormal  distributions
associated to the surfaces $K(x_j,\xi_j;t_0,s_0)$
 that 
do not intersect $\dot\gamma_{x_5,\xi_5}$
  Next we consider $ \ell=3$.

Note that as $b_j$, $j=1,2,3,$ are light-like, they are linearly independent
only when at least two of them are parallel. As this is not possible by
our assumption that $\gamma_{x_j,\xi_j}(\R)$ do not coincide, we can
assume that $b_j$, $j=1,2,3,$ are linearly independent.
In the case where $\B_j=\Phi_0$, $S_1={\bf Q}$, and $W_k$ 
are as above. Similar computations
to the above as in (\ref{U1term})  and (\ref{Q U1U2term}) give  
\beq \nonumber
T^{(3),\beta}_\tau&=&\bra u_\tau,U_3\,\cdotp {\bf Q}(U_2\,\cdotp U_1)\cet_{L^2(\hattuM _0)},
\\ \label{eq: 3rd order term}
&=&\tau^{3}
\int_{\R^{7}}e^{i\tau (\theta_1z^1+\theta_2z^2+\theta_3z^3+\varphi(z))}
c_1(z,\tau \theta_1,\tau \theta_2)\cdotp\\ \nonumber
& &\quad 
\cdotp  a_3(z,\tau \theta_3)a_5(z,\tau )\,d\theta_1d\theta_2d\theta_3dz_1dz_2dz_3dz_4 .
\eeq
%
Assume that the phase function of this  integral has a critical point $(\theta^0,q_0)$.
Then  $z^1_0=z^2_0=z^3_0=0$ and
$\theta_j^0+\p_j\varphi(z_0)=0$ for $j=1,2,3$,
%
and $\p_4\varphi(z_0)=0$.
If at the critical point $\im \varphi(z_0)=0$,  then there is $t_0=t_0(x_5,\xi_5)\in \R$ such that 
 $z_0=\gamma_{x_5,\xi_5}(t_0)\in K_{123}$ 
 and $\zeta=d\varphi(z_0)^\flat=c\dot \gamma_{x_5,\xi_5}(t_0)$.
 Then $\p_4\varphi(z_0)=0$ implies that $\zeta^\flat\in N^*_{z_0}K_{123}$.
This can not hold when $s_0$ is small enough as
we have $x_5\not \in {\mathcal Y}(((x_{j},\xi_{j}))_{j=1}^4;t_0)$.
Hence, $T^{(3),\beta}_\tau=O(\tau^{-N})$ for all $N>0$.
Similar analysis with different $\B_j$ and $S_j$ show  the claim  in the case $\ell=3$.
\hfill \Box \medskip
}

 \begin{proposition}\label{lem:analytic limits C}
 Let the assumptions of Proposition \ref{lem:analytic limits A} be valid.
Moreover, assume that
$(x_5,\xi_5)$ comes from the 4-intersection
of  rays corresponding to $(\vec x,\vec \xi)$ and $q$ is the corresponding
intersection point. 
Then the point $x_5$ has a neighborhood $V$ so that 
$\M^4$ in $V$ satisfies $\M^4|_V\in  \I(V; \Lambda^+_q)$.
%
\end{proposition}

\noindent{\bf Proof.}
 Let us fix $(x_j,\xi_j)$, $j\leq 4$, $s_0$, and the
 waves $u_j\in \I(K_j)$, $j\leq 4$.

Let us consider the same condition (LI) that was used in the proof of Prop. \ref{lem:analytic limits A}. 
First we observe that if (LI) is not valid, we see using the proof of Prop. \ref{lem:analytic limits A}, (see the end of the proof where the
case when (LI) is not valid is considered) that
$\M^{(4)}$ is $C^\infty$ smooth in $\V\setminus (\Y\cup \bigcup_{j=1}^4 K_j)$.
Thus to prove the claim of the proposition we can assume that (LI) is valid.

Let us decompose  $\F^{(4)}$, given by
 (\ref{w1-3 solutions}) and (\ref{M4 terms})-(\ref{tilde M4 terms}) 
 as  $\F^{(4)}=\F_{1}^{(4)}+\F_{2}^{(4)}$
 where $\F^{(4)}_{p}$ is defined similarly to
 $\F^{(4)}$ in  (\ref{w1-3 solutions}) and 
 (\ref{M4 terms})-(\ref{tilde M4 terms})
by modifying these formulas so that 
the operator $S_1^\beta$ is replaced by 
$S_{1,p}^\beta$, where $S_{1,p}^\beta={\bf Q}_p$, when
 $S_{1}^\beta={\bf Q}$, and $S_{1,p}^\beta=(2-p)I$, when $S_{1}^\beta=I$.
 Here, the operators ${\bf Q}_p$  are defined as above using the parameters
  $\e_2$ and $\e_1$ defined below.
 
%

 Using formulae (\ref{w1-3 solutions}), (\ref{U1term}), (\ref{Q U1U2term}), and (\ref{QU4Ut term})
 we see that near near $q$ in the $Y$ coordinates $\M_1^{(4)}={\bf Q} \F_1^{(4)}$ 
 can be calculated using that
 \beq\label{f4-formula}
 \F_1^{(4)}(y)=\int_{\R^4} e^{iy^j\theta_j}b(y,\theta)\,d\theta,
 \eeq
 where $K_j$ in local coordinates is given by $\{y^j=0\}$ and
  $b(y,\theta)$  is a finite sum of terms  that are products of some of the following
 terms: at most one  product type symbol $c_l(y,\theta_j,\theta_k)\in S(W_0;\R\times (\R\setminus \{0\}))$ (they appear in the terms  (\ref{T-type source})-(\ref{tilde T-type source})  where the $S_j^\beta$ operators are ${\bf Q}$ 
 and do not appear if these operators are the identity),
 and one ore more term which is either
 the symbols $a_j(y,\theta_j)\in S^n(W_0;\R)$,
or the functions $q_1(y,y,\omega_\beta(\theta))$,
  cf.\  (\ref{eq; T tau asympt modified}),
 where $\omega_\beta(\theta)$ is equal to some of the vectors
$(\theta_1,\theta_2,\theta_3,0)$,
 $(\theta_1,\theta_2,0,\theta_4)$, $(\theta_1,0,\theta_3,\theta_4)$, or
 $(0,\theta_1,\theta_3,\theta_4)$, depending on the permutation $\sigma$.

Let us consider next the source $F_\tau$ is determined by the functions $(p,h)$ in (\ref{Ftau source}). 
 Then using the  method of  stationary phase gives  the asymptotics, c.f. (\ref{eq; T tau asympt modified}),
 \ba
 \bra u_\tau,\F_1^{(4)}\cet \hspace{-1mm}\sim \hspace{-1mm}\tau^8\hspace{-1mm}\int_{\R^8} e^{i\tau (\varphi(y)+y^j\theta_j)}(a_5(y,\tau),b(y,\tau \theta))_{\hat G}d\theta dy
\sim \sum_{k=m}^\infty  
 s_k(p,h) \tau^{-k}
 \ea
where $\hat G$ is a Riemannian metric
of the fiber of $\B^L$ at $y$, that is isomorphic to $\R^{10+L}$, and
the critical point of the phase function is $y=0$ and $\theta=-d\varphi(0)$.
  As we saw above, we have that
 when $\e_2>0$ is small enough then for $p=2$ we have
$\bra u_\tau,\F_p^{(4)}\cet=O(\tau^{-N})$ for all $N>0.$
%

Let us choose sufficiently
small $\e_3>0$ and choose a function $\chi(\theta)\in C^\infty(\R^4)$ that 
vanishes in a $\e_3$-neighborhood (in the $\hat g^+$ metric) of
$\mathcal A_q$, 
\beq\label{set Yq}
\mathcal A_q&:=&
N_q^*K_{123}
\cup N_q^*K_{134}\cup N_q^*K_{124}\cup N_q^*K_{234}
\eeq
and is equal to 1 outside the $(2\e_3)$-neighborhood of this set.

Let $\phi\in C^\infty_0(W_1)$ be a function that is one near $q$.
Also, let  \ba\hbox{$b_0(y,\theta)=\phi(y)\chi(\theta)b(y,\theta)$}\ea be a classical symbol,
$p=\sum_{j=1}^4 p_j$,
and let  $
\F^{(4),0}(y)
\in \I^{p-4}(q) 
$ be the
 conormal distribution 
 that is given by the formula (\ref{f4-formula}) with
$b(y,\theta)$ being replaced by 
 $b_0(y,\theta)$.

When $\e_3$ is small enough (depending on the point $x_5$), we see that 
$F_\tau$ is determined by functions $(p,h)$ 
and the corresponding 
gaussian beams $u_\tau$
propagating on the geodesic $\gamma_{x_5,\xi_5}(\R)$ such
that the geodesic
passes through 
$x_5\in V$, we have 
\ba
\bra u_\tau,\F^{(4),0}\cet\sim \sum_{k=m}^\infty s_k(p,h) \tau^{-k},
\ea
that is, we have $\bra u_\tau,\F^{(4),0}\cet-\bra u_\tau,\F^{(4)}\cet=O(\tau^{-N}$
for all $N$.
When $\gamma_{x_5,\xi_5}(\R)$  does not pass through $q$,
we have that  $\bra u_\tau,\F^{(4)}\cet$ and $\bra u_\tau,\F^{(4),0}\cet$ 
are both of order $O(\tau^{-N})$ for all $N>0$.

Let $V\subset \V((\vec x,\vec \xi),t_0)\setminus \bigcup_{j=1}^4\gamma_{x_j,\xi_j}([0,\infty))$, see   (\ref{eq: summary of assumptions 2}), be an open set.
By varying  the source $F_\tau$,  defined in  (\ref{Ftau source}),  
 we see, by multiplying the solution  with a smooth
 cut of function and  using  Corollary 1.4 in \cite{Delort} in local coordinates,
 or \cite{MelinS},
%
 we see 
 that the function
$\M^4-{\bf Q}\F^{(4),0}$ has no wave front set in $T^*(V)$ and it is
thus $C^\infty$-smooth function in $V$.

As by \cite{GU1}, ${\bf Q}:\I^{p-4}(\{q\})\to \I^{p-4-3/2,-1/2}(N^*(\{q\}),\Lambda_q^+)$, the above implies that
\beq\label{eq: conormal}
\M^4|_{V\setminus {\mathcal Y}}\in \I^{p-4-3/2}(V\setminus  {\mathcal Y};\Lambda_q^+),
\eeq
where $  {\mathcal Y}=  {\mathcal Y}((\vec x,\vec \xi),t_0,s_0)$. When $x_5$ if fixed, choosing $s_0$ to
be small enough, we obtain the claim.
\hfill \Box \medskip

Next we will show that the function $\mathcal G$ is not identically vanishing.

%

\generalizations
{{\bf Remark 3.6.A.} 
Similar considerations to the 
proof of  Prop.\ \ref{lem:analytic limits C}, can be used to show that $\M^3|_{(U\setminus \cup_j K_j)\cap I^-(x_5)}$
is a lagrangian distribution. As an example of this, let us analyze the term
\ba
m={\bf Q}(\B_3 u_3\,\cdotp v_{12}),\quad \hbox{where }v_{12}={\bf Q}(\B_2 u_1\,\cdotp \B_2 u_3),
\ea
where $u_j\in \I^{p_j}(K_j)$ are solutions
of the linear wave equation.
Let $S=(\bigcup_j N^*K_j)\cup(\bigcup _{j,k}N^*K_{jk})$.
Next, we do computations in local coordinates $X:V\to \R^n$,  $x^j=X^j(x)$, such that $V\subset  I^-(x_6)$ 
and $K_j\cap V =\{x\in V;\ x^j=0\}$. 
As $v_{12}=v_{12}^1+v_{12}^2\in \I(K_1,K_{12})+\I(K_2,K_{12})$, we can write the function
$f=\B_3 u_3\,\cdotp v_{12}$  in
$V$ as an oscillatory integral
\beq\label{def; f}
f(x)=\int_{\R^3}e^{i(\theta_1x^1+\theta_2x^2+\theta_3x^3)}b(x,\theta_1,\theta_2,\theta_3)d\theta_1 d\theta_2d\theta_3.
\eeq
Here, $b(x,\theta_1,\theta_2,\theta_3)\in S_{cl}^{p_1+p_2+p_3-2}(V;\R^3)$ away from
the set where some  some of the variables $\theta_j$ vanish and $\xi=\sum_{j=1}^3 \theta_jdx^j$
is light-like.
Let $\eta(x,\theta_1,\theta_2,\theta_3)=\phi(x,\sum_{j=1}^3 \theta_jdx^j),$
where $\phi$ is  a cut-off function
that
vanishes  in a conic neighborhood $S$ 
and is one outside a suitable conic neighborhood of $S$.
Then 
 $c(x,\theta_1,\theta_2,\theta_3)=\eta(x,\theta_1,\theta_2,\theta_3)b(x,\theta_1,\theta_2,\theta_3)$ is a classical symbol. Writing 
$f=f_1+f_2$ where $f_1$ and $f_2$ have the representation 
(\ref{def; f}) with the symbol $b$ is replaced by $c$ and $b-c$, correspondingly.
Then $f_1\in \I^{p_1+p_2+p_3-2}(K_{123})= \I^{p_1+p_2+p_3-5/2}(N^*K_{123})$ and 
we see that ${\bf Q}f_1\in
\I^{p_1+p_2+p_3-4}(Y;\tilde \Lambda),$  is a lagrangian distribution on  $Y=I^-(x_5)\setminus ( \cup_j K_j)$,
where $\tilde \Lambda$ is the flow-out of $N^*K_{123}$,
and ${\bf Q}f_2$ is a $C^\infty$-smooth in  $Y$.
Hence, 
\beq\label{3-singularity order}
\M^3|_{Y}\in \I^{\tilde p-10}(Y;\tilde \Lambda),\quad \tilde p=p_{1}+p_{2}+p_{3}.
\eeq
}

\subsubsection{WKB computations and the indicator functions in the Minkowski space}

To show that the function (\ref{definition of G})  is not identically vanishing, we will consider 
waves in Minkowski space.

In this section $\hat g_{jk}=\diag(-1,1,1,1)$ denotes the metric
in the standard coordinates of the Minkowski space $\R^4$. Below we call
the principal symbols of the linearized waves the polarizations to emphasize their physical meaning.
To show that $\mathcal  G(\bsequence,{\bf w})$ is non-vanishing, recall that ${\bf w}=(w_j)_{j=1}^5$
where for $j\leq 4$ the polarizations $w_j=(v_j,v_j^{\prime})$,
represented as a pair of metric and scalar field polarizations, the metric part of the polarization $v_j$
has to satisfy 4 linear conditions. 
Because of this, below we study the case when all polarizations of the matter fields at $q$ vanish,
that is, $v^{\prime}_j=0$ for all $j$, and $v_j$ satisfies 4 conditions. 
In this case, in Minkowski space the function $\mathcal G(v,0,\bsequence)$
can
be analyzed by assuming that there are no matter fields, which we do next.
Later we return to the case  of general  polarizations.
Next, we denote  the $g$-components amplitudes by $v_j=v^{(j)}$.
%


Instead of the indicator function $\Theta^{(4)}({\bf v},0,\bsequence)$ given in (\ref{indicator})
para\-metr\-iz\-ed by $({\bf v},0,\bsequence)$ we will consider the function parametrized by
the variables $({\bf v},\bsequence)$, where  $\bsequence=(b^{(j)})_{j=1}^5$, $b^{(j)}\in \R^4$,
are as before but ${\bf v}=(v^{(j)})_{j=1}^5$ are the amplitudes of the linear 
waves taking
values in symmetric $4\times 4$ matrices.
In addition, we assume that  the waves $u_j(x)$, $j=1,2,3,4$, solving the linear wave equation
in the Minkowski space,  are of the form
\ba
 u_j(x)= v^{(j)}\, \bigg(b^{(j)}_px^p \bigg)^{a}_+,\quad t^a_+=|t|^aH(t),
\ea
where $b^{(j)}_pdx^p$, $p=1,2,3,4$ are  {four linearly independent} light-like co-vectors of $\R^4$, \HOX{April 10, 2013: The text below is new and needs to be checked.}
$a>0$ and $v^{(j)}$ 
are constant $4\times 4$ matrices. 
{We also assume that $b^{(5)}$ is not in the linear span
of any three vectors $b^{(j)}$, $j=1,2,3,4$.} 
In the following,
we denote 
$b^{(j)}\,\cdotp x:=b^{(j)}_p x^p.$
Let us next consider the wave produced by interaction of two plane wave type solutions
in the  Minkowski space.\hiddenfootnote{We start with the observation
 that
\ba
& &(\p_0^2-\p_1^2-\p_2^2-\p_3^2)\bigg(  (x_0-x_1)^a_\pm\,\cdotp (x_0-x_2)^c_\pm   \bigg)\\
& &=2ac\,
(x_0-x_1)^{a-1}_\pm\,\cdotp (x_0-x_2)^{c-1}_\pm ,  \\
& &(\p_0^2-\p_1^2-\p_2^2-\p_3^2)\bigg(  (x_0-x_1)^a\,\cdotp (x_0-x_2)^c_\pm   \bigg)\\
& &=2ac\,
(x_0-x_1)^{a-1}_\pm\,\cdotp (x_0-x_2)^{c-1}.
\ea
More generally,  for $\tilde w=(0,\tilde w^1,\tilde w^2,\tilde w^3)\in \R^4$, $\tilde w\,\cdotp \tilde w=1$ and  $w=(0,w^1,w^2,w^3)\in \R^4$, $w\,\cdotp w=1$,
we have
\ba
& &(\p_0^2-\p_1^2-\p_2^2-\p_3^2)\bigg(  (x_0-x_1)^a_+\,\cdotp (x_0-\tilde w\,\cdotp x)^c_+   \bigg)=
\\&=&2ac\,
(x_0-x_1)^{a-1}_+\,\cdotp (x_0-\tilde w\,\cdotp x^{\prime})^{c-1}_+ + \\ 
& &-2ac\,
(x_0-x_1)^{a-1}_+\,\cdotp \tilde w_1(x_0-\tilde w\,\cdotp x^{\prime})^{c-1}_+
\ea
and
\ba
& &(\p_0^2-\p_1^2-\p_2^2-\p_3^2)\bigg(  (x_0-w\,\cdotp x^{\prime})^a_+\,\cdotp (x_0-\tilde w\,\cdotp x^{\prime})^c_+   \bigg)=
\\&=&2ac\,
(x_0-w\,\cdotp x^{\prime})^{a-1}_+\,\cdotp (x_0-\tilde w\,\cdotp x^{\prime})^{c-1}_+ + \\ 
& &+\sum_{j=1}^3 2ac\tilde w_jw_j
(x_0-w\,\cdotp x^{\prime})^{a-1}_+\,\cdotp \tilde w_1(x_0-\tilde w\,\cdotp x^{\prime})^{c-1}_+
\\&=&2ac(1-\tilde w\,\cdotp w)\,
(x_0-w\,\cdotp x^{\prime})^{a-1}_+\,\cdotp (x_0-\tilde w\,\cdotp x^{\prime})^{c-1}_+.
\ea}

Next, we will consider an operator ${\bf Q}_0$, which is an algebraic inverse of $\square_{\hat g}$ in Minkowski space, that 
we will define for certain products of Heaviside functions and polynomials.
Let $b^{(1)}=(1,p^1,p^2,p^3)$ and $b^{(2)}=(1,q^1,q^2,q^3)$ be light like vectors.
We use the notations
\ba
& &v^{a_1,a_2}(x;{b^{(1)},b^{(2)}})=(b^{(1)}\,\cdotp x)^{a_1}_+\,\cdotp (b^{(2)}\,\cdotp x)^{a_2}_+
\ea
and define 
\ba
{\bf Q}_0(v^{a_1,a_2}(x;{b^{(1)},b^{(2)}}))=\frac 1{2(a_1+1)(a_2+1)\,\hat g( b^{(1)},b^{(2)})} v^{a_1+1,a_2+1}(x;{b^{(1)},b^{(2)}}).
\ea
Then $\square_{\hat g}({\bf Q}_0(v^{a_1,a_2}(x;{b^{(1)},b^{(2)}})))=v^{a_1,a_2}(x;{b^{(1)},b^{(2)}})$.
{
Note that the function $v^{a_1,a_2}(x;{b^{(1)},b^{(2)}})$ is a product of two
plane waves and it is supported in the causal future of the space-like set $K_{12}$.
This implies that 
\ba
{\bf Q}_0(v^{a_1,a_2}(x;{b^{(1)},b^{(2)}}))={\bf Q}(v^{a_1,a_2}(x;{b^{(1)},b^{(2)}}))
\ea
where ${\bf Q}$ is the causal inverse of the wave operator $\square$
in the Minkowski space.

%
%
%

Similarly to the above
we denote
\ba
& &v^{a,0}_{\tau} (x;b^{(4)},b^{(5)})= u_4(x)\, u^\tau_0(x),
\quad 
u_4(x)= (b^{(4)}\,\cdotp x)^a_+,\quad  u^\tau_0(x)=e^{i\tau\,b^{(5)}\cdotp x} .
\ea
Then
$\square_{\hat g} ( v^{a,0}_{\tau} (x;b^{(4)},b^{(5)}))=2a\,\hat g(b^{(4)},b^{(5)})i\tau \, v^{a-1,0}_{\tau} (x;b^{(4)},b^{(5)})
$
and hence  we 
define 
\beq\label{Q0 def}\\
\nonumber
{\bf Q}_0( v^{a,0}_{\tau} (x;b^{(4)},b^{(5)}))=\frac 1{2i(a+1)\,\hat g(b^{(4)},b^{(5)})\tau}v^{a+1,0}_{\tau} (x;b^{(4)},b^{(5)}).
\eeq
Later, we will consider the relation between ${\bf Q}_0$ and the causal inverse ${\bf Q}$.
}
%
%


 
Next we prove that the indicator function
  $\mathcal G(v,\bsequence)$ in (\ref{definition of G}) does not vanish identically
  by showing that it coincides with  the {\it formal}
  indicator function $\mathcal G^{({\bf m})}(v,\bsequence)$,
 which is a real-analytic function that does not vanish identically
%
%
%

{
Below, let $x=(x^0,x^1,x^3,x^4)$ be the standard coordintes
in the Minkowski space and let $z=(z^j)_{j=1}^4$ be
light-like coordinates  $z^j=b^{(z)}\,\cdotp x$. We denote $x=X(z)$ and
$z=Z(x)$. Also, let
$P$ be a vector such that   $b^{(5)}\,\cdotp x=P\,\cdotp z$.

Let $h\in C^\infty_0(\R^4)$ be a function that has value 1 in a neigbhorhood
of $x=0$,  $T_0>1$ be such that $\supp(h)$ is contained
in the set  $\{x;\ x^0<T_0\}$.

Let $\chi=\chi(x^0)\in C^\infty(\R)$ be zero for $x^0>T_0+1$ and one for  $x^0<T_0$.}

We define  the (Minkowski) indicator function (c.f.\ (\ref{definition of G}))
\ba
\mathcal G^{({\bf m})}(v,{\bf b})=
\lim_{\tau\to\infty} \tau^{m}(\sum_{\b\leq n_1}
\sum_{\sigma\in \Sigma(4)} 
T^{({\bf m}),\b}_{\tau,\sigma}+\tilde T^{({\bf m}),\b}_{\tau,\sigma}),
\ea 
where the super-index $({\bf m})$ refers
to the word "Minkowski".
Above, 
 $\sigma:\{1,2,3,4\}\to \{1,2,3,4\}$ runs over all permutations
of indexes of the waves $u_j$,
where ${\bf b}=(b^{(1)},b^{(2)},\dots,b^{(5)})$ and
$T^{({\bf m}),\b}_{\tau,\sigma}$ and $\tilde T^{({\bf m}),\b}_{\tau,\sigma}$ are counterparts
of the functions $T^{(4),\b}_{\tau}$ and $\tilde T^{(4),\b}_{\tau}$,
see (\ref{T-type source})-(\ref{tilde T-type source}),
obtained by replacing the pieces of the spherical waves and
the gaussian beam by plane waves by replacing 
the parametrix ${\bf Q}$ with a formal parametrix ${\bf Q}_0$
and including in the obtained formula
 a smooth cut off function   $h\in C^\infty_0(M)$ which is one near 
 the intersection point of  the waves, and permutating indexes,
that is,
\beq
\label{Term type 1}
& &\hspace{- .5cm} T^{({\bf m}),\b}_{\tau,\sigma}=\bra S^0_2(u^\tau \,\cdotp \B_{4}u_{\sigma(4)}), h\,\cdotp \B_{3}u_{\sigma(3)}\,\cdotp 
S_1^0(\B_{2}u_{\sigma(2)}\,\cdotp \B_{1}u_{\sigma(1)})\cet_{L^2(\R^4)},\hspace{-1cm} \\
\label{Term type 2}
& &\hspace{- .5cm} \tilde T^{({\bf m}),\b}_{\tau,\sigma}=\bra u^\tau,h\,\cdotp S_2^0(\B_4u_{\sigma(4)} \,\cdotp \B_{3}u_{\sigma(3)})\,\cdotp 
S^0_1(\B_2u_{\sigma(2)}\,\cdotp \B_1u_{\sigma(1)})\cet_{L^2(\R^4)},\hspace{-1cm} 
\eeq
where $u_j =v_{(j)} (b^{(j)}\,\cdotp x)^a_+$, $j=1,2,3,4$ and
\ba
u^\tau(x) =\chi(x^0)v_{(5)} \exp(i\tau b^{(5)}\,\cdotp x).
\ea
 Moreover, $\B_j=\B_{j,\beta}$ and finally,
 $S_j^0=S^0_{j,\beta}\in \{{\bf Q}_0,I\}$.
We note that here that the algebraic inverse ${\bf Q}_0$ is used to 
replace both the causal parametrix ${\bf Q}$ and the anti-causal parametrix ${\bf Q}^*$, and
the commutator terms do not appear at all.

\HOX{ T17: Task was to show that $\mathcal G^{({\bf m})}({\bf w},\bsequence)=\mathcal G({\bf v},\bsequence)$. This is now done, but arguments need to be checked.}

Let us now consider the orders of the differential  operators appearing above.
The orders $k_j=ord(\B^\beta_j)$ 
of the differential operators $\B^\beta_j$, 
defined in (\ref{extra notation}), 
 depend on 
$\vec S_\beta^0=(S_{1,\beta}^0,S_{2,\beta}^{0})$ 
as follows: 
When $\beta$ is such that $\vec S_\beta^0=(  {\bf Q}_0,  {\bf Q}_0)$,  for the terms $ T^{({\bf m}),\b}_{\tau,\sigma}$
 we have 
\beq\label{k: s for (Q_0,Q_0) and T}
& &k_1+k_2+k_3+k_4\leq 6,\quad  k_3+k_4\leq 4,\quad
k_4\leq 2
\eeq
and for the terms $\tilde T^{({\bf m}),\b}_{\tau,\sigma}$ we have 
\beq\label{k: s for (Q_0,Q_0) and tilde T} 
k_1+k_2+k_3+k_4 \leq 6,\quad
k_1+k_2\leq 4,\quad
k_3+k_4\leq 4.\hspace{-1.5cm}
\eeq

When  $\beta$ is such that $\vec S_\beta^0=(I,  Q_0)$
 we have  for terms $T^{({\bf m}),\b}_{\tau,\sigma}$
\beq\label{k: s for (I,Q_0) and T} 
k_1+k_2+k_3+k_4 \leq 4,\quad k_4\leq 2,
\eeq
 and for terms $\tilde T^{({\bf m}),\b}_{\tau,\sigma}$ we have 
\beq\label{k: s for (I,Q_0) and tilde T}  k_1+k_2+k_3+k_4 \leq 4,\quad
k_1+k_2\leq 2.
\eeq

When  $\beta$ is such that $\vec S_\beta^0=(  {\bf Q}_0,I)$, 
both for the terms $T^{({\bf m}),\b}_{\tau,\sigma}$ 
and $\tilde  T^{({\bf m}),\b}_{\tau,\sigma}$
we have 
\beq\label{k: s for (Q_0,I) and tilde T} 
& &k_1+k_2+k_3+k_4\leq 4,\quad
k_3+k_4\leq 2.
\eeq
Finally, when  $\beta$ is such that $\vec S_\beta^0=(I,I)$,
 for the terms   $T^{({\bf m}),\b}_{\tau,\sigma}$ and $\tilde T^{({\bf m}),\b}_{\tau,\sigma}$ we have
$k_1+k_2+k_3+k_4 \leq 2$.

%

{
\begin{lemma}  When
 $b^{(j)}$, $j=1,2,3,4$ are   linearly independent light-like co-vectors and light-like
 co-vector $b^{(5)}$ is not in the linear span
of any three vectors $b^{(j)}$, $j=1,2,3,4$ we have
$\mathcal G({\bf w},{\bf b})=\mathcal G^{({\bf m})}({\bf v},{\bf b})$
when  $w_{(j)}=(v_{(j)},0)\in \R^{10}\times \R^L$.
\end{lemma}

{\bf Proof.}
Let us start by considering the relation of ${\bf Q}_0$ with the causal inverse ${\bf Q}$
in (\ref{Q0 def}).
Let
\ba
\nonumber
w_{\tau,0}&=&{\bf Q}_0(v^{a,0}_{\tau} (\,\cdotp;b^{(4)},b^{(5)}))\\
&=&\int_{\R}e^{i\theta_4z^4+ i\tau P\cdotp z}
a_4(z,\theta_4)d\theta_4,\\
w_\tau&=&{\bf Q}^*(J ),\\
J&=&
 u_4 \,\cdotp (\chi\,\cdotp u^\tau_0),\ea
where 
\ba
J(z)&=&\chi(X^0(z))\, \square w_{\tau,0}\\
&=&\int_{\R}e^{i\theta_4z^4+ i\tau P\cdotp z}
(\tau b_1(z,\theta_4)+
b_2(z,\theta_4))a_5(z,\tau)d\theta_4,
\ea
where $a_5(z,\tau)=1$.
Then 
\ba
& &\square(w_{\tau}-\chi w_{\tau,0})=J_1,\quad \hbox{for } x\in \R^4,\\
& &J_1=[\square,\chi]w_{\tau,0}\\
&&\quad\ =\int_{\R}e^{i\theta_4z^4+ i\tau P\cdotp z}
(\tau b_3(z,\theta_4)+
b_4(z,\theta_4))d\theta_4,
\ea
where $b_3(z,\theta_4)$ and $b_4(z,\theta_4)$ are
supported in the domain $ T_0< X^0(z)<T_0+1$ and
$w_{\tau}-\chi w_{\tau,0}$ is supported in domain $X^0(z)<T_0+1$ .
Thus 
\ba
w_{\tau}=\chi w_{\tau,0}+{\bf Q}^*J_1.
\ea
Here, 
we can write 
\beq\label{modified J}
& &J_1(z)= u_4^{(1)}(z) u^{\tau,(1)}(z)+
u_4^{(2)}(z) u^{\tau,(2)}(z)
,\quad\hbox{where}\\ \nonumber
& &u_4^{(1)}(z)= \int_{\R}e^{i\theta_4z^4}
b_3(z,\theta_4)d\theta_4,\\ \nonumber
& & u^{\tau,(1)}(z)=\tau u^\tau(z),\\ \nonumber
& &u_4^{(2)}(z)= \int_{\R}e^{i\theta_4z^4}
b_4(z,\theta_4)d\theta_4,\\ \nonumber
& & u^{\tau,(2)}(z)=u^\tau(z).
\eeq

Let us now substitute this in to the above  microlocal computations
done in the proof of Prop.\ \ref{lem:analytic limits A}. 
%
%

Recall that $b^{(j)}$, $j=1,2,3,4$, are four linearly independent co-vectors.
This means that a condition analogous to (LI) in the proof of Prop.\ \ref{lem:analytic limits A}
is satisfied, and that $b^{(5)}$ is not in the space spanned by any of three  of the
co-vectors $b^{(j)}$, $j=1,2,3,4$.
Also, observe that hyperplanes $K_j=\{x\in \R^4;\ b^{(j)}\,\cdotp x=0\}$
intersect at origin of $\R^4$. Thus, we see that
the arguments in  the proof of Prop.\ \ref{lem:analytic limits A}
are valid mutatis mutandis if the phase function of the gaussian beam 
$\varphi(x)$ is replaced
by the phase function of the plane wave, $b^{(5)}\,\cdotp x$, 
and  the geodesic $\gamma_{x_5,\xi_5}$, on
which the gaussian beam propagates, is replaced by the whole space $\R^4$.
In particular, as $b^{(5)}$ is not in the space spanned by any of three  of those
co-vectors, the case (A1) in the proof of Prop.\ \ref{lem:analytic limits A} cannot occur.
In particular, we see that
 the  leading order asymptotics of the terms
$T^{\b}_{\tau,\sigma}$ and $\tilde T^{\b}_{\tau,\sigma}$ does not change
as these asymptotics are  obtained using
the method of stationary phase for the integral 
(\ref{eq; T tau asympt modified}) and the other analogous integrals
at the critical point $z=0$.
\HOX{The explanation needs to be clarified.}
In other words, we can replace the gaussian beam by a plane wave in our considerations
similar to those in the proof of Prop.\ \ref{lem:analytic limits A}.

Using (\ref{modified J})  and the fact that $b_3(z,\theta_4)$ and $b_4(z,\theta_4)$ 
vanish near $z=0$, we see that if $u_4$ and $u^{\tau}$ are replaced by
 $u_4^{(j)}$ and $u^{\tau,(j)}$, respectively, where $j\in \{1,2\}$ and
we can do similar computations based on the method of stationary phase  as are done in the 
proof of Proposition \ref{lem:analytic limits A}. Then both terms
$T^{\b}_{\tau,\sigma}$ and $\tilde T^{\b}_{\tau,\sigma}$
have asymptotics  $O(\tau^{-N})$ for all $N>0$ as $\tau\to \infty$.
In other words, in the proof of Prop.\ \ref{lem:analytic limits A} the term
$w_{\tau}={\bf Q}^*( u_4u^\tau )$ can be replaced by $\chi w_{\tau,0}$
without changing the leading order asymptotics. 
This shows
that $\mathcal G({\bf w},{\bf b})=\mathcal G^{({\bf m})}({\bf v},{\bf b})$,
where $w_{(j)}=(v_{(j)},0)\in \R^{10}\times \R^L$.  
\hfill \Box \medskip
}

%
%


\begin{proposition}\label{singularities in Minkowski space}
 Let $X$ be the set of $({\bf b},v^{(2)},v^{(3)},v^{(4)})$,
where ${\bf b}$ is a  5-tuple of
light-like covectors ${\bf b}=(b^{(1)},
b^{(2)},b^{(3)},b^{(4)},b^{(5)})$ and 
$v^{(j)}\in \R^{10}$, $j=2,3,4$ are the polarizations that satisfy the equation (\ref{divergence condition for symbol})
with respect to $b^{(j)}$,
i.e.,  the divergence condition for the principal symbols.
Also,  let $a\in \Z_+$ be  large enough.
Then for $({\bf b},v^{(2)},v^{(3)},v^{(4)})$ in a  generic (i.e. open and dense)
subset of $X$   
 there exist linearly independent vectors  $v^{(5)}_q$, $q=1,2,3,4,5,6,$
so that if
 $v^{(5)}\in \hbox{span}(\{v^{(5)}_p;\ p=1,2,3,4,5,6\})$  is non-zero, then 
 there exists 
 a vector $v^{(1)}$ for which  
the pair 
$(b^{(1)},v^{(1)})$
satisfies the equation (\ref{divergence condition for symbol})
 and
 $\mathcal G^{({\bf m})}({\bf v},{\bf b})\not =0$
 with ${\bf v}=(v^{(1)},v^{(2)},v^{(3)},v^{(4)},v^{(5)})$.

%
%
%
\end{proposition}
%


\noindent{\bf Proof.} 
To show that the coefficient $\mathcal G^{({\bf m})}({\bf v},{\bf b})$  of the 
leading order term in the asymptotics  
is non-zero, we consider a special case when the direction vectors of the
intersecting plane waves  in the  Minkowski space are 
the linearly independent light-like vectors of the form
\ba
b^{(5)}=(1,1,0,0),\quad b^{(j)}=(1,1-\frac 12\rhoepsilon_j^2,\rhoepsilon_j+O(\rhoepsilon_j^3),\rhoepsilon_j^3),\quad j=1,2,3,4,
\ea
where $\rhoepsilon_j>0$ are small parameters for which
\beq\label{b distances}
& &\|b^{(5)}-b^{(j)}\|_{(\R^4,\hat g^+)}= \rhoepsilon_{j}(1+o(\rhoepsilon_{j})),\quad j=1,2,3,4.
\eeq
With an appropriate choice of $O(\rhoepsilon_k^3)$ above, the vectors $b^{(k)}$, $k\leq 5$
are  light-like and
\ba
\,\hat g(b^{(5)},b^{(j)})&=&-1+(1-\frac 12 \rhoepsilon_j^2)=-\frac 12 \rhoepsilon_j^2,\\
\,\hat g( b^{(k)},b^{(j)})&=&
-\frac 12 \rhoepsilon_k^2-\frac 12 \rhoepsilon_j^2+O(\rhoepsilon_k\rhoepsilon_j).
\ea
Below, we denote $\omega_{kj}=\,\hat g( b^{(k)},b^{(j)})$. 
We consider the case when the orders of $\rhoepsilon_j$
are   
(Note here the "unordered" numbering 4-2-3-1)
\beq\label{eq: ordering of epsilons}
\rhoepsilon_4=\rhoepsilon_2^{100},\ \rhoepsilon_2=\rhoepsilon_3^{100},\hbox{ and }\rhoepsilon_3=\rhoepsilon_1^{100}.  
\eeq
{
Note that when $\rho_1$ is small enough, $b^{(5)}$ is not a linear combination
of any three vectors $b^{(j)}$, $j=1,2,3,4$.}

The coefficient $\mathcal G^{({\bf m})}$ of the leading order asymptotics
is computed by analyzing the leading order terms of 
all 4th order interaction terms, similar to those 
given in (\ref{Term type 1}) and (\ref{Term type 2}).
We will start by analysing the most important terms $ T^{({\bf m}),\beta}_\tau$
of the type (\ref{Term type 1})
%
%
when $\beta$ is such that $\vec S_\b=({\bf Q}_0,{\bf Q}_0)$.
%
When $k_j=k_j^\b$ is the order of $\B_j$, 
and we denote $\vec k_\b=(k_1^\b,k_1^\b,k_3^\b,k_4^\b)$,
we see that 
\beq\label{first asymptotical computation}
 & &T^{({\bf m}),\beta}_\tau
=\bra {\bf Q}_0(\B_4 u_4\,\cdotp u^\tau), h\,\cdotp \B_3u_3\,\cdotp {\bf Q}_0(\B_2u_2\,\cdotp \B_1 u_1)\cet\\
\nonumber
& &\hspace{-1cm}=C\frac {\P_\beta} {\omega_{45}
\tau} \frac 1{\omega_{12}}
\bra v^{a-k_4+1,0}_{\tau} (\,\cdotp ;b^{(4)},b^{(5)}), h\,\cdotp u_3\,\cdotp v^{a-k_2+1,a-k_1+1} (\,\cdotp ;b^{(2)},b^{(1)})\cet
\\
\nonumber&&\hspace{-1cm} =C\frac {\P_\beta}{\omega_{45}\tau} \frac 1{\omega_{12}} 
\int_{\R^4}
(b^{(4)}\,\cdotp x)^{a-k_4+1}_+e^{i\tau(b^{(5)}\,\cdotp x)} 
h(x)(b^{(3)}\,\cdotp x)^{a-k_3}_+\cdotp\\
\nonumber& &\quad\quad\quad\cdotp
(b^{(2)}\,\cdotp x)^{a-k_2+1}_+
(b^{(1)}\,\cdotp x)^{a-k_1+1}_+\,dx,
\eeq
where $\P=\P_\beta$ is a polarization factor involving the coefficients of $\B_j$, 
the directions $b^{(j)}$, and the polarization $v_{(j)}$. Moreover, $C=C_a$ is a generic constant
 depending on $a$ and $\beta$ but not on  $b^{(j)}$ or $v_{(j)}$.

We will analyze the polarization factors later, but as a sidetrack, 
let us already now explain the nature of the polarization
term 
when $\beta=\beta_1$,
see (\ref{beta0 index}).  Observe that this term appear only when  we analyze
 the term
 $\bra F_\tau,{\bf Q}(A[u_4,{\bf Q}(A[u_3,{\bf Q}(A[u_2,u_1])])])\cet$ where
 all operators $A[v,w]$ are of the type $A_2[v,w]=\hat g^{np}\hat g^{mq}v_{nm}\p_p\p_q w_{jk}$,
 cf.\ (\ref{eq: tilde M1}) and (\ref{eq: tilde M2}).
 Due to this, we have the polarization factor 
 \beq\label{pre-eq: def of D}
\P_{\beta_1} =(v_{(4)}^{rs}b^{(1)}_rb^{(1)}_s)(v_{(3)}^{pq}b^{(1)}_pb^{(1)}_q)(v_{(2)}^{nm}b^{(1)}_nb^{(1)}_m)\D,
\eeq
 where $v_{(\ell)}^{nm}=\hat g^{nj}\hat g^{mk}v^{(\ell)}_{jk}$ and
\beq\label{eq: def of D}
 \D =\hat g_{nj}\hat g_{mk}v_{(5)}^{nm}v_{(1)}^{jk}.
 \eeq
We will postpone the analysis of  the polarization factors
$\P_{\beta}$ in  $ T^{({\bf m}),\beta}_\tau$  with $\beta\not=\beta_1$ later.

Let us now return back to the computation (\ref{first asymptotical computation})\hiddenfootnote{
In the computations below, we analyze the following terms:

  1. The $T^{{\bf Q}, {\bf Q}}$ term in formula (\ref{first asymptotical computation}) comes from the 2nd term in (\ref{w4 solutions})
 
  2. $\tilde T^{{\bf Q}, {\bf Q}}$ term on p. 55 that comes from the 1st term in (\ref{w4 solutions}) inside
  the large brackets,
  
  3. $T^{I, {\bf Q}}$ term on p. 56  that comes from the 3rd term in (\ref{w4 solutions}) inside
  the large brackets,

  4. $\tilde T^{I, {\bf Q}}$ term on p. 56  that comes from the 4th term in (\ref{w4 solutions}) inside
  the large brackets,

  5. $T^{{\bf Q}, I}$ term on p. 56  that comes also from the 4th term in (\ref{w4 solutions}) inside
  the large brackets
  if we use
another permutation of $\{1, 2, 3, 4\}$,

 6. $T^{{\bf Q}, I}$ term on p. 56  that comes from the 5th term in (\ref{w4 solutions}) inside
  the large brackets.
}.
We next use in $\R^4$ the coordinates $y=(y^1,y^2,y^3,y^4)^t$
where $y^j=b^{(j)}_kx^k$, i.e., and let $A\in \R^{4\times 4}$
be the matrix for which $y=A^{-1}x$.
Let ${\bf p}=(A^{-1})^tb^{(5)}$.
 {In the $y$-coordinates,  $b^{(j)}=dy^j$ for $j\leq 4$ and
 $b^{(5)}=\sum_{j=1}^4{\bf p}_jdy^j$ and
 \ba
 {\bf p}_j
 =\hat g(b^{(5)},dy^j)=\hat g(b^{(5)},b^{(j)})= \omega_{j5}=-\frac 12 \rhoepsilon_j^2.
 \ea
Then   $b^{(5)}\,\cdotp x={\bf p}\,\cdotp y$.
 We use the  notation
$
{\bf p}_j=\omega_{j5}=-\frac 12 \rhoepsilon_j^2,
$ 
that is, we denote the same object with several symbols,
to clarify the steps we do in the computations.

%
%
%
 
 }

 Then \HOX{ In the final version, we could 
substitute the form of $\det (A)$ 
 in all terms. Now it is not done for clarity.}
 $\det (A)=8\rhoepsilon_1^{-3}\rhoepsilon_2^{-2}\rhoepsilon_3^{-1}(1+O(\rhoepsilon_1))$ and
\ba
 T^{({\bf m}),\beta}_\tau
=\frac {C\P_\beta}{ \omega_{45}\tau} \frac {\det (A)}{ \omega_{12}} \int_{(\R_+)^4}
e^{i\tau {\bf p}\cdotp y} 
h(Ay)y_4^{a-k_4+1}y_3^{a-k_3}y_2^{a-k_2+1}
y_1^{a-k_1+1}\,dy.
\ea
Using repeated integration by parts 
we see that
%
%
%
\beq\label{t 4 formula}\\ \nonumber
 T^{({\bf m}),\beta}_\tau&=& {C\det (A)\,{\P_\beta}}
\frac {(i\tau)^{-(12+4a-|\vec k_\b|)}(1+O(\tau^{-1}))}{ \rhoepsilon_4^{2(a-k_4+1+2)}\rhoepsilon_3^{2(a-k_3+1)}\rhoepsilon_2^{2(a-k_2+2)}\rhoepsilon_1^{2(a-k_1+1+2)}}\,.
\eeq
{\newtext Note that here and below $O(\tau^{-1})$ may depend also on $\rhoepsilon_j$, that is,
we have  $|O(\tau^{-1})|\leq C(\rhoepsilon_1,\rhoepsilon_2,\rhoepsilon_3,\rhoepsilon_4)\tau^{-1}.$}

To show that $ \mathcal G^{({\bf m})}({\bf v},{\bf b})$ is  non-vanishing 
we need to estimate $\P_{\beta_1}$ from below. In doing this we encounter 
the difficulty that $\P_{\beta_1}$ can go to zero, and moreover,
simple computations show that as  the pairs
$(b^{(j)},v^{(j)})$
satisfies the divergence condtion (\ref{divergence condition for symbol})
we have $ v_{(r)}^{ns} b^{(j)}_nb^{(j)}_s=O(\rhoepsilon_r+\rhoepsilon_j).$
%
%
%
%
However, to show that $ \mathcal G^{({\bf m})}({\bf v},{\bf b})$ is  non-vanishing 
 for some  ${\bf v}$ we consider a 
  particular choice of polarizations $v^{(r)}$,
 namely  
\beq\label{chosen polarization}
 v^{(r)}_{mk}= b^{(r)}_mb^{(r)}_k,\quad \hbox{for $r=2,3,4,$ but not for $r=1,5$}
\eeq
so that for $r=2,3,4$, we have
\ba
\hat g^{nm}b^{(r)}_nv^{(r)}_{mk}=0,\quad
\hat g^{mk}v^{(r)}_{mk}=0,\quad
\hat g^{nm}b^{(r)}_nv^{(r)}_{mk}-\frac 12 (\hat g^{mk}v^{(r)}_{mk})b^{(r)}_k=0.
\ea
Note that for this choice of $ v^{(r)}$ the linearized divergence conditions hold.
Moreover, for this choice of $ v^{(r)}$ we see that for $\rhoepsilon_j\leq \rhoepsilon_r^{100}$
\beq\label{vbb formula}
v_{(r)}^{ns} b^{(j)}_nb^{(j)}_s
 = \hat g(b^{(r)},b^{(j)})\, \hat g(b^{(r)},b^{(j)})
 = \rhoepsilon_r^4+O (\rhoepsilon_r^5). 
\eeq


In particular, when $\beta=\beta_1$, so that
 $k_{\b_1}=(6,0,0,0)$ and  
 the  polarizations are given by (\ref{chosen polarization}), we have
\ba
\P_{\beta_1}= (\D +O(\rhoepsilon_1))\rhoepsilon_1^4\cdotp\rhoepsilon_1^4\cdotp\rhoepsilon_1^4,
\ea
where $\D$ is the inner product of $v^{(1)}$ and $v^{(5)}$ given in (\ref{eq: def of D}). 
Then\hiddenfootnote{Indeed, when $\beta=\beta_1$ we have
\ba
T^{\beta_1,0}_\tau
\ba
T^{(1)}_\tau
&=&\frac {c_1^{\prime}}{4 } \cdotp {(i\tau)^{-(12 +4a-|\vec k_\b|)}(1+O(\tau^{-1}))}\vec \e^{-2\vec a}
 (\omega_{45}\omega_{12})^{-1} \rhoepsilon_4^{2(k_4-1+1)}\rhoepsilon_2^{2(k_2-1+1)}\rhoepsilon_3^{2(k_3+1)}\rhoepsilon_1^{2(k_1-1+1)}{\P_\beta}\\
 &=&{c_1^{\prime}}
 {(i\tau)^{-(12 +4a-|\vec k_\b|)}(1+O(\tau^{-1}))}\vec \e^{-\vec a-\vec 1}
 (\omega_{45}\omega_{12})^{-1}\rhoepsilon_4^{-2}\rhoepsilon_2^{-2}\rhoepsilon_3^{0}\rhoepsilon_1^{10}{\P_\beta}
 \\
 &=&{c_1^{\prime}}
 {(i\tau)^{-(12 +4a-|\vec k_\b|)}(1+O(\tau^{-1}))}\vec \e^{-2(\vec a+\vec 1)}
\rhoepsilon_4^{-4}\rhoepsilon_2^{-2}\rhoepsilon_3^{0}\rhoepsilon_1^{20}\D.
\ea
} 
the term $T^{({\bf m}),\beta_1}_\tau$, which later turns out to
have the strongest asymptotics in our considerations, has the asymptotics
\beq\label{leading term}
& &\hspace{-1cm}T^{({\bf m}),\beta_1}_\tau=\L_\tau,\quad\hbox{where}\\
\nonumber
& &\hspace{-1cm}\L_\tau=
C\det (A)
{(i\tau)^{-(6 +4a)}(1+O(\tau^{-1}))}\vec\rhoepsilon^{-2(\vec a+\vec 1)}
\rhoepsilon_4^{-4}\rhoepsilon_2^{-2}\rhoepsilon_3^{0}\rhoepsilon_1^{20}\D,\hspace{-1.5cm}
\eeq 
where $\vec\rhoepsilon=(\rho_1,\rho_2,\rho_3,\rho_4)$,
$\vec a=(a,a,a,a)$, and $\vec 1=(1,1,1,1)$.
To compare different terms,
we express
$\rhoepsilon_j$ in powers of $\rhoepsilon_1$ as explained in formula (\ref{eq: ordering of epsilons}), 
that is, we write $\rhoepsilon_4^{n_4}\rhoepsilon_2^{n_2}\rhoepsilon_3^{n_3}\rhoepsilon_1^{n_1}=\rhoepsilon_1^m$
with $m=100^3n_4+100^2n_2+100n_3+n_1$.
In particular, we will below write
\ba
& &\L_\tau=C_{\beta_1}(\vec\rhoepsilon)\,\tau^{n_0}
(1+O(\tau^{-1}))\quad\hbox{as $\tau\to \infty$ for each fixed $\vec \e$, and}\\ 
& &C_{\beta_1}(\vec\e)=c^{\prime}_{\beta_1}\,\rhoepsilon_1^{m_0}(1+o(\rhoepsilon_1))
\quad\hbox{as  } \rhoepsilon_1\to 0. 
\ea
 Below we will show that $c^{\prime}_{\beta_1}$
does not vanish for generic $(\vec x,\vec \xi)$ and $(x_5,\xi_5)$  and polarizations ${\bf v}$.
We will  consider below $ \beta\not= \beta_1$
and show that also  these terms have the asymptotics 
\ba
& &T^{({\bf m}),\beta}_\tau=C_{\beta}(\vec\rhoepsilon)\,\tau^{n}
(1+O(\tau^{-1}))\quad\hbox{as $\tau\to \infty$ for each fixed $\vec \e$, and}\\ 
& &C_{\beta}(\vec\e)=c^{\prime}_{\beta}\,\rhoepsilon_1^{m}(1+o(\rhoepsilon_1))
\quad\hbox{as  } \rhoepsilon_1\to 0. 
\ea
When\hiddenfootnote{
Alternatively, the ordering could be formulated as follows: 
We say that
\beq \label{Slava 1}
C' \tau^{-c'_\tau} \rhoepsilon_4^{c'_4} \rhoepsilon_3^{ c'_3} \rhoepsilon_2^{ c'_2} \rhoepsilon_1^{ c'_1}
\prec C \tau^{-c_\tau} \rhoepsilon_4^{c_4} \rhoepsilon_3^{ c_3} \rhoepsilon_2^{ c_2} \rhoepsilon_1^{ c_1} 
\eeq
if $C \neq 0$ and some of the following conditions is valid:
\begin{itemize}
\item [(i)] if $c_\tau < c'_\tau$;
\item [(ii)] if $c_\tau = c'_\tau$ but $c_4 <c_4'$;
\item [(iii)] if $c_\tau = c'_\tau$ and $c_4 = c_4'$ but $c_2 <c_2'$;
\item [(iv)] if $c_\tau = c'_\tau$ and $c_4 = c_4',\,c_2 =c_2' $ but $c_3 <c_3'$;
\item [(v)] if $c_\tau = c'_\tau$ and $c_4 = c_4',\,c_3 =c_3',\, c_2 =c_2' $ but $c_1 <c_1'$.
\end{itemize}
Note that this is an ordering not just a partial ordering.} we have that 
either $n\leq n_0$ and $m<m_0$,
or $n<n_0$, we say that 
$T^{({\bf m}),\beta}_\tau$ has weaker asymptotics than
$T^{({\bf m}),\beta_1}_\tau$ and denote
$T^{({\bf m}),\beta}_\tau \prec \L_\tau.$

{

As we consider here the asymptotic of five
small parameters $\tau^{-1}$ and $\e_i,$ $i=1,2,3,4,$ and compare in which order
we make them tend to $0$, let us explain the above ordering
in detail. Above, we have chosen the order: first $\tau^{-1}$, then $\e_4$, $\e_2$, $\e_3$ and finally,
 $\e_1$. 
 In correspondence with this choice we can introduce an ordering on all monomials $c \tau^{-c_\tau} \e_4^{n_4} \e_3^{n_3} \e_2^{n_2} \e_1^{n_1}$. Namely, we
say that
\beq \label{SL1}
C^\prime  \tau^{-n^\prime _\tau} \e_4^{n^\prime _4} \e_3^{ n^\prime _3} \e_2^{ n^\prime _2} \e_1^{ n^\prime _1}\prec C \tau^{-n_\tau} \e_4^{n_4} \e_3^{ n_3} \e_2^{ n_2} \e_1^{ n_1} 
\eeq
if $C \neq 0$ and one of the following holds
\begin{itemize}
\item [(i)] if $n_\tau < n^\prime _\tau$;
\item [(ii)] if $n_\tau = n^\prime _\tau$ but $n_4 <n_4^\prime $;
\item [(iii)] if $n_\tau = n^\prime _\tau$ and $n_4 = n_4^\prime $ but $n_2 <n_2^\prime $;
\item [(iv)] if $n_\tau = n^\prime _\tau$ and $n_4 = n_4^\prime ,\,n_2 =n_2^\prime  $ but $n_3 <n_3^\prime $;
\item [(v)] if $n_\tau = n^\prime _\tau$ and $n_4 = n_4^\prime ,\,n_3 =n_3^\prime ,\, n_2 =n_2^\prime  $ but $n_1 <n_1^\prime $.
\end{itemize}
Note that this ordering is not  a partial ordering.

Then, we can analyze terms $T^{({\bf m}),\beta}_{\tau,\sigma}$ and $
\tilde T^{({\bf m}),\beta}_{\tau,\sigma}$ in the formula for 
\beq \label{SL2}
\Theta_\tau^{(4)}=
\Theta_{\tau,\vec\e}^{(4)}=\sum_{\beta\in J_\ell}\sum_{\sigma\in \Sigma(\ell)}
\left(T^{({\bf m}),\beta}_{\tau,\sigma} + 
\tilde T^{({\bf m}),\beta}_{\tau,\sigma} \right).
\eeq
Note that here the terms, in which the permutation $\sigma$ is either the identical permutation $id$ or
the permutation  $\sigma_0=(2,1,3,4)$, are the same.

\medskip

{\bf Remark 3.3.} We can find the leading order asymptotics of the strongest terms in the decomposition
(\ref{SL2}) using the  following algorithm.
 First, let us multiply $\Theta_{\tau,\vec\e}^{(4)}$ by $\tau^{\hat n_\tau}$,
where $\hat n_\tau= \min_\beta n_\tau(\beta)$. Taking then $\tau \to \infty$ will give non-zero contribution from
only those terms $T^{({\bf m}),\beta}_{\tau,\sigma}$ and 
$\tilde T^{({\bf m}),\beta}_{\tau,\sigma} $ where $n_\tau(\beta)=\hat n_\tau$.
This corresponds to step $(i)$ above.
Multiplying next by $\e_4^{-\hat n_4}$, where $\hat n_4= \min_\beta n_4(\beta)$ under
the condition that $n_\tau(\beta)=\hat n_\tau$ and taking $\e_4 \to 0$ corresponds 
to selecting terms $T^{({\bf m}),\beta}_{\tau,\sigma} $ with $n_\tau(\beta)=\hat n_\tau$ and $n_4(\beta)=\hat n_4$
and terms $ \tilde T^{({\bf m}),\beta}_{\tau,\sigma} $ with  $ \, \tilde n_\tau(\beta)=
\hat n_\tau$ and $\tilde n_4(\beta)=\hat n_4$.
This corresponds to step $(ii)$.
Continuing this process we  obtain a scalar value that gives 
the leading order asymptotics of the strongest terms in the decomposition
(\ref{SL2}).\medskip

The next results  tells  what are the strongest terms in (\ref{SL2}).

\begin{proposition} \label{SL:order} In 
(\ref{SL2}), the strongest term are $T^{({\bf m}),\beta_1}_{\tau}=T^{({\bf m}),\beta_1}_{\tau,id}$
and $T^{({\bf m}),\beta_1}_{\tau,\sigma_0}$
in the sense that for all $(\beta,\sigma)\not \not\in \{ (\beta_1,id),(\beta_1,\sigma_0)\}$ we have  $T^{({\bf m}),\beta}_{\tau,\sigma}\prec 
T^{({\bf m}),\beta_1}_{\tau,id}$. 
\end{proposition}

{\bf Proof.}} When  $\vec S_\b=(S_1,S_2)=({\bf Q}_0,{\bf Q}_0)$,
similar computations to the above ones yield
\ba
\tilde  T^{({\bf m}),\beta}_\tau
&=&\bra  u^\tau,h\,\cdotp  {\bf Q}_0(\B_4u_4\,\cdotp \B_3u_3)\,\cdotp {\bf Q}_0(\B_2u_2\,\cdotp \B_1u_1)\cet\\
&=&C\det (A) {{\P_\beta}} \frac {(i\tau)^{-(12+4a-|\vec k_\b|)}(1+O(\tau^{-1}))}{ \rhoepsilon_4^{2(a-k_4+2)}\rhoepsilon_3^{2(a-k_3+1+2)}\rhoepsilon_2^{2(a-k_2+2)}\rhoepsilon_1^{2(a-k_1+1+2)}}\,.
\ea

Let us next consider the case when  $\vec S_\b=(S_1,S_2)=(I,{\bf Q}_0)$.
Again, the computations similar to the above ones show that
\ba
 T^{({\bf m}),\beta}_\tau
&=&\bra {\bf Q}_0(u^\tau\,\cdotp \B_4 u_4 ), h\,\cdotp \B_3u_3\,\cdotp I(\B_2u_2\,\cdotp \B_1u_1)\cet\\
&=& {i\,C{\P_\beta} \det (A)}\frac {(i\tau)^{-(10+4a-|\vec k_\b|)}(1+O(\tau^{-1}))}{
 \rhoepsilon_4^{2(a-k_4+1+2)}\rhoepsilon_3^{2(a-k_3+1)}\rhoepsilon_2^{2(a-k_2+1)}\rhoepsilon_1^{2(a-k_1+1)}}\,
\ea
and\hiddenfootnote{In fact, when $(S_1,S_2)=(I,{\bf Q}_0)$, the term $\tilde  T^{({\bf m}),\beta}_\tau$ is similar to the 
corresponding term in the above case $(S_1,S_2)=({\bf Q}_0,I)$
up to a permutation of indexes, and thus there is
 no need to analyze these terms separately.} 
\ba
& &\hspace{-1cm}\tilde  T^{({\bf m}),\beta}_\tau
=\bra  u^\tau, h\,\cdotp {\bf Q}_0(\B_4u_4\,\cdotp \B_3\,u_3)\,\cdotp I(\B_2\,u_2\,\cdotp \B_1\, u_1)\cet\\
&=&{{\P_\beta} C\det (A)} 
\frac {(i\tau)^{-(10+4a-|\vec k_\b|)}(1+O(\tau^{-1}))}
{\rhoepsilon_4^{2(a-k_4+2)}\rhoepsilon_3^{2(a-k_3+1+2)}\rhoepsilon_2^{2(a-k_2+1)}\rhoepsilon_1^{2(a-k_1+1)}}.
\ea
When $\vec S_\b=(S_1,S_2)=({\bf Q}_0,I)$ we have $\tilde  T^{({\bf m}),\beta}_\tau
=T^{({\bf m}),\beta}_\tau
$ and
\ba
 T^{({\bf m}),\beta}_\tau
&=& \bra I(u^\tau \,\cdotp  \B_4 u_4), h\,\cdotp  \B_3 u_3\,\cdotp {\bf Q}_0(\B_2 u_2\,\cdotp \B_1 u_1)\cet,\\
&=& {i{\P_\beta}\,\det (A)C}
\frac {(i\tau)^{-(10+4a-|\vec k_\b|)} (1-O(\tau^{-1}))}{\rhoepsilon_4^{2(a-k_4+1)}\rhoepsilon_3^{2(a-k_3+1)}\rhoepsilon_2^{2(a-k_2+2)}\rhoepsilon_1^{2(a-k_1+1+2)}}\, ,
\ea
and finally when $\vec S_\b=(S_1,S_2)=(I,I)$ 
\ba
\tilde  T^{({\bf m}),\beta}_\tau
&=&\bra  u^\tau, h\,\cdotp I(\B_4u_4\,\cdotp \B_3u_3)\,\cdotp I(\B_2u_2\,\cdotp \B_1 u_1)\cet\\
&=&{{\P_\beta} C_a\det (A)}
\frac {(i\tau)^{-(8+4a-|\vec k_\b|)}(1+O(\tau^{-1}))}{\rhoepsilon_4^{2(a-k_4+1)}\rhoepsilon_3^{2(a-k_3+1)}\rhoepsilon_2^{2(a-k_2+1)}\rhoepsilon_1^{2(a-k_1+1)}}.
\ea

Next we consider all $\beta$ such that  $\vec S^\beta=({\bf Q}_0,{\bf Q}_0)$
but
$\b\not=\beta_1$. Then
\ba
 \tilde T^{({\bf m}),\beta}_\tau
=C\det (A)
(i\tau)^{-(6 +4a)}(1+O(\tau^{-1}))
\vec\rhoepsilon^{\,-2(\vec a+\vec 1)+2\vec k_\b}\rhoepsilon_4^{-2}\rhoepsilon_2^{-2}\rhoepsilon_3^{-4}\rhoepsilon_1^{-4}\cdotp{\P_{\beta}}
\ea where
$\vec k_\b$ is as in (\ref{k: s for (Q_0,Q_0) and T}). 
Note that for $\b=\beta_1$ we have 
$\P_{\beta_1}= (\D +O(\rhoepsilon_1))\rhoepsilon_1^4\cdotp\rhoepsilon_1^4\cdotp\rhoepsilon_1^4$
while $\b\not=\beta_1$ we just use an estimate
 ${\P_\beta}=O(1)$.
Then we see that
$\tilde  T^{({\bf m}),\beta}_\tau\prec \L_\tau$.

When  $\beta$ is such that $\vec S^\beta=({\bf Q}_0,I)$, we see that 
\ba
 T^{({\bf m}),\beta}_\tau
=C\det (A)
(i\tau)^{-(10 +4a-|\vec k_\b|)}(1+O(\tau^{-1}))
\vec\rhoepsilon^{\,-2(\vec a+\vec 1)+2\vec k_\b}
\rhoepsilon_4^{0}\rhoepsilon_2^{-2}\rhoepsilon_3^{0}\rhoepsilon_1^{-4}{\P_\beta},\\
\tilde  T^{({\bf m}),\beta}_\tau
=C\det (A)
(i\tau)^{-(10 +4a-|\vec k_\b|)}(1+O(\tau^{-1}))
\vec\rhoepsilon^{\,-2(\vec a+\vec 1)+2\vec k_\b}
\rhoepsilon_4^{-2}\rhoepsilon_2^{0}\rhoepsilon_3^{-4}\rhoepsilon_1^{0}{\P_\beta}\
\ea where ${\P_\beta}=O(1)$ and 
$\vec k_\b$ is as in (\ref{k: s for (Q_0,I) and tilde T}) 
and hence $ T^{({\bf m}),\beta}_\tau\prec \L_\tau$ and 
$\tilde  T^{({\bf m}),\beta}_\tau\prec \L_\tau$.

When  $\beta$ is such that $\vec S^\beta=(I,{\bf Q}_0)$ or 
$\vec S^\beta=(I,I)$, using inequalities  of the type (\ref{k: s for (I,Q_0) and T}) and
 (\ref{k: s for (I,Q_0) and tilde T}) in appropriate cases,
 we see that $T^{({\bf m}),\beta}_\tau\prec \L_\tau$ 
 and $\tilde  T^{({\bf m}),\beta}_\tau\prec \L_\tau$.

%
%
%
%
%


The above shows that all terms 
$T^{({\bf m}),\beta}_\tau$ and $\tilde  T^{({\bf m}),\beta}_\tau$ with
maximal allowed $k$.s
have asymptotics with the same power of $\tau$ but their $\vec \rhoepsilon$ asymptotics vary, and when the
asymptotic orders of $\rhoepsilon_j$
are given as explained in after (\ref{eq: ordering of epsilons}), there
is only one term, namely $\L_\tau=T^{({\bf m}),\beta_1}_\tau$, that has the strongest 
order asymptotics given in (\ref{leading term}).

Next we analyze the effect of the permutation $\sigma:
\{1,2,3,4\}\to \{1,2,3,4\}$  of the indexes $j$ of the waves $u_{j}$. We assume below that the permutation 
 $\sigma$ is not the identity map.
%

 Recall that in the computation (\ref{first asymptotical computation})
 there appears a term $\omega_{45}^{-1}\sim \rhoepsilon_4^{-2}$. 
As this term does not appear 
  in the computations of the terms $\tilde  T^{({\bf m}),\beta}_{\tau,\sigma}$,
 we see that  $\tilde  T^{({\bf m}),\beta}_{\tau,\sigma}\prec \L_\tau$. 
Similarly, if
  $\sigma$  is such that $\sigma(4)\not =4$, the
 term $\omega_{45}^{-1}$ does not appear in the computation of
$ T^{({\bf m}),\beta}_{\tau,\sigma}$
 and hence
$ T^{({\bf m}),\beta_1}_{\tau,\sigma}\prec  \L_\tau$. 
Next we consider the permutations
 for which  $\sigma(4) =4$.

Next we consider $\sigma$ that is either
$\sigma=(3,2,1,4)$ or
$\sigma=(2,3,1,4)$.
 These terms are very similar
and thus we analyze the case when $\sigma=(3,2,1,4)$.
First we consider the case when $\beta=\beta_2$ is such that $\vec S^{\beta_2}=({\bf Q}_0,{\bf Q}_0)$,
 $\vec k_{\beta_2}=(2,0,4,0)$  
This term appears in the analysis of
the term  $A^{(1)}[  u_{\sigma(4)},{\bf Q}(A^{(2)}[ u_{\sigma(3)},{\bf Q}(A^{(3)}[ u_{\sigma(2)}, u_{\sigma(1)}])])]$ when $(A^{(1)},A^{(2)},A^{(3)})=(A_2,A_1,A_2)$, see (\ref{A alpha decomposition2}). 
By a permutation of the indexes in (\ref{first asymptotical computation}) we obtain 
the formula 
\beq\label{Formula sic!}
T^{({\bf m}),{\beta_2}}_{\tau,\sigma}&=&
c_1^{\prime}\det (A){(i\tau)^{-(6 +4a)}(1+O(\tau^{-1}))}\vec \rhoepsilon^{-2(\vec a+\vec 1)}
\\ \nonumber
& & \cdotp (\omega_{45}\omega_{32})^{-1}  \rhoepsilon_4^{2(k_4-1)}\rhoepsilon_1^{2k_3}\rhoepsilon_2^{2(k_2-1)}\rhoepsilon_3^{2(k_1-1)}
 \P_{\beta_2},\\
\nonumber
\tilde \P_{{\beta_2}}&=&(v^{pq}_{(4)}b^{(1)}_pb^{(1)}_q)
(v^{rs}_{(3)}b^{(1)}_rb^{(1)}_s)
(v^{nm}_{(2)}b^{(3)}_nb^{(3)}_m)\D. 
\hspace{-1.5cm}
\eeq
Hence, in the case when we use the polarizations (\ref{chosen polarization}), we obtain
\ba
  T^{({\bf m}),\beta_2}_{\tau,\sigma}=
c_1{(i\tau)^{-(6 +4a)}(1+O(\tau^{-1}))}\vec \rhoepsilon^{-2(\vec a+\vec 1)}
 \rhoepsilon_4^{-4}\rhoepsilon_2^{-2}\rhoepsilon_3^{0+4}\rhoepsilon_1^{6+8}\D.
 \ea
Comparing the power of $\rhoepsilon_3$ in the above expression,
we see that in this case $ T^{({\bf m}),\beta_2}_{\tau,\sigma}\prec \L_\tau$. 
When $\sigma=(3,2,1,4)$, we see in  a straightforward way  also for 
 other $\beta$  for which  $\vec S^{\beta}=({\bf Q}_0,{\bf Q}_0)$, that  $ T^{({\bf m}),\beta}_{\tau,\sigma}\prec \L_\tau$.


When $\sigma=(1,3,2,4)$,
we see that for all $\beta$ with  $|\vec k_\beta|=6$,
\ba
  T^{({\bf m}),\beta}_{\tau,\sigma}\hspace{-1mm}&=&\hspace{-1mm}
C\det (A){(i\tau)^{-(6 +4a)}(1+O(\frac 1\tau))}\vec \rhoepsilon^{-2(\vec a-\vec k+\vec 1)}\omega_{45}^{-1}\omega_{13}^{-1}
 \rhoepsilon_4^{-2}\rhoepsilon_2^{0}\rhoepsilon_3^{-2}\rhoepsilon_1^{-2}\P_{\beta}\\
 \hspace{-1mm}
 &=&\hspace{-1mm}
C\det (A){(i\tau)^{-(6 +4a)}(1+O(\frac 1\tau))}\vec \rhoepsilon^{-2(\vec a-\vec k+\vec 1)}
 \rhoepsilon_4^{-4}\rhoepsilon_2^{0}\rhoepsilon_3^{-2}\rhoepsilon_1^{-4}\P_{\beta}. \ea
 Here,  $\P_{\beta}=O(1)$.
Thus when $\sigma=(1,3,2,4)$, by comparing the powers of $\rhoepsilon_2$ we see that $ T^{({\bf m}),\beta}_{\tau,\sigma}\prec \L_\tau$. The same holds in the case when
$|\vec k_\beta|<6$.
The case when  $\sigma=(3,1,2,4)$ is similar to $\sigma=(1,3,2,4)$. 
This proves Proposition \ref{SL:order}\hfill\Box\medskip

\hiddenfootnote{
THIS WAS
Next, in the  case  when $\sigma=(3,1,2,4),$ we need to analyze the term
\ba
 & &\bra u^\tau ,A_2[u_4,{\bf Q}(A_2[u_2,{\bf Q}(A_1[u_1,u_3])])]\cet\\
&=& \bra u^\tau, u^{pq}_4\,{\bf Q} ((\p_p\p_q\p_r\p_s \p_j\p_k u_1)\,{\bf Q}(u^{rs}_2\, u^{jk}_3))\cet
+\hbox{similar terms}.
\ea
The first term on the right hand side, that we index by  $\beta=\beta_2$,
corresponds to $k_{\beta_2}=(0,0,6,0)$. Then we see that $\P_{\beta_2}=\P_{\beta_0}$ and
\ba
  T^{({\bf m}),\beta_2}_{\tau,\sigma}\hspace{-1mm}&=&\hspace{-1mm}
C\det (A){(i\tau)^{-(6 +4a)}(1+O(\frac 1\tau))}\vec \rhoepsilon^{-2(\vec a+\vec 1)}\omega_{45}^{-1}\omega_{32}^{-1}
 \rhoepsilon_4^{-2}\rhoepsilon_2^{0}\rhoepsilon_3^{-2}\rhoepsilon_1^{4}\P_{\beta_2}\\
 \hspace{-1mm}&=&\hspace{-1mm}
C\det (A){(i\tau)^{-(6 +4a)}(1+O(\frac 1\tau))}\vec \rhoepsilon^{-2(\vec a+\vec 1)}
 \rhoepsilon_4^{-4}\rhoepsilon_2^{0}\rhoepsilon_3^{-4}\rhoepsilon_1^{4}\P_{\beta_2}. \ea
Comparing the powers of $\rhoepsilon_2$ we see that $ T^{({\bf m}),\beta_2}_{\tau,\sigma}\prec \L_\tau$.
When $\sigma=(3,1,2,4)$ we see in a straightforward manner for also values of 
$\beta\not=\beta_2$ 
 that   $ T^{({\bf m}),\beta}_{\tau,\sigma}\prec \L_\tau$.}

Summarizing; we have analyzed the terms $ T^{({\bf m}),\beta}_{\tau,\sigma}$ 
corresponding to any $\b$ and all $\sigma$
except $\sigma=\sigma_0=(2,1,3,4)$.
Clearly, the sum  $ \sum_\b T^{({\bf m}),\beta}_{\tau,\sigma_0}$ 
is equal to the sum $ \sum_\b T^{({\bf m}),\beta}_{\tau,id}$.
Thus, when
the asymptotic orders of $\rhoepsilon_j$
are given  in (\ref{eq: ordering of epsilons})
and
the polarizations satisfy (\ref{chosen polarization}), we have
\beq\nonumber
\mathcal G^{({\bf m})}(v,{\bf b})&=&\lim_{\tau\to \infty}
\sum_{\b,\sigma} \frac {T^{({\bf m}),\beta}_{\tau,\sigma}}{(i\tau)^{(6 +4a)}}
=\lim_{\tau\to \infty}
 \frac {2 T^{({\bf m}),\beta_1}_\tau
(1+O(\rhoepsilon_1))} {(i\tau)^{(6 +4a)}}\\
\label{eq: summary}
&=&2c_1\det (A)
(1+O(\rhoepsilon_1))\,\vec \rhoepsilon^{-2(\vec a+\vec 1)}
\rhoepsilon_4^{-4}\rhoepsilon_2^{-2}\rhoepsilon_3^{0}\rhoepsilon_1^{20}\D.\hspace{-1.5cm}
\eeq
Notice that here $\rhoepsilon_j$ depend only on $b^{(k)}$, $k=1,2,3,4,5$.

Let $Y=\hbox{sym}(\R^{4\times 4})$ 
and consider the quadratic form $B:(v,w)\mapsto  \hat g_{nj}\hat g_{mk} v^{nm}w^{jk}$
as a  inner product in $Y$.
Then
$ \D =B(v^{(5)},v^{(1)})$. 

 
Let $L(b^{(j)})$ denote the subspace of dimension $6$ of the symmetric
 matrices $v\in Y$ that satisfy 
 equation (\ref{divergence condition for symbol}) with covector $b^{(j)}$.
%
%

%

%
%

 Let $\L$ be the real analytic manifold consisting  of  $\eta=({\bf b}, {\underline v},V^{(1)},V^{(5)})$,
where ${\bf b}=(b^{(1)},b^{(2)},b^{(3)} ,b^{(4)},
b^{(5)})$ is a sequence of  light-like vectors and ${\underline v} =(v^{(2)},v^{(3)},v^{(4)})$ 
satisfy $v_{(j)}\in L(b^{(j)})$ for all $j=2,3,4$, 
and $V^{(1)}=(v_p^{(1)})_{p=1}^6$ be basis of $L(b^{(1)})$ and 
 $V^{(5)}=(v_p^{(5)})_{p=1}^6$ be vectors in $Y$ such 
 that $B(v_p^{(5)},v_q^{(1)})=\delta_{pq}$ for $p\leq q$.
 
 We define for $\eta\in \L$ 
\ba
\kappa(\eta):=\det\bigg(\mathcal G^{({\bf m})}({\bf v}_{(p,q)},{\bf b})\bigg) _{p,q=1}^6,\ \hbox{where}\ {\bf v}_{(p,q)}=(v^{(1)}_p,v^{(2)},v^{(3)},v^{(4)},v_q^{(5)}).
\ea
Then $\kappa(\eta)$ is a real-analytic function on $\L$.
 
 Let us next  consider linearly independent light-like vectors,
 $\hat {\bf b}= (\hat b^{(1)},\hat b^{(2)},\hat b^{(3)} ,\hat b^{(4)},\hat b^{(5)})$ satisfying
 (\ref{b distances}) with $\vec \rhoepsilon$ 
 given in (\ref{eq: ordering of epsilons})
 with some small $\rhoepsilon_1>0$ 
 and let the polarizations $\hat {\underline v}=(\hat v^{(2)},\hat v^{(3)},\hat v^{(4)})$ be  such that $\hat v^{(j)}\in  L(b^{(j)})$, $j=2,3,4$, are those given
 by (\ref{chosen polarization}),
  and
 $\hat V^{(1)}=(\hat v_p^{(1)})_{p=1}^6$  be a basis of $ L(b^{(1)})$.
 Let $\hat V^{(5)}=(\hat v_p^{(5)})_{p=1}^6$ be vectors in $Y$ such 
 that $B(\hat v_p^{(5)},\hat v_q^{(1)})=\delta_{pq}$ for $p\leq q$.
When $\rhoepsilon_1>0$ is small enough,  formula (\ref{eq: summary}) 
 yields that $\kappa(\hat \eta)\not =0$ for   $\hat \eta=(\hat {\bf b},\hat {\underline v},\hat V^{(1)},\hat V^{(5)})$.
 As $\kappa(\eta)$ is a real-analytic function on $\L$,
 we see that $\kappa(\eta)$ is non-vanishing on a generic
 subset of the component of $\L$ containing $\hat \eta$.
 Note that for any ${\bf b}$ there is $\eta=({\bf b},{\underline v},V^{(1)},V^{(5)})$
 that is in this component.
 
 Consider next $\eta=({\bf b},{\underline v},V^{(1)},V^{(5)})$
 that is in the component of $\L$ containing $\hat \eta$.
 As  $v^{(5)}\mapsto \mathcal G^{({\bf m})}({\underline v},{\bf b})$ is linear
 and thus $\kappa(\eta)$ can be considered as an alternative $6$-multilinear
 form of $V^{(5)}$.
 Thus if $\tilde v^{(5)}=\sum_{p=1}^6
 a_pv_p^{(5)}\not =0$ is such that 
$G^{({\bf m})}(v^{(1)},{\underline v},\tilde v^{(5)},{\bf b})=0$ for all
$v^{(1)}\in  L(b^{(1)})$,
we see that $\kappa(\eta)=0$. As the image of a open and dense set
in the projection $({\bf b},{\underline v},V^{(1)},V^{(5)})\mapsto ({\bf b},{\underline v})$
is open and dense, we conclude that for an open and dense 
set of pairs $ ({\bf b},{\underline v})$ there is $V^{(5)}$
so that for all $\tilde v^{(5)}\in \hbox{span}(V^{(5)})$   there is
$v^{(1)}\in  L(b^{(1)})$ such that $G^{({\bf m})}(v^{(1)},{\underline v},\tilde v^{(5)},{\bf b})\not=0$.
%
\hfill \Box \medskip

{

\section{Observations in normal coordinates}

\label{sec: normal coordinates}
\subsection{Detection of singularities}
\HOX{Improve the text on motivation.}
Above we have considered the singularities of the metric $g$ in the wave gauge
coordinates, that is, we have used the coordinates of manifold $M_0$ where 
the metric $g$ solves the $\hat g$-reduced Einstein equations. As the wave gauge coordinates may  also be non-smooth,
we do not know if the observed singularities are caused by the metric or
the coordinates. Because of this, we next consider the metric in normal coordinates.

%

%



Let  $(g^{\vec \e},\phi^{\vec \e})$ be the solution 
 of the  $\hat g$-reduced Einstein equation (\ref{eq: adaptive model}) with
 the source ${\bf f}_{\vec \e}$ given in (\ref{eq: f vec e sources}). 
 We emphasize that  $g^{\vec \e}$ is the metric in the
$\hat g$-wave gauge coordinates.

Let $(z,\eta)\in \U_{(z_0,\eta_0)}(\hat h)$ and denote by $\mu_{\vec \e}([-1,1])=\mu_{g^{\vec \e},z,\eta}([-1,1])$ the freely falling observers, i.e. time-like geodesic, 
on $(\hattuM _0,g^{\vec \e})$
having the same initial data as $\mu_{\hat g,z,\eta}$, see Sect. \ref{sec:notations 1}. When $\tilde r\in (-1,1]$, we call the family 
$\mu_{\vec \e}([-1,\tilde r])$, where $\vec \e=(\e_j)_{j=1}^4$, $\e_j\geq 0$ the observation
geodesics and note that $\vec \e$ may be here the  zero-vector.
For given $(z,\eta)$, let us choose  $(Z_j({-1}))_{j=1}^4$ to be linearly
independent  vectors at $\mu_{\hat g,z,\eta}(-1)$ such that $Z_1({-1})=\dot \mu_{\hat g,z,\eta}(-1)$.
For  $s\in [-1, 1]$, let 
$Z_{j,\vec \e}(s)$ be the parallel translation of $Z_j({-1})$
along $\mu_{\vec \e}$.
Also, assume that $\mu_{\vec \e}([-1,1])\subset U_{\vec\e}$ for $|\vec\e|$ small enough
and denote $p_{\vec \e}=\mu_{\vec \e}(\tilde r)$ where $ \tilde r<1$.

%
%
%
%

Let then $\Psi_{\vec \e}$ denote  normal coordinates
of $(\hattuM _0,g^{\vec \e})$ defined using the center $p_{\vec \e}$
and the frame $Z_{j,\vec \e}$, $j=1,2,3,4.$
Below, we denote $ \mu_0=\mu_{\vec\e}|_{\vec \e=0}$
and $Z_{j,0}= Z_{j,\vec \e}|_{\vec \e=0}.$ We say that these normal
coordinates are associated to the observation geodesics $\mu_{\vec \e}$,  see Fig.\ 15. 


Next we consider the metric $g^{\vec \e}$ in the normal coordinates and
 study when $\p_{\vec \e}^4 (\Psi_{\vec \e})_*g^{\vec \e}|_{\vec \e=0}$
is smooth. 
Below, we denote  $g_{\vec \e}= g^{\vec \e}$ and
$U_{\vec \e}=U_{g^{\vec \e}}$. Recall that $\hat U=U_{\hat g}$. In the next Lemma, we 
consider observations in  normal coordinates, see   Fig.\ 15. 

%


\begin{figure}[htbp] \label{Fig-14}
\begin{center}

\psfrag{1}{$q$}
\psfrag{2}{\hspace{-2mm}$x$}
\includegraphics[width=7.5cm]{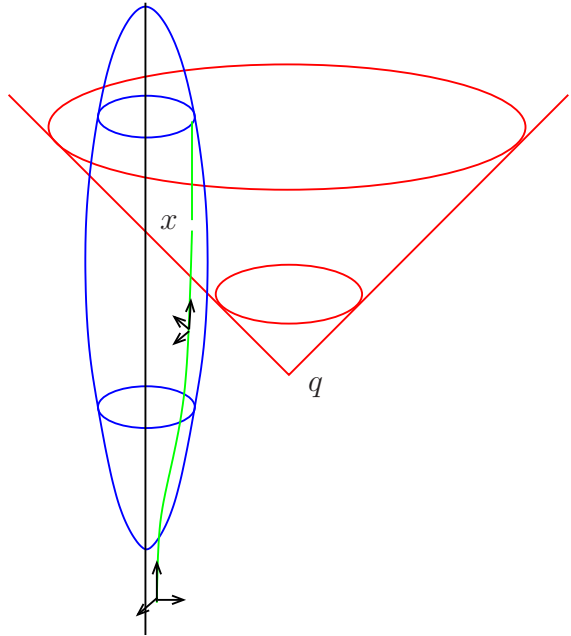}
\end{center}
\caption {A schematic figure where the space-time is represented as the  3-dimensional set $\R^{1+2}$. The future light cone $\mathcal L^+_{\hat g}(q)$ corresponding to the point $q$ is  shown
as a red cone. The green curve is the geodesic
 $\mu_0=\mu_{\hat g,z,\eta}$.
This geodesic intersect the future light cone $\mathcal L^+_{\hat g}(q)$
at the point $x$.
The black vectors are the frame $(Z_j)$ that is obtained using parallel translation along
the geodesic $\mu_0$. Near the intersection point $x$ we use the normal coordinates
centered at $x$ and associated to the frame obtained via parallel translation. 
}
 \end{figure}

\begin{lemma}\label{lem: sing detection in wave gauge coordinates 1}
Let $u^{\vec \e}=(g^{\vec \e},\phi^{\vec \e})$ be the solution of the reduced Einstein equations (\ref{eq: adaptive model})
with the source ${\bf f}_{\vec \e}$ given in (\ref{eq: f vec e sources}), and $ \mu_{\vec \e}([-1,\tilde r])$, $-1<\tilde r<1$ be  the observation geodesics in $U_{\vec \e}$. 
Let us consider at the point $p_{\vec \e}=\mu_{\vec \e}(\tilde r)$ the 
frame $Z_{j,\vec \e}=Z_{j,\vec \e}(\tilde r)$ and
let $\Psi_{\vec\e}: W_{\vec \e}\to  \Psi_{\vec \e}(W_{\vec \e})\subset
\R^4$ be  the normal coordinates  centered at  $p_{\vec \e}$ and associated to the frame $(Z_{j,\vec \e})_{j=1}^4$. Denote $W=W_0$.



   Let
 $S\subset \hat U$ be a smooth 3-dimensional surface such that $p_0=p_{\vec \e}|_{\vec \e=0}\in S$ 
and
\ba
& & g^{(\alpha)}=\p_{{\vec \e}}^\alpha g^{\vec \e}|_{\vec \e=0},\quad
\phi_\ell^{(\alpha)}=\p_{{\vec \e}}^\alpha\phi_\ell^{\vec \e}|_{\vec \e=0},\quad\hbox{for } |\alpha|\leq 4,\ \alpha\in \{0,1\}^4,
\ea
and assume that $ g^{(\alpha)}$ and $\phi_l^{(\alpha)}$ are in $C^\infty(W)$ for $|\a|\leq 3$ and
 $g^{(\alpha_0)}_{pq}|_W\in\I^{m_0}(W\cap S)$ and $ \phi_l^{(\alpha_0)}|_W\in\I^{m_0}(W\cap S)$
for $\a_0=(1,1,1,1)$.  \HOX{Slava: Consider using the notation $\I^{m_0}(W;W\cap S)= \I^{m_0}_{W,\hat g}(W\cap S)$.}
\smallskip

(i) Assume that   $S\cap W$ is empty. Then
$\p_{{\vec \e}}^4 ((\Psi_{\vec\e})_*g_{\vec \e})|_{\vec \e=0}$ and
 $\p_{{\vec \e}}^4 ((\Psi_{\vec\e})_*\phi_\ell ^{\vec \e})|_{\vec \e=0}$
 are  $C^\infty$-smooth 
in $\Psi_0(W)$.

\smallskip

(ii)
Assume that  
 $ \mu_{0}([-1,1])$ intersets $S$
 transversally  at $p_0$.
 Consider the conditions

 \smallskip

 \noindent (a) There is a 2-contravariant tensor field $v$ that is a smooth section of $TW\otimes TW$  such that
$v(x)\in T_xS\otimes T_xS$ for $x\in S$  
and the principal symbol of
$ \bra v,g^{(\alpha_0)}\cet\in \I^{m_0}(W;W\cap S)$
is non-vanishing at $p_0$.
 \smallskip

 \noindent (b) The principal symbol of $ \phi_\ell^{(\alpha_0)}$ is non-vanishing
at $p_0$ for some $\ell=1,2,\dots,L$.
 \smallskip

If (a) or (b) holds, then either  $\p_{{\vec \e}}^4 ((\Psi_{\vec\e})_*g_{\vec \e})|_{\vec \e=0}$ or
 $\p_{{\vec \e}}^4 ((\Psi_{\vec\e})_*\phi_\ell ^{\vec \e})|_{\vec \e=0}$ is not $C^\infty$-smooth 
in $\Psi_0(W)$.

\end{lemma}

\HOX{  
Note that we could change notation $\p^4$ to $\p^{\a_0}$ as Slava suggested.}

{\bf Proof.}  (i) Using the metric $g_{\vec \e}$ in the 
$\hat g$-wave gauge coordinates we can
compute the $\Psi_{\vec \e}$-coordinates and thus
find  $(\Psi_{\vec \e})_*g_{\vec \e}$. As $g^{(\alpha)}$, and thus $\p^\a_{\vec \e}\Psi_{\vec \e}|_{\vec \e=0}$
 for $|\a|\leq 3$
 and $g^{(\alpha_0)}$ and $\p^{\a_0}_{\vec \e}\Psi_{\vec \e}|_{\vec \e=0}$ are smooth in 
 $\Psi_0(W)$,
 we see that $\p^4_{\vec \e}(\Psi_{\vec \e})_*g_{\vec \e}|_{\vec \e=0}$
is smooth in  $\Psi_0(W)$. This proves (i).

 (ii)  Let $g^{\vec \e}$ denote the metric in the $\hat g$ wave gauge
coordinates on $\hattuM _0$.
%
%
%
%
Denote $\gamma_{\vec \e}(t)=\mu_{\vec \e}(t+\tilde r)$.
Let $X:W_0\to V_0\subset \R^4$, $X(y)=(X^j(y))_{j=1}^4$ be such local coordinates in $W_0$
that $X(p_0)=0$ and $X(S\cap W_0)=\{(x^1,x^2,x^3,x^4)\in V_0;\ x^1=0\}$  
and $y(t)=X( \gamma_0(t))=(t,0,0,0)$. 
 Note that the   coordinates $X$ are
independent of ${\vec \e}$.

Let $\hat Y_j(t)$ be $\hat g$-parallel vector fields on $\gamma_0(t)$
such that $\hat Y_j(0)= \p/\p X^j$ are the coordinate vector fields corresponding to the coordinates $X$. Let $b_j^k\in \R$ be such that 
$\hat Y_j(0)=b_j^kZ_{0,k}(\tilde r)$ and define
$Y^{\vec \e}_j(t)=b_j^kZ_{\vec \e,k}(t+\tilde r)$.

To do computations
in local coordinates, let us denote
\ba
\tilde g^{(\alpha)}=\p_{{\vec \e}}^\alpha (X_*g^{\vec \e})|_{\vec \e=0},\quad\tilde \phi_\ell^{(\alpha)}=\p_{{\vec \e}}^\alpha (X_*\phi_\ell^{\vec \e})|_{\vec \e=0},
\quad\hbox{for } |\alpha|\leq 4,\ \alpha\in \{0,1\}^4.
\ea
We note that as the $X$-coordinates are independent of $\vec \e$,
we have $\p_{{\vec \e}}^\alpha (X_*g^{\vec \e})|_{\vec \e=0}=X_*(\p_{{\vec \e}}^\alpha (g^{\vec \e})|_{\vec \e=0})$.
In the local $X$ coordinates, let $v(x)=v^{pq}(x)\frac\p{\p x^p}\frac\p{\p x^q}$  be such that  $v^{pq}(x)=0$ if $(p,q)\not \in\{2,3,4\}^2$ and $x\in X(S\cap W_0)$.

%

Let    $R^{\vec \e}$  be the curvature tensor  of $g^{\vec \e}$
and  define the functions
\ba
h^{\vec\e}_{mk}(t)&=&g^{\vec \e}(R^{\vec \e}(\dot \gamma_{\vec \e}(t),Y^{\vec \e}_m(t))\dot \gamma_{\vec \e}(t),Y^{\vec \e}_k(t)),\\
J_v(t)&=&
\p_{{\vec \e}}^4( \hat v^{mq}\,
h^{\vec \e}_{mq}(t))|_{{\vec \e}=0},
\ea
where $\hat v^{mq}=v^{mq}(0)\in \R^{4\times 4}$. Here we can consider $\hat v^{mq}$
as a constant matrix or alternatively, a tensor field whose representation in the $X$ coordinates
is given by a constant matrix.

%
%
%
%
Observe that  $J_v(t)$ is a function defined
on the geodesic $\gamma_0(t)$, $t\in I$ that is parametrized along arc length
and thus it has a coordinate invariant definition.
As all $\vec\e$ derivatives of order 3 or less 
of $g_{\vec \e}$ are smooth, we see that if $\p_{{\vec \e}}^4 ((\Psi_{\vec\e})_*g_{\vec \e})|_{\vec \e=0}$ would
be smooth near $0\in \R^4$, then also $\p_{{\vec \e}}^4 ((\Psi_{\vec\e})_*R^{\vec \e})|_{\vec \e=0}$
and $\p_{{\vec \e}}^4 ((\Psi_{\vec\e})_*\dot \gamma_{\vec \e})|_{\vec \e=0}$
and  $\p_{{\vec \e}}^4 ((\Psi_{\vec\e})_* Y_k^{\vec \e})|_{\vec \e=0}$ would be smooth,
and thus also the
function $J_v(t)$ would be smooth near $t=0$. Hence, to show that the $\p_{\vec\e}^4$-derivatives of the metric tensor  in the normal
coordinates are not smooth, 
it is enough to show that  for some values of $\hat v^{mq}$ the function
$J_v(t)$ is non-smooth  at $t=0$.

The curvature tensor and thus $h^{\vec\e}_{mk}$ can be written  in the $X$
coordiantes as a sum of terms which are products of  the $x$-derivatives
of $g_{\vec \e}$ up to order 2, its inverse matrix $g^{-1}_{\vec \e}$, evaluated at
$\gamma_{\vec \e}(t)$, and the vector fields  $\dot \gamma_{\vec \e}(t)$ and 
$Y_k^{\vec \e}(t)$.
When we apply the product rule in the differentation and the chain rule
(when we compute  e.g.\ $\p_{\vec\e}(g^{jk}_{\vec \e}(\gamma^{\vec \e}(t))$), the first derivative $\p_{\e_1}$ operates either to
 $\p^\beta_x g_{\vec \e}$, with $|\beta|\leq 2$, or  $\dot \gamma_{\vec \e}(t)$, or
 $Y_k^{\vec \e}(t)$, or due to chain rule, it produces the $x$-derivatives
   of   $\hat g$  multiplied by $\p_{\e_1}\gamma_{\vec \e}(t)$.
 As in $W$ all $\e$-derivatives of the metric tensor $g_{\vec \e}$ 
 up to order 3 are smooth,
 we see that all other $\e$-derivatives, namely  $\p_{\e_j}$, $j=2,3,4$
 have to operate on the same term on which the $\p_{\e_1}$ derivative
 operated or otherwise, the produced term is $C^\infty$-smooth in the $X$-coordinates.
 Thus we need to consider only terms where all four $\e$-derivatives operate
on the same term. 

 Below,  $\Gamma_{\vec \e}$ are
 the connection coefficients corresponding to 
$g_{\vec \e}$. 
We will work in the $X$ coordinates and  \HOX{Choose a better notation for $\tilde R$.}
denote $\tilde R(x)=\p_{\vec \e}^4 X_*(R_{\vec \e})(x)|_{{\vec \e}=0}$,
and $\tilde \Gamma^j_{nk}$ and $\tilde \gamma^j$ are  the analogous
4th order $\e$-derivatives, and denote $\tilde g=\tilde g^{(\alpha_0)}$. 
  For simplicity we also denote
   $X_*\hat g$ and $ X_*\hat \phi_l$  by $\hat g$ and $\hat \phi_l$, respectively.

We analyze the functions of $t\in I=(-t_1,t_1)$, e.g., $a(t)$, where  $t_1>0$ is small. We say
that $a(t)$ is of order $n$ if $a(\,\cdotp)\in \I^{n}(\{0\}).$

When  $a(t)$ solves an ordinary differential equation (ODE) of the type $\p_ta(t)+K(t)a(t)=b(t)$, $a(0)=a_0$
where $K(t),b(t)\in \I^{n}(\{0\})$, we say that $a(t)$ solves an ODE involving
$K(t)$. When $n<-1$, this implies, due to \cite{GU1} and bootstrap arguments, that
then $a\in \I^{n-1}(\{0\})$, i.e., $a$ is one order smoother than $b$ and $K$.

When $t_0$ is small enough,  $S\subset \hattuM _0$ intersects $\gamma_0=
\gamma_{\vec \e}|_{\vec \e=0}$ only at the point $p_0$, and the
intersection there is transversal. Then, we see that the restrictions of conormal distributions in
 $\I(\hattuM _0;S)$ on $\gamma_0$ are conormal distributions associated to the submanifold
 $\{p_0\}\subset \gamma_0$. Thus 
by the assumptions of the theorem, $\big(\p_x^\beta \tilde g_{jk}\big)(\gamma(t))\in \I^{m_0+\beta_1}(\{0\})$
when $\beta=(\beta_1,\beta_2,\beta_3,\beta_4)$.
As
 $\hat g$ and $\hat \phi_l$ are $C^\infty$-smooth in $W_0$
 and  the geodesic $\gamma_0$
intersects $S$ transversally,
we see that $ (\tilde \gamma(t),\p_t \tilde  \gamma(t))=\p^4_{\vec \e} (\gamma^{\vec \e}(t),\p_t \gamma^{\vec \e}(t))|_{\vec \e=0}$ and 
$\tilde Y_k=\p^4_{\vec \e} Y_k^{\vec \e}|_{\vec \e=0}$  in the 
$X$-coordinates
are solutions of ODEs (the latter one obtained
by differentiating, with respect to $\vec\e$,  the equation $\nabla^{\vec \e}_{\dot\gamma^{\vec \e}}Y_k^{\vec \e}=0$) with coefficients depending on
 the Christoffel symbols, i.e., on the derivatives 
 $\p^4_{\vec \e}(\p_j  g_{pq}^{\vec \e})|_{{\vec \e}=0}\in  \I^{m_0+1}(\{0\})$.
  Thus $(\tilde \gamma(t),
 \p_t \tilde \gamma(t))$ and $\tilde Y_k$ 
 are in $ \I^{m_0}(\{0\})$. 
 


As the curvature $R_{\vec \e}$ depends on the 2nd order derivatives of the metric,
the above analysis shows that $\tilde R|_{\gamma_0(t)}\in  \I^{m_0+2}(\{0\})$. Thus in the $X$ coordinates
$\p_{\vec \e }^4(h^{\vec \e}_{mk}(t))|_{{\vec \e}=0}\in \I^{m_0+2}(\{0\})$ can be written as
\ba
\p_{\vec \e }^4(h^{\vec \e}_{mk}(t))|_{{\vec \e}=0}
&=& \hat g( \tilde  R(\dot  \gamma_0(t),\hat Y_m(t))
\dot \gamma_0(t),\hat Y_k(t))+\hbox{smoother terms}\\
&=& 
 \hat g_{kq}
\tilde  R^q_{11m}+\hbox{smoother terms}\\
&=& \hat g_{qk}
(\frac \p{\p x^1}\tilde  \Gamma^q_{1m}-
\frac \p{\p x^m}\tilde  \Gamma^q_{11})
+\hbox{smoother terms}\\
&=&\frac 12
\bigg(\frac \p{\p x^1}(\frac  {\p \tilde  g_{km}}{\p x^1}+\frac  {\p \tilde  g_{k1}}{\p x^m}-\frac  {\p \tilde  g_{1m}}{\p x^k})\\
& &-\frac \p{\p x^m}(\frac  {\p \tilde  g_{k1}}{\p x^1}+\frac  {\p \tilde  g_{k1}}{\p x^1}-\frac  {\p \tilde  g_{11}}{\p x^k})\bigg)
+\hbox{smoother terms},
\ea
where all "smoother terms" are in  $\I^{m_0+1}(\{0\})$.
%

%
Later, we will use the fact that $\p/\p x^1$
 raises
the order of the singularity by one.

Consider next the case (a).
Assume next that   for given $(k,m)\in\{2,3,4\}^2$, the principal symbol of $ \tilde g^{(\alpha_0)}_{km}$
is non-vanishing at $0=X(p_0)$. Let $v$ be such a tensor field that $v^{mk}(0)= v^{km}(0)\not=0$
and $ v^{in}(0)=0$ when $(i,n)\not\in \{(k,m),(m,k)\}$.
Then
the above yields (in the formula below, we do not sum over $k,m$)
\ba
J_v(t)=e_{km}v^{km}\p_{\vec \e}^4(h^{\vec \e}_{mk}(t))|_{\vec \e=0}
=\frac {e_{km}}2\hat v^{km} 
\bigg(\frac \p{\p x^1}\frac  {\p \tilde  g_{km}}{\p x^1}\bigg)
+\hbox{smoother terms},
\ea
where $e_{km}= 2-\delta_{km}$. 
Thus the principal symbol of $J_v(t)$  in $\I^{m_0}(\{0\})$ is non-vanishing and
 $J_v(t)$ is not a smooth function. Thus in this case
$\p_{{\vec \e}}^4 ((\Psi_{\vec\e})_*g_{\vec \e})|_{\vec \e=0}$  is not
smooth.


Next, we consider   the case (b).
Assume that there is $\ell$ such that  the principal
symbol of the field $\tilde \phi_\ell$ is non-vanishing. 
As $\p_t\tilde \gamma(t))\in \I^{m_0}(\{0\})$, we see that $\tilde \gamma(t))\in \I^{m_0+1}(\{0\})$,
Then as $\phi_\ell^{\vec \e}$ are scalar fields, 
\beq\label{j-smoothness}
& & j_\ell(t)=
\p_{{\vec \e}}^4 \bigg(
 \phi^{\vec \e}_\ell(\gamma_{\vec \e}(t)))\bigg)
\bigg|_{\vec \e=0}\\ \nonumber
&&=\tilde  \phi_\ell(\hat \gamma(t))+
\tilde \gamma^j(t) \frac {\p\hat \phi_\ell}{\p x^j}(\hat \gamma(t))+
\hbox{smoother terms},
\eeq
where $ j_\ell \in  \I^{m_0}(\{0\})$ and the smoother terms are in $ \I^{m_0+1}(\{0\})$.
Again, if both $\p_{{\vec \e}}^4 ((\Psi_{\vec\e})_*g_{\vec \e})|_{\vec \e=0}$ and
 $\p_{{\vec \e}}^4 ((\Psi_{\vec\e})_*\phi_\ell ^{\vec \e})|_{\vec \e=0}$ are smooth, 
we see that $j_\ell(t)$ is  smooth, too. Thus to show the claim it
is enough to show that $j_\ell(t)$ is not smooth at $t=0$.  

In (\ref{j-smoothness}), $\p_1\tilde  \phi_\ell$ has the order $(m_0+1)$.
As we saw above, $\tilde \gamma(t)$ and $\p_t \tilde \gamma(t)$ have the order $m_0$. Thus $j_\ell(t)$ is not smooth which proves the claim.  
%
\hfill \Box \medskip

\begin{figure}[htbp] \label{Fig-14B}
\begin{center}

\psfrag{1}{$q$}
\psfrag{2}{\hspace{-2mm}$y$}
\includegraphics[width=5.5cm]{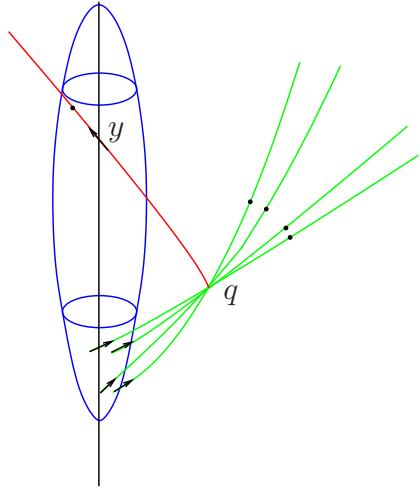}
\end{center}
\caption { A schematic figure where the space-time is represented as the  3-dimensional set $\R^{1+2}$. The light-like geodesic emanating from the point $q$  is  shown
as a red curve. The point $q$ is the intersection of
light-like geodesics corresponding to the starting points
and directions $(\vec x,\vec \xi)=((x_j,\xi_j))_{j=1}^4$. A light like geodesic
starting from $q$ passes through the  point $y$ and has the direction $\eta$ at $y$. 
The black points are the first conjugate points on the geodesics
$\gamma_{x_j(t_0),\xi_j(t_0)}([0,\infty))$, $j=1,2,3,4,$ and $\gamma_{q,\zeta}([0,\infty))$.
The figure shows the case when the interaction condition (I) is satisfied for $y\in \hat U$ with light-like
vectors $(\vec x,\vec \xi)$.
}
 \end{figure}

Next we use the above result to detect singularities
in normal coordinates. We will consider
the condition that an intersection point $q$ exists and
the light cone of $q$ intersects $y$: We say
that the {\it interaction condition} (I) is satisfied for $y\in \hat U$ with light-like
vectors $(\vec x,\vec \xi)=((x_j,\xi_j))_{j=1}^4$ and $t_0\geq0$,
 if
 \medskip
}

(I) There exists
 $q\in \bigcap_{j=1}^4\hat \gamma_{x_j(t_0),\xi_j(t_0)}((0,{\bf t}_j))$, ${\bf t}_j=\rho(x_j(t_0),\xi_j(t_0))$,
 $\zeta\in L^+_q(\hattuM _0,\hat g)$ and $t\geq 0$
 such that $y=\hat \gamma_{q,\zeta}(t)$.

 \medskip
 
 \noindent
 where $(x_j(h),\xi_j(h))$ are given in (\ref{eq: x(h) notation}) and
 the function $\rho$ is defined in (\ref{eq: max time}), see Fig.\ 16. 
 When (I) holds, we sometimes say that it holds for $y$ with
  parameters  $q$, $\zeta$, $t$, $t_0$, and $\eta=\p_t\hat \gamma_{q,\zeta}(t)$.

Let $\W_j(s)=\W_j(s;x_j,\xi_j)$ be the $s$-neighborhood of $(x_j,\xi_j)$ in $T\hattuM _0$ 
in the Sasaki metric corresponding to $\hat g^+$.

As earlier, let $  \mu_{\vec \e}([-1,1])$ be a family of observation  geodesics in $(U_{\vec \e},g_{\vec \e})$,
determined by a geodesic $\mu_0([-1,1])\subset U_{\hat g}$ of $(\hattuM _0,\hat g)$  
and $\Psi_{\vec \e}:W_{\vec \e}=B_{g^+_{\vec\e}}(p_{\vec\e},R_1)\to \R^4$ be the normal coordinates 
at the point $p_{\vec\e}=\mu_{\vec \e}(\tilde r)$, $-1\leq \tilde r<1$
associated to the frame  $(Z_{j,\vec\e }(\tilde r))_{j=1}^4$ obtained by the parallel translation along the geodesic $\mu_{\vec \e}$.

Next we investigate when some solution 
$u^{\vec \e}$ corresponding to
  $(\vec x,\vec \xi)$ has observable singularities  in the normal
coordinates $\Psi_{\vec \e}$ determined by  the observation geodesics  $\mu_{\vec \e}([-1,1])$,
that have the center $y=\mu_{\vec \e}(\tilde r)$ and the frame obtained by parallel translation along $\mu_{\vec \e}([-1,1])$.
Below, we denote $\mu_{\vec \e}=\mu_0$ when $\vec \e=\vec 0$ and say that
$\Psi_{\vec \e}$ are the normal coordinates associated 
to $\mu_{\vec \e}([-1,1])$ and $y$.

Using such normal coordinates, we define that point $y\in \hat U$, 
  satisfy the singularity {\it detection condition} ({D}) with light-like directions $(\vec x,\vec \xi)$
  and $t_0,\hat s>0$ 
 if 
 \medskip  

(D) For any $s,s_0\in (0,\hat s)$  there are $(x_j^{\prime},\xi_j^{\prime})\in \W_j(s;x_j,\xi_j)$, $j=1,2,3,4$,
and ${\bf f}_j\in {\mathcal I_{S}}^{n-3/2}(Y((x_j^{\prime},\xi_j^{\prime});t_0,s_0))$, see (\ref{eq: Z-sources}) and (\ref{eq: f vec e sources}), \HOX{Improve notation ${\mathcal I_{S}}^{n-3/2}(Y(((x_j^{\prime},\xi_j^{\prime});t_0,s_0)))$.}
and a family of observation geodesics $\mu_{\vec \e}([-1,1])$ with $y=\mu_{0}(\tilde r)$, such that when  $u_{\vec \e}$ of is the solution  of (\ref{eq: adaptive model}) 
with the source ${\bf f}_{\vec \e}=\sum_{j=1}^4
\e_j{\bf f}_j$
and $\Psi_{\vec \e}$ are the  the normal
coordinates corresponding to the center $\mu_{\vec \e}(\tilde r)$ and the frame obtained by parallel translation along $\mu_{\vec \e}([-1,1])$, then
the function 
$\p_{\vec \e}^4((\Psi_{\vec \e})_*u_{\vec \e})|_{\vec \e=0}$ is not
$C^\infty$-smooth in any neighborhood of $0=\Psi_0(y)$. 
 \medskip

Below we use the $\mathcal Y((\vec x,\vec \xi);t_0,\hat s)$ set defined in (\ref{K and X sets}).

 \medskip


\begin{lemma}\label{lem: sing detection in in normal coordinates 2}
Let $(\vec x,\vec \xi)$,   and ${\bf t}_j$ with $j=1,2,3,4$, $t_0>0$, and $x_6
\in \hat U$
satisfy (\ref{eq: summary of assumptions 1})-(\ref{eq: summary of assumptions 2}).
Let $t_0,\hat s>0$ and assume  that   $y\in  \V((\vec x,\vec \xi),t_0)\cap U_{\hat g}$ satisfies $y\not \in 
 {\mathcal Y}((\vec x,\vec \xi);t_0,\hat s)\cup\bigcup_{j=1}^4
 \gamma_{x_j,\xi_j}(\R) $. 
Then

(i) If $y$ does not satisfy condition (I)  with  $(\vec x,\vec \xi)$ and  $t_0$ then 
$y$ does not satisfy
the condition ({D}) with  $(\vec x,\vec \xi)$ and  $t_0,\hat s>0$.

(ii) Assume  
$y\in \hat U$ satisfies the condition (I) with  $(\vec x,\vec \xi)$ and  $t_0$ and parameters $q,\zeta$, and  $0<t<\rho(q,\zeta)$.
Then $y$ satisfies condition ({D}) with  $(\vec x,\vec \xi)$, $t_0$, and  any sufficiently
small $\hat s>0$.
   \end{lemma}

   \noindent

   {\bf Proof.}  
      (i)   
   If $y\not \in 
{\mathcal Y}((\vec x,\vec \xi);t_0,\hat s)\cup\bigcup_{j=1}^4
 \gamma_{x_j,\xi_j}(\R) $, the same condition holds also for $(\vec x^\prime,\vec \xi^\prime)$ close 
 to ($ \vec x,\vec \xi)$.
   Thus
Prop.\ \ref{lem:analytic limits A}, Prop.\ \ref{lem:analytic limits C},  
and Lemma \ref{lem: sing detection in wave gauge coordinates 1} imply that (i) holds.

(ii) Let  $\mu_{\vec \e}([-1,1])$ be  a family of observation geodesics
such that $y=\mu_{0}(\tilde r)$ and $\Psi_{\vec \e}$ be the normal coordinates associated 
to $\mu_{\vec \e}([-1,1])$ and $y$.

 Let $R_1>0$ be such that 
$B_{\hat g^+}(y,R_1)\subset\V((\vec x,\vec \xi),t_0)\cap \hat U$, see  (\ref{eq: summary of assumptions 2}).
  Our aim is next to show that there is a source ${\bf f}_{\vec \e}$
so that 
$\p_{\vec \e}^4((\Psi_{\vec \e})_*u_{\vec \e})|_{\vec \e=0}$ is not
$C^\infty$-smooth in $B_{\hat g^+}(y,R_1)$ for any $R_1>0$.
By making $R_1$ smaller if necessary, we see that when 
$\hat s>0$ is small enough, we have that 
 if $p$ is a cut point  of $\gamma_{x_j,\xi_j}([t_0,\infty))$ for some $j\leq 4$,
then $B_{\hat g^+}(y,R_1)\cap J^+_{\hat g}(p)=\emptyset$
and $B_{\hat g^+}(y,R_1)\cap {\mathcal Y}((\vec x,\vec \xi);t_0,\hat s)\cup\bigcup_{j=1}^4
 \gamma_{x_j,\xi_j}(\R) =\emptyset$.

Let $\eta=\p_t \hat \gamma_{q,\zeta}(t)$ and
 denote
$(y,\eta)=(x_5,\xi_5)$.
Let $t_j\in \R$ be such that  $\hat \gamma_{x_j,\xi_j}(t_j)=q$
and denote $b_j=\p_t \hat\gamma_{x_j,\xi_j}(t_j)$, $j=1,2,3,4,5$. Note that then $b_5=\zeta$, $x_5=y$,
and $t_5=t$.

By Propositions\ \ref{lem:analytic limits A} and  \ref{singularities in Minkowski space},  arbitrarily near to $b_j\in L^+_q\hattuM _0$ there are  $b_j^{\prime}\in L^+_q\hattuM _0$  
 and polarizations $v^{(j)}\in \R^{10+L}$, $j=2,3,4$, i.e., principal symbols at $(q,b_j^{\prime})$,
and linearly independent polarizations $v^{(5)}_p\in \R^{10+L}$, $p=1,2,3,4,5,6,$
and $v^{(1)}_r\in \R^{10+L}$, $r=1,2,3,4,5,6$,
having the following properties:
\medskip

(a) All $j\leq 4$, $v^{(j)}$, $j=2,3,4$, and $v^{(1)}_r$, $r=1,2,\dots,6$  are such that
their metric components (i.e., $g$-components of $(g,\phi)$) satisfy the divergence 
conditions for the symbols  (\ref{divergence condition for symbol})
with the covector $\xi$ being $b_j^{\prime}$ and $b_1^{\prime}$,  respectively.
\smallskip

(b)
If
 $v^{(5)}\in X_5=\hbox{span}(\{v^{(5)}_p;\ p=1,2,3,\dots,6\})\setminus\{0\}$ then 
 there exists
 a vector $v^{(1)}\in X_1=\hbox{span}(\{v^{(1)}_r;\ r=1,2,3,\dots,6\})$
such that for  ${\bf v}=(v^{(1)},v^{(2)},v^{(3)},v^{(4)},v^{(5)})$
 and ${\bf b}^{\prime}=(b_j^{\prime})_{j=1}^5$
 we have
 $\mathcal G({\bf v},{\bf b}^{\prime})\not =0$.
 \medskip

  Let $s\in (0,\hat s)$ and 
$x^{\prime}_j= \hat\gamma_{q,b^{\prime}_j}(-t_j)$ and
$\xi^{\prime}_j=\p_t \hat\gamma_{q,b^{\prime}_j}(-t_j)$, $j=1,2,\dots,4$,
and $\xi^{\prime}_5=\p_t \hat\gamma_{q,b^{\prime}_5}(t_5)$.
As the function $\rho$ is lower semi-continuous, we can assume that the $b_j^{\prime}$ are above chosen to be so  close to $b_j$ that $t_j>\rho(q,b_j^\prime)$,
 $(x_j^{\prime},\xi_j^{\prime})\in \W_j(s;x_j,\xi_j)$
 and  $x_5^\prime \in \V((\vec x,\vec \xi),t_0)\cap U_{\hat g}$, see   (\ref{eq: summary of assumptions 2}).
 We denote $(\vec x^{\prime},\vec \xi^{\prime})=((x_j^{\prime},\xi_j^{\prime}))_{j=1}^4$.

  As  $\rho(q,\zeta)>t$, and the function $\rho$ is lower semi-continuous,
  we can also assume that $b_5^\prime$ is so close to $\zeta$
  that $\rho(q,b_5^\prime)>t$ and $x_5^\prime  \in B_{\hat g^+}(y,R_1)$.
  By
assuming that $V\subset B_{\hat g^+}(y,R_1)$ is a sufficiently small neighborhood of $x_5^\prime$,
we have that $S:=\L_{\hat g}^+(q)\cap V$ is a smooth 3-submanifold.

 Next, consider the parametrix ${\bf {\bf Q}}_{\hat g}^*$ corresponding
 to  the linear wave equation with a reversed causality and
 the gaussian beam $u_\tau={\bf Q}_{\hat g}^* F_\tau$, produced by a source $F_\tau$ and function $h$
 given in (\ref{Ftau source}).
When $h(x_5^{\prime})=\overline w$ and  $w$ 
 is  the principal symbol (i.e.\  the polarization) of $u_\tau$ at $(q,b_5^{\prime})$,
   we can use the techniques of \cite{KKL,Ralston},
   see also \cite{Babich,Katchalov} \HOX{20.4.2013: Here was a change 
to citation to \cite{Ralston}  showing that the propagation
of polarization for gaussian beam.}
to obtain an analogous result to Lemma
\ref{lem: lagrangian} for the propagation of singularities  along the geodesic
$\hat \gamma_{q,b_5^{\prime}}([0,-t_5])$, and see that
$w=(R_{(5)})^*\overline w$, where  $R_{(5)}$ is 
a bijective linear map similar to map $R_{(5)}(q,b_5^{\prime};x_5^{\prime},\xi_5^{\prime})$ considered in the formula (\ref{eq: R propagation}), that is obtained by
solving a system of linear ordinary differential equations along the geodesic
connecting $q$ to $x_5^\prime$.  \observation{We note
that a formula similar to (\ref{new eq: R propagation}) holds for 
map $R_{(5)}(q,b_5^{\prime};x_5^{\prime},\xi_5^{\prime})$ in the case
when we consider the Einstein equation  (\ref{eq: adaptive model}).}
 Let us also denote
$R_{(1)}=R({q,b_1^{\prime}};x_1^{\prime},\xi_1^{\prime})$
and recall that the map $R_{(1)}$  is bijective, too.

%
%
Consider  next ${\bf b}^{\prime}$ and $w^{(2)},w^{(3)},w^{(4)}$ as parameters and let
\ba
\W_5&=&\{(R_{(1)})^{-1}w^{(1)},(R_{(5)}^*)^{-1}w^{(5)})\ ;\ w^{(5)}\in X_5,\ w^{(1)}\in X_1, 
  \\
& &\quad \hbox {and $\mathcal G({\bf w},{\bf b}^{\prime})\not =0$ for }{\bf w}=(w^{(1)},w^{(2)},w^{(3)},w^{(4)},w^{(5)})\}.
\ea
Let $W_5=\pi_2 (\W_5)\cup \{0\}$, where $\pi_2:(\overline w_1,\overline w_5)\mapsto \overline w_1$.
When $\overline w^{(5)}\in W_5$ is non-zero,
there is $\overline w^{(1)}\in X_1$
such that $(\overline w^{(1)},\overline w^{(5)})\in \W_5$.
Then, by  Lemma \ref{lem: lagrangian}, there 
are ${\bf f}_j\in {\mathcal I_{S}}^{n-3/2}(Y((x_j^{\prime},\xi_j^{\prime});t_0,s_0))$
that have the principal symbols
$\overline w^{(j)}=R({q,b_j^{\prime}};x_j^{\prime},\xi_j^{\prime})^{-1}w^{(j)}$
at $(x_j^{\prime},\xi_j^{\prime})$, $j\leq 4$.
Moreover, let $u_{\vec \e}=(g^{\vec \e},\phi^{\vec \e})$ be the solution
corresponding to 
${\bf f}_{\vec \e}=\sum_{j=1}^4 \e_j {\bf f}_j $
 and
$u_\tau={\bf Q}_{\hat g}^*F_\tau$ 
is
a past propagating gaussian beam sent in the direction $(x_5^{\prime},-\xi_5^{\prime})$,
defined in (\ref{Ftau source}) with functions  $F_\tau$ and  $h(x)$ such that $h(x_5^\prime)=\overline w^{(5)}$.
%
%
Then  we see for $\M^{(4)}=\p_{\vec \e}^4u_{\vec \e}|_{\vec \e=0}$  that
 the inner product $\bra F_\tau,\M^{(4)}\cet_{L^2(\hat U)}$ is not of the order $O(\tau^{-N})$
for all $N>0$.  Next, let us continue $\overline w^{(5)}$ from the point $x_5^\prime$ 
to a smooth section of $\B^L$.  Then
the above implies that the function $x\mapsto \bra \overline w^{(5)}(x),\M^{(4)}(x)\cet_{\B^L}$ is
not smooth in any neighborhood of $x_5^\prime$.

Roughly speaking, this means that in the wave gauge
coordinates  
$\M^{(4)}$ has wave front set
at $(x_5^{\prime},\xi_5^{\prime})$ with a polarization  that is
not perpendicular to $\overline w^{(5)}$.

Let us consider a family of   observation geodesics $\mu_{\vec \e}^\prime([0,1])$
such that
$p_{\vec\e}=\mu_{\vec \e}^\prime(\tilde r)$ satisfies $p_0=x_5^\prime$
and that for some $r_-,r_+\in (-1,1)$  satisfying $r_-<\tilde r<r_+$,  the curve $\hat \gamma(r)=\mu_0(r)$, $r\in [r_-,r_+]$, 
is a causal geodesic that intersects $S$ only at $x_5^\prime$, the intersection
is  transversal, and
 $\hat \gamma\subset V$. Let $\Psi_{\vec \e}$ be the normal coordinates 
 at $p_{\vec\e}$ associated to the observation geodesics $\mu_{\vec \e}^\prime([-1,1])$
centered at the point $\mu_{\vec \e}^\prime(\tilde r).$

Let $\X=\hbox{symm}(T_{x_5^{\prime}}\hattuM _0\otimes T_{x_5^{\prime}}\hattuM _0)+\R^L$
be the linear space of dimension $(10+L)$.
By the property (b) above,  $W_5\subset \X$ is a linear subspace containing  $X_5$
that has  dimension  $6$ and the
dimension of $V_5=\hbox{symm}(T_{x_5^{\prime}}S\otimes T_{x_5^{\prime}}S)+\R^L
\subset X$
is $(6+L)$.  Thus the dimension of $\V_5\cap W_5$
is at least two. In particular, it contains a non-zero vector.

Next, 
let  $\overline w^{(5)}\in (W_5\cap V_5)\setminus\{0\}$ and $\overline w^{(1)}\in (R_{(1)})^{-1}X_1$ be such
that $(\overline w^{(1)},\overline w^{(5)})\in \W_5$.
By Lemma \ref{lem: sing detection in wave gauge coordinates 1},
we then see that
$\p_{\vec \e}^4((\Psi_{\vec \e})_*u_{\vec \e})|_{\vec \e=0}$
is not $C^\infty$-smooth in any neighborhood of $\Psi_0(x_5^{\prime})$.
This means that there are 
sources  ${\bf f}_j\in {\mathcal I_{S}}^{n-3/2}(Y((x_j^{\prime},\xi_j^{\prime});t_0,s_0))$ with
polarizations $(\overline w^{(1)},\overline w^{(2)},\overline w^{(3)},\overline w^{(4)})$
that cause  singularities in $\p_{\vec \e}^4((\Psi_{\vec \e})_*u_{\vec \e})|_{\vec \e=0}$ near $
\Psi_0(x_5^{\prime})$, that is, the singularities that one can observe 
in the normal 
coordinates. As $R_1$ is arbitrarily small so that $x_5^{\prime}$ can be assumed
to be in an arbitrary neighborhood of $y$, we see using the above and the
normal coordinates associated a family of observation geodesics  $\mu_{\vec \e}([-1,1])$ and  $y=\mu_0(\tilde r)$
 that condition 
({D}) is valid for $y$. This proves (ii).
 \hfill \Box \medskip

\observation{

{\bf Remark 4.1.}
  When we consider the Einstein equation  (\ref{eq: adaptive model}),
we can use in
the above proof the fact that
  the formula  (\ref{new eq: R propagation}) holds for the
  maps $R_{(1)}=R(q,b_1^{\prime};x_1^{\prime},\xi_1^{\prime})$
 and $R_{(5)}(q,b_5^{\prime};x_5^{\prime},\xi_5^{\prime})$.
Then, we  can assume that the principal symbols of the $\phi$-components of the sources
${\bf f}_j\in {\mathcal I_{S}}^{n-3/2}(Y((x_j^{\prime},\xi_j^{\prime});t_0,s_0))$
vanish, that is, the leading order singularities of the sources 
are in the $g$-component. This is due to the fact that then
the principal symbols of the $\phi$-components the waves $u_\tau$ and $u_j$
vanish at the point $q$ and thus our earlier
considerations in the proof
of Proposition \ref{singularities in Minkowski space} on the 4:th order interaction of the waves are valid.  Moreover, this implies that in the condition (D) we
can require that the $\phi$-components of the sources
${\bf f}_j\in {\mathcal I_{S}}^{n-3/2}(Y((x_j^{\prime},\xi_j^{\prime});t_0,s_0))$
vanish identically and that we observe only the $g$-component
of the wave $\mathcal M^{(4)}$ 
and still the claim of Lemma \ref{lem: sing detection in in normal coordinates 2}
will be valid.\medskip
}


Above
we considered the  solution $u_{\vec \e}$ and a source ${\bf f}_{\vec \e}$ satisyfing the
conditions (\ref{eq: source causality condition})  in  the wave guide coordinates.
However, we do not know the wave guide coordinates in the set $(U_g,g)$ and thus
thus we do not know which element of the 
data set ${\cal D}(\hat g,\hat \phi,\e)$, see (\ref{eq: main data}), corresponds to the source 
 ${\bf f}_{\vec \e}$. This problem is solved in the following lemma.

\begin{figure}[htbp] \label{Fig-14d}
\begin{center}

\psfrag{1}{$\hat U$}
\psfrag{2}{$\Sigma$}
\psfrag{3}{$J_{\hat g}(\hat p^-,\hat p^+)$}
\includegraphics[width=5.5cm]{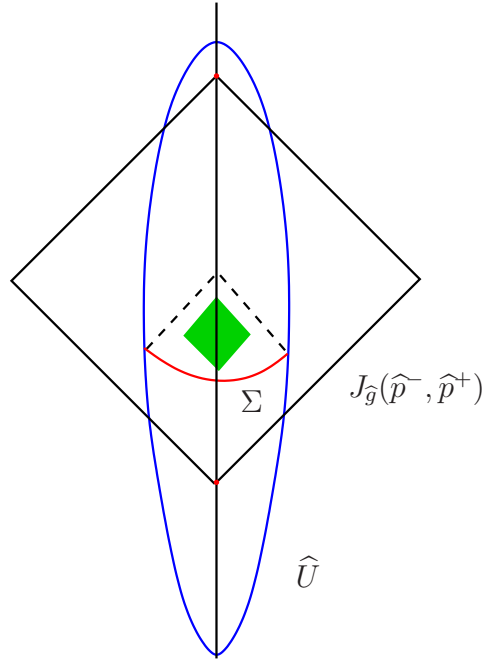}
\end{center}
\caption { A schematic figure where the space-time is represented as the  3-dimensional set $\R^{1+1}$. 
The figure shows the objects used in the proof of Lemma \ref{lem: sources in wave guide}.
Figure displays the subsets of $(M_0,\hat g)$. The green diamond in the figure is the set $V$ where the source is supported. The dashed line shows the set from which all causal curves intersect the 
surface $\Sigma\subset \hat U$. 
}
 \end{figure}


\begin{lemma}\label{lem: sources in wave guide}


%

Assume that we are  given ${\cal D}(\hat g,\hat \phi,\e)$ where  $\e>0$ is small enough.
Let $0<r_3<r_2<r_1$, where $r_1$ is the parameter used to define $W_g$ and $s_-+r_1<s_1<s_+$.
When $\e_1,\e_2>0$ are small enough the following holds:

(i)  Assume that we are  given 
$(U^\prime,g^\prime,\phi^\prime,F^\prime)$  such that the equivalence
class $[(U^\prime,g^\prime,\phi^\prime,F^\prime)]$ is in 
${\cal D}(\hat g,\hat \phi,\e)$, and moreover, we have 
$\mathcal N^{(16)}_{\hat g}(g^\prime)<\e_1$,
$\mathcal N^{(16)}(F^\prime)<\e_1$,
and 
 \beq\label{eq support 1}
K^\prime\subset  I_{g^\prime}( \mu_{g^\prime}(s_1-r_2),\mu_{g^\prime}(s_1))
 \quad
 \hbox{where } K^\prime=\supp(F^\prime).
 \eeq 
 
Then we can determine a source $F\in C^\infty_0(W_{\hat g})$  such  that
$(U^\prime,g^\prime,\phi^\prime,F^\prime)\in [(U_g,g|_{U_g},\phi|_{U_g},F|_{U_g})]$,
where $(g,\phi)$ is the solution of 
the $\hat g$-reduced
Einstein equations  (\ref{eq: adaptive model})
with the source  $F$.

(ii) Let  $K\subset \hat U$ be a compact set such that
\beq\label{eq support 1 B}
K\subset  I_{\hat g}( \mu_{\hat g}(s_1-r_3),\mu_{\hat g}(s_1)).
\eeq
 When
  $F\in C^\infty_0(W_{\hat g})$ satisfies 
$\supp(F)\subset K$ and $\mathcal N^{(16)}(F)<\e_2$, 
we can find the element
 $[(U_g,g|_{U_g},\phi|_{U_g},F|_{U_g})]$ in  ${\cal D}(\hat g,\hat \phi,\e)$,
 where $(g,\phi)$ is the solution of 
the $\hat g$-reduced
Einstein equations  (\ref{eq: adaptive model})
with the source  $F$.
Moreover, we can find $(\Psi_{\mu})_*g$ and $(\Psi_{\mu})_*\phi$
where $\Psi_{\mu}$ are normal coordinates associated
to a given  geodesic $\mu=\mu_{g,z,\eta}([-1,1])\subset U_g$,
that is, these normal coordinates are centered at the end
point of $\mu$ and are associated to the frame obtained
by parallel translation along  $\mu$.

\end{lemma}
     \noindent
   {\bf Proof.} 
   (i) 
   As $[(U^\prime,g^\prime,\phi^\prime,F^\prime)]\in
{\cal D}(\hat g,\hat \phi,\e)$, there exists 
a  source  $F$ on $M_0$ and
a  solution  $(g,\phi)$ of 
the $\hat g$-reduced
Einstein equations  (\ref{eq: adaptive model}) on $M_0$
with the source  $F$
such that
$(U_g,g|_{U_g},\phi|_{U_g},F|_{U_g})\in [(U^\prime,g^\prime,\phi^\prime,F^\prime)]$.

By definition, $(U_g,g|_{U_g},\phi|_{U_g},F|_{U_g})\in [(U^\prime,g^\prime,\phi^\prime,F^\prime)]$. This
implies that 
there exists a diffeomorphic  isometry \HOX{We replaced $I_{U^\prime,g^\prime}$ by $I_{g^\prime}$. Check if we need to add an explanation.}
$f:(U^\prime,g^\prime) \to (U_g, g)$. 
Let us next denote the causal domain in $(U^\prime,g^\prime)$ by $I^+_{g^\prime}(p^\prime)$ and
 $I^-_{g^\prime}(p^\prime)$  etc., where $p^\prime\in U^\prime$.
Let $V^\prime=I_{g^\prime}( \mu_{g^\prime}(s_1-r_2),\mu_{g^\prime}(s_1))$ 
and denote  $V=f(V^\prime)$, see Fig. 17. 
For clarity,  we denote the causal domains in $(M,g)$ by $I^+_{g}(p)$ and
 $I^-_{g}(p)$,  etc., where $p\in M$.

Let 
   $\Sigma^\prime \subset U^\prime$ be  a smooth space-like 3-dimensional surface such that
$J^+_{g^\prime}(K^\prime)\cap \Sigma^\prime=\emptyset$ and
for all $x\in  \hbox{cl}(V^\prime)$ any non-extendable  past-directed causal curve in $(U^\prime,g^\prime)$  starting from $x$
intersects $\Sigma^\prime$. Then the space-like surface $\Sigma=f(\Sigma^\prime)\subset \hat U$ is such that
$(J^-_{g^\prime}(\Sigma^\prime),g^\prime)$ is isometric to $(J^-_{U_g,g}(\Sigma),g)$
and also to  $(J^-_{\hat g}(\Sigma),\hat g)$. 
Next we  identify these sets as well $\Sigma^\prime$ and $\Sigma$. 
Then the restriction $f_1=f|_{J^-_{g^\prime}(V^\prime)\cap J^+_{g^\prime}(\Sigma)}$,
that is, the isometry
 $f_1:(J^-_{g^\prime}(V^\prime)\cap J^+_{g^\prime}(\Sigma),g^\prime)
\to (J^-_{g}( V)\cap J^+_{g}(\Sigma)$, g),
 is the identity map in a neighborhood of   $ \Sigma$.

Assume next that $\e_1$ and $r_2$ are  so small that
$f(J^-_{g^\prime}(V^\prime)\cap J^+_{g^\prime}(\Sigma))\subset 
 \hat U$.
 Considering
 the particular case when $(U^\prime,g^\prime,\phi^\prime,F^\prime)$  is equal 
to $(U_g,g, \phi, F)$ and $f_1$ is the identity map, we see 
using the harmonicity condition (\ref{Harmonicity condition}) that
 $f_1$ is also a $(g^\prime,\hat g)$-wave map  
  $f_1:J^-_{g^\prime}(V^\prime)\cap J^+_{g^\prime}(\Sigma)
\to \hat U$ 
 that is an identity map near  $ \Sigma$.  
Thus, as the given data contain the pairs $(U^\prime,g^\prime)$ and $(\hat U,\hat g)$ and we have
fixed above the surface $\Sigma$,
 we  can  determine $f_1$
by solving a Cauchy problem for the wave map equation. Hence,  we can find $F=(f_1)_*F^\prime$ in 
$f_1(J^-_{g^\prime}(V^\prime)\cap J^+_{g^\prime}(\Sigma))$. As $F$ vanishes outside this this set, we can find $F$ in $\hat U$ by extending the obtained
function by zero. This proves (i).

(ii) 
Assume that $\e_2\in (0,\e)$ is so small that the Einstein equation  (\ref{eq: adaptive model})
 have solution $(g,\phi)$ for  all $F\in \H$,
where 
\ba
\H=\{F\in C^\infty_0(W_{\hat g});\  
K=\supp(F)\hbox{ satisfies (\ref{eq support 1 B}) and }
\mathcal N^{(16)}_{\hat g}(F)<\e_2\}.\ea 
When $\e_2>0$ and $r_3>0$ are  small enough,
for any $F\in \H$ is such that the unique solution $(g,\phi)$ 
and $F$  satisfy $\mathcal N^{(16)}_{\hat g}(g)<\e_1$
and (\ref{eq support 1}). Then,
  there is a unique equivalence class $[(U_g,g|_{U_g},\phi|_{U_g},F|_{U_g})]\in {\cal D}(\hat g,\hat \phi,\e)$
and thus the map $\M:\H
\to 
{\cal D}(\hat g,\hat \phi,\e)$ is injective. 
Observe that for given $[(U_g,g|_{U_g},\phi|_{U_g},F|_{U_g})]\in {\cal D}(\hat g,\hat \phi,\e)$
we can verify if it belongs in $\H$. 
By (i) we can construct  the
inverse of the map $\M:\H\to \M(\H)$, and thus we can construct also the map $\M:\H\to \M(\H)$.
 Hence,
when $F\in \H$, we can determine $[(U_g,g|_{U_g},\phi|_{U_g},F|_{U_g})]$
 and represent $g$ and $\phi$ in the normal coordinates given in the claim. This yields (ii).
 \hfill \Box \medskip

By Lemma \ref{lem: sources in wave guide},
for the smooth sources ${\bf f}_{\vec \e}$ satisfying conditions (\ref{eq: source causality condition}) and $|\vec\e|$ small enough
we can find  $(\Psi_{\vec \e})_*u_{\vec \e}$ where $\Psi_{\vec \e}$ are the normal
coordinates associated to an observation geodesic on
$(M_0,g_{\vec \e})$, where $u_{\vec \e}=(g_{\vec \e},\phi_{\vec \e})$. 
The non-smooth sources $F$, given on $\hat U$ e.g.\ in the Fermi coordinates
of $(\hat U,\hat g)$,  for which  $\mathcal N^{(16)}(F)<\e$, can clearly been approximated by smooth sources.
Thus when  $(\vec x,\vec \xi)$ 
and the observation geodesic  $\mu_0([-1,1])$  are given as well as  the data set 
${\cal D}(\hat g,\hat \phi,\e)$,  
 for the sources ${\bf f}_{\vec \e}$ with 
${\bf f}_j\in {\mathcal I_{S}}^{n-3/2}(Y((x_j^{\prime},\xi_j^{\prime});t_0,s_0))$, where $(-n)$ is large enough and 
$(x_j^{\prime},\xi_j^{\prime})\in \W
(s;x_j,\xi_j)$,
we can compute the derivatives 
$\p_{\vec \e}^4((\Psi_{\vec \e})_*u_{\vec \e})|_{\vec \e=0}$.
Hence, using the  observation geodesics  $\mu_{\vec \e}([-1,1])$ and $\mu_0(\tilde r)=y$ we can  check
if the condition ({D}) is valid  for $y$  with the given 
 $(\vec x,\vec \xi)$, $s_0$ and $t_0$  or not.


\section{Determination of light observation sets for Einstein equations}\label{subsection combining}

In this section we use only the metric $\hat g$ and denote often $\hat g=g$, $\hat \gamma=\gamma$, $U=U_{\hat g}$, etc.  
Below, let $\hbox{cl}(A)$ denote the closure of the set $A$.

Our next  
aim is to handle the technical problem that  
in the set $ {\mathcal Y}((\vec x,\vec\xi))$, see (\ref{K and X sets}) we  
have not analyzed if we observe singularities or not. To/aim  
 this end,  
we define next the  
 sets $\S_{reg} ((\vec x,\vec \xi),t_0)$ of points near which we observe   
singularities in a 3-dimensional set.

%

\begin{definition}\label{def: Sreg}  
Let $(\vec x,\vec \xi)=((x_j,\xi_j))_{j=1}^4$ be a collection of light-like  
vectors with $x_j\in U_{g}$ and $t_0>0$. We define to  
 $\S((\vec x,\vec \xi),t_0)$   
be the set of those $y\in U_{g}$    
that satisfies   
the property  ({D}) with $(\vec x,\vec \xi)$ and $t_0$ and some $\hat s>0$.  
\HOX{Add citation to (D) and number it.}  
  
Moreover, let   
 $\S_{reg} ((\vec x,\vec \xi),t_0)$ be the set  
of the points  $y_0\in \S((\vec x,\vec \xi),t_0)$ having a neighborhood $W\subset U_{g}$  
such that the intersection $W\cap \S((\vec x,\vec \xi),t_0)$  
is a non-empty smooth 3-dimensional submanifold.  
We denote  
\ba  
& &{\Sclo}(\vec x,\vec \xi),t_0)=\hbox{cl}\,(\S_{reg} ((\vec x,\vec \xi),t_0))\cap U_g.  
\ea  
\end{definition}

Note that the  data  ${\cal D}(\hat g,\hat \phi,\e)$ determines the sets ${\Sclo}((\vec x,\vec \xi),t_0)$.  
Our goal is to show that ${\Sclo}((\vec x,\vec \xi),t_0)$ coincides with the intersection   
of the light cone $\L^+_{g}(q)$ and $U_{g}$ where $q$ is the intersection  
point of the geodesics corresponding to $(\vec x,\vec \xi)$, see Fig.\ 15 and 16.

\HOX{Change here $y$ to $\hat x$, as $y$ is the observation point.}  
Let us next motivate the analysis we do below:   
We will consider how to create an artificial point source using interaction  
of spherical waves propagating along light-like geodesics $\gamma_{x_j,\xxi_j}(\R_+)$  
where $(\vec x,\vec \xxi)=((x_j,\xxi_j))_{j=1}^4$ are perturbations  
of a light-like $(y,\zeta)$, $y\in \hat \mu$. 
We will use the fact that for all $q\in J(p^-,p^+)\setminus \hat \mu$ there  
is a light-like geodesic $\gamma_{y,\zeta}([0,t])$ from $y=\hat \mu(f^-_{\hat \mu}(q))$  
to $q$ with $t\leq \rho(y,\zeta)$.  
We will next  
show that   
when we choose $(x_j,\xxi_j)$ to be suitable   
perturbations of $\p_t \hat \gamma_{y,\zeta}(t_0)$, $t_0>0$,  
it is possible that all geodesics $\hat \gamma_{x_j,\xxi_j}(\R_+)$ intersect at $q$ before   
their first cut points,  
that is, $\hat \gamma_{x_j,\xxi_j}(r_j)=q$, $r_j<\rho(x_j,\xxi_j)$.   
We note that we can not analyze the interaction of the waves  
if the geodesics intersect after the cut points as then  
 the spherical  
waves can have caustics. Such interactions of wave caustics can, in principle,  
cause propagating singularities. Thus the sets ${\Sclo}((\vec x,\vec \xi),t_0)$  
contain singularities propagating along the light cone $\L^+_{g}(q)$  
and in addition that they may contain singularities produced by caustics that we do not know how to analyze (that could  
be called "messy waves"). Fortunately, near an open and dense set of geodesics  
$\mu_{z,\eta}$ the nice singularities propagating along the light cone $\L^+_{g}(q)$  
arrive before the "messy waves". This is the reason why we consider below the  
first observed singularities on geodesics  
$\mu_{z,\eta}$.  
 Let us now return to the rigorous analysis.  
 
Below in this section we fix $t_0$ to have the value  
 ${t_0}=4\kappa_1$, cf.\ Lemma \ref{lem: detect conjugate 0}. Recall the notation that  
\ba  
& &(x({t_0}),\xxi({t_0}))=(\gamma_{x,\xxi}({t_0}),\dot \gamma_{x,\xxi}({t_0})),\\  
& &(\vec x({t_0}),\vec \xxi({t_0}))=((x_j({t_0}),\xxi_j({t_0})))_{j=1}^4.  
\ea

\HOX{It seems that the claim "open set of directions" can be removed.}  
\begin{lemma} \label{lem: existence of vec x vec zeta}    
Let $\vartheta>0$ be arbitrary, $q\in J(p^-,p^+)\setminus {\hat \mu}$ and   
let $y={\hat \mu}(f_{\hat \mu}^-(q))$  
and $\zzeta\in L^+_y\hattuM _0$, $\|\zeta\|_{g^+}=1$ be such that $\gamma_{y,\zzeta}([0,r_1])$, $r_1>t_0=4\kappa_1$ is   
a longest causal (in fact, light-like) geodesic connecting $y$ to $q$.   
Then there exists  a set $\mathcal G$ of 4-tuples   of light-like vectors $(\vec x,\vec \xxi)=((x_j,\xxi_j))_{j=1}^4$   
 such that the points  $x_j(t_0)=\gamma_{x_j,\xxi_j}(t_0)$  
 and the directions $\xxi_j(t_0)=\dot \gamma_{x_j,\xxi_j}(t_0)$ have the following properties:  
\begin{itemize}  
\item[(i)] $x_j\in U_{g}$, $x_j\not \in J^+(x_k)$ for $j\not =k$,  
  
\item[(ii)] $d_{g^+}  
((x_l,\xxi_l),(y,\zeta))<\vartheta$ for $l \leq 4$,  
  
\item[(iii)] $q=\gamma_{x_j,\xxi_j}(r_j)$ and $\rho(x_j(t_0),\xxi_j(t_0))+t_0>r_j$,  
  
\item[(iv)] when $(\vec x,\vec \xxi)$ run through the set  $\mathcal G$,  
the directions  
 $ (\dot \gamma_{x_j,\xxi_j}(r_j))_{j=1}^4$ form an open  
 set in $(L^+_qM_0)^4$.  
 \end{itemize}  
   
In addition, $\mathcal G$ contains elements $(\vec x,\vec \xxi)$ for which $(x_1,\xxi_1)=  
(y,\zzeta)$.  
\end{lemma}  
  
{\bf Proof.}   
Let  $\eta=\dot \gamma_{y,\zzeta}(r_1)\in L^+_q\hattuM _0$.  
Let us choose light-like directions  $\eta_j\in T_q\hattuM _0$, $j=1,2,3,4,$   
close to $\eta$ so that $\eta_j$ and $\eta_k$ are not parallel for $j\not  =k$.  
In particular, it is possible (but not necessary) that   
$\eta_1=\eta$. Let ${\bf t}:M\to \R$ be a time-function on $M$ that can be  
used to identify $M$ and $\R\times N$. Moreover, let  us choose  
 $T\in ({\bf t}(y),{\bf t}(\gamma_{y,\zeta}(r_0-t_0)))$ and  for $j=1,2,3,4$, let $s_j>0$   
be such that   
${\bf t}(\gamma_{q,\eta_j}(-s_j))=T$. Choosing  first $T$ to be sufficiently close ${\bf t}(y)$  
and then all $\eta_j$, $j=1,2,3,4$ to be sufficiently  
close to $\eta$     
 and defining  
$x_j=\gamma_{q,\eta_j}(-s_j)$ and $\xxi_j=\dot\gamma_{q,\eta_j}(-s_j)$  
we obtain the pairs $(x_j,\xxi_j)$ satisfying the properties stated in the claim.  
As vectors $\eta_j$ can be varied in sufficiently small open sets   
so that the properties stated in claim stay valid, we obtain   
that  claim  concerning the open set of light-like directions.  
  
The last claim follows from the fact that $\eta_1$ may be equal to $\eta$ and $T={\bf t}(y)$.  
 \hfill \Box \medskip

Next we  analyze the   
set ${\Sclo}((\vec x,\vec \xi),t_0)=\hbox{cl}\,(\S_{reg} ((\vec x,\vec \xi),t_0))$.  
We recall that the set $\be_U(q)$ is the points on $\mu_{z,\eta}$  
on which the light cone $\L_g^+(q)$ is observed at the earliest time, see Def.\ \ref{def: earliest element}.

\begin{lemma}\label{lem: sing detection in in normal coordinates 3}  
Let $(\vec x,\vec \xi)$,   and ${\bf t}_j=\rho(x_j(t_0),\xxi_j(t_0))$ with $j=1,2,3,4$, $t_0=4\kappa_1$, and $x_6\in U_{g}$  
and   
satisfy (\ref{eq: summary of assumptions 1})-(\ref{eq: summary of assumptions 2})  
and assume that $\vartheta_0$ in  (\ref{eq: summary of assumptions 1}) and   
Lemma \ref{lem: detect conjugate 0}   is so small  
that  for all $j\leq 4$,   
$x_j\in  I_{g}( \mu_{g}(s-r_1),\mu_{ g}(s))=  
 I_{g}^+( \mu_{g}(s-r_1))\cap I_{g}^-(\mu_{ g}(s))$ with  some $s\in (s_-+r_1,s_+)$,  
where $r_1$ appears in Lemma \ref{lem: sources in wave guide}.  
Let $ \V((\vec x,\vec \xi),t_0)$ be the set defined in  (\ref{eq: summary of assumptions 2})  
and consider $y\in \V((\vec x,\vec \xi),t_0)\cap U_{g}$.  
%
%
 Then \HOX{Add formula number and cite to (I).}  
\smallskip  
  
(i) Recall that if $y=\gamma_{q,\zeta}(t)$ with $t\leq \rho(q,\zeta)$, then $y\in \be_U(q)$.  
Assume that $y$ satisfies the condition (I) with $(\vec x,\vec \xi)$ and $t_0$ with parameters   
$q,\zeta$, and $t$ such that $0\leq t\leq \rho(q,\zeta)$,  
that is, $y\in \be_U(q)$.  
 Then   
 $y\in {\Sclo}((\vec x,\vec \xi),t_0)$.    
\smallskip   
  
(ii) Assume  $y$ does not satisfy condition (I) with $(\vec x,\vec \xi)$ and $t_0$.  
Then  $y$ has a neighborhood  
that does not intersect $\S_{reg}((\vec x,\vec \xi),t_0)$.   
In particular, $y\not \in {\Sclo}((\vec x,\vec \xi),t_0)$.   
\smallskip   
  
(iii) If ${\Sclo}((\vec x,\vec \xi),t_0)\cap \V((\vec x,\vec \xi),t_0)\not =\emptyset$ then   
the geodesics corresponding to $(\vec x,\vec \xi)=((x_j,\xi_j))_{j=1}^4$  
intersect (see Def.\ \ref{def:  4-intersection of rays}) and there is  
a unique point $q\in  \V((\vec x,\vec \xi),t_0)$ where the   
 the intersection takes place.   \end{lemma}  
    
  The point $q$ considered in Lemma \ref{lem: sing detection in in normal coordinates 3} (iii) is the earliest   
point in the set  
 $\cap_{j=1}^4 \gamma_{x_j,\xi_j}([t_0,\infty))$ and we denote  it by  
$q=Q((\vec x,\vec \xi),t_0)$.  
If such intersection point does not exists, we define $Q((\vec x,\vec \xi),t_0)$  
to be the empty set.  
    
   \noindent  
   {\bf Proof.}    
   (i) Clearly, if $t=0$ so that $q\in U_{g}$, we have $q \in  {\Sclo}((\vec x,\vec \xi),t_0)$.

   Assume next that $0<t< \rho(q,\zeta)$.  
    Let $r_2>0$ and $V(y,r_2)$ be the set of points $y^\prime\in B_{g^+}(y,r_2)$   
 such that there exists a family of observation geodesics  $\mu_{\vec \e}([-1,1])$  
and $y^\prime=\mu_{\vec \e}(\tilde r)$  
 satisfying condition ({D}).  
Let   
$\Gamma=\exp_q^{g}( \X((\vec x,\vec \xi),t_0)$,  
see (\ref{K and X sets}),  
be the set of points on the light cone on which singularities  
caused by 3-interactions may appear. Then the closure of $ \L^+_{g}(q)\setminus \Gamma$ is $ \L^+_{g}(q)$.  
 As $t<\rho(q,\zeta)$, we see that    
 when $r_2>0$ is small enough,  
   $B_{g^+}(y,r_2)\cap \L^+_{g}(q)$ is a smooth  3-submanifold  
  having non-empty intersection with any neighborhood of $y$ and  
  $ B_{g^+}(y,r_2)\cap \Gamma$ is a submanifold of dimension 2.  
   Moreover,  we see using  Lemma  \ref{lem: sing detection in in normal coordinates 2}  
   that for   
 $r_2>0$ small enough that  
 $ B_{g^+}(y,r_2)\cap  (\L^+_{g}(q)\setminus  
  \Gamma)\subset V(y,r_2)$ and  
  $V(y,r_2)\subset \L^+_{g}(q)$.  
    Thus $(\L^+_{g}(q)\setminus \Gamma)\cap U_{g} \subset S_{reg}\subset \L^+_{g}(q)\cap U_{g}$, where the complement of $\Gamma$ in  $\L^+_{g}(q)$ is  dense in $\L^+_{g}(q)$.  
  The claim (i) follows from this in the case when $t< \rho(q,\zeta)$.  
    
  Assume next that   
  $y=\gamma_{q,\zeta}(t)$ where $\zeta$  
is light-like vector and $t=\rho(q,\zeta)$.  
Let $t_j<t$ be such that $t_j\to t$ as $j\to \infty$. Then  
 the $t_j<\rho(q,\zeta)$  
and the above    
 yields  
 $y_j=\gamma_{q,\zeta}(t_j)\in {\Sclo}((\vec x,\vec \xi),t_0)$. As $y_j\to y$  
 as $j\to \infty$ we  
 see that $y\in {\Sclo}((\vec x,\vec \xi),t_0)$. This proves (i).  
    
  \HOX{Slava said that (ii) follows directly from   
   Lemma  \ref{lem: sing detection in in normal coordinates 2}  
so in the final version of the paper the proof can be omitted.}  
     
  (ii) It follows from the assumption that either there are no intersection point $q$  
   or   
  that $y\not \in \L^+_{g}(q)$.   
  If $y\not \in A=\Gamma \cup \bigcup_{j=1}^4  
 \gamma_{x_j,\xi_j}([t_0,\infty)) $,  
 it follows from Lemma  \ref{lem: sing detection in in normal coordinates 2}  
 that  $y$ has a neighborhood $V$ that  does not intersect the set  
 $S ((\vec x,\vec \xi),t_0)$.  
If  $y\in A$,  
 we see, using the fact that $y$ does not satisfy condition (I) with $(\vec x,\vec \xi)$ and $t_0$, that $y$  has a neighborhood $V$   
  such that $V\setminus A$ does not intersect   
 $S((\vec x,\vec \xi),t_0)$. The Hausdorff dimension of the set $A\cap V$ is  less  
 or equal to  
2. Thus $V\cap S_{reg}((\vec x,\vec \xi),t_0)=\emptyset $. This proves (ii).

(iii) Using the conditions posed for  $(\vec x,\vec \xi)$,  $t_0$ and ${\bf t}_j$ with $j=1,2,3,4,$ and $x_6$  
we see that geodesics $\gamma_{x_j,\xi_j}([t_0,\infty)$ can intersect only  
once in $\V((\vec x,\vec \xi),t_0) $. Moreover, if there is no such intersection, the  
above shows that  ${\Sclo}((\vec x,\vec \xi),t_0)\cap \V((\vec x,\vec \xi),t_0) =\emptyset$, which proves the claim (iii).  
\hfill \Box \medskip

{Recall that ${\cal D}(\hat g,\hat \phi,\e)$  
 determines  the sets  
${\Sclo}((\vec x,\vec \xi),t_0)$. 
Below, denote  
\ba  
{\Sclo}_{z,\eta}((\vec x,\vec \xi),t_0)=\pointear_{z,\eta}({\Sclo}((\vec x,\vec \xi),t_0)),  
\ea  
(see (\ref{Earliest element sets}) and Def.\ \ref{def: Sreg}).  

Next we summarize the above results in a form that is convenient below.

\begin{lemma}\label{lem: finding earliest observ. in J(x6)}   
Let $(\vec x,\vec \xi)$,   and ${\bf t}_j=\rho(x_j(t_0),\xxi_j(t_0))$ with $j=1,2,3,4$, $t_0=4\kappa_1$, and $x_6\in U_{g}$  
and   
conditions (\ref{eq: summary of assumptions 1})-(\ref{eq: summary of assumptions 2}) are some satisfied  
with  some sufficiently small $\vartheta_0$.

 Let $(z,\eta)\in \U_{z_0,\eta_0}$. Then we have  
\smallskip  
  
\noindent  
(i) If $Q((\vec x,\vec \xi),t_0)=\emptyset$ or $Q((\vec x,\vec \xi),t_0)\not \in J^-(x_6)$  
then   
\beq\label{beU formula emptyset}  
{\Sclo}_{z,\eta}((\vec x,\vec \xi),t_0)\cap J^-(x_6)=\emptyset,  
\eeq  
  
\smallskip  
  
\noindent  
(ii) If $Q((\vec x,\vec \xi),t_0)\not =\emptyset$ and $q=Q((\vec x,\vec \xi),t_0)\in J^-(x_6)$  
then   
\beq\label{beU formula1}  
\quad{\Sclo}_{z,\eta}((\vec x,\vec \xi),t_0)\cap J^-(x_6)=  
\left \{ \begin{array}{cl}  
\{y_{z,\eta}\},&\hbox{if }y_{z,\eta}\in  J^-(x_6),\\  
\emptyset ,&\hbox{if }y_{z,\eta}\not \in  J^-(x_6),\end{array}\right.\hspace{-1cm}  
\eeq  
where $y_{z,\eta}=\mu_{z,\eta}(f^+_{\mu(z,\eta)}(q))$, i.e.,  
$\be_{z,\eta}(q)=\{y_{z,\eta}\}$.  
\end{lemma}

Note that this lemma combined with Theorem \ref{main thm} prove Theorem  
\ref{main thm Einstein} in the case when there are  
no cut points on the manifold $(M,g)$. Later we consider general manifolds with cut points.  
  
{\bf Proof.}    
 The claim (i) follows by applying  Lemma \ref{lem: sing detection in in normal coordinates 3} (ii) for all points $y\in J^-(x_6)\cap U_{g}$.  
  
Let us next consider the claim (ii). Denote $\V=\V((\vec x,\vec \xi),t_0)$.  
Let $q=Q((\vec x,\vec \xi),t_0)$. Then $q\in J^-(x_6)\subset \V$.   
  
Let $y\in U_{g}  
\cap \V$.  
If    
$y\not \in \L_{g}(q)$, Lemma \ref{lem: sing detection in in normal coordinates 3} yields  
  $y \not \in {\Sclo}((\vec x,\vec \xi),t_0)$.  
On the other hand, if $y\in \be_U(q)$, Lemma \ref{lem: sing detection in in normal coordinates 3} yields  
 that $y\in {\Sclo}((\vec x,\vec \xi),t_0)$.  
 Thus  $\be_U(q)\cap \V\subset {\Sclo}((\vec x,\vec \xi),t_0)\subset \L_{g}(q)$ yielding  
 \beq\label{beU formula2}  
 \be_{z,\eta}(q)\cap \V =\pointear_{z,\eta}({\Sclo}((\vec x,\vec \xi),t_0))\cap \V.  
\eeq  
If $y_{z,\eta}\in J^-(x_6)$ then $ J^-(x_6)\subset \V$ implies that  
 $\be_{z,\eta}(q)\cap \V=\{y_{z,\eta}\}$.  
 On the other hand, if   $y_{z,\eta}\not \in J^-(x_6)$, then   
 $ \pointear_{z,\eta}(F)\cap  J^-(x_6)= \be_{z,\eta}(q)\cap J^-(x_6)=\emptyset.$  
This yields (ii).  
    \hfill \Box \medskip  
  
 Below, let $t_0=4\kappa_1$, cf.\ Lemma \ref{lem: detect conjugate 0}  and  
 $\K_{t_0}\subset \hat U$ be the set of  points  
 $x=\gamma_{y,\zeta}(r)$ where $y=\hat \mu(s)$, $s\in [s^-,s^+]$,  
 $\zeta\in L^+_yM_0$ satisfies $g^+(\zeta,\zeta)=1$ and $r\in [0,2{t_0})$.  
 Recall that $\U_{z_0,\eta_0}=\U_{z_0,\eta_0}(\hat h)$  
 was defined using the parameter $\hat h$. We see that if $t_0^\prime >0$  
 and  $\hat h^\prime\in (0, \hat h)$ are small enough and $(z,\eta)\in \U_{z_0,\eta_0}(\hat h^\prime)$,  
 then the \HOX{The explanation that  $\be_U(\K_{t_0})$ is assumed to be known should be improved}  
 longest geodesic from $x\in \overline \K_{t_0^\prime}$ to the point  
 $\be_{z,\eta}(x)$ is contained in $\U_g$,  
 and hence we can determine the point $\be_{z,\eta}(x)$ for such $x$ and $(z,\eta)$.  
 Let us replace the parameters $\hat h$ and $t_0$ by $\hat h^\prime$  
 and $t_0^\prime$, correspondingly in our considerations below. Then we  
 may  assume  
 that in addition to the data given in the original formulation of the problem,  
 we are given also the set $\be_U(\K_{t_0})$. Next we do this.  
   
For technical reasons, we will next replace $U$ by   
$V=U\cap I^-(\hat p^+)$ and consider the sets   
$\be_V(q)$, $q\in J^+(\hat p^-)\cap I^-(\hat p^+)$.  
   
  \HOX{Simplification for final paper: Returning geodesics can be seen using linearized waves}  
Next we consider step-by-step construction of the set $\be_V(J^+(\hat p^-)\cap I^-(\hat p^+))$ by constructing  
 $\be_V(J^+(y)\cap I^-(\hat p^+))$ with $y={\hat \mu}(s)$ and decreasing $s$ in small steps.   
{The difficulty we encounter here is that we do not know how the spherical waves propagating  
along geodesics   
interact after the geodesics have a cut point.   

 Our aim is construct sets ${\mathcal Z}_k=\be_V(I^-(\hat p^+)\cap J^+(y_k))\setminus   
\be_V(I^-(\hat p^+)\cap J^+(y_{k-1}))$, with ${\mathcal Z}_0=\emptyset$ and $y_k={\hat \mu}(s_k)$,  
 $s_{k}<s_{k-1}$, $s_0,\dots,s_K\in [s_-,s_+]$ with $s_K=s_-$ and $s_0$   
 being  close to $s_+$.  
 The union  
 of ${\mathcal Z}_k,$ $k=1,2,\dots, K$ is the set  
 $\be_V(J^+(\hat p^-)\cap I^-(\hat p^+))$. The idea of this construction is to choose $s_j$  
so that when ${\mathcal Z}_{k-1}$ is constructed we obtain the sets  
 ${\mathcal Z}_k$ as a union of  sets  
  $\be_V(\gamma_{y,\zzeta}([0,r(y,\zzeta))))$ where $y\in {\hat \mu}([s_k,s_{k-1}))$, $\zzeta$ is light-like,  
  and  
 $r(y,\zzeta)$ is chosen so that the geodesic $\gamma_{y,\zzeta}([{t_0},r(y,\zzeta)))$  
 does not contain cut points and does not intersect $J^+(y_{k-1})$ but  
 still lies in $I^-(\hat p^+)$.}

The construction  is the following:  
\medskip  
  
Let below $\kappa_1,\kappa_2,\kappa_3$ be constants  given in Lemma \ref{lem: detect conjugate 0}.  
{Let  $s_0\in [s_-,s_+]$   
be so close to $s_+$  
 that $J^+(\hat \mu(s_0))\cap I^-(\hat p^+)\subset \K_{t_0}$. Then the given data  
  determines $\be_V(J^+(\hat \mu(s_0))\cap I^-(\hat p^+))$.  
Moreover, let  
 $s_{k}<s_{k-1}$, $s_0,\dots,s_K\in [s_-,s_+]$ be such that  
 $s_{j+1}>s_j-\kappa_3$ and $s_K=s_-$, and denote  
 $y_k=\hat \mu(s_k)$, see Fig.\ 18 for  the points $y_1$ and $y_2$. 

Above,   $\be_V(J^+(\hat \mu(s_0))\cap I^-(\hat p^+))$  was  
determined from the data. Next we use induction: We    
consider  
 $s_1\in (s_-,s_+)$ and assume  
 we are given  $\be_V(J^+(y_1)\cap I^-(\hat p^+))$ with $y_1={\hat \mu}(s_1)$.} Let us then consider  $s_2\in (s_1-\kappa_3,s_1)$.  
Our next aim is to find the light observation points $\be_V(J^+(y_2)\cap I^-(\hat p^+))$ with $y_2={\hat \mu}(s_2)$. To this end   
we need to make the following definitions  (see Fig.\ 18, 19, and 20). 
\medskip  
  
\begin{figure}[htbp]  \label{Fig-16}  
\begin{center}  
  
\psfrag{1}{$y_2$}  
\psfrag{2}{$p_2$}  
\psfrag{3}{$J_{g}(\hat p^-,\hat p^+)$}  
\psfrag{4}{$y_2^\prime$}  
\psfrag{5}{}  
\psfrag{6}{$y_1$}  
\psfrag{6}{$y_1$}  
\psfrag{7}{$\tilde p$}  
\psfrag{8}{$$}  
\includegraphics[width=7.5cm]{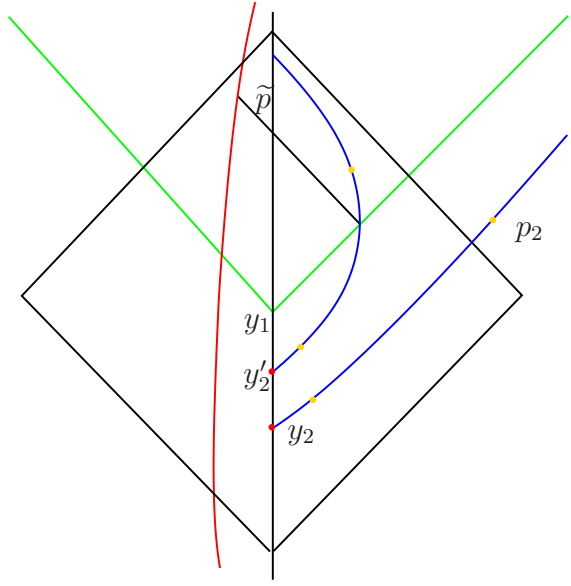}  
\end{center}  
\caption{A schematic figure where the space-time is represented as the  2-dimensional set $\R^{1+1}$.   
In section 5 we consider geodesic $\gamma_{y_2,\zeta_2}([0,\infty))$ sent from $y_2=\hat \mu(s_2)$.   
We consider two cases corresponding to geodesics  $\gamma_{y_2,\zeta_2}([0,\infty))$   
and $\gamma_{y_2^\prime,\zeta_2^\prime}([0,\infty))$.  
For $t_0>0$, the first cut point $p_2$, denoted by a golden point in the figure, of the geodesic $\gamma_{y_2,\zeta_2}([t_0,\infty))$  
is either outside the set $J_{g}(\hat p^-,\hat p^+)$ or is in the  
set $J^+(y_1)$, denoted by the green boundary, where  $y_1=\hat \mu(s_1)$.   
At the point $\tilde p=\mu_{z,\eta}((\mathbb{S} (y ,\zeta ,z,\eta,s_1)))$  we observe the first time  on the geodesic $\mu_{z,\eta}\cap I^-(\hat p)$  that  
the geodesic $\gamma_{y_2,\zeta_2}([0,\infty))$ has entered   
in the set $J^+(y_1)$.  The red curve is $\mu_{z,\eta}$  
}  
 \end{figure}

\begin{figure}[htbp] \label{Fig-17}  
\begin{center}  
  
\psfrag{1}{$y_2$}  
\psfrag{2}{$p_2$}  
\psfrag{3}{$J_{g}(\hat p^-,\hat p^+)$}  
\psfrag{4}{$y_2^\prime$}  
\psfrag{5}{}  
\psfrag{6}{$y_1$}  
\psfrag{7}{}  
\psfrag{8}{$\tilde  p$}  
\psfrag{A}{\hspace{-3mm}$\mu_{z,\eta}$}  
\psfrag{B}{$\mu_{z_0,\eta_0}$}  
\includegraphics[width=7cm]{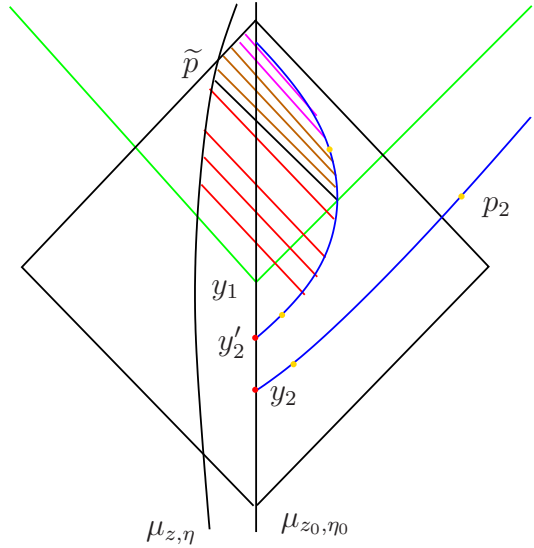}  
\end{center}  
  
\caption{A schematic figure where the space-time is represented as the  2-dimensional set $\R^{1+1}$.   
We consider  
the geodesics emanating from a point $x(r)=\gamma_{(y_2^\prime,\zeta_2^\prime)}(r)$,  
When $r$ is smallest value for which $x(r)\in J^+(y_1)$, a light-like geodesic  (black line segment)  
emanating from $x(r)$ is observed at the point $\tilde p\in \mu_{z,\eta}  
\cap  
I^-(\hat p^+)$.  
Then $\tilde p= \mu_{z,\eta}(\mathbb{S} (y ,\zeta ,z,\eta,s_1))$.  
When $r$ is small enough so that $x(r)\not \in J^+(y_1)$,  
the light-geodesics (red line segments) can be observed on $\mu_{z,\eta}$ in the set  $J^-(\tilde p)$.  
Moreover, when  $r$  is such that $x(r)\in J^+(y_1)$, the light-geodesics can be observed at $\mu_{z,\eta}$ in the set  $J^+(\tilde p)$. The golden point is the cut point  on $\gamma_{y_2^\prime,\zeta_2^\prime}([t_0,\infty))$ and the singularities on the light-like geodesics   
starting before this point (brown line segments) can be   
analyzed, but after the cut point the singularities on the light-like geodesics (magenta line segments)  
are not analyzed in this paper.  
}

 \end{figure}

\HOX{Note: We can simplify the construction by first showing that $m(z,\eta)$  
is continuous. Note that $m(z,\eta)=f^-_{\mu(z,\eta)}(\hat p^+)$.}  
{Let \beq\label{max s value}  
m (z,\eta)=\sup \{s\leq 1;\ \mu_{z,\eta}(s)\in I^-(\hat p^+)\},\quad (z,\eta)\in \U_{z_0,\eta_0}.  
\hspace{-1cm}  
\eeq  
Note that we can determine $m (z,\eta)$ using $U_{g}$.  
 }

\begin{definition}\label{def: hat s 1}  
Let $s_-\leq s_2<s_1\leq s_+$ satisfy  $s_1<s_2+\kappa_3$, let  
$s\in [s_2,s_1)$,  
 $y_j={\hat \mu}(s_j)$, $j=1,2$, $y={\hat \mu}(s)$ and $\zeta\in L^+_{y}U$, $\|\zeta\|_{g^+}=1$.  
  Moreover, recall that here $g=\hat g$ and let $(z,\eta)\in \W_{z_0,\eta_0}$ and $\mu=\mu(z,\eta)=\mu(g,z,\eta)=  
  \mu_{g,z,\eta}$  
and recall that $\hat \mu=\mu(g,z_0,\eta_0)$.  
  
%
%
%
%
%
When  
 $\gamma_{y ,\zeta }(\R_+)$ intersects $J^+(\hat \mu(s_1))\cap {I}^-(\hat p^+)$  
{we define  
\beq\label{special S formula}  
\mathbb{S} (y ,\zeta ,z,\eta,s_1)=\min(m(z,\eta),f^+_{\mu(z,\eta)}(q_0)),\quad \eeq  
where $ q_0=\gamma_{y ,\zeta }(r_0)$ and   
$r_0>0$ is the smallest $r\geq 0$ such that $\gamma_{y ,\zeta }(r)  
\in J^+(\hat \mu(s_1))$. In the case when  
$\gamma_{y ,\zeta }(\R_+)$ does not intersect $J^+(\hat \mu(s_1))\cap  {I}^-(\hat p^+)$,  
we define $\mathbb{S} (y ,\zeta ,z,\eta,s_1)=m (z,\eta)$.} \end{definition}  
  
In other words, if  
$\mathbb{S} (y ,\zeta ,z,\eta,s_1)<m(z,\eta)$, then it is  
 the first time when $\gamma_{y ,\zeta }(t)$ is observed on $\mu({z,\eta})$  
to enter the set $J^+(\hat \mu(s_1))\cap  {I}^-(\hat p^+)$.

\begin{definition}\label{def: hat s}  
Let   
$s\in [s_2,s_1)$, $y={\hat \mu}(s)$ and $\zeta\in L^+_{y}U$, $\|\zeta\|_{g^+}=1$.  
  
For $\vartheta\in (0,\vartheta_0)$, let  
$\qP^{(0)}_\vartheta (y ,\zeta )$  be the set of $(\vec x,\vec \xxi)=((x_j,\xxi_j))_{j=1}^4$  such that   
\begin{itemize}  
\item [(i)] $(x_j,\xxi_j)\in L^+U_{g}$, $j=1,2,3,4$,  
\item [(ii)]    
$(x_1,\xxi_1)=(y ,\zeta )$  
and  
 $(x_k,\xxi_k)$, $k=2,3,4$ are in $\vartheta$-neighborhood of   
 $(y ,\zeta )$  in $(T\hattuM _0,g^+)$ and  
$x_l\not \in J_{g}^+(x_j)$ for $ j,l=1,2,3,4$,  
\end{itemize}  
  
Also, we define  
$\qP_{\vartheta} (y ,\zeta ,z,\eta,s_1)$  to be the set of   
4-tuples compatible with earlier observations near $\mu_{z,\eta}$, that is,  
the set $\qP_{\vartheta} (y ,\zeta ,z,\eta,s_1)$ consist of all such  
$(\vec x,\vec \xxi)\in \qP^{(0)}_{\vartheta} (y ,\zeta )$  for which we have  
  
\begin{itemize}  
  
\item [(iii)]  
there is an open neighborhood $\mathcal W_0\subset \U_{z_0,\eta_0}$ of $(z,\eta)$  
such  that   
 \ba  
& &{\Sclo}_{z_1,\eta_1}((\vec x,\vec \xi),t_0)\in \be_{z_1,\eta_1}(J^+(y_1)\cap I^-(\hat p^+)),  
\quad \hbox{for all }(z_1,\eta_1)\in \W_0,  
\ea  
\hspace{-.7cm}and   
 \ba  
& &{\Sclo}_{z_1,\eta_1}((\vec x,\vec \xi),t_0)\in I^-(\hat p^+),  
\quad \hbox{for all }(z_1,\eta_1)\in \W_0.  
\ea  
  
\end{itemize}  
We denote $\qP^{(0)} (y ,\zeta )=\qP^{(0)}_{\vartheta_0} (y ,\zeta )$  
and $\qP(y ,\zeta ,z,\eta,s_1)=\qP_{\vartheta_0} (y ,\zeta ,z,\eta,s_1)$.  
\HOX{Notations  $\qP^{(0)} $ and $\qP$ are clumsy and should be changed.}

\end{definition}

Below, let $\kappa_4=\min(\kappa_2,\kappa_0)$ cf.\ Lemma \ref{lem: detect conjugate 0}.  
{
Let   $R_0(y,\zeta,s_1)\in [0,\infty]$ be the smallest value for which  
 $q_0=\gamma_{y,\zeta}(r_0)\in \p J^+(\hat\mu (s_1))$, or $R_0(y,\zeta,s_1)=\infty$  
 if no such value of $r$ exists.  
 By Lemma \ref{lem: detect conjugate 0} (iv), $\gamma_{y ,\zzeta }([0,t_0])$ does not   
intersect $J^+(y_1)\cap I^-(\hat p^+)$. Thus,  if   
$\gamma_{y ,\zzeta }([0,\infty))$  
intersects $J^+(y_1)\cap I^-(\hat p^+)$ if and only if  
$\gamma_{y ,\zzeta }([{t_0},\infty))$ intersects it.

\begin{lemma}\label{conjugatepoints are nice}  
  
Assume that $y=\hat \mu(s)$, $y_1=\hat \mu(s_1)$ with $-1<s_1-\kappa_3\leq s<s_1<s^+$, $\zeta\in L^+_yM$, $\|\zeta\|_{g^+}=1$. Then   for all  $\delta>0$ there  
  is $\vartheta_1(y,\zeta,s_1,\delta)>0$   
   such that if $0<\vartheta<\vartheta_1(y,\zeta,s_1,\delta)$ and  
 $(\vec x,\vec \xi)\in \qP^{(0)}_\vartheta(y,\zeta)$ and for some $j=1,2,3,4,$ we have ${\bf t}_j=  
 \rho(x_j(t_0),\xi_j(t_0))<\mathcal T(x_j({t_0}),\xxi_j({t_0}))$, then the cut point    
 $p_j=\gamma_{x_j(t_0),\xi_j(t_0)}({\bf t}_j)$ of   
 the geodesics $ \gamma_{x_j(t_0),\xi_j(t_0)}([0,\infty))$ satisfies either  
 \smallskip  
   
(i)  $p_j\not \in J^-(\hat p^+)$,  
 \smallskip  
  
\noindent  
or  
 \smallskip  
  
\HOX{Slava commented: these considerations can be simplified using a new parameter $\tilde \kappa_4$.}  
  
(ii) it  holds that  $r_0=R_0(y,\zeta,s_1)<\infty$,  
$q_0=\gamma_{y,\zeta}(r_0)$ satisfies  
$q_0\in J^-(\hat \mu(s_{+2})) $  
and   
 $r_1=r_0+\kappa_4<t_0+\rho(y(t_0),\zeta(t_0))$ and   
 $q_1=\gamma_{y,\zeta}(r_1)$ satisfy  
   $q_1\in  I^-(\hat \mu(s_{+3}))$, cf. (\ref{kaava D}) and text above it,  
 and  $f^+_{\mu(z,\eta)}(p_j)>f^+_{\mu(z,\eta)}(q_1)-\delta$.

  \end{lemma}  
  
}  
\noindent  
{\bf Proof.}  
We start with a Lipschitz estimate that we will need below.  
Consider $r>0$ such that  
$\gamma_{y,\zeta}(r)\in J^-(\hat \mu(s_{+3}))$ and    
let $(\vec x,\vec\xi)\in \mathcal R^{(0)}(y ,\zeta )$.  
Using the Lipschitz estimate (\ref{Lip estimate kiire}) for the   
geodesic flow in $J(\hat \mu(-1),\hat \mu(1))$,  
we see that there are $L_0=L_0(y,\zeta,r)>0$ and  
$\vartheta_2(y,\zeta,r)>0$ such that if  
$0<\vartheta<\vartheta_2(y,\zeta,r)$   
 and $d_{g^+}((x_j,\xi_j),(y ,\zeta ))<\vartheta$ then   
\beq\label{Li estim 2}  
d_{g^+}(\gamma_{x_j,\xi_j}(r), \gamma_{y ,\zeta }(r))<  
L_0\vartheta.  
\eeq

Next, let    $t_2=\rho(y ,\zzeta )$,  
$t_3=\rho(y(t_0),\zeta(t_0))+t_0$, and $(\vec x,\vec\xi)\in \mathcal R^{(0)}(y ,\zeta )$.  
 Let $q_2=\gamma_{y ,\zzeta }(t_2)$  
and $q_3=\gamma_{y ,\zzeta }(t_3)$.   
By Lemma \ref{lem: detect conjugate 0} (i)-(ii), we see  
that if  $t_2\geq R_1+\kappa_0$ then  
$\rho(x_j(t_0),\xi_j(t_0))+t_0> R_1$ and hence   
$p_j\not \in J^-(\hat \mu(s_{+2}))$. Thus it is enough  
to consider the case when $t_2<R_1+\kappa_0$  
  
We see using   
Lemma \ref{lem: detect conjugate 0}  (ii) that  
$t_3\geq t_2+3\kappa_2$. Note that  
$q_3=\gamma_{y({t_0}),\zeta({t_0})}(\rho (y({t_0}),\zeta({t_0})))$   
is the first cut point on $\gamma_{y({t_0}),\zeta({t_0})}([0,\infty))$.   
If this cut point satisfies  
$q_3 \in J^-(\hat \mu(s_{+2}))$, then by Lemma \ref{lem: detect conjugate 0} (iii),  
we have  
$f^-_{\hat \mu}(q_3)>f^-_{\hat \mu}(q_2)+3\kappa_3\geq s_2+3\kappa_3$,  
implying  $q_3\in I^+(y_1)$ as $s_1<s_2+2\kappa_3$. Thus the   
first cut point $q_3$ on $\gamma_{y({t_0}),\zeta({t_0})}([0,\infty))$ satisfies (see Fig.\ 18) 
\beq\label{eq: estimate for cut point}  
q_3\in I^+(y_1)\quad\hbox{or}\quad q_3\not \in J^-(\hat \mu(s_{+2})).  
\eeq  
  
Below, we consider two cases separately:  
\smallskip  
  
In case (a):  $q_3 \not \in J^-(\hat \mu(s_{+2}))$.  
%
\smallskip

In case (b):  $q_3 \in J^-(\hat \mu(s_{+2}))$.   
\smallskip

Assume next the alternative (b) is valid.  
  
 Then  
by (\ref{eq: estimate for cut point}) we have $q_3\in I^+(y_1)$ and hence  
$\gamma_{y ,\zeta }([0,\infty))$ intersects $J^+(y_1)\cap J^-(\hat \mu(s_{+2}))$.  
There exists    
 the smallest number  $r_0>0$ such that  
 $\gamma_{y,\zeta}(r_0)\in \p J^+(y_1)$.  
 We define $q_0=\gamma_{y,\zeta}(r_0)\in \p J^+(y_1)$. Note  
 that then $q_0\in J^-(\hat \mu(s_{+2}))$.

  
Let $r_1=r_0+\kappa_4$. Then by   
 Lemma \ref{lem: detect conjugate 0}  (i), we have   
$q_1=\gamma_{y,\zeta}(r_1)\in I^-(\hat \mu(s_{+3}))$.

Using Lemma \ref{lem: detect conjugate 0} (iii) for the point $q_0$, we see that  
as $f^-_{\hat \mu}(q_0)=s_1<s+2\kappa_3$, we have  
$r_0<\rho(y ,\zeta )+\kappa_2$.  
On the other hand,  by Lemma \ref{lem: detect conjugate 0} (ii),   
$\rho(x_j({t_0}),\xi_j({t_0}))+t_0>\rho(y ,\zeta )+3\kappa_2$,  
and hence $\rho(x_j({t_0}),\xi_j({t_0}))+t_0>r_0+2\kappa_2\geq r_1+\kappa_2$.

%
{
Let $S^*=f^+_{\mu(z,\eta)}(q_1)$ and $\delta>0$.  
As the functions $(x,\xi)\mapsto \gamma_{x,\xi}(r_1)$ and  $q\mapsto f^+_{\mu(z,\eta)}(q)$ are continuous and  
$f^+_{\mu(z,\eta)}(q_1)>S^*-\delta$, we see that  
when there is $\vartheta_1(y,\zeta,s_1,r_1,\delta)$   
such that if $0<\vartheta<\vartheta_1(y,\zeta,s_1,r_1,\delta)$  
then  
$q_j^\prime=\gamma_{x_j,\xxi_j}(r_1)$ satisfies  
$f^+_{\mu(z,\eta)}(q_j^\prime)>S^*-\delta$.  
As $q_j^\prime<p_j$, we have   
$f^+_{\mu(z,\eta)}(p_j)>S^*-\delta$.  
}

Let us next consider the case when the alternative  (a) is valid, that is,  
$q_3 \not \in J^-(\hat \mu(s_{+2}))$.  
  
Let $r_4=r_4(y,\zeta,s_2)>0$ be the smallest number  such that  
 $\gamma_{y,\zeta}(r_4)\in \p J^-(\hat \mu (s_{+2}))$  
 and define $q_4=\gamma_{y,\zeta}(r_4)$. Then using this  
and the definition of $q_4$, we have $t_3  
=\rho(y(t_0),\zeta(t_0))+t_0>r_4$.  
As $(x,\xi)\mapsto \rho(x,\xi)$ is a lower-semicontinuous function,  
 there exists $\vartheta_3(y,\zeta,s_1)>0$ such that   
if $0<\vartheta<\vartheta_3(y,\zeta,s_1)$ then  
$\rho(x_j(t_0),\xi_j(t_0))+t_0>r_4$.  
  
Observe that $q_4\not \in J(\hat \mu(s_{-2}),\hat p^+)=J^+(\hat \mu(s_{-2}))\cap J^-(\hat p^+)$  
and thus $h_2:=d_{g^+}(q_4,J(\hat \mu(s_{-2}),\hat p^+))>0$.  
Using (\ref{Li estim 2})  
we see that there exists $\vartheta_1(y,\zeta,s_1)\in (0,\vartheta_3(y,\zeta,s_1))$ such that   
if $0<\vartheta<\vartheta_1(y,\zeta,s_1)$  and  
%
$(\vec x,\vec\xi)\in \mathcal R^{(0)}_\vartheta(y ,\zeta )$, then $x_j\in J^+( \hat \mu(s_{-2}))$ and  
\ba  
d_{g^+}(\gamma_{x_j,\xi_j}(r_4),q_4)<\frac 12 h_2,  
\ea   
and hence  
\ba  
d_{g^+}(\gamma_{x_j,\xi_j}(r_4),J(\hat \mu(s_{-2})),\hat p^+))>\frac 12 h_2.  
\ea  
Assume next that  $0<\vartheta<\vartheta_1(y,\zeta,r_4)$.  
Then,    
as $x_j\in J^+(\hat \mu(s_{-2}))$,  we have  $\gamma_{x_j,\xi_j}(r_4)\in J^+(\hat \mu(s_{-2}))$  
and hence  
$\gamma_{x_j,\xi_j}(r_4)\not \in  J^-(\hat p^+)$.  
Then  
 $\rho(x_j(t_0),\xi_j(t_0))+t_0>r_4$ and $\gamma_{x_j,\xi_j}(r_4)\not  
 \in J^-(\hat p^+)$ imply that  
 $p_j=\gamma_{x_j({t_0}),\xxi_j({t_0})}(\rho(x_j({t_0}),\xxi_j({t_0})))\not \in J^-(\hat p^+)$. This proves the claim.  
  \hfill \Box \medskip

Next our aim is show that we can determine the function $\mathbb{S} (y ,\zzeta ,z,\eta,s_1)$.   
To this end we have to take care of the difficulty that the geodesic $\gamma_{y ,\zeta }$  
can exit $U_{g}$, later return to it and intersects  the geodesic  $\hat \mu$ or a  
geodesic $\mu_{z,\eta}$. This  
happen for instance in the Lorentzian manifold $\R\times \mathbb S^3$.  
First we consider the case when the geodesic $\gamma_{y ,\zeta }$  
does not intersect e.g.\ a geodesic  $ \mu_{z,\eta}$.

{

\begin{lemma}\label{lem: cut points are not bad}  
Let $s_-\leq s_2\leq s<s_1\leq s_+$ satisfy  $s_1<s_2+\kappa_3$, let  
 $y_j={\hat \mu}(s_j)$, $j=1,2$,  $y={\hat \mu}(s)$ and $\zeta \in L^+_{y }U$,  
  $\|\zeta\|_{g^+}=1$.   
Assume that   
$(z,\eta)\in \U_{z_0,\eta_0}$   
is such that  
 $\gamma_{y,\zeta}([t_0,\infty))$ does not intersect  
$\mu_{z,\eta}([-1,1])$.   
  
There is $\vartheta_3(y,\zeta,s_1,z,\eta)>0$ such that   
if  $0<\vartheta<\vartheta_3(y,\zeta,s_1,z,\eta)$ and   
$(\vec x,\vec \xxi)\in \qP^{(0)}_\vartheta(y ,\zeta)$ the following holds:  
  
If  
$\rho(x_j({t_0}),\xxi_j({t_0}))<\mathcal T(x_j({t_0}),\xxi_j({t_0}))$ for some $j=1,2,3,4$,  
then  
the cut point  $p_j=\gamma_{x_j({t_0}),\xxi_j({t_0})}(\rho(x_j({t_0}),\xxi_j({t_0})))$  
satisfies  either $p_j\not \in J^-(\hat p^+)$ or   
$f^+_{\mu(z,\eta)}(p_j)>\mathbb{S} (y ,\zeta ,z,\eta,s_1)$.  
  
\end{lemma}

\noindent{\bf Proof.} We denote below $S=\mathbb{S} (y ,\zeta ,z,\eta,s_1)$.  
The fact that $\gamma_{y,\zeta}([t_0,\infty))$ does not intersect  
$\mu_{z,\eta}([-1,1])$ implies   
that when $\vartheta>0$ is small enough, we have that  
$(\vec x,\vec \xxi)\in \qP^{(0)}_\vartheta(y ,\zeta)$ that  
$\gamma_{x_j,\xi_j}([t_0,\infty))$ does not intersect  
$\mu_{z,\eta}([-1,1])$.  
 Assume next that $\vartheta$ is so small  
that this is valid.

Using a short-cut argument  
 that if $t_2>t_1$ then  for all $j=1,2,3,4,$  
 \beq\label{eq: strict ineq for f+}  
 f^+_{\mu(z, \eta)}(\gamma_{x_j,\xi_j}(t_2))>  
 f^+_{\mu(z,\eta)}(\gamma_{x_j,\xi_j}(t_1)).  
 \eeq

%
%
  
Assume that $\rho(x_j({t_0}),\xxi_j({t_0}))<\mathcal T(x_j({t_0}),\xxi_j({t_0}))$ for some $j=1,2,3,4$. We can assume that $p_j\in J^-(\hat p^+)$ as otherwise the claim is proven.  
  
By  
 Lemma \ref{conjugatepoints are nice} (ii), then    
  $r_0=R(y,\zeta,s_1)<\infty$, $\gamma_{y ,\zzeta }([{t_0},\infty))$ intersects $J^+(y_1)\cap I^-(\hat \mu(s_{+2})))$,  
and  
  $q_0=\gamma_{y,\zeta}(r_0)\in \p J^+(y_1)$.  
Let $r_1=r_0+\kappa_4$ so that $q_1=\gamma_{y,\zeta}(r_1)\in  I^-(\hat \mu(s_{+3}))$.  
Using (\ref{eq: strict ineq for f+}) we see the  
$S_1=f^+_{\mu(z,\eta)}(q_1)>S$, and define  
$\delta=(S_1-S)/2>0$.  
%
%
By Lemma \ref{conjugatepoints are nice},   when  
 $0<\vartheta<\vartheta_1(y,\zeta,s_1,\delta)$, we have $\rho(x_j({t_0}),\xxi_j({t_0}))+t_0>r_1$  
and $f^+_{\mu(z,\eta)}(p_j)>S_1-\delta>S$.  
  \hfill \Box \medskip  
}

\begin{definition}  
Let  $T (y ,\zeta ,z,\eta,s_1)$ be the infimum of $s\in [-1,  
m (z,\eta)]$   
for which for every $\vartheta\in (0,\vartheta_0)$ there exists   
$(\vec x,\vec \xxi)\in \qP_\vartheta(y ,\zeta ,z,\eta,s_1)$  
such that  
 $\mu_{z,\eta}(s)=\pointear_{z,\eta}({\Sclo}((\vec x,\vec \xxi),{t_0}))$  
if such values of $s$ exist, and otherwise, let $T(y ,\zeta ,z,\eta,s_1)=  
{m (z,\eta)}$.  
\end{definition}

\begin{figure}[htbp] \label{Fig-18} 
\psfrag{1}{$x_1$}  
\psfrag{2}{$x_2$}  
\psfrag{3}{$q$}  
\psfrag{4}{$p$}  
\psfrag{5}{}  
\psfrag{6}{$y_1$}  
\psfrag{7}{$x$}  
\psfrag{8}{$z$}  
\begin{center}  
\includegraphics[width=6.5cm]{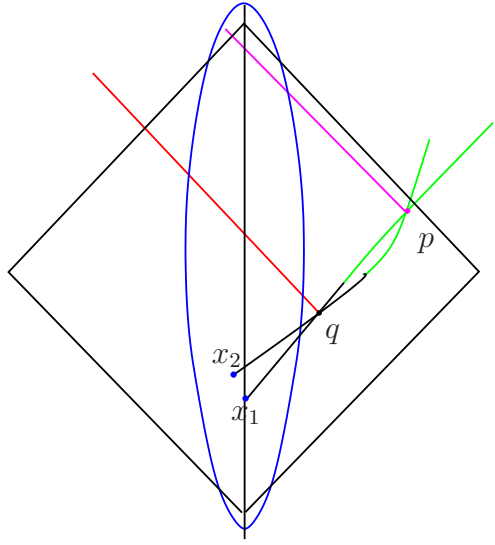}  
\end{center}  
\caption{  
A schematic figure where the space-time is represented as the  2-dimensional set $\R^{1+1}$.   
In section 5 we consider geodesics $\gamma_{x_j,\xi_j}([0,\infty))$, $j=1,2,3,4,$ that all intersect for the first time at a point  
$q$. In the figure we consider only geodesics $\gamma_{x_j,\xi_j}([0,\infty))$, $j=1,2$.  
When we consider geodesics $\gamma_{x_j,\xi_j}([t_0,\infty))$, $j=1,2$ with $t_0>0$, they may  
have cut points at $x^{cut}_j=\gamma_{x_j,\xi_j}({\bf t}_j)$.  
In the figure $\gamma_{x_j,\xi_j}([0,{\bf t}_j))$ are colored by black and the geodesics   
 $\gamma_{x_j,\xi_j}([{\bf t}_j,\infty))$ are colored by green. In the figure  
 the geodesics intersect at the point $q$ before the cut point and for  
the second time at $p$  
 after the cut point. We can analyze the singularities caused by the spherical waves   
 that interact at $q$ but not the interaction of waves after the cut points of the geodesics.  
 It may be that e.g. the intersection at $p$ causes new singularities to appear and we observe  
 those in $\hat U$. As we cannot ianalyze these singularities, we consider these singularities  
 as "messy waves". However, the "nice" singularities caused by the interaction at $q$ propagate along  
 the future light cone of the point $q$ arrive to $\hat U$ so that on a dense and open subset of geodesics   
 $\mu_{g,z,\eta}$  
 these "nice" singularities are observed before the "messy waves". Due to this we consider the first singularities  
 observed on  geodesics $\mu_{g,z,\eta}$.  
}  
 \end{figure}

\begin{lemma}\label{lem: Determination of hat s}  
Let $s_-\leq s_2\leq s<s_1\leq s_+$ satisfy  $s_1<s_2+\kappa_3$, let  
 $y_j={\hat \mu}(s_j)$, $j=1,2$,  $y={\hat \mu}(s)$ and $\zeta \in L^+_{y }U$,  
  $\|\zeta\|_{g^+}=1$, and   
$(z,\eta)\in \U_{z_0,\eta_0}$.   
  
  
Then  
\smallskip  
  
{\mattibf  (i)   
We always have $T (y ,\zeta ,z,\eta,s_1)\geq \mathbb{S} (y ,\zeta ,z,\eta,s_1)$.}  
\smallskip  
  
 (ii) Assume that   
$(z,\eta)\in \U_{z_0,\eta_0}$   
is such that  
 $\gamma_{y,\zeta}([t_0,\infty))$ does not intersect  
$\mu_{z,\eta}([-1,1])$.   
Then $T (y ,\zeta ,z,\eta,s_1)=\mathbb{S} (y ,\zeta ,z,\eta,s_1)$.  
  
(iii) For all $(z,\eta)\in \U_{z_0,\eta_0}$ we have   
\ba  
\mathbb{S} (y ,\zeta ,z,\eta,s_1)=\min(m(z,\eta),\liminf_{(z_1,\eta_1)\to (z,\eta)}  
T (y ,\zeta ,z_1,\eta_1,s_1)).  
\ea

\end{lemma}

\noindent{\bf Proof.}  
In the proof we continue to use the notations given in Def.\ \ref{def: hat s}  
and denote $T=T (y ,\zeta ,z,\eta,s_1)$. The objects used in the proof  
are shown in Figs.\ 19 and  20. 

(i)  
{Let $s^+_1<\mathbb{S} (y ,\zeta ,z,\eta,s_1)$  
and  $x_6=\mu_{z,\eta}(s_1^+)$.  
Note that as then $s^+_1<m (z,\eta)$. As $\mu_{z,\eta}(m (z,\eta))\in J^-(\hat p^+)$, we have  $x_6\ll \hat p^+$.}  
  
{
Let $\delta\in (0,\mathbb{S} (y ,\zeta ,z,\eta,s_1)-s^+_1)$.  
Assume that $(\vec x,\vec\xi)\in \qP^{(0)}_\vartheta(y,\zeta)$   
with some  
$0<\vartheta<\vartheta_1(y,\zeta,s_1,r_1,\delta)$, where  
$\vartheta_1(y,\zeta,s_1,r_1,\delta)$ is defined in  
Lemma \ref{conjugatepoints are nice}.    
  
%
  
Assume that for some $j=1,2,3,4$ it holds that $\rho(x_j({t_0}),\xxi_j({t_0}))<\mathcal T(x_j({t_0}),\xxi_j({t_0}))$ and consider the cut point  
$p_j=\gamma_{x_j({t_0}),\xxi_j({t_0})}(\rho(x_j({t_0}),\xxi_j({t_0})))$.  
By Lemma \ref{conjugatepoints are nice}, $p_j$  
 either satisfies  
$p_j\not\in J^-(\hat p^+)$, or alternatively, $r_0=R(y,\zeta,s_1)<\infty$ and $q_0=\gamma_{y,\zeta}(r_0)\in \p J^+(y_1)$ and for $r_1=r_0+\kappa_4$ we have  
 $q_1=\gamma_{y,\zeta}(r_1)\in  I^-(\hat \mu(s_{+3}))$ and    
 \ba  
 f^+_{\mu(z,\eta)}(q_1)\geq f^+_{\mu(z,\eta)}(q_0)=\mathbb{S} (y ,\zeta ,z,\eta,s_1)  
 \ea   
 and  
 \ba  
 f^+_{\mu(z,\eta)}(p_j)\geq f^+_{\mu(z,\eta)}(q_1)-\delta\geq \mathbb{S} (y ,\zeta ,z,\eta,s_1)-  
 \delta>s^+_1.  
 \ea   
Notice that  in the latter case $f^+_{\mu(z,\eta)}(p_j)>s^+_1$.  
}

%
%

Thus in all cases the cut points $p_j$   
of geodesics $\gamma_{x_j(t_0),\xxi_j(t_0)}$ satisfy $p_j\not \in J^-(x_6)$.  
Assume below that $\vartheta$ is so small enough that the above holds.

Next, assume  that the geodesics corresponding to $(\vec x,\vec\xi)\in \qP^{(0)}_\vartheta(y ,\zeta)$  intersect  
at some point $q\in J^-(x_6)$. Then  
$\be_{z,\eta}(q)\cap J^-(x_6)
=\{\mu_{z,\eta}(\tilde s)\}$ with $\tilde s=f^+_{\mu(z,\eta}(q)\leq s_1^+$.  
Hence Lemma \ref{lem: finding earliest observ. in J(x6)}  
yields that  ${\Sclo}_{z,\eta}((\vec x,\vec \xi),t_0)\cap J^-(x_6)
$ is  
equal to $\{\mu_{z,\eta}(\tilde s)\}$.  
However,  for all $q^\prime \in \gamma_{y ,\zeta }\cap  (J^+(y_1)\cap  
I^-(\hat p^+))$,  
we have $\be_{z,\eta}(q^\prime )=\{\mu_{z,\eta}(s^\prime)\}$ with  
$s^\prime\geq \mathbb{S} (y ,\zeta ,z,\eta,s_1).$ As   
$\tilde s\leq s^+_1<\mathbb{S} (y ,\zeta ,z,\eta,s_1)$, we see that   
the condition  (iii) in Def.\ \ref{def: hat s} can not be satisfied for $(\vec x,\vec\xi)$ and   
hence $(\vec x,\vec \xxi)\not \in \qP_\vartheta(y ,\zeta ,z,\eta,s_1)$.   

Assume next that  
$(\vec x,\vec \xxi)\in \qP_\vartheta(y ,\zeta ,z,\eta,s_1)$.    
The above yields that either  
 the geodesics corresponding to $(\vec x,\vec\xi)$ do not  intersect, that is,  
 $Q((\vec x,\vec\xi),t_0)=\emptyset$ or  
the intersection point $q=Q((\vec x,\vec\xi),t_0)$ is not in $J^-(x_6)$. Then we see using  
 Lemma \ref{lem: finding earliest observ. in J(x6)} that   
 ${\Sclo}((\vec x,\vec \xi),t_0)\cap J^-(x_6)=\emptyset $.  
 As ${\Sclo}((\vec x,\vec \xi),t_0)$ is closed, we have   
 either $e_{z,\eta}({\Sclo}((\vec x,\vec \xi),t_0))=\emptyset $ or  
 $\mu_{z,\eta}(s_1^+)=x_6\ll e_{z,\eta}({\Sclo}((\vec x,\vec \xi),t_0)) $.  
%
%
By definition of $T$, this gives that $T\geq s_1^+$.  
As above $s_1^+<\mathbb{S} (y ,\zeta ,z,\eta,s_1)$ was arbitrary, this yields $T\geq\mathbb{S} (y ,\zeta ,z,\eta,s_1)$.

%
%
%
%
%
     

  
(ii) Denote below $S=\mathbb{S} (y ,\zeta ,z,\eta,s_1)$.   
%
{
If $S=m(z,\eta)$  
and  the claim holds trivially. Thus,  
let us assume below that  $S<m(z,\eta)$.  
By definition (\ref{special S formula}), this implies that   
$r_0=R(y,\zeta,s_1)<\infty$ and $q_0=\gamma_{y,\zeta}(r_0)\in \p J^+(y_1)$  
is such that $f^+_{\mu(z,\eta)}(q_0)=S$.

%
%
%
%
  
 By Lemma \ref{lem: cut points are not bad},  
if $\rho(x_j({t_0}),\xxi_j({t_0}))<\mathcal T(x_j({t_0}),\xxi_j({t_0}))$ for some $j=1,2,3,4$,  
and $\vartheta$ is small enough,  
the cut point  $p_j=\gamma_{x_j({t_0}),\xxi_j({t_0})}(\rho(x_j({t_0}),\xxi_j({t_0})))$  
satisfies  either $p_j\not \in J^-(\hat p^+)$ or   
 $f^+_{\mu(z,\eta)}(p_j)>S$.  
  
}  
  

%
%
%
%
%
%
%
%

 {  
As $S<m(z,\eta)$ and above holds for all $j\in \{1,2,3,4\}$, we see that there is $s_{*}\in (S,m(z,\eta))$  
 such that $   
p_j=\gamma_{x_{j}({t_0}),\xxi_{j}({t_0})}(\rho(x_j({t_0}),\xxi_j({t_0})))\not \in J^-(\mu(s_{*}))$.  
%
}  
%
%
We observe that then $(\vec x({t_0}),\vec \xxi({t_0}))$   
satisfies conditions (\ref{eq: summary of assumptions 1})-(\ref{eq: summary of assumptions 2})  
with $x_6= \mu(s_*)$ and $q_0\in I^-(x_6)$.

 Let  
  $s^\prime\in (S,s_*)$. Note that  then $\mu_{z,\mu}(s^\prime)\in  I^-(x_6)$.  
{As $f^+_{\mu(z,\eta)}(q_0)=S$,  
the functions $(z_1,\eta_1) \mapsto f^+_{\mu(z_1,\eta_1)}(q_0)$  and  $(z_1,\eta_1) \mapsto\mu_{z_1,\eta_1}(s^\prime)$ are continuous, we see that   
 $(z,\eta)$ has a neighborhood $\W_1\subset \U_{z_0,\eta_0}$ such that  
for all $(z_1,\eta_1)\in \W_1$  we have  
$f^+_{\mu(z_1,\eta_1)}(q_0)<s^\prime$ and $\mu_{z_1,\eta_1}(s^\prime)\subset   
I^-(x_6)$. In particular then $\be_{z_1,\eta_1}(q_0)=\{\mu_{z_1,\eta_1}(f^+_{\mu(z_1,\eta_1)}(q_0))\}\subset I^-(x_6)$.  
Then  
by  Lemma \ref{lem: finding earliest observ. in J(x6)},   
${\Sclo}_{z_1^\prime,\eta_1^\prime}((\vec x,\vec \xi),t_0))\cap  J^-(x_6)$  
coincides with the set $\be_{z_1^\prime,\eta_1^\prime}(q_0)\cap J^-(x_6)$  
for all $(z_1^\prime,\eta_1^\prime)\in \U_{z_0,\eta_0}$.  
This and the above imply that   
%
${\Sclo}_{z_1,\eta_1}(((\vec x,\vec \xi),t_0)=\be_{z_1,\eta_1}(q_0)$ when $({z_1,\eta_1})  
\in \W_1$.  
As $q_0\in  J^+(y_1)\cap I^-(\hat p^+)$, we have $(\vec x,\vec \xxi)=((x_j,\xxi_j))_{j=1}^4\in \qP_\vartheta (y ,\zeta ,z,\eta,s_1)$.}

The above shows that   
there is  $(\vec x,\vec \xxi)\in \qP_\vartheta (y ,\zzeta ,z,\eta,s_1)$  
{\mattibf such that  the geodesics corresponding to $(\vec x,\vec \xxi)$ intersect  
in $q_0$ and  
$ \mu_{z,\eta}(s)={\Sclo}_{z,\eta}(((\vec x({t_0}),\vec \xxi({t_0})),{t_0})$  
with $s=S=f^+_{ \mu(z,\eta)}(q_0)$.  
  
Hence  
$T(y ,\zzeta ,z,\eta,s_1)\leq S=\mathbb{S} (y ,\zzeta ,z,\eta,s_1)$ and  
 the claim (ii) follows from (i). }

(iii)  Using standard results of differential topology, we see that there is  
an open and dense set $\W\subset \U_{z_0,\eta_0}$  
such that if  $(z,\eta)\in \W$ then   
$\mu_{z,\eta}([-1,1])$  does not intersect  $\gamma_{y,\zeta}([t_0,\infty))$.  
  
 Let $q_0=\gamma_{y,\zeta}(r_0)\in  
   \p J^+(y_1)$  
   be the point defined in the proof of (ii).  
By  Lemma \ref{B: lemma} (v),  
the function $(z,\eta)\mapsto f^+_{z,\eta}(q_0)$ is continuous  
and as $M_0\setminus I^-(\hat p^+)$ is closed,   
similarly to the first part of the proof of Lemma \ref{B: lemma} (v), we   
see that the function $(z,\eta)\mapsto m(z,\eta)$ is lower-semicontinuous.  
\HOX{We could simplify the final construction using continuity of $m(z,\eta)$.}  
{As   
$\mathbb{S} (y ,\zzeta ,z,\eta,s_1)=\min(m(z,\eta),f^+_{z,\eta}(q_0))$  
the claim (iii) follows  easily  
from (ii).}  
 \hfill \Box \medskip

Next we reconstruct $\be_V(q)$ when $q$ runs over a geodesic segment.   
  
\begin{lemma}\label{lem: Determination of be-sets}  
Let $s_-\leq s_2\leq s<s_1\leq s_+$ satisfy  $s_1<s_2+\kappa_3$, let  
 $y_j={\hat \mu}(s_j)$, $j=1,2$,  $y={\hat \mu}(s)$, and $\zeta \in L^+_{y }U$,  
   $\|\zeta\|_{g^+}=1$.  
When we are given the data set  ${\cal D}(\hat g,\hat \phi,\e)$, 
we can determine the collection $\{\be_V(q) 
;\ q\in (\gamma_{y ,\zzeta }([{t_0},\infty))\cap  I^-(\hat p^+) )\setminus J^+(y_1)\}$, where   
$V=U\cap I^-(\hat p^+)$.  
\end{lemma}  
  
  
{\bf Proof.}   
{In the proof, we consider  
$y ,\zzeta ,s_1$, and $t_0$ as fixed parameters and do not always  
indicate the dependency on the other parameters on those.  
  
%
   Let  us denote by $\mathcal W=\W(y ,\zeta ) \subset \U_{z_0,\eta_0}$ the open and dense    
set of those $(z,\eta)\in \U_{z_0,\eta_0}$ for which  
$\mu_{z,\eta}([-1,1])\cap \gamma_{y,\zeta}([t_0,\infty))=\emptyset$.  
  
Next, let  $(z,\eta)\in \W$ and denote $S=\mathbb{S} (y ,\zzeta ,z,\eta,s_1)$.   
Let $0<\vartheta<\vartheta_3(y,\zeta,s_1,z,\eta)$ and  
 consider $(\vec x,\vec \xxi)\in \qP_\vartheta^{(0)}(y ,\zeta)$.   
%
%
%
%
  
We define  
$x^\prime=x^\prime(z,\eta):=\mu_{z,\eta}(\mathbb{S} (y ,\zzeta ,z,\eta,s_1))$.  
Then $x^\prime  
\in \p J^+(q_0)$, $x^\prime\in J^-(\hat p^+)$.  


{
As  $0<\vartheta<\vartheta_3(y,\zeta,s_1,z,\eta)$, using  
Lemma \ref{lem: cut points are not bad} we see that  
if  
$\rho(x_j({t_0}),\xxi_j({t_0}))<\mathcal T(x_j({t_0}),\xxi_j({t_0}))$ for  some $j=1,2,3,4$, then $p_j=\gamma_{x_j({t_0}),\xxi_j({t_0})}(\rho(x_j({t_0}),\xxi_j({t_0})))$  
 either satisfies  
$p_j\not\in J^-(\hat p^+)$ or $f^+_{\mu(z,\eta)}(p_j)>S$.  
In both cases the cut points $p_j$   
of geodesics $\gamma_{x_j(t_0),\xxi_j(t_0)}$ satisfy $p_j\not \in J^-(x^\prime)$.  
  
%
%
%
%
%
%
Recall that for $(\vec x,\vec \xxi)$  
we denote the first intersection point of the geodesics $\gamma_{x_j,\xi_j}$  by $q=Q(\vec x,\vec \xxi)$ if  
such intersection point exists and otherwise we define  
$Q(\vec x,\vec \xxi)=\emptyset$. 

With the above definitions all cut points satisfy $p_j\not \in J^-(x^\prime)$  for all $j\leq 4$.  
By definition of $x^\prime$ and  Lemma \ref{lem: finding earliest observ. in J(x6)},   
we see  using $x^\prime$ as the point $x_6$, that  
if the geodesic corresponding to $(\vec x,\vec \xxi)$ intersect  
at some point $q=Q(\vec x,\vec \xxi) \ll x^\prime$ we have   
 ${\Sclo}_{z,\eta}((\vec x,\vec \xxi),{t_0})=\be_{z,\eta}(q) \ll x^\prime,$ and otherwise,  
${\Sclo}_{z,\eta}((\vec x,\vec \xxi),{t_0})\cap I^-(x^\prime)=\emptyset$.

{\mmmbf   
  
Consider next a general $(z,\eta)\in \U_{z_0,\eta_0},$ that is, we  
do not anymore require that $(z,\eta)\in \W(y,\eta)$.  
{We say that a sequence $((\vec x^{(\ell)},\vec\xi^{(\ell)}))_{\ell=1}^\infty$  
is a $\A_{z,\eta}(y ,\zeta )$ sequence  
and denote $((\vec x^{(\ell)},\vec\xi^{(\ell)}))_{\ell=1}^\infty\in \A_{z,\eta}(y ,\zeta )$,  
 if $(\vec x^{(\ell)},\vec\xi^{(\ell)})  
\in \mathcal R^{(0)}_{\vartheta(\ell)}(y ,\zeta )$ with $\vartheta(\ell)=1/\ell$  
and there is $\ell_0>0$ such that either there is  $p\in \mu_{z,\eta}$ such that for all $\ell\geq \ell_0$  
\ba  
{\Sclo}_{z,\eta}((\vec x^{(\ell)},\vec \xi^{(\ell)}),{t_0})  
\cap I^-(x^\prime(z,\eta)) =\{p\}  
\ea  
or alternatively,  for all $\ell\geq \ell_0$  
\ba  
{\Sclo}_{z,\eta}((\vec x^{(\ell)},\vec \xi^{(\ell)}),{t_0})\cap I^-(x^\prime(z,\eta))=\emptyset.  
\ea  
Note that as $x^\prime(z,\eta)=\mu_{z,\eta}(\mathbb{S} (y ,\zzeta ,z,\eta,s_1))\leq \hat p^+$,  
we have $I^-(x^\prime(z,\eta)) \subset I^-(\hat p^+)$ and we can therefore have   
\ba  
& &{\Sclo}_{z,\eta}((\vec x^{(\ell)},\vec \xi^{(\ell)}),{t_0})\cap I^-(x^\prime(z,\eta))\\  
& &={\Sclo}_{z,\eta}((\vec x^{(\ell)},\vec \xi^{(\ell)}),{t_0})\cap   
\mu_{z,\eta}((-1,\mathbb{S} (y ,\zzeta ,z,\eta,s_1))).  
\ea  
 We denote then   
 \ba  
 {\bf p}_{z,\eta}(((\vec x^{(\ell)},\vec\xi^{(\ell)}))_{\ell=1}^\infty)  
={\Sclo}_{z,\eta}((\vec x^{(\ell_1)},\vec \xi^{(\ell_1)}),{t_0})\cap I^-(x^\prime(z,\eta))  
\ea  
where $\ell_1$ above is chosen so that the right hand side does  
not change if $\ell_1$ is replaced with any larger value. Note that  
we can use any $\ell_1\geq \ell_0$. Note that here   
$ {\bf p}_{z,\eta}(((\vec x^{(\ell)},\vec\xi^{(\ell)}))_{\ell=1}^\infty)$   
is in fact a set-valued function; its value can  
either be one point, or an empty set.  
  
We say also that   $((\vec x^{(\ell)},\vec\xi^{(\ell)}))_{\ell=1}^\infty$  
is a $\A(y ,\zeta )$ sequence if there is open and dense  
set $\W^\prime\subset \U_{z_0,\eta_0}$ and a non-empty open   
set $\W^\prime_0\subset \W^\prime$  
 such that $((\vec x^{(\ell)},\vec\xi^{(\ell)}))_{\ell=1}^\infty$  
is a $\A_{z,\eta}(y ,\zeta )$ sequence  for all $(z,\eta)\in \W^\prime$  
and the set  
${\bf p}_{z,\eta}(((\vec x^{(\ell)},\vec\xi^{(\ell)}))_{\ell=1}^\infty)$ is non-empty  
for all $(z,\eta)\in \W^\prime_0$.  
  
Next, let us denote the $\A(y ,\zeta )$ sequences by $\mathcal X=((\vec x^{(\ell)},\vec\xi^{(\ell)}))_{\ell=1}^\infty$.  
  
Observe that ${\cal D}(\hat g,\hat \phi,\e)$  
 determines  the points $x^\prime(z,\eta)$,  
 the  sets $\A(y ,\zeta )$,   
  and ${\bf p}_{z,\eta}(\mathcal X)$  
for $\mathcal X \in  \A(y ,\zeta )$.

Let  $({z,\eta})\in \W$  
and consider   
a $\A_{z,\eta}(y ,\zeta )$  sequence  
 $\X=((\vec x^{(\ell)},\vec\xi^{(\ell)}))_{\ell=1}^\infty$.  
 Below,   
we use  $x_6=x^\prime(z,\eta)$. We consider two cases separately:  
   
\smallskip

Case (a): Assume  that for all $\ell$ large enough  the geodesics corresponding to $(\vec x^{(\ell)},\vec \xxi^{(\ell)})$  
intersect at a point $q_\ell\in (I^-(\hat p^+)\cap I^-(x^\prime(z,\eta)))\setminus J^+(y_1)$.  
As $\be_{z,\eta}(q_\ell)\ll x_6$,  
 Lemma \ref{lem: finding earliest observ. in J(x6)}  
implies $\be_{z,\eta}(q_\ell)= {\Sclo}_{z,\eta}((\vec x^{(\ell)},\vec \xi^{(\ell)}),{t_0})$  
for $\ell$ is large enough. By definition of $\X$ we have $\be_{z,\eta}(q_\ell)  
={\bf p}_{z,\eta}(\mathcal X)$ for all $\ell$ large enough.  
Recall that when $\be_{z,\eta}(q_\ell)$  is given for all  $(z,\eta)\in \W^\prime_0$ we can determine by  
Lemma \ref{lemma: global injectivity} the point  $q_\ell$  
uniquely. Thus we see that the points $q_\ell$ have to coincide  when $\ell$ is large enough.  
Next we denote the point $q_\ell$, when $\ell$ is large enough, by $q$.  
  
\smallskip  
  
Case (b): Assume that there are arbitrarily large $\ell$ such that  the intersection point $q_\ell$ of  
the geodesics corresponding to $(\vec x^{(\ell)},\vec \xxi^{(\ell)})$  
is not in $(I^-(\hat p^+)\cap I^-(x^\prime(z,\eta)))\setminus J^+(y_1)$ or it  does not  
exists. Then, when  
$\ell$ is large enough, then the set ${\Sclo}_{z,\eta}(\vec x^{(\ell)},\vec \xi^{(\ell)}),{t_0})$ does not intersect $I^-(\hat p^+)\cap I^-(x^\prime(z,\eta))$.  
  
\smallskip

By definition, for all $\X=((\vec x^{(\ell)},\vec\xi^{(\ell)}))_{\ell=1}^\infty\in \A(y ,\zeta )$ there is a non-empty open set  
$\W_0^\prime$ such that for large enough $\ell$ the intersection of  
${\Sclo}_{z,\eta}(\vec x^{(\ell)},\vec \xi^{(\ell)}),{t_0})$ and  
 $I^-(\hat p^+)\cap I^-(x^\prime(z,\eta))$ is non-empty for all $(z,\eta)\in \W_0^\prime$.  
Thus, as $\W_0^\prime\cap \W$ is non-empty, we see that the case (b) is not possible.  
 Thus for   
 all $\X\in \A(y ,\zeta )$ the case (a) has to hold and we have a well-defined intersection point  
 $q=Q(\vec x^{(\ell)},\vec \xxi^{(\ell)})$ where  
  the geodesics corresponding to $(\vec x^{(\ell)},\vec \xxi^{(\ell)})$  
intersect  when $\ell$ is large enough.  
Below we denote  
 it by $Q(\X)=q$. This point has to be on the geodesic   
 $\gamma_{y,\eta}$, and by the above considerations,  
 we see that it has to in the set  $\gamma_{y ,\zzeta }([{t_0},\infty))\cap (I^-(\hat p^+)\setminus J^+(y_1))$. On the other hand,  
  by Lemma \ref{lem: existence of vec x vec zeta}, for all  $q\in \gamma_{y ,\zzeta }([{t_0},\infty))\cap (I^-(\hat p^+)\setminus J^+(y_1))$ and all $\vartheta>0$  
there is $(\vec x,\vec \xxi)\in \mathcal R_{\vartheta}(y ,\zeta )$ such that the   
geodesics corresponding to $(\vec x,\vec \xxi)$  
intersect at $q$. As $q\in I^-(\hat p^+)$, we see that there  
is an open set $\W_0^\prime\subset \W$ such that for  
all $(z,\eta)\in \W_0^\prime$ we have  
$\mu_{z,\eta}((-1,m(z,\eta)))\cap \L^+_g(q)\not=\emptyset$ and thus  
 $f^+_{\mu(z,\eta)}(q)<m(z,\eta)$.  
  Combining these observations we conclude   
that the set $\{Q(\X);\ \X\in \A(y ,\zeta )\}$ coincides with the set   
$\gamma_{y ,\zzeta }([{t_0},\infty))\cap (I^-(\hat p^+)\setminus J^+(y_1))$.

}

When $\W^\prime$  is equal to  the set  $\W$ we have  
\beq\label{final condition 1}  
& & {\bf p}_{z,\eta}(\X)\cap (I^-(\hat p^+)\cap I^-(x^\prime(z,\eta)))=  
\hspace{-1cm}\\  
& & \nonumber \be_{z,\eta}(Q(\X)) \cap (I^-(\hat p^+)\cap I^-(x^\prime(z,\eta)))\   
 \hbox{for all }({z,\eta})\in \W^\prime. \hspace{-1cm}  
\eeq  
%


Note that  then the equation  
 (\ref{final condition 1})  is valid in particular when $\W^\prime=\W$.  
Unfortunately,  we have not above determined the set $\W$ and thus can not assumed it to be known. However, for any open and dense set  $\W^\prime$  
 the intersection of $\W^\prime\cap \W$ is an open and  dense subset of $\W$  
 and thus    
 \ba  
 \{{\bf p}_{z,\eta}(\mathcal X);\ (z,\eta)\in \W^\prime\cap \W\}=  
 \bigcup_{ (z,\eta)\in \W^\prime\cap \W} \mu_{z,\eta}((-1,m(z,\eta)))\cap  
\be_V(q)  
\ea  
 is a dense subset of  
$\be_V(q)\cap I^-(\hat p^+)$. Hence for all open and dense sets $\W^\prime$,  
  \ba  
  \be_V(q)\cap I^-(\hat p^+)\subset   
 \hbox{cl}(\{{\bf p}_{z,\eta}(\mathcal X);\ (y,\eta)\in \W^\prime \})\cap I^-(\hat p^+)\ea  
 and the equality holds when $\W^\prime$ is  the set $\W$.  
Using this  we see that if we take  
 the intersection of the all sets $\hbox{cl}(\{{\bf p}_{z,\eta}(\mathcal X);\ (y,\eta)\in \W^\prime \})\cap I^-(\hat p^+)$  
 where $\W^\prime \subset \U_{z_0,\eta_0}$   
 is an open and dense subset, we obtain  
the set $\be_{V}(q)$ for $q=Q(\X)$, where $V=U_{g}\cap I^-(\hat p^+)$.

Doing this construction for all for $\mathcal X\in \A(y ,\zeta )$,  
we determine the set   
$\be_{V}((\gamma_{y ,\zzeta }([{t_0},\infty))\cap I^-(\hat p^+))\setminus J^+(y_1)).$   
\hfill \Box \medskip

The above shows that  
the given data  ${\cal D}(\hat g,\hat \phi,\e)$ determine   
the  
collection   
$\{\be_V(q);  
\ q\in \gamma_{y ,\zzeta }([{t_0},\infty))\cap   (I^-(\hat p^+) \setminus J^+(y_1))\}$  
for all $y={\hat \mu}(s)$ and $\zeta   
\in L^+_{y }M_0$, $\|\zeta\|_{g^+}=1$, where $s \in [s_1-\kappa_3,s_1)$. Taking  
the union of all such collections and   
$\be_V(J^+({\hat \mu}(s_1))\cap I^-(\hat p^+))$  
and  $\be_V(\K_{t_0}\cap J^+({\hat \mu}(s )))$,  
we obtain the set $\be_V( J^+({\hat \mu}(s ))\cap I^-(\hat p^+))$.  
Iterating this construction for $s_2,s_3,\dots,s_K$,  with $s_{k+1}\in (s_k-\kappa_3,s_k)$,  
we can find $\be_V( I(\hat \mu (s^\prime),\hat p^+))$  
for all $s_-<s^\prime<s_+$.  
  
Let now for $s_-<s^\prime<s^{\prime\prime}<s_+$.  
Observe that   
 $\be_V (q)\in \be_V( J^+({\hat \mu}(s ))\cap I^-(\hat p^+))$  
 satisfies  $\be_V (q)\in \be_V( I(\hat \mu (s^\prime),\hat \mu (s^{\prime\prime})))$  
 if and only if $\be_{z_0,\eta_0} (q)\ll \hat \mu (s^{\prime\prime}).$  
 Thus we can find the   
sets $\be_V( I(\hat \mu (s^\prime),\hat \mu (s^{\prime\prime})))$  
for all $s_-<s^\prime<s^{\prime\prime}<s_+$.  
   
By Theorem \ref{main thm}  we can reconstruct the manifold  
 $I(\hat \mu (s^\prime),\hat \mu (s^{\prime\prime}))$  
for all $s_-<s^\prime<s^{\prime\prime}<s_+.$ Glueing these  
constructions together we obtain  
 $ I(\hat p^-,\hat p^+)$  
and the conformal structure on it.    
This proves   
Theorem \ref{main thm Einstein}. \hfill \Box \medskip

\observation{   
{\bf Remark 5.1.}   
When the adaptive source functions $\mathcal S_\ell$ are the  
ones given in Appendix $B$ we can consider an  
improvement of   Theorem \ref{main thm Einstein}:  
When we are given not  ${\cal D}(\hat g,\hat \phi,\e)$ but assume  
that we know ${\cal D}_{0}(\hat g,\hat \phi,\e)$,  
we can use Remarks  
3.1, 3.2, 3.3, 4.1, and 4.2 and the proofs of this section  
 to show that  $ I(\hat p^-,\hat p^+)$  
and the conformal structure on it can be reconstructed.  
This yields similar result to Theorem \ref{main thm Einstein}.  
}\medskip

{\bf Proof.}  (of Thm.\ \ref{alternative main thm Einstein}) As noted in Appendix B, by assuming that Condition A is valid and   
by making the parameter $\hat h$ determining  
$U_{g}$ smaller, there are  adaptive source function, or  
$\mathcal S$-functions,  given in the  
formula (\ref{S sigma formulas}) in Appendix B  
with $L\geq 5$ and  $K=L\,\cdotp (L!)+1$.  
 for which  
 the Assumption S is valid. Next we consider these $\mathcal S$-functions.  
  
Let us denote by $(g^\prime,\phi^\prime)$  
the solutions of (\ref{eq: adaptive model with no source})   
with some sources $(\F_1,\F_2)$ supported compactly in $U_{g^\prime}$.  
We want to find out when the observations $(U_{g^\prime},  
g^\prime|_{U_{g^\prime}},\phi^\prime|_{U_{g^\prime}},\F_1,\F_2)$,
representing an equivalence class in ${\cal D}^{alt}(\hat g,\hat \phi,\e)$ that 
satisfy the conservation law $\nabla^g_j({\bf T}^{jk}(g,\phi)+\F_1^{jk})=0$, are equivalent  
to the observations of some solution $(g,\phi)$ of the model (\ref{eq: adaptive model})  
when the $\mathcal S$-functions are  given by the  
formula (\ref{S sigma formulas}) in Appendix B.  
To do this, assume that we are given  the restrictions of $\hat g$ and  
the solution of (\ref{eq: adaptive model with no source}), denoted  $g^\prime,\phi^\prime$ in $U_{g^\prime}$  and  $(\F_1,\F_2)$.  
Then we can  compute   
\ba  
& &R^\prime:=\F_1,\\  
& &S_\ell=(\F_2)_\ell,\quad \ell=1,2,\dots,L,\\  
& &Z^\prime:=\sum_{\ell=1}^L S_\ell \phi_\ell^\prime  
\ea  
and find $P^\prime:=R^\prime-Z^\prime g^\prime.$  
%
After we have found these functions, we test if equations   
(\ref{meson number}) and (\ref{meson number2}) in the Appendix  
 B hold for  $g=g^\prime,$ $\phi=\phi^\prime,$  
  $P=P^\prime$,  $Z=Z^\prime$, $R=R^\prime$,   
 $Q_{L+1}=Z$ and some functions $(Q_\ell)_{\ell=1}^L$ that are compactly  
 supported in $U_g=U_{g^\prime}$. If this is the case  
the functions  
 $g,\phi$  in $U_{g}$  and $(P,Q)$, supported in $U_{g}$,   
 are restrictions of the solutions of   
  (\ref{eq: adaptive model})   
with the $\mathcal S$-functions  given by the  
formula (\ref{S sigma formulas}) in Appendix B  
with some $P$ and $Q$.   
  
Clearly, all solutions  
(\ref{eq: adaptive model})   
with the $\mathcal S$-functions  (\ref{S sigma formulas})  correspond to the   
solutions of  (\ref{eq: adaptive model with no source})  with some  
$(\F_1,\F_2)$,

Summarizing, for a given $(U_{g^\prime},  
g^\prime|_{U_{g^\prime}},\phi^\prime|_{U_{g^\prime}},\F_1,\F_2)$ that
represents an equivalence class in ${\cal D}^{alt}(\hat g,\hat \phi,\e)$ and
satisfies the conservation law $\nabla^g_j({\bf T}^{jk}(g,\phi)+\F_1^{jk})=0$,    
we can find out if there exists an element $(U_{g^\prime},  
g^\prime|_{U_{g^\prime}},\phi^\prime|_{U_{g^\prime}},F)\in   
{\cal D}(\hat g,\hat \phi,\e)$ with some $F$ supported in $U_{g^\prime}$, and if so, we can find this element.  
Hence, when we are given   
 the collection  
${\cal D}^{alt}(\hat g,\hat \phi,\e)$,  we can choose from it all elements that  
correspond to some element of ${\cal D}(\hat g,\hat \phi,\e)$,  
and thus we can find ${\cal D}(\hat g,\hat \phi,\e)$. The claim follows then from   
 Theorem \ref{main thm Einstein}.   
%
%
\hfill \Box \medskip

{\bf Proof.}\  (of Corollary \ref{coro of main thm original Einstein})  
{In the above proof of  Theorem \ref{main thm Einstein}, we used the  
 assumption that $\hat Q=0$ and $\hat P=0$ to obtain  
 equations (\ref{w1-3 solutions}). In the   
setting of Cor.\ \ref{coro of main thm original Einstein}  
 where the background source fields  
 $\hat Q$ and $\hat P$  
  are not zero,  
we need to assume   
in the computations related to sources (\ref{eq: f vec e sources})  
that there are  
neighborhoods $V_j$ of the geodesics $\gamma_j$ that satisfy  
$\supp({\bf f}_j)\subset V_j$,  
the linearized waves $u_j={\bf Q}_{\hat g}{\bf f}_j$ satisfy  singsupp$(u_j)\subset V_j$,  
and     
$V_i\cap V_j\cap (\supp(\hat Q)\cup \supp(\hat P))=\emptyset$ for $i\not =j$.}  
To this end, we have first consider  measurements for the linearized waves   
and check for given $(\vec x,\vec\xi)$ that no two geodesic $\gamma_{x_j,\xi_j}(\R_+)$  
do intersect at $U_{\hat g}$ and restrict all considerations for such $(\vec x,\vec\xi)$.  
Notice that such $(\vec x,\vec\xi)$ form an open and dense set in   
$(TU_{\hat g})^4$. If then the width $\hat s$ of the used  
spherical waves is chosen to be small enough, we see that condition   
$V_i\cap V_j\cap (\supp(\hat Q)\cup \supp(\hat P))=\emptyset$ is satisfied.  
  
The above restriction causes only minor modifications in the above\HOX{We have to check  
this text and that the above proof with minor modifications works}  
proof and thus, mutatis mutandis, we see that we can determine  
the conformal  
type of the  metric in all relatively compact subsets 
$I_{\hat g}(\hat\mu(s^\prime),\hat\mu(s^{\prime\prime}))\setminus (\supp(\hat Q)\cup \supp(\hat P))$, $s_-<s^\prime<s^{\prime\prime}<s_+$,
of $I_{\hat g}(\hat p^-,\hat p^+)\setminus (\supp(\hat Q)\cup \supp(\hat P))$.  
By glueing these  manifolds and $U_{\hat g}$ together, we find the conformal  
type of the  metric in $I_{\hat g}(\hat p^-,\hat p^+)$.  
After this the claim follows from Corollary \ref{coro of main thm original Einstein pre}.  

\hfill \Box \medskip  
%


\generalizations{
 \section{Generalizations and outlook}\label{sec: Generalization}

\HOX{T23. This section needs to be checked and made less messy. - Matti }
{
In this section we present a sketch of the analysis how the above results can
be applied for the single measurement inverse problem, that is, how we can  obtain an approximation
of the space-time with a single measurement. Recall that  
 Theorem \ref{main thm Einstein} concerns the case when all possible
measurements near the freely falling observer $\mu$ are known. In principle this implies that the perfect
space time cloaking with a smooth globally hyperbolic metric is not possible, that is, no two
space-times having different conformal structure appear the same in
all measurements. However, one can ask if 
one can make an approximative image of the space-time knowing
 only one measurement.  In general, in many
inverse problems several measurements can be packed together to 
one measurement.
For instance, for the wave equation with a  time-independent simple metric
this is done in \cite{HLO}. Similarly, Theorem \ref{main thm Einstein} and its proof
make it possible to do approximate reconstructions as discuss below.

\subsection{Pre-compact collection of Lorentzian manifolds}
For $m\in \Z_+$ and $\Lambda_{m}\in \R_+$, let $\M_{m,\Lambda}$
be the set of pointed  globally hyperbolic Lorentzian manifolds $(M,g,x_0)$ such that
the corresponding Riemannian manifold $(M,g^+)$
has the property that in the closure of the ball $B_{g^+}(x_0,m)$ with center $x_0$
and radius $r=m$
the covariant derivatives of the curvature tensor $R=R(g^+)$ of the Riemannian
metric $g^+$ satisfy $\|\nabla^j R\|_{g^+}\leq \Lambda$
for all $j\leq m$, the  injectivity radius $i_0(x)$ of $(M,g^+)$ satisfies
$i_0(x)>\Lambda^{-1}$ for all $x\in \hbox{cl}(B_{g^+}(x_0,m))$. 
Let us choose some numbers $\Lambda_m$ for all $m\in \Z_+$,
and let 
\ba
\M=\bigcap_{m\in \Z_+}\M_{m,\Lambda_m}
\ea
be the set endowed with the initial topology for which all
identity maps $\iota_m:\M\to \M_{m,\Lambda_m}$ are continuous.
Using the Cheeger-Gromov theory and the Cantor diagonalization
procedure, we see that the set $\B_R$ of the balls  $B_{g^+}(x_0,R)$ of radius $R$ of the pointed manifolds  
 $(M,g,x_0)\in \M$ form a precompact set in the  topology given by the Gromov-Hausdorff metric
 and the closure of $\B_R$, denoted $\hbox{cl}_{GH}(\B_R)$, consists of smooth manifolds,
 see \cite{And1,Cheeger}, see  also the applications for inverse problems in  \cite{AKKLT}.
 Indeed, we see that the balls of radius $R$ of the pointed manifolds in $\M_{m,\Lambda_m}$ are
 precompact for all $m$ and thus  using a Cantor diagonalization argument
 we see claimed is precompactness.
  
Let us fix $R>0$ and define  $\mathcal N$ to be the set 
$((M,g),(x_0,\xi_0))$  such that  $(M,g,x_0)$ is a pointed  Lorentzian manifold for which  $B_{g^+}(x_0,R)\in \hbox{cl}_{GH}(\B_R)$,
and a time-ike vector $(x_0,\xi_0)\in TM$ is such that the  
 freely falling observer $\hat \mu$ in $(M,g,x_0)$, such that $\hat \mu(-1)=x_0$ and 
 $\p_s\hat \mu(-1)=\xi_0$ satisfies
$U_{g}\cup J_{g}(\hat \mu(-1-\e_0),\hat \mu(1+\e_0))\subset B_{g^+}(x_0,R)$, with fixed parameters 
$\e_0>0$ and  $\hat h>0$ (that determines $U_g$), 
and finally that $\hat \mu([-1,1])$ is such that no light-like geodesic has cut points
in $J_{g}(\hat \mu(-1),\hat \mu(1))$.

\subsubsection{Geometric preparations for single source measurement}
Let $S^{g^+}(y,r)$ denote the sphere of $T_yM$ of radius $r$ with respect to the metric
$g^+$ and  $S^{g^+}(y)= S^{g^+}(y,1)$.
For  $(y,\xi)\in L^+M$ we define a spherical surfaces
\ba
& &\L((y,\xi);\rho)=\{\exp_y(t\theta);\ \theta\in L_y^+M,\ 
 \|\theta-\frac 1{\|\xi\|_{g^+}}\xi\|_{g^+}<\rho,\ t>0\},\\
 & &\L^c((y,\xi);\rho)=\{\exp_y(t\theta);\ \theta\in L_y^+M,\  
 \|\theta-\frac 1{\|\xi\|_{g^+}}\xi\|_{g^+}\geq \rho,\ t>0\}.
  \ea
  We say that $y$ is the center point of these surfaces.

In particular,
$((M,g),(x_0,\xi_0))\in \mathcal N$ implies that
no light-like geodesic starting from $\mu_{z,\eta}$ can intersect $\mu_{z,\eta}$ again.

 Using compactness of $\mathcal N$, \HOX{We need to check (\ref{consequence of compactness})
 and existence of $\rho>0$  and add parameter $r_0$ below. Idea below is first
 to choose $\e$, i.e. the accuracy on which we want to find $(M,g)$, then
 suitable $L$, then perturbation so that light-like geodesics do not intersect
 $y_a$, then $r_0$ and then $\rho$.}
for any $r_0>0$ we can choose $\rho=\rho(r_0)>0$ such that 
for all  $((M,g),(x_0,\xi_0))\in \mathcal N$  and any $y_j\in J_{g}(\hat \mu(s_-),\hat \mu(s_+)) $, $j=1,2,3,4$ and $\xi_j\in L^+_{y_j}M$ such that 
\beq\label{far away}
d_{\hat g^+}(y_j,
\gamma_{y_k,\xi_k}(\R_+))\geq r_0,\quad \hbox{ when $j\not   =k$},
\eeq
the  set
\beq\label{consequence of compactness}
\bigcap_{j=1}^4 \L((y_j,\xi_j);\rho)\cap J_{g}(\hat \mu(s_-),\hat \mu(s_+))\hbox{ contains at most one point.}\hspace{-2.1cm}
 \eeq
 Indeed, if there are no such $\rho$ we find a sequence of manifolds $(M_k,g_k)$
 and $(\vec y_k,\vec \xi_k)\in (TM_k)^4$ such that the sets
 (\ref{consequence of compactness}) contain at least two points with the parameter $\rho_k$
 and $\rho_k\to 0$ as $k\to \infty$. Considering suitable subsequences of $(M_k,g_k,x_k)$
 that converge to $(M,g,x)$ 
 and $(\vec y_k,\vec \xi_k)$ that converge to $(\vec y,\vec \xi)$ we 
 see that some  geodesics in $(M,g,x)$ have cut points in $J_{g}(\hat \mu(s_-),\hat \mu(s_+))$
 that is not possible due to the definition of $\mathcal N$. Hence required parameter $\rho>0$ exists.

Next, for $((M,g),(x_0,\xi_0))\in  \hbox{cl}( \mathcal N)$ we denote 
\ba
J((M,g),(x_0,\xi_0))=J_{g}(\hat \mu(s_-),\hat \mu(s_+))
\ea
and 
\ba
\mathcal J=\{(J((M,g),(x_0,\xi_0)),[g],x_0,\xi_0);\ ((M,g),(x_0,\xi_0))\in \hbox{cl}( \mathcal N)\},
\ea
where $[g]$ denotes the conformal class of $g$. Equivalence classes $[g]$
are closed sets and we use Hausdorff distance for to define
the distance for the equivalence classes $[g_1]$ and $[g_2]$ using representatives
of the equivalence classes for which the determinant of the metric tensor in the Fermi
coordinates are equal to the determinant of the metric tensor of the Euclidean
metric in the Fermi coordinates associated to a line. We call such representative of 
the equivalence class $[g]$ a Fermi-normalized metric and denote it by $g^{(n)}$.
The corresponding Riemannian metric is denoted by $g^{(n),+}$

If $p\in J_{g}(\hat \mu(s_-),\hat \mu(s_+))$, then for given $(z,\eta)$ there
is unique point $y\in \mu_{z,\eta}(s)$
for which there exists a light like geodesic
 $\gamma_{y,\xi}([0,t])$ that connects $y$ to $p$. Moreover, then
 $s=f^-_{z,\eta}(p)$. Let $(z_j,\eta_j)\in \U_{z_0,\eta_0}$, $j=1,2,\dots$,
 be dense set 
 $(z_j,\eta_j)\not= (z_0,\eta_0)$. We consider these points to be defined so that they have fixed 
 representation in the Fermi coordinates. \HOX{Add details on the points in the Fermi coordinates}
 Using compactness of $J_{ g}(\hat \mu(s_-),\hat \mu(s_+))$,
  we see that  there is $m_0=m_0(M,g,x_0,\xi_0)$ such that
 for all $p\in J_{ g}(\hat \mu(s_-),\hat \mu(s_+))$ there are four pairs $(z_{j_k},\eta_{j_k})$
 with $j_k\leq m_0$ such
 that $q\mapsto ( f^-_{z_{j_k},\eta_{j_k}}(q))_{k=1}^4$ forms smooth coordinates
 in a neighborhood of $p$.

Let choose dense set 
 of numbers $s_j^k\in (-1,1)$, $j,k\in \Z_+$ and let
 $y_{j,k}=\gamma_{z_j,\eta_j}(s_j^k)$. For each
 $y_{j^\prime,k^\prime}$ we choose a dense
 set of directions, $\xi_{j^\prime,k^\prime,n}$, $n=1,2,\dots,$ of $S^{g^+}(y_{j^\prime,k^\prime})\cap L^+_{y_{j^\prime,k^\prime}}M$. 
 After this, let us use the pairs $(y_{j^\prime,k^\prime},\xi_{j^\prime,k^\prime,n})$ to define
the pairs $(y_{j^\prime,k^\prime}(t_{j^\prime,k^\prime}),\xi_{j^\prime,k^\prime,n}(t_{j^\prime,k^\prime}))$ with some 
 $t_{j^\prime,k^\prime}>0$, and 
 re-enumerate the obtained pairs 
 $(y_{j^\prime,k^\prime}(t_{j^\prime,k^\prime}),\xi_{j^\prime,k^\prime,n}(t_{j^\prime,k^\prime}))\in L^+M$,
 $j^\prime,k^\prime,n\in \Z_+$ as 
$(y_a,\xi_a)$, $a\in \Z_+$. 
 
 Let $\A$ be the set of 4-tuples $\vec a=(a_\ell)_{\ell=1}^4\in \Z_+^4$ that contain
 four pairwise non-equal positive integers.  
 Let
 \ba
 q(\vec a)=\bigcap_{\ell=1}^4 \L((y_{a_\ell},\xi_{a_\ell});\rho)\cap J_{ g}(\hat \mu(s_-),\hat \mu(s_+)).
 \ea
 Note that $q(\vec a)$ contains then only one point or is empty
 if (\ref{far away}) is satisfied with some $r_0$ and $\rho<\rho(r_0)$. Also,
 let 
 \ba
 & & \Psi^{(1)}_k(\vec a)=\prod_{\ell,\kappa\in \{1,2,3,4\},\ \ell\not=\kappa } \phi(\frac 1k \dist_{g^{(n),+}}(y_{a_\kappa},
\gamma_{y_{a_\ell},\xi_{a_\ell}}(\R_+))),\\
& & \Psi^{(2)}_k(\vec a)=\prod_{\ell=1}^4 \phi(\frac 1k \dist_{g^{(n),+}}(q(\vec a), 
 \L^c((y_{a_\ell},\xi_{a_\ell});\rho)),\\
& & \Psi_k(\vec a)=\Psi^{(1)}_k(\vec a)\,\cdotp \Psi^{(2)}_k(\vec a), 
  \ea
where $\phi\in C(\R\cup\{\infty\})$ is a function with $\phi(t)=0$ for $t\leq 1$ and $\phi(t)=1$
for $t>2$ and $t=\infty$ and $\dist_{g^{(n),+}}$ is the distance
with respect to the metric $g^{(n),+}$ that is conformal to the metric $g^+$
and which determinant in the Fermi coordinates is normalized as above.

Using the existence
 of $m_0(M,g,x_0,\xi_0)$ and local coordinates associated to four pairs $(z_{j_k},\eta_{j_k})$,
we see that $\{q(\vec a);\ \vec a\in \A,\ q(\vec a)\not=\emptyset\}$ is a dense set of $J_{ g}(\hat \mu(s_-),\hat \mu(s_+))$.

 Let us then consider the continuous map $G:\mathcal J\to \R^{\Z_+\times \Z_+}$,
 \ba
G: (J((M,g),(x_0,\xi_0)),[g],x_0,\xi_0)\mapsto (\Psi_k(\vec a)\,f^+_{z_j,\eta_j}(q(\vec a)))_{\vec a\in \A,\ j\in \Z_+,\ k\in \Z_+}.
 \ea
  \HOX{Check that $G$ is continuous}
In the case when the set $q(\vec a)$ is empty, we define 
 $f^+_{z_j,\eta_j}(q(\vec j(n),\vec k(n))$ to be equal to $\sup\{s;\ \mu_{z_j,\eta_j}(s)\in 
 J_{ g}(\hat \mu(s_-),\hat \mu(s_+))\}$. In fact, we can define this value
 to be arbitrary because of the "cut-off" function $\Psi_k(\vec a)$.
 
  Note that the sets $\P_V(q)$ are  closed. We see using the arguments of the main text that the map $G$ is injective in $\mathcal J$. 
As 
 $\mathcal J$ is compact and $G$ is continuous, we see for $G_N:\mathcal J\to \R^{\Z_+\times \Z_+}$,
 \ba
G_N: (J((M,g),(x_0,\xi_0)),[g],x_0,\xi_0)\mapsto (\Psi_k(\vec a)\,f^+_{z_j,\eta_j}(q(\vec a)))_{|\vec a|\leq N,\ j\leq N,\ k\leq N}
 \ea
that for any $\e>0$ there is  $N$ 
 such that $G_N(g)$ determines the conformal type $(J((M,g),(x_0,\xi_0)),[g])$ up to $\e$-error in the  metric  in $\mathcal J$. 
 Here, for on $ \mathcal J $  we use
 the distance function defined using 
 $((B_1,[g_1]),(B_2,[g_2]))\mapsto
 d_{GH}((B_1,g_1^{(n),+}),(B_2,g_2^{(n),+}))$.
  \HOX{We need to define the topology of $\mathcal J$ carefully}

 \subsection{Single source measurement}
 Let us now consider spherical waves propagating
 on the spherical surfaces  $\L((y_{a_\ell},\xi_{a_\ell});\rho)$, $|\ell |\leq L$.
 When $N$ is large, these waves interact and produce 3-wave and 4-wave interactions.
 These waves can be produced by sources ${\bf f}_{a_\ell}$ supported in an arbitrarily
 small neighborhoods  $W_{a_\ell}\subset U_{\hat g}$ of the "conic" points $y_{a_\ell}$ of the spherical surfaces  $\L((y_{a_\ell},\xi_{a_\ell});\rho)$.

  We see that for generic set  (i.e.\ in an open and dense set) of  manifolds in $\mathcal N$ it
holds that   there are no points $y_{a_\ell}$  \HOX{This "generic" statement needs to be checked}
 and $y_{a_\ell^\prime}$ with ${a_\ell}\not= {a_\ell^\prime}$ that
 can be connected with a light-like geodesic. This is due to
 the fact that this condition holds in an open set and we see that if for any $L\in \Z_+$
 we can order the points $y_{a(\ell)}$, $\ell=1,2,\dots,L$ in causal order, such that $y_{a(\ell)}\not \in J^+(y_{a(\ell+1)})$, and then making a small perturbation to the metric near
each point $y_{a(\ell)}$, in causal order, we can  choose such perturbations
that the point $y_{a(\ell)}$ is not on the surface $\cup_{j<\ell}\L^+(y_{a(j)})$. 
{Below we will assume that neighborhoods \HOX{This generic statement needs to be checked}
 $W_{a_\ell}$ of $y_{a_\ell}$ of the  points $y_{a_\ell}$ are so small
 that $W_{a_\ell}\subset \hat U$ does not intersect any $\L^+(y_{a(j)})$, where
$j\not =\ell$ and $j\leq L$.

In the main text we used sources ${\bf f}_{a_\ell}$ that are elements of  the space
$\I(Y_{a_\ell})$  and produce waves $u_{a_\ell}\in \I (N^*Y_{a_\ell},N^*K_{a_\ell})$,
where $K_{a_\ell}\subset \L^+(y_{a(\ell)})$ are 3-dimensional subsets  and 
 $Y_{a_\ell}= \L^+(y_{a(\ell)})\cap {\bf t}^{-1}(T_{a(\ell)})$ are 2-dimensional sub-manifolds.
 As the 
products of  distributions associated to two lagrangian manifolds are difficult to analyze, we modify the construction by 
replacing $K_a$ by surfaces $K_a\subset \L^+_{\hat g}(y_a)$ and 
  ${\bf f}_{a_\ell}\in \I(Y_{a_\ell})$ by sources $\tilde {\bf f}_{a_\ell}\in \I(K_{a_\ell})$,
supported in the  neighborhood  $W_{a_\ell}$ of $y_{a_\ell}$.
As then $\tilde {\bf f}_{a_\ell}\in \I(N^*K_{a_\ell})$ and $N^*K_{a_\ell}$ is invariant in the future
bicharacteristic flow of $\square_{\hat g}$, the sources $\tilde {\bf f}_{a_\ell}$
 produce by \cite[Prop. 2.1]{GU1} waves  $ u_{a_\ell}\in \I (N^*K_{a_\ell})$,
where the principals symbols are not anymore evaluations of the principal
symbol of the source at some point but integrals of the principal symbol along a bicharacteristic.

The sources $\tilde {\bf f}_{a_\ell}\in \I (K_{a_\ell})$  can be chosen so that there are 
 neighborhoods $V_{a_\ell}$ of the surfaces  $\L((y_{a_\ell},\xi_{a_\ell});\rho)$  that satisfy
$W_{a_\ell}\subset V_{a_\ell}$,
the linearized waves $ u_{a_\ell}={\bf Q}_{ \hat g}\tilde {\bf f}_{a_\ell}$ satisfy  singsupp$( u_{a_\ell})\subset V_{a_\ell}$,
and   
\beq\label{eq: V ehto}
\hbox{$V_{a_\kappa}\cap W _{a_\ell}=\emptyset$  for all $\ell,\kappa\leq L$,
 $\ell\not =\kappa$.}
\eeq
We also assume that
$W_{a_\kappa}$  are such that $J^+_{\hat g}(W_{a_\ell})\cap
J^-_{\hat g}( W_{a_\ell})\subset  U_{\hat g}$ and that the wave map coordinates
corresponding to the Euclidean metric are well defined in
this set. Note that the above properties of $W_{a_\kappa}$ can be obtained
by choosing all surfaces $K_{a_\kappa}$, $\kappa\leq L$ with a the same parameter
$s_0$ that is small enough. Below, we assume that $s_0$ is chosen in such a way. 

Consider now the non-linear interactions produced by the source
\beq\label{Fe source}
F_\e(x)=\e\left(\sum_{\ell=1}^N\tilde {\bf f}_{a_\ell}(x)\right),
\eeq
where $\e>0$ is small.  
This source produces the wave
\beq\label{nl wave}
u(x)=u^{(4)}_\e(x)+O(\e^5),\quad u^{(4)}_\e(x)=\sum_{j=1}^4 \e^j\M^{(j)}(x),
\eeq
where we consider $O(\e^5)$ as a term which is so small that it can be considered to be negligible. \HOX{Maybe we could use $\vec \e$ with $\e_j=10^{-100-10j}$ etc and justuse the considerations of main text}
Let us next assume that we can measure  the singularities of $\M^{(j)}$ for all $j\leq 4$.
Let $S_j$ be the intersection of $\hat U$ and singular support of  $\M^{(j)}$, $j\leq 4$.

\subsubsection{Analysis of 1st, 2nd, and 3rd order waves}
As in the source $F_\e$ we use only one parameter $\e$ instead of four parameters
$\vec\e=(\e_1,\e_2,\e_3,\e_4)$ we need to consider self-interaction terms,
that is, in our computations the permutations $\sigma:\{1,2,3,4\}\to \{1,2,3,4\}$
need to be replaced by general (non-injective) maps $\sigma:\{1,2,3,4\}\to \{1,2,3,4\}$.
This does not cause difficulties in analyzing the terms $\M^{(j)}$ for all $j\leq 3$
as because of condition (\ref{eq: V ehto}) the source terms do not cause
self-interactions with terms corresponding to different Lagrangian manifolds
$\Lambda_{a_\ell}=N^*K_{a_\ell}$ and the interactions of linearized waves
corresponding to two such manifolds produce  waves that are in
$\I(N^*K_{a_\ell},N^*K_{a_\ell,a_\kappa})$.

For generic sources singsupp$(\M^{(1)}) $ is equal to $S_1=\cup_{\ell\leq L} K_{a_\ell}\cap \hat U$ and 
singsupp$(\M^{(2)})$  is equal to $S_2=\cup_{\ell,\kappa\leq L} K_{a_\ell,a_\kappa}\cap \hat U
\subset S_1$. 

Let us then consider singularities of $\M^{(3)}$ using similar methods that we used to analyze the fourth
order terms. Again, observe that we need to consider also the terms for 
which the map $\sigma:\{1,2,3\}\to \{1,2,3\}$ is not bijective. For such terms,
e.g., if $\sigma(1)=r_1$, $\sigma(2)=r_2$ and $\sigma(3)=r_3$ with $r_1=r_3$,   
we see using \cite[Lemmas 1.2 and 1.3]{GU1} 
 \cite[Thm.\ 1.3.6]{Duistermaat}  
 that e.g. $v\in u_{r_1}\cdotp {\bf Q}(u_{r_1}u_{r_2})$ has
 wave front set that is subset of $N^*K_{r_1r_2}$, implying
 that wave front set of ${\bf Q}v$ is subset of  $N^*K_{r_1r_2}$.
 Using this we see that outside $S_2$ only  the terms for which map $\sigma:\{1,2,3\}\to \{1,2,3\}$
 is bijective cause singularities. Let us next show that this kind of singularities can be observed.

To this end, let us consider the term
 $\bra F_\tau,{\bf Q}(A[u_3,{\bf Q}(A[u_2,u_1])])\cet$ where
 all operators $A[v,w]$ are of the type $A_2[v,w]=\hat g^{np}\hat g^{mq}v_{nm}\p_p\p_q w_{jk}$,
 cf.\ (\ref{eq: tilde M1}) and (\ref{eq: tilde M2}), we obtain  formulas similar to
  (\ref{pre-eq: def of D}) and   (\ref{eq: def of D}), where the first factor in 
   (\ref{pre-eq: def of D}) is removed. Next we show that using this and the similar analysis we did
   for the fourth order terms with the WKB analysis, we can see that the 
   principal symbol of  $\M^{(3)}$ in the normal coordinates does not vanish
   for generic interacting spherical waves. Indeed, the indicator function
   $\Theta^{(3)}_\tau$ can be obtained as a sum of terms similar to $T^{(3),\beta}_\tau$ given in
    (\ref{eq: 3rd order term}). For these terms we obtain in the Minkowski space formulas
    analogous to (\ref{first asymptotical computation}) and (\ref{t 4 formula}). Let us consider  the case when
     $b^{(j)}$, $j\leq 5$ are  light-like co-vectors such that
    $b^{(5)}\in \hbox{span}(b^{(1)},b^{(2)},b^{(3)})$ and vector $b^{(4)}$ is linearly independent of $b^{(1)},b^{(2)},b^{(3)}$  
 so that
    $y=A^{-1}x$ are coordinates such  that
   ${\bf p}_4=(A^{-1})^tb^{(4)}=0$. The particular example of such co-vectors that we consider, are
   \ba
& &b^{(5)}=(1,1,0,0),\quad b^{(1)}=(1,1-\frac 12\rhoepsilon_1^2+O(\rhoepsilon_1^3),\rhoepsilon_1,\rhoepsilon_1^3),\\
& &b^{(3)}=(1,-1,0,0),\quad b^{(2)}=(1,1-\frac 12\rhoepsilon_2^2+O(\rhoepsilon_2^3),\rhoepsilon_2,\rhoepsilon_2^{201})
\ea
 and $ \rhoepsilon_1=\rhoepsilon^{100}$, $\rhoepsilon_2=\rho$, and $b^{(4)}$ is some fixed vector.
 Then $b^{(5)}=\sum_{j=1}^3 a_j b^{(j)}$,
 where $a_1=O(1)$, $a_2\sim -a_1\rho^{99}$, and $a_3=O( \rho^{2})$ as $\rho\to 0$.   For
 $T^{(3),\beta}_\tau$ given in
    (\ref{eq: 3rd order term}), we have then 
   \HOX{WKB needs to be checked.}
   \ba
T^{(3),\beta}_\tau&=&
\bra  u^\tau, h\,\cdotp \B_3u_3\,\cdotp {\bf Q}_0(\B_2u_2\,\cdotp \B_1 u_1)\cet\\
& &\hspace{-1cm}=C\frac {\P_\beta^\prime} {\omega_{12}}
\bra u^\tau, h\,\cdotp u_3\,\cdotp v^{a-k_2+1,a-k_1+1} (\,\cdotp ;b^{(2)},b^{(1)})\cet
\\
\nonumber&&\hspace{-1cm} =C\frac {\P^\prime_\beta}{\omega_{12}} 
\int_{\R^4}
e^{i\tau(b^{(5)}\cdotp x)} 
h(x)(b^{(3)}\,\cdotp x)^{a-k_3}_+\cdotp\\
\nonumber& &\quad\quad\quad\cdotp
(b^{(2)}\,\cdotp x)^{a-k_2+1}_+
(b^{(1)}\,\cdotp x)^{a-k_1+1}_+\,dx,\\
& &\hspace{-1cm}
=\frac {C\P^\prime_\beta\det (A^\prime)}{ \omega_{12}} \int_{(\R_+)^3\times \R}
e^{i\tau {\bf p}\cdotp y} 
h(Ay)y_3^{a-k_3}y_2^{a-k_2+1}
y_1^{a-k_1+1}\,dy
\\
\nonumber&&\hspace{-1cm} =
 {C\det (A^\prime)\,{\P^\prime_\beta}}
\frac {(i\tau)^{-(10+3a-|\vec k_\b|)}(1+O(\tau^{-1}))}{{\bf p}_3^{(a-k_3+1)}{\bf p}_2^{(a-k_2+2)}{\bf p}_1^{(a-k_1+2)}\omega_{12}}\,,
\ea
   where ${\bf p}_j=g(b^{(5)},b^{(j)})=-\frac 12 \rhoepsilon_j^2+O(\rhoepsilon_j^3)$, $j\leq 3$ and $\omega_{12}=g(b^{(1)},b^{(2)})=
   -\frac 12 \rhoepsilon_2^2+O(\rhoepsilon_2^3)$. 
Moreover, with the appropriate choice of polarizations   $v_{(\ell)}^{nm}=\hat g^{nj}\hat g^{mk}v^{(\ell)}_{jk}$ given in  (\ref{chosen polarization}) for $j=2,3,4$, c.f. (\ref{vbb formula}), we have 
\beq\label{pre-eq: def of D prime}
\P_{\beta_1}^\prime =(v_{(3)}^{pq}b^{(1)}_pb^{(1)}_q)(v_{(2)}^{nm}b^{(1)}_nb^{(1)}_m)\D,
\eeq
 and
\ba
 \D =\hat g_{nj}\hat g_{mk}v_{(5)}^{nm}v_{(1)}^{jk},
 \ea 
 where we note that $v_{(5)}^{nm}$ does not need to satisfy any divergence type conditions.
We see that the term $T^{(3),\beta_0}_\tau$, where $k_1=k_1^{\beta_0}=4$, dominates the other terms 
$T^{(3),\beta}_\tau$,when 
$\rho\to 0$. Thus, using the real analyticity of the leading order
coefficient in the asymptotic of  $\Theta^{(3)}_\tau$, we see similarly to the main text that the principal symbol of $\M^{(3)}$ in the normal
coordinates does not vanish in the generic case on the manifolds
$\Lambda_{a_1,a_2,a_3}^{(3)}$ that are the flow out of the light-like
directions of $N^*K_{a_1,a_2,a_3}$.
     
   \HOX{We need to check here arguments!!}

   The above analysis implies that for generic manifold
   and for generic sources singsupp$(\M^{(3)}) $ is equal to 
   \ba
   S_3=\hat U\cap (\bigcup_{(\vec x_j,\vec \xi_j)=(y_{a_\ell},\zeta_{a_\ell}),\ \ell\leq L}
    \mathcal Y((\vec x,\vec \xi),t_0,s_0)).
   \ea
   In the generic
   case we can thus assume the set $S_3$ has the above form and that the given set
  singsupp$(\M^{(3)}) $ determines it.
  
\subsubsection{Analysis of the 4th order terms}

Let us  next consider singularities of $\M^{(4)}$ outside the set $S_1\cup S_2\cup S_3$.
Outside
these surfaces we look singularities of  $\M^{(4)}$ using the indicator functions  $\Theta^{(4)}_\tau$,  constructed using gaussian
beams.
Assume that $x_5\in \hat U\setminus (S_1\cup S_2\cup S_3)$.
As pointed out before we need to consider general (non-injective) maps $\sigma:\{1,2,3,4\}\to \{1,2,3,4\}$.
The case when $\sigma$ is a permutation is considered in the main text.
Let us thus assume that $\sigma$ is not a permutation.

We need to consider cases where the linearized waves $u_1,u_2,u_3$ are replaced by
$u_{r_1},u_{r_2},u_{r_3}$, where $r_1,r_2,r_3\in \{1,2,3,4\}$ are arbitrary,
and $u_4$, kept the same (changing it would be just
renumeration of indexes.) Also, we can assume that $\{r_1,r_2,r_3\}$
is not the set  $\{1,2,3\}$, as this case we have already analyzed.

We start by  considering the $\tilde T$-functions.
We need to consider terms similar to $\tilde T$ in the main text containing e.g.\
$\bra u_\tau,{\bf Q}(u_{r_1}u_{r_2})\cdotp {\bf Q}(u_{r_3}u_{r_4})\cet$,
when some $r_j$ are the same. Assume next that
e.g. $r_1$ and $r_3$ are the same.
We use that fact that  $v\in \I(\Lambda_1, \Lambda_2)$,
we have
WF$(v)\subset \Lambda_1\cup \Lambda_2$. Moreover,
by \cite[Thm.\ 1.3.6]{Duistermaat},
\ba
\hbox{WF}(v\,\cdotp w)\subset 
\hbox{WF}(v)\cup \hbox{WF}(w)
\cup\{(x,\xi+\eta);\ (x,\xi)\in \hbox{WF}(v),\ (x,\eta)\in \hbox{WF}(w)\}
\ea
and thus we see that
the wave front set of term ${\bf Q}(u_{r_1}u_{r_2})\cdotp {\bf Q}(u_{r_1}u_{r_4})$
is in the union of the sets $N^*K_{r_j}$,
$N^*(K_{r_1}\cap K_{r_2})$, and 
$N^*(K_{r_1}\cap K_{r_2}\cap K_{r_4})$. The elements in
the first two sets do not propagate in the bicharacteristic flow
and the last one propagates along $\mathcal Y((\vec x,\vec \xi),t_0,s_0)$.
Observe that in the generic case $\mathcal Y((\vec x,\vec \xi),t_0,s_0)\subset S_3$.
The other  $\tilde T$-functions can be analyzed similarly.
Thus in the generic case the  $\tilde T$-functions do not cause observable 
singularities in the complement of $S_3$.

Consider next the $T$-functions.
We need to consider terms similar to $T$ in the main text containing e.g.\
$\bra u_\tau,u_{r_1}\,\cdotp {\bf Q}(u_{r_2}\,\cdotp {\bf Q}(u_{r_3}\,\cdotp u_{r_4}))\cet$,
when some $r_j$ are the same. 
 Again we decompose ${\bf Q}={\bf Q}_1+{\bf Q}_2$ and consider cases $p=1,2$ as in the main text. In the case $p=2$ we need to consider
critical points of modified phase functions 
\ba
\nonumber \Psi^{\vec r}_2(z,y,\theta)
&=&\theta_1z^{r_1}+\theta_2z^{r_2}+\theta_3z^{r_3}+\theta_4y^4+\varphi(y).
\ea
Under the assumption that $x_5\not \in K_4$, we see similarly
as in the main text that there are no critical points. Now
\ba
d_{z,y,\theta}\Psi^{\vec r}_2&=&(\omega^{\vec r};
r+d_y\psi_4(y,\theta_4),
z^{r_1},z^{r_2},z^{r_3},d_{ \theta_4}\psi_4(y,\theta_4))
\ea
where $\omega^{\vec r}=(\omega^{\vec r}_j)_{j=1}^4$,
$\omega^{\vec r}_j=\sum_{n=1}^3\theta_n\delta_{n\in R(j)}$, $R(j)=\{n\in \{1,2,3\}; r_n=j\}.$
Note that $(z,\omega^{\vec r})\in \bigcap_{n\in \{r_1,r_2,r_3\}}N^*K_{n}$ and some of $r_j$
may be the same. 
We see as in the main text that as $x_5\not \in {\mathcal Y},$
alternative (A1) can not hold.
Note that if (A2) holds then $\omega^{\vec r}\in \X$, and we see as in 
the main text that this is not possible.  Thus the term with $p=2$ causes no observable singularities
near $x_5$.

In the case $p=1$, let $(z,\theta,y,\xi)$ be a critical point of
 the phase function
  \ba
 \Psi^{\vec r}_3(z,\theta,y,\xi)=\theta_1z^{r_1}+\theta_2z^{r_2}+\theta_3z^{r_3}+(y-z)\,\cdotp \xi+\theta_4y^4+\varphi(y).
 \ea
Then
 \beq\label{eq: critical points 1}
  & &\p_{\theta_j}\Psi_3=0,\ j=1,2,3\quad\hbox{yield}\quad z\in K_{{r_1}}\cap K_{r_2}\cap K_{r_3},\\
  \nonumber
& &\p_{\theta_4}\Psi_3=0\quad\hbox{yields}\quad y\in K_4,\\
  \nonumber
 & &\p_z\Psi_3=0\quad\hbox{yields}\quad\xi=\omega^{\vec r},\\
   \nonumber
& &\p_{\xi}\Psi_3=0\quad\hbox{yields}\quad y=z,\\
  \nonumber
 & &\p_{y}\Psi_3=0\quad\hbox{yields}\quad \xi=
 -\p_y\varphi(y)-\theta_4 w.
  \eeq
  The critical points we need to consider for the asymptotics satisfy also
\beq\label{eq: imag.condition copu} \\ \nonumber
\im \varphi(y)=0,\quad\hbox{so that }y\in \gamma_{x_5,\xi_5},\ \im d\varphi(y)=0,\
\re d\varphi(y)\in L^{*,+}_y\hattuM _0.\hspace{-1cm}
\eeq
Now, as we assume that $\{r_1,r_2,r_3\}$
is not the set  $\{1,2,3\}$, some of $r_j$ is equal to $4$ or two of these are the same.   For simplicity, assume that
$r_1$ is either 4 or alternatively, $r_1$ coincides with $r_2$ or $r_3$. As $w\in N^*K_4$, these imply that 
$\p_y\varphi(y)\in N^*K_{r_2}\cap  N^*K_{r_3}\cap N^*K_4\subset \mathcal X((\vec x,\vec\xi);t_0)$.
Hence the assumption that $x_5\not \in  {\mathcal Y}$ implies that
terms with $p=1$ do not cause any singularities near $x_5$. Thus only the singularities of
$T$-type terms for which $\sigma$ is bijection,
that were analyze in the main text, are visible in $\hat U\setminus (S_1\cup S_2\cup S_3)$.
The singularities could be analyzed  in much more careful way
if we could analyze products of $v\in \I(K_1,K_{13})$ and $w\in \I(K_2,K_{23})$
and show that such product is a sum of conormal distributions.\hiddenfootnote{HERE ARE SOME EXTRA COMPUTATIONS:
By applying method of stationary phase in $\xi$ and $z$ variable we obtain
a integral similar to (\ref{eq; T tau asympt modified}),
 \beq \label{eq; T tau asympt modified2}
& &T_{\tau,1,1}^{(4),\beta}
=c\tau^{4}
\int_{\R^{8}}e^{i(\theta_1y^{r_1}+\theta_2y^{r_1}+\theta_3y^{r_1}+\theta_4y^{r_1})+i\tau \varphi(y)}
c_1(y, \theta_1, \theta_2)\cdotp\\ \nonumber
& &
\cdotp a_3(y, \theta_3)
q_{1,1}(y,y,\omega^{\vec r}_\beta(\theta))a_4(y,\theta_4) a_5(y,\tau )\,d\theta_1d\theta_2d\theta_3d\theta_4dy.
\eeq
where $\omega^{\vec r}_\beta$ is similar to $\omega^{\vec r}$  is defined with 
$R_\beta(j)=\{n\in \{1,2,3\}; r_n=\sigma_\beta(j)\}$ is  similar to $R(j)$.

 Similarly, to the analysis related to the function $\chi(\theta)$ that we vanish
 in a neighborhood of the set $\mathcal A_q$ given in (\ref{set Yq}),
 we can define a fucntion 
 $\tilde \chi(\theta)\in C^\infty(\R^4)$ that 
vanishes in a $\e_3$-neighborhood $\W$ (in the $\hat g^+$ metric) of
$\tilde {\mathcal A_q}$,
\ba
\tilde {\mathcal A_q}&:=&
\bigcup_{1\leq j<k\leq 4}N_q^*K_{jk}.
\ea
The idea of the sets $N_q^*K_{jk}$ is that in the complement of  their  neighborhoods at least three of the coordinates $(\theta_1,\theta_2,\theta_3,\theta_4)$
 are comparable to $|\vec\theta|$, that is, there is $j_0$ such that $|\theta_k|\geq c_j|\vec \theta|$ for $k\not =j_0$.
 Indeed, e.g.\ $(q;(\theta_1,\theta_2,\theta_3,\theta_4))\not \in N_q^*K_{jk}$ implies 
 that $(\theta_j,\theta_k)\not =0$. As this holds for all $(j,k)$, we see
 that $\theta\not \in \W$ and $|\theta|_{\R^4}=1$, only one coordinate of  $\theta$ can vanish. This makes the symbol
$\tilde \chi(\theta)b(y,\theta)$ a sum of a product type symbols $S^{p,l}(M_0;\R^3\times \R\setminus \{0\})$.

Let  \ba\hbox{$\tilde b_0(y,\theta)=\phi(y)\tilde \chi(\theta)b(y,\theta)$}\ea be a  product type symbol determining 
a  distribution 
$
\tilde \F^{(4),0}$ 
 that is given by the formula (\ref{f4-formula}) with
$b(y,\theta)$ being replaced by 
 $b_0(y,\theta)$.
We  can consider
 $\tilde \F^{(4),0}$ as a sum of in three components: 
 \ba
 \tilde \F^{(4),0}=\tilde \F^{(4),0,1}+\tilde \F^{(4),0,2}+\tilde \F^{(4),0,2}\ea
  where 
 \ba
 \tilde \F^{(4),0,1}&\in& \sum_{ \{k_1,k_2,k_3,k_4\}=\{1,2,3,4\}} \I^{p(\vec k)-4+1/2}(K_{k_2,k_3,k_4})\cap  \I(N^*K_{k_2,k_3,k_4}\setminus N^*\{q\})\\
& &\subset 
\sum_{ \{k_1,k_2,k_3,k_4\}=\{1,2,3,4\}} \I^{p(\vec k)-4+1/2-3/2+1}(N^* K_{k_2,k_3,k_4})\
,
\ea
where $p(\vec k)=p_{k_2}+p_{k_3}+p_{k_4}$ and
\ba
\tilde \F^{(4),0,2}&\in& \I^{p-4}(\{q\}\cap  \I(N^*\{q\}\setminus (\cup_{\vec k} N^*K_{k_2,k_3,k_4}))\\
& &\subset \I^{p-4-2+1}(N^*\{q\})
\ea
 and  $\tilde \F^{(4),0,2}$ is microlocally supported in some (arbitrarily) small
 neighborhood $\W$ of  $\cup (N^*K_{k_2,k_3,k_4}\cap N^*\{q\})$. 
  Let $\Lambda^{(3)}$ be the flowout of $\cup_{\vec k} N^*K_{k_2,k_3,k_4}$.
We see that $\M^{(4)}$ is Lagrangian distribution
 with known order in $I^-(x_6)\setminus  {\mathcal Y}((\vec x,\vec \xi),t_0)$, that is outside
 the set $\cap_{s_0>0}{\mathcal Y}((\vec x,\vec \xi),t_0,s_0)$ that has the Hausdorff dimension
 two. Indeed,
\ba
& &{\bf Q}_{\hat g}\tilde \F^{(4),0,1}\in
\sum_{ \{k_1,k_2,k_3,k_4\}=\{1,2,3,4\}} \I^{p(\vec k)-4-3/2}
\I(\hat U\setminus K_{k_2,k_3,k_4};
\Lambda^{(3)}),\\
& &{\bf Q}_{\hat g}\tilde \F^{(4),0,2}\in 
\I^{p-5-3/2}(\hat U\setminus\{q\};\Lambda_q),\\
  & &\hbox{WF}({\bf Q}_{\hat g}\tilde \F^{(4),0,3})\subset  \V .\ea
where $\V\subset T^*M_0$ is a neighborhood of $ 
T^*\mathcal Y((\vec x,\vec \xi),t_0)$ that is
the flow out of the set $\W$ in the bicharacteristic flow.
This makes it possible to analyze $
 \bra u_\tau,\F_1^{(4)}\cet$ and see that ${\bf Q}\F_1^{(4)}$ is Lagrangian distribution
 in $\hat U\setminus {\mathcal Y}((\vec x,\vec\xi),t_0)$.

To study the singularities of $u^{(4)}_\e(x)$, we consider next interaction of four spherical waves, 
corresponding the $\vec a=(a_\ell)_{\ell=1}^4 \in \A$, $a_\ell\leq N$ that are singular on the surfaces
$K_{a_\ell}=\L((y_{a_\ell},\xi_{a_\ell});\rho)$.
  Assume now that $\bigcap_{\ell=1}^4 \K_{a_\ell}$ contains a point $q$ and
  let $K_{a_1,a_2}=\bigcap_{\ell=1}^2 \K_{a_\ell}$ and $K_{a_1,a_2,a_3}=\bigcap_{\ell=1}^3 \K_{a_\ell}$ . 
 Then, we see that   $\Lambda_q^+$ and
  the flowout $\Lambda_{a_1,a_2,a_3}$ of $N^*K_{a_1,a_2,a_3}$ in the canonical
  relation of $\square_g$ are in $T^*U_g$ four dimensional Lagrangian
  manifolds that intersect only on a three dimensional set.
  Assume now that the source ${\bf f}_{a_\ell}$ have such orders 
  that $u_{a_\ell}={\bf Q}_g{\bf f}_{a_\ell}\in \I^{p_{a_\ell}}(K_{a_\ell})$ outside support of ${\bf f}_{a_\ell}$.

The linearized term $\M^{(1)}$ in (\ref{nl wave}) has the wave front set that is
a subset of  $S_1=\bigcup_{a\leq N}N^*K_a$,
the term  $\M^{(2)}$  has wave front set is a subset of  $S_1\cup S_2$ where  $S_2=\bigcup_{a_1<a_2\leq N}N^*K_{a_1,a_2}$. Let us note that when analyze the term $\M^4$, we have
to consider the term $\M^{(1)}\,\cdotp\M^{(1)}$ containing terms
$u_{a_1}\,\cdotp u_{a_1}$, where $u_{a_1}\in \I^{p}(K_{a_1})$. 
Denote $b(x^1)=u_{a_1}(x^1,x^\prime)$ and assume that $\hat b(\xi_1)\sim c_0 |\xi_1|^p$ as $\xi_1\to \infty$,
we see that 

\ba
(\hat b*\hat b)(\xi_1)=\int_\R \hat b(\xi_1-t)\hat b(t)dt
\sim 
|\xi_1|^p
\int_\R c_0(1-|\xi_1|^{-1}t)^p\hat b(t)dt
\sim c_0 b(0)|\xi_1|^p,
\ea

****GENERALIZATION***

Let  $B_1\in \I^{p_1}(\{0\})$ and $B_2\in \I^{p_2}(\{0\})$,
Let us write symbol $b_1(x,\xi)$ of $B_1$ as a sum of positive homogeneous term and
a compactly $\xi$-supported term $b_1(x,\xi)=b_1^h(x,\xi)+b_1^c(x,\xi)$ and similar
lower order terms.
Then $B_1^c,B_2^c\in C^\infty$, and
 $B_1^h\,\cdotp B_2^c\in \I^{p_1}(\{0\})$,  $B_1^c\,\cdotp B_2^h\in \I^{p_2}(\{0\})$,
  $B_1^c\,\cdotp B_2^c\in C^\infty$, and symbol of  $B_1^h\,\cdotp B_2^h$ is 
  obtained via substitution $t=|\xi_1|s$,
\ba
(b_1*b_2)(x,\xi_1)&=&\int_\R b_1(x,\xi_1-t) b_2(x,t)dt\\
&=&|\xi_1|^{p_1+p_2+1}
\int_\R b_1(x,1-s) b_2(x,s)ds
\ea
and thus  $B_1^h\,\cdotp B_2^h\in \I^{p_1+p_2+1}(\{0\})$.
Hence,  $B_1\,\cdotp B_2\in \I^{r}(\{0\})$, $r=\max(p_1,p_2)$.

****GENER. 2***

Let  $B_1\in \I^{p_1}(\{0\})$ and $B_2\in \I^{p_2}(\{0\})$,
Let us estimate the symbol $b_1(x,\xi)$ of $B_1$
by
\ba
|b_1(x,\xi)|\leq b_1^h(x,\xi)+b_1^c(x,\xi)
\ea 
where $b_1^c(x,\xi)$ is a compactly $\xi$-supported  symbol and
$b_1^h(x,\xi)$ is $p_1$-homogeneous. Let $B_1^c$ be the conormal distribution
with symbol $b_1^c(x,\xi)$ and  $B_1^h=B_1-B_1^c$, and introduce similar notations for $B_2$.
Then $B_1^c,B_2^c\in C^\infty$, and
 $B_1^h\,\cdotp B_2^c\in \I^{p_1}(\{0\})$,  $B_1^c\,\cdotp B_2^h\in \I^{p_2}(\{0\})$,
  $B_1^c\,\cdotp B_2^c\in C^\infty$, and symbol $c^h$ of  $B_1^h\,\cdotp B_2^h$ can be estimated
using substitution $t=|\xi_1|s$,
\ba
|c^h(x,\xi_1)|&\leq &\int_\R b_1(x,\xi_1-t) b_2(x,t)dt\\
&\leq &|\xi_1|^{p_1+p_2+1}
\int_\R b_1(x,1-s) b_2(x,s)ds
\ea
and thus  $B_1^h\,\cdotp B_2^h\in \I^{p_1+p_2+1}(\{0\})$.
Hence,  $B_1\,\cdotp B_2\in \I^{r}(\{0\})$, $r=\max(p_1,p_2)$.

****GENER. 3***

We need also to analyze the products of 
  $B_1\in \I^{p_1,l_1}(K_1,K_{13})$ and $B_2\in \I^{p_2,l_2}(K_2,K_{13})$.
  When we analyze this in the case when $|\theta_2|>c|\theta|$ and $|\theta_3|>c|\theta|$ 
  and  $|\theta_1|>c|\theta|$, this is equivalent to analyzing product of 
  $ \I(K_{13})$ and  $ \I(K_{12})$.
  The case where only  $|\theta_2|>c|\theta|$ and $|\theta_3|>c|\theta|$ but $\theta_1$ may be
  small  is difficult. Note that it is possible that we do not need to analyze
  this to understand light-like directions on $N^*K_{123}$. Maybe we should
  first find $S_3$ from $\e^3$ terms.

***

Plan: 

1. It seems likely that in generic case we can see $S_1$, $S_2$ and $S_3$ from different $\e$-order of singularities
The difficulty why we need to consider different $\e$ order is that we need
to consider  the product of distributions 
  $v_1\in \I(K_{13})$ and  $ v_1\in\I(K_{12})$ and we can not find the order of ${\bf Q}(v_1v_2)$.

2. Outside the surfaces  $S_1$, $S_2$ and $S_3$ we can detect $S_4$ and can separated different components.
For this we use the product theorem for the wave front sets of products
of   $v_1\in \I(K_{13})$ and  $ v_1\in\I(K_{12})$ to see that the wave front set of ${\bf Q}(v_1v_2)$ is contained
in $\Lambda^{(3)}$. 

***

and thus $u_{a_1}\,\cdotp u_{a_1}\in \I^{p}(K_{a_1})$ has a non-vanishing principal symbol
at $(x,\xi)\in N^*K_{a_1}$ if and only $u_{a_1}$ has a non-vanishing principal symbol at $(x,\xi)$ and 
a non-vanishing value at $x$.

*** ABOVE: We had problem: We can not analyze terms $\M^{(2)}\cdotp \M^{(2)}$
appearing in $\M_4$. The problem can be solved by considering solutions
outside Hausdroff 2-dimensional set $\Z$.
}
\footnote{We could use Piriou,	A.: Calcul	symbolique non	 lineare pour une onde conormale simple. Ann. Inst. Four. 38, 173-188 (1988), or \cite{GU1}, Lemma 1.3. Look also
Bougrini, H.; Piriou, A.; Varenne, J. P. Propagation et interaction des symboles principaux pour les ondes conormales semi-lin'eaires. (French) [Propagation and interaction of principal symbols for semilinear conormal waves] Comm. Partial Differential Equations 23 (1998), no. 1-2, 333--370.}

\subsection{Separation of singularities from different point sources}
Summarizing, above we have seen that in the generic case  
the terms  $\M^{(j)}$, $j\leq 3$   have wave front sets which union is $S_1\cup S_2\cup S_3$ 
and,  $\M^{(4)}$  has wave front set is a subset of   $S_1\cup S_2\cup S_3\cup S_4$ where  $S_4=\bigcup_{a_1<a_2<a_3<a_4\leq N}\Lambda^+_{q(a_1,a_2,a_3,a_4)}$.
Now we notice that for a generic  (i.e.\ in an open and dense set) $((M,g),(x_0,\xi_0))\in \mathcal N$,
we see that when $\vec a$  and $\vec b$ run over any finite subset of  $\A$,
\HOX{Check the the claims on "generic manifolds"}
the sets $\Lambda^+_{q(b_1,b_2,b_3,b_4)}$ intersect
the sets $N^*K_a$, $N^*K_{a_1,a_2}$, $\Lambda_{a_1,a_2,a_3}^{(3)}$, and 
the other sets 
$\Lambda^+_{q(a_1,a_2,a_3,a_4)}$, 
on at most 3-dimensional sets. Indeed,
 by making a small perturbation to the metric in a small neighborhood of
  $q(a_1,a_2,a_3,a_4)$ the sets  $\Lambda^+_{q(a_1,a_2,a_3,a_4)}$ can be perturbed so
  that the above non-intersection condition becomes valid, and the non-intersection condition
is clearly valid in an open set of metric tensors.
Thus the set of manifolds for which the stated condition holds is  open and dense.
Let us then choose the sources $\tilde {\bf f}_{a_\ell}$ to have generic polarizations.
We see that  in a generic set of polarizations the principal symbol
 of the 4th order interaction term $\M^{(4)}$ is non-vanishing on an open and dense set
 of $S_4$. Recall also that   in a generic situation the singular support of $\M^{(3)}$
 determines $S_3$.

 Thus in a generic situation
we can observe an open and dense subset of the surface $S_4$. Let us now consider how
we can decide which points on the union  of the surfaces 
$\Lambda^+_{q(a_1,a_2,a_3,a_4)}$ belong to the same surface.
Recall that the order of singularity (at points that do not belong to the 3-dimensional intersections)
on $\Lambda^+_{q(a_1,a_2,a_3,a_4)}$ is $(p_{a_1}+p_{a_2}+p_{a_3}+p_{a_4})-4-3/2$, see (\ref{eq: conormal}).
Let us now assume that $p_j$, $j\leq N$ are such that all the numbers
\ba
& &(p_{a_1}+p_{a_2}+p_{a_3}+p_{a_4})-1/2-n_1,
\ea
are different,
where  $n_1$ and $n_2$ are arbitrary integers (corresponding to different orders of the
classical symbols of terms) and \HOX{Gunther, do you know if the symbols of of waves $\M_j$ classical, assuming
that the sources have classical symbols, and if so, what would be a reference for this?}
the indexes $a_1,a_2,a_3,a_4\leq N$
satisfy $a_j\not =a_k$ for $j\not=k$. Let $\Sigma$ be the set of points $(x,\xi)\in T^*U_g$, $x\not \in S_3$ on the 
wave front set of  $\M^{(4)}$ such that the point $x$
has a neighborhood $V$ where $\M^{(4)}$  is a conormal
distribution associated smooth surface $S\subset V$ such that
$\M^{(4)}|_V\in \I^p(N^*S)$, where $p=(p_{a_1}+p_{a_2}+p_{a_3}+p_{a_4})-4-3/2$
but $p$ can not be replaced by any smaller number. Then for
generic set of manifolds  in $\mathcal N$
the closure of the set $\Sigma$ coincides with 
 $\Lambda^+_{q(a_1,a_2,a_3,a_4)}$. \HOX{Check the construction of  $\Lambda^+_{q(a_1,a_2,a_3,a_4)}$}
This implies
that we can determine $\Lambda^+_{q(a_1,a_2,a_3,a_4)}\cap T^*U_g$ for  all
$\vec a=(a_1,a_2,a_3,a_4)$ and thus
we can find the whole smooth surface $\Lambda^+_{q(a_1,a_2,a_3,a_4)}\cap T^*U_g$.
This implies that we can find  
$G_N(J((M,g),(x_0,\xi_0)),[g],x_0,\xi_0)$. This determines the manifold
$(J((M,g),(x_0,\xi_0)),[g])$ up to a small error.

   Summarizing the above  sketch of construction
   means the following:
   \medskip
   
   {\bf Summary:} {\it For any $\e>0$ one can $L$ such that the following holds
   in a generic set of manifolds
    $((M,g),(x_0,\xi_0))\in \mathcal N$.
    Let   $(y_{a_\ell},\xi_{a_\ell})$, $|\ell |\leq L$ be constructed above and
    $F_\e$ be a source (\ref{Fe source}) that sends spherical waves to directions
       $(y_{a_\ell},\xi_{a_\ell})$, $|\ell |\leq L $ with generic polarizations  $w_{a_\ell}$  and
       satisfies (\ref{eq: V ehto}). Then,
       by measuring singularities of terms 
   $\M^{(j)}$, $j\leq 4$ of the     
       the wave  $u^{(4)}_\e(x)$ on $U_g$
    produced by the source  $F_\e$,
   we can find the manifold  
    $J_{g}(\hat \mu(s_-),\hat \mu(s_+))$ and the conformal type 
    of the metric $[g]$ on it up to an error $\e$.}
    \medskip
   
 In other words,  when we consider the above measurement
 with single source and approximate the $O(\e^5)$ terms by zero,
 the observations of the singularities of the terms with different $\e^j$ magnitudes $j\leq 4$, 
 make it possible in a generic case to construct
a  discrete approximation of the light observation point
sets, $\{\mathcal P_V(q_k);\ k\leq K\}.$
Such data, under appropriate conditions, could be used to determine a discrete approximation 
for the conformal structure of the space-time in an analogous manner used
in \cite{AKKLT,K2L}. However, making the approach described above to a detailed proof is outside
the scope of this paper  will be considered
elsewhere.}
\medskip
}}


\subsection{A discussion on an application for a dark matter related example}


%

The determination of the conformal class of the Lorentzian metric considered above  
can be done also for a model that is related to dark matter and energy
\cite{darkenergy}. 
We note that by reconstructing the conformal class of the metric tensor $g$ 
in
area of space that contain dark matter but not usual "observable matter" tells how the 
"dark matter" would change the path of light rays
that would travel in this area, even the path of light rays on which
we can not do direct measurements.
 
Let us consider a model
where the fields $\phi_\ell$, $\ell\leq L-1$ can be observed  and correspond to "usual" matter. The field $\phi_L$ could correspond
to "dark" matter.
We write $\phi=(\phi^\prime,\phi_L)$, where $\phi^\prime=(\phi_\ell)_{\ell=1}^{L-1}$
and $Q=(Q^\prime,Q_L)$, where $Q^\prime=(Q_\ell)_{\ell=1}^{L-1}$.

Moreover, we assume that  in the model (\ref{eq: adaptive model})  the  adaptive source functions 
${\mathcal S}_\ell(\phi, \nabla^g
 \phi,Q,\nabla^g Q,P,\nabla^g P,g)$
are such that
\beq\label{eq: adaptive model dark model}
& & {\mathcal S}_\ell(\phi, \nabla^g
 \phi,Q,\nabla^g Q,P,\nabla^g P,g)=
\tilde  {\mathcal S}_\ell(\phi^\prime , \nabla^g
 \phi^\prime ,Q^\prime ,\nabla^g Q^\prime ,P,\nabla^g P,g),\hspace{-1cm}
 \\& & \nonumber\quad\hbox{for }\ell\leq L-1,\\
 &  & \nonumber
 {\mathcal S}_L(\phi, \nabla^g
 \phi,Q,\nabla^g Q,P,\nabla^g P,g)=0.
 \eeq
\generalizations{Also, we assume that the terms (\ref{B-interaction}) modeling interaction of matter fields,
are of the form
 \ba
B_\ell(\phi)=\sum_{\kappa=1}^{L-1} a^\kappa_\ell\phi_\kappa(x)
+ \sum_{\kappa,\alpha,\beta =1}^{L-1}b^{\kappa,\a,\beta}_\ell \phi_\kappa (x)
\phi_\alpha (x)\phi_\beta (x)\hspace{-1cm}\hspace{-1cm}
\ea  
for $\ell\leq L-1$ and $B_L(\phi)=c_L \phi_L(x)^3$.}
We consider the model (\ref{eq: adaptive model})
assuming that  $Q_L=0$ and that we
 can observe only the $g$ and $\phi^\prime$ components of the waves in $U_g$.
 Moreover, we assume that the Assumption $B$ is valid with permutations $\sigma$
for which $L\not \in \{\sigma(j);\ j=1,2,3,4,5\}$.
Analyzing the wave equation in the $\hat g$ wave map coordinates, we see
that the $\phi_L$ component of linearized waves $(\dot g,\dot \phi)$ satisfies
 \beq\label{eq: lin u}
 (\hat g^{jk}\p_j\p_k-\hat g^{pq}\hat \Gamma_{pq}^j\p_j+m)\dot \phi_L=
-(\dot g^{jk}\p_j\p_k+\dot g^{pq}\hat \Gamma_{pq}^j\p_j)\hat \phi_L.
\eeq
 Thus, using the notations
introduced earlier in the paper,
we consider $(\vec x,\vec \xi)$ such that $\gamma_{x_j,\xi_j}([t_0,\infty))$
intersect at a point $q$ and that there is a light-like geodesic from $q$
to the point $y\in U_{\hat g}$. Assume
that $K_j\subset \L_{\hat g}^+(x_j)$  are such that  $\gamma_{x_j,\xi_j}([t_0,\infty))\subset K_j$ and consider 
linearized waves  $\dot u^{(j)}=(\dot g^{(j)},\dot \phi^{(j)})\in \I(K_j)$, $j=1,2,3,4$.
We showed earlier that in a generic case the interaction of the waves at $q$ can
be observed in  normal coordinates at $y$ at least in two polarizations,
that is, in two components of $\dot u^{(j)}$.
If we do not observe the $ \phi_L$ component of the waves,
then the interaction of the waves can
be observed in  normal coordinates at $y$ at least in one polarization.
We also showed that there are some principal symbols
of the waves $\dot u^{(j)}$ at $q$ that produce such observable singularities
and these principal symbols are such that the $\phi$-components
of the principal symbols are zero. In particular,
the $\phi_L$-component the principal symbols of  $\dot u^{(j)}$ is zero at $(q,\eta)\in N^*K_j$.
As the linearized equation (\ref{eq: lin u}) contains no derivatives of $\dot g$, we see then that
the $\phi_L$-component the principal symbols of  $\dot u^{(j)}$  is zero 
also at $(x_j,\xi_j^\flat)$ that is on the same bicharacteristic as $(q,\eta)$. A linearized wave  $\dot u^{(j)}$ with such
principal symbol can be produced
with sources $Q$ for which $Q_L=0$, that is, without source terms in the "dark" matter
component. 
Because of the above considerations, we see that 
 with a slight modification of the proof, 
we can find the conformal type of the metric.

\section*{Appendix A: Reduced Einstein equation}

In this section we review known results on Einstein equations and wave maps.

\HOX{The appendices need to be shortened  in the final version of the paper.}

\subsection*{A.1.\ Summary of the used notations}


Let us recall some definitions given in Introduction, in the Subsection
\ref{subsec: Gloabal hyperbolicity}.
Let $(\hattuM ,\hat g)$ be a $C^\infty$-smooth globally hyperbolic Lorentzian
manifold  and  $\tilde g$ be a $C^\infty$-smooth globally hyperbolic
metric on $M$ such that $\hat g<\tilde g$. 

Recall that there is an isometry $\Phi:(M,\tilde g)\to (\R\times N,\tilde  h)$,
where $N$
is a 3-dimensional manifold and the metric $\tilde  h$ can be written as
$\tilde  h=-\beta(t,y) dt ^2+\overline h(t,y)$ where $\beta:\R\times  N\to (0,\infty)$ is a smooth function and 
$\overline h(t,\cdotp)$ is a Riemannian metric on $ N$ depending smoothly
on $t\in \R$. 
As in the main text identify these isometric manifolds and denote $\hattuM =\R\times N$.
%
Also, for $t\in \R$, recall that $\hattuM (t)=(-\infty,t)\times N$. We use parameters $t_1>t_0>0$
and denote $\hattuM _j=\hattuM (t_j)$, $j\in \{0,1\}$.
We use  the time-like geodesic    $\hat \mu =\mu_{\hat g}$, $\mu_{\hat g}:[-1,1]\to M_0$ on $(M_0,\hat g)$
and  the set $\K_j:=J^+_{\tilde g}(\hat p^-)\cap M_j$ with $\hat p^-=
\hat \mu(s_-)\in (0,t_0)\times N$, $s_-\in (-1,1)$
and recall that  
$J^+_{\tilde g}(\hat p^-)\cap M_j$ is compact and 
there exists $\e_0>0$ such that if $\|g-\hat g\|_{C^0_b(\hattuM _1;\hat g^+)}<\e_0$, then
$g|_{\K_1}<\tilde g|_{\K_1}$,
and in particular, we have  $J^+_{g}(p)\cap \hattuM _1\subset \K_1$
for all $p\in \K_1$.

Let us use local coordinates on $\hattuM _1$ and denote by
$\nabla_k=\nabla_{X_k}$  the covariant derivative with respect to the metric $g$
to the direction $X_k=\frac \p{\p x^p}$
 and by
 $\hat \nabla_k=\hat \nabla_{X_k}$ the covariant derivative with respect to the metric $\hat g$
to the direction $X_k$.

\subsection*{A.2.\ Reduced Ricci and Einstein tensors}
Following \cite{FM} we recall that 
\beq\label{q-formula2copy}
\Ric_{\mu\nu}(g)&=& \Ric_{\mu\nu}^{(h)}(g)
+\frac 12 (g_{\mu q}\frac{\p \Gamma^q}{\p x^{\nu}}+g_{\nu q}\frac{\p \Gamma^q}{\p x^{\mu}})
\eeq
where  $\Gamma^q=g^{mn}\Gamma^q_{mn}$,
\beq\label{q-formula2copyB}
& &\hspace{-1cm}\Ric_{\mu\nu}^{(h)}(g)=
-\frac 12 g^{pq}\frac{\p^2 g_{\mu\nu}}{\p x^p\p x^q}+ P_{\mu\nu},
\\ \nonumber
& &\hspace{-2cm}P_{\mu\nu}=  
g^{ab}g_{ps}\Gamma^p_{\mu b} \Gamma^s_{\nu a}+ 
\frac 12(\frac{\p g_{\mu\nu }}{\p x^a}\Gamma^a  
+ \nonumber
g_{\nu l}  \Gamma^l _{ab}g^{a q}g^{bd}  \frac{\p g_{qd}}{\p x^\mu}+
g_{\mu l} \Gamma^l _{ab}g^{a q}g^{bd}  \frac{\p g_{qd}}{\p x^\nu}).\hspace{-2cm}
\eeq
Note that $P_{\mu\nu}$ is a polynomial of $g_{jk}$ and $g^{jk}$ and first derivatives of $g_{jk}$.
The harmonic  Einstein tensor is  
\beq\label{harmonic Ein}
\Ein^{(h)}_{jk}(g)=
\Ric_{jk}^{(h)}(g)-\frac 12 g^{pq}\Ric_{pq}^{(h)}(g)\, g_{jk}.
\eeq 
The  harmonic  Einstein tensor is extensively used to study
Einstein equations in local coordinates where one can use the Minkowski
space $\R^4$ as the background space. To do global constructions
with a background space $(M,\hat g)$ one uses  the  reduced Einstein tensor.
The $\hat g$-reduced Einstein tensor 
 $\Ein_{\hat g} (g)$ and the $\hat g$-reduced  Ricci tensor
  $\Ric_{\hat g} (g)$  
are given 
by
\beq\label{Reduced Einstein tensor}
& &(\Ein_{\hat g} (g))_{pq}=(\Ric_{\hat g} (g))_{pq}-\frac 12 (g^{jk}(\Ric_{\hat g} g)_{jk})g_{pq},\\& &
\label{Reduced Ric tensor}
(\Ric_{\hat g} (g))_{pq}=\Ric_{pq} g-\frac 12 (g_{pn} \hat \nabla _q\hat F^n+ g_{qn} \hat \nabla _p\hat F^n)
\eeq
where $\hat F^n$ are the harmonicity functions given by
\beq\label{Harmonicity condition}
\hat F^n=\Gamma^n-\hat \Gamma^n,\quad\hbox{where }
\Gamma^n=g^{jk}\Gamma^n_{jk},\quad
\hat \Gamma^n=g^{jk}\hat \Gamma^n_{jk},
\eeq
where $\Gamma^n_{jk}$ and $\hat \Gamma^n_{jk}$ are the Christoffel symbols
for $g$ and $\hat g$, correspondingly.
Note that $\hat \Gamma^n$ depends also on $g^{jk}$.
As $\Gamma^n_{jk}-\hat \Gamma^n_{jk}$ is the difference of two connection
coefficients, it is a tensor. Thus $\hat F^n$ is tensor (actually, a vector field), implying that both
$(\Ric_{\hat g} (g))_{jk}$ and $(\Ein_{\hat g} (g))_{jk}$ are  2-covariant tensors.
Observe that  the $\hat g$-reduced Einstein tensor is sum of the harmonic Einstein tensor 
and a term that is a zeroth order in $g$, 
\beq\label{hat g reduced einstein and reduced einstein}
(\Ein_{\hat g} (g))_{\mu\nu}=\Ein^{(h)}_{\mu\nu} (g)+\frac 12 (g_{\mu q}\frac{\p \hat \Gamma^q}{\p x^{\nu}}+g_{\nu q}\frac{\p \hat \Gamma^q}{\p x^{\mu}}).
\eeq

%

%

%

\subsection*{A.3.\ Wave maps and reduced Einstein equations}
Let us consider the manifold $M_1=(-\infty,t_1)\times N$ with a $C^m $-smooth metric $g^{\prime}$,
$m\geq 8$,
 which
is a perturbation of the metric $\hat g$ and satisfies the Einstein equation
\beq\label{Einstein on M_0}
\Ein(g^{\prime})=T^{\prime}\quad \hbox{on }M_1,
\eeq
or equivalently, 
\ba
\Ric(g^{\prime})=\rho^{\prime},\quad \rho^{\prime}_{jk}=T_{jk}^{\prime}-\frac 12 ((g^{\prime})^{  nm}T^{\prime}_{nm})g^{\prime}_{jk}\quad \hbox{on }M_1.
\ea
Assume also that  $g^{\prime}=\hat g$ in the domain $A$,
where $A=M_1\setminus\K_1$ and  $\|g^{\prime}-\hat g\|_{C^2_b(M_1,\hat g^+)}<\e_0$,
so that $(M_1,g^{\prime})$ is globally hyperbolic. Note that then $T^{\prime}=\hat T$ in the set $A$ 
and that the metric $g^{\prime}$ coincides with $\hat g$
in particular in the set $M^-=\R_-\times N$

We recall next the considerations of \cite{ChBook}.
Let us  consider the Cauchy problem for the wave map
$f:(M_1,g^\prime)\to (M,\hat g)$, namely
\beq
\label{C-problem 1}& &\square_{g^{\prime},\hat g} f=0\quad\hbox{in } M_1,\\
\label{C-problem 2}& &f=Id,\quad \hbox{in  }\R_-\times N,
\eeq
where $M_1=(-\infty,t_1)\times N\subset M$. In (\ref{C-problem 1}), 
$\square_{g^{\prime},\hat g} f=g^\prime\,\cdotp\hat \nabla^2 f$ is the wave map operator, where $\hat \nabla$
is the covariant derivative of a map $(M_1,g^\prime)\to (M,\hat g)$, see \cite[formula (7.32)]{ChBook}.
In  local coordinates
$X:V\to \R^4$ of $V\subset M_1$, denoted
$X(z)=(x^j(z))_{j=1}^4$ and $Y:W\to \R^4$ of $W\subset \hattuM $, denoted
$Y(z)=(y^A(z))_{A=1}^4$, 
the wave map $f:M_1\to \hattuM $ has representation $Y(f(X^{-1}(x)))=(f^A(x))_{A=1}^4$ and
the wave map operator in equation (\ref{C-problem 1}) is given by
\beq\label{wave maps2}
& &(\square_{g^{\prime},\hat g} f)^A(x)=
(g^{\prime})^{  jk}(x)\bigg(\frac \p{\p x^j}\frac \p{\p x^k} f^A(x) -\Gamma^{{\prime} n}_{jk}(x)\frac \p{\p x^n}f^A(x)
\\
& &\quad \quad\quad\quad\quad\quad\quad\quad\nonumber
+\hat \Gamma^A_{BC}(f(x))
\,\frac \p{\p x^j}f^B(x)\,\frac \p{\p x^k}f^C(x)\bigg)\eeq
where $\hat \Gamma^A_{BC}$ denotes the Christoffel symbols of metric $\hat g$
and  $\Gamma^{ {\prime} j}_{kl}$ are the Christoffel symbols of metric $g^{\prime}$.
When (\ref{C-problem 1})
is satisfied, we say that  $f$ is wave map with respect to the pair $(g^{\prime},\hat g)$.
The important property of the wave maps is that, if 
$f$ is wave map with respect to the pair $(g^{\prime},\hat g)$ and
$g=f_*g^{\prime}$ 
then, as follows from (\ref{wave maps2}), the identity map $Id:x\mapsto x$ is a wave map with respect to the pair $(g,\hat g)$
and, the wave map equation for the identity map is equivalent to (cf.\ \cite[p.\ 162]{ChBook})
\beq\label{wave equation id}
\Gamma^n=\hat\Gamma^n,\quad \hbox{where }\Gamma^n=g^{jk}\Gamma^n_{jk},
\quad \hat \Gamma^n=g^{jk}\hat\Gamma^n_{jk}
\eeq
%
where the Christoffel symbols $\hat\Gamma^n_{jk}$ of the metric $\hat g$ are smooth functions.

As  $g=g^{\prime}$ outside a compact set $\K_1\subset (0,t_1)\times N$,
we see that   this Cauchy problem is equivalent to the same equation
restricted to the
set $(-\infty,t_1)\times B_0$, where $B_0\subset N$ is such an open
relatively compact set that $\K_1\subset (0,t_1]\times B_0$
with the boundary condition $f=Id$ on  $(0,t_1]\times \p B_0$.
Then using results of \cite {HKM},  that can be applied for equations on manifold  as is
done  in Appendix C, and combined with the Sobolev embedding theorem, 
we see\footnote{See also:  
Thm.\ 4.2 in App.\ III  of
of \cite{ ChBook}, and its proof for the estimates
for the time on which the solution exists.} that 
 there is $\e_1>0$ such that if $\|g^\prime-\hat g\|_{C^{m}_b(\hattuM _1;\hat g^+)}<\e_1$,  then there  is a map  $f:M_1 \to M$ satisfying
the Cauchy problem (\ref{C-problem 1})-(\ref{C-problem 2})
and the solution depends continuously, 
 in
$C^{m-5}_b([0,t_1]\times N,g^+)$, on the metric $g^{\prime}$.
 Moreover, by the uniqueness of the wave map, we have
$f|_{M_1\setminus \K_1}=id$ so that
 $f(\K_1)\cap M_0\subset \K_0$.

 As
 the inverse function of the wave map $f$ depends continuously,
 in
$C^{m-5}_b([0,t_1]\times N,g^+)$, on the metric $g^{\prime}$
we can also assume that $\e_1$ is  so small that
$\hattuM _0\subset f(M_1)$.

Denote next $g:=f_*g^{\prime}$, $T:=f_*T^{\prime}$, and $\rho:=f_*\rho^{\prime}$ 
and define $\hat \rho=\hat T-\frac 12 (\hbox{Tr}\, \hat T)\hat g.$ 
Then $g$ is $C^{m-6}$-smooth and the equation (\ref{Einstein on M_0}) implies 
\beq\label{Einstein on M1}
\Ein(g)=T\quad \hbox{on }\hattuM _0.
\eeq 
As $f$ is a wave map, $g$ satisfies (\ref{wave equation id}) and thus  by the
definition of the reduced Einstein tensor,  (\ref{Reduced Einstein tensor}), we have
\ba
\Ein_{pq}(g)=
(\Ein_{\hat g} (g))_{pq}\quad \hbox{on } \hattuM _0.
\ea
This and  (\ref{Einstein on M1}) yield the $\hat g$-reduced Einstein equation
\beq\label{hat g reduced einstein equations}
(\Ein_{\hat g} (g))_{pq}=T_{pq}\quad \hbox{on } \hattuM _0.
\eeq
This equation is useful for our considerations as it is a quasilinear,
hyperbolic equation on $\hattuM _0$. Recall that 
 $g$ coincides with
$\hat g$ in $M_0\setminus \K_0$. The unique solvability of this 
Cauchy problem is studied in e.g.\ \cite[Thm.\ 4.6 and 4.13]{ ChBook} and \cite {HKM} 
and in Appendix C below.


\subsection*{A.4.\ Relation of the reduced Einstein equations and for the original Einstein equation}

The metric $g$ which solves the $\hat g$-reduced Einstein
equation $\Ein_{\hat g} (g)=T$ is a solution of the original
Einstein equation $\Ein (g)=T$ if the harmonicity
functions $\hat F^n$
vanish identically. Next we recall the result that 
the harmonicity functions vanish on $\hattuM _0$
when 
\beq\label{eq: good system}
& &(\Ein_{\hat g}(g))_{jk}=T_{jk},\quad \hbox{on }M_0,\\
\nonumber & &\nabla_pT^{pq}=0,\quad \hbox{on }M_0,\\
\nonumber & &g=\hat g,\quad \hbox{on }M_0\setminus \K_0.
\eeq
To see this, let
us  denote $\Ein_{jk}(g)=S_{jk}$, $S^{jk}=g^{jn}g^{km}S_{nm}$,
and $T^{jk}=g^{jn}g^{km}T_{nm}$. 
Following the standard arguments,
see \cite{ ChBook}, we see from (\ref{Reduced Einstein tensor}) that in local coordinates
\ba
S_{jk}-(\Ein_{\hat g}(g))_{jk}=\frac12(g_{jn}\hat \nabla_k \hat F^n+g_{kn}\hat \nabla_j\hat F^n-g_{jk}\hat \nabla_n \hat  F^n).
\ea
Using equations (\ref{eq: good system}), the Bianchi identity  $\nabla_pS^{pq}=0$,
and the basic property of Lorentzian connection,
$\nabla_kg^{nm}=0$,
we obtain 
\ba
0&=&
2\nabla_p(S^{pq}-T^{pq})\\
&=&\nabla_p(g^{qk}\hat  \nabla_k F^p+g^{pm}\hat\nabla_m \hat F^q-  g^{pq}\hat  \nabla_n \hat F^n)
\\
&=&g^{pm}\nabla_p \hat \nabla_m \hat F^q+(g^{qk} \nabla_p \hat \nabla_k \hat F^p-  g^{qp} \nabla_p \hat \nabla_n \hat F^n)\\
&=&g^{pm}\nabla_p \hat \nabla_m \hat F^q+W^q(\hat F)
\ea
where $\hat F=(\hat F^q)_{q=1}^4$ and
the operator 
\ba
W:(\hat F^q)_{q=1}^4\mapsto (g^{qk} (\nabla_p \hat \nabla_k \hat F^p- \nabla_k \hat \nabla_p \hat F^p))_{q=1}^4
\ea
is a linear first order differential operator which coefficients are polynomial functions of
$\hat g_{jk}$, $\hat g^{jk}$, $g_{jk}$, $g^{jk}$ and their first derivatives.

Thus the harmonicity functions $\hat F^q$ satisfy on $\hattuM _0$ the hyperbolic initial 
value problem 
\ba
& &g^{pm}\nabla_p \hat \nabla_m \hat F^q+W^q(\hat F)=0,\quad\hbox{on }M_0,\\
& & \hat F^q=0,\quad \hbox{on }M_0\setminus \K_0,
\ea
and as this initial 
Cauchy problem is uniquely solved by \cite[Thm.\ 4.6 and 4.13]{ ChBook}
or \cite {HKM}, we see that $\hat F^q=0$ on $\hattuM _0$. Thus
equations (\ref{eq: good system})
yield that
Einstein equation  $\Ein(g)=T$ holds on $\hattuM _0$.

\section*{Appendix B: 
An example  satisfying Assumption S}
 
Next we give an example of functions ${\mathcal S}_\ell(\phi, \nabla
 \phi,Q^{\prime},\nabla^g P,Q_K,\nabla Q_K,g)$ in the model (\ref{eq: adaptive model})
 for which Assumption S is valid.

 Let $L\geq 5$, $g$ be a $C^2$-smooth metric and $\phi=(\phi_\ell)_{\ell=1}^L$
 be $C^2$-smooth functions on $\hat U\subset M$.
   Let us fix 
  a symmetric (0,2)-tensor $P$ and a scalar function $Z$
  that are $C^2$-smooth and compactly supported in  $\hat U$.
Let $[P_{jk}(x)]_{j,k=1}^4$ be  the coefficients of 
  $P$ in local coordinates. 
 

To obtain adaptive source functions satisfying 
 the assumption S, 
 let us start implications of the conservation law (\ref {conservation law0}).
To this end, consider $C^2$-smooth functions 
 $S_\ell(x)$ on $\hat U$.
  The conservation law (\ref {conservation law0})
 gives for all $j=1,2,3,4$ equations (see \cite[Sect. 6.4.1]{ChBook}) 
 \ba
0&=&\frac 12\nabla_p^g (g^{pk}T_{jk})\\
&=& \sum_{\ell=1}^L(g^{pk} \nabla^g_p\p_k\phi_\ell )\,\p_j\phi_\ell 
-(m_\ell \phi_\ell \p_p \phi_\ell) \delta_j^p +\frac 12 \nabla^g_p (g^{pk}g_{jk}S_\ell\phi_\ell+g^{pk}P_{jk})
\\
&=& \sum_{\ell=1}^L(g^{pk} \nabla^g_p\p_k\phi_\ell -m_\ell \phi_\ell ) \,\p_j \phi_\ell +\frac 12 \nabla^g_p (g^{pk}g_{jk}S_\ell\phi_\ell+g^{pk}P_{jk})
\\
&=&\sum_{\ell=1}^L S_\ell \,\p_j\phi_\ell+\frac 12 \nabla^g_p (g^{pk}g_{jk}S_\ell \phi_\ell)+\frac 12 g^{pk} \nabla^g_p P_{jk}\\
&=&\left (\sum_{\ell=1}^L S_\ell \,\p_j\phi_\ell\right)+\frac 12\p_j\left(
\sum_{\ell=1}^L S_\ell \phi_\ell \right)+\frac 12 g^{pk}\nabla_p^g P_{jk}.
\ea
Recall that $S_\ell$ should satisfy
\beq\label{meson number}
\sum_{\ell=1}^L S_\ell \phi_\ell=Z.
\eeq

Then, the  conservation law  (\ref {conservation law0}) holds if we have
\beq\label{meson number2}& &
\sum_{\ell=1}^L S_\ell \,\p_j\phi_\ell
=-\frac 12 g^{pk}\nabla^g_p R_{jk},\quad
 R_{jk}=(P_{jk}+g_{jk}Z),
\eeq
for $j=1,2,3,4.$

Equations (\ref{meson number}) and (\ref{meson number2}) give together  five point-wise equations
 for the  functions $S_1,\dots,S_L$.

%
%
%
Next we denote the set of permutations $\sigma:\{1,2,\dots,L\}\to \{1,2,\dots,L\}$
by $\Sigma(L)$.
Next we assume Condition A, that is, that at any $x\in \hbox{cl}(\hat U)$ there is
a permutation $\sigma:\{1,2,\dots,L\}\to \{1,2,\dots,L\}$ such that the 
$5\times 5$ matrix $( B_{jk}^\sigma(\hat \phi(x),\nabla \hat \phi(x)))_{j,k\leq 5}$
is invertible, where
  \ba
( B_{jk}^\sigma(\phi(x),\nabla \phi(x)))_{j,k\leq 5}=\left(\begin{array}{c}
( \,\p_j  \phi_{\sigma(\ell)}(x))_{j\leq 4,\ \ell\leq 5}\\
(   \phi_{\sigma(\ell)}(x)) _{\ell\leq 5}\end{array}\right).
 \ea
Let $V_\sigma\subset \hbox{cl}(\hat U)$ be the set where $( B_{jk}^\sigma(\hat \phi(x),\nabla \hat \phi(x)))_{j,k\leq 5}$
is invertible.

Below, let us use  $K=L\,\cdotp(L!)+1$ and identify the set 
index set $\{1,2,\dots, K-1\}$ with the
the set $\Sigma(L)\times \{1,2,\dots,L\}$. We  consider
a $\R^K$ valued function  $Q(x)=(Q^\prime(x),Q_K(x))$,
where  $$Q^\prime=(Q_{\sigma,\ell})_{\ell\in \{1,2,\dots,L\},\ \sigma\in \Sigma(L)}.$$
 Note that we have introduced the following renumbering:
 identify the set 
index set $\{1,2,\dots,\tilde K-1\}$ with the
the set $\Sigma(L)\times \{1,2,\dots,L\}$ using a bijective  map $j\mapsto (\sigma(j),l_j)$.

Also, below $R_{jk}=P_{jk}+g_{jk}Z$ and  we set
\beq
& &Q_K=Z.
\label{S R Z equations 3bbb}
\eeq

Our next aim is to consider first a fixed permutation $\sigma$ and point $x\in V_\sigma$,
and construct scalar functions
${\mathcal S_{\sigma,\ell}}(Q^\prime,Q_K,R,g,\phi)$, $\ell=1,2,\dots,L$ that satisfy
\beq\label{S R Z equations 1bbb}
& &\sum_{\ell=1}^5 {\mathcal S_{\sigma,\ell}}(Q^\prime,Q_K,R,g,\phi)
 \,\p_j\phi_{\sigma(\ell)}=-\frac 12 g^{pk}\nabla^g_p 
R_{jk} -
\sum_{\ell=6}^{L} Q_{\sigma,\ell} \,\p_j\phi_{\sigma(\ell)},\hspace{-1cm} \\
& &\sum_{\ell=1}^5 {\mathcal S_{\sigma,\ell}}(Q^\prime,Q_K,R,g,\phi)\, \phi_{\sigma(\ell)}=Q_K-\sum_{\ell=6}^L Q_{\sigma,\ell} \phi_{\sigma(\ell)}.\hspace{-1cm}
\label{S R Z equations 2bbb}
\eeq
Recall that for $x\in V_\sigma$
 the matrix $\B^\sigma(\phi,\nabla \phi)=(B^\sigma_{jk}(\phi(x),\nabla \phi(x)))_{j,k=1}^5$ is invertible.

%
%
Let $(Y_\sigma(\phi,\nabla \phi))(x)=(B^\sigma(\phi,\nabla \phi))^{-1}$ for
$x\in V_\sigma$, and zero for $x\not \in V_\sigma$.
Then we define 
${\mathcal S_{\sigma,\ell}}={\mathcal S_{\sigma,\ell}}(Q^\prime,Q_K,R,g,\phi)$, $\ell=1,2,\dots,L,$ to be
\beq\label{S sigma formulas}
& &(S_{\sigma,\ell})_{\ell\leq 5}= Y(\phi,\nabla \phi) \left(\begin{array}{c}
(-\frac 12 g^{pk}\nabla^g_p 
R_{jk} -
\sum_{\ell=6}^L Q_{\sigma,\ell} \,\p_j\phi_{\sigma(\ell)})_{j\leq 4}
\\
Q_K-\sum_{\ell=6}^L Q_{\sigma,\ell}\phi_{\sigma(\ell)} \end{array}\right)
\hspace{-1cm}\\ \nonumber
& &  \quad   \quad \quad    \quad =Y(\phi,\nabla \phi) \left(\begin{array}{c}
(-\frac 12 g^{pk}\nabla^g_p 
R_{jk} 
\\
Q_{K}\end{array}\right)+
Q_{\sigma,\ell},
\hspace{-1cm}\\
& & \nonumber
(S_{\sigma,\ell})_{\ell\geq 6}= (Q_{\sigma,\ell})_{\ell\geq 6}.
\eeq
Observe that $Q_{\tilde \sigma,\ell}$, with $\tilde \sigma\not=\sigma$, do not appear 
in the formula (\ref{S sigma formulas}).

Above, $Z$ and $P$ were fixed. Let us now choose  $Q^\prime$ to
be  arbitrary $C^1$-functions. Recall that  $R_{jk}=P_{jk}+g_{jk}Z$.
Then, we see that if we denote (note that we later will change the 
meaning of symbols $ S_\ell$)
\ba
 Q_K&=&Z,\\
 S_{\sigma,\ell}&=&{\mathcal S}_{\sigma,\ell}(\phi, \nabla
 \phi,Q^{\prime},Q_{K},\nabla Q_{K},\nabla^g P,g),
 \ea
we see that 
 the equations
(\ref{meson number}) and (\ref{meson number2}) 
are satisfied in
$x\in V_\sigma$ when functions $S_\ell$ are replaced
by $S_{\sigma,\ell}$. Moreover,
we see that
 the derivative of 
$ {\mathcal S}_\sigma=( {\mathcal S}_{\sigma,\ell})_{\ell=1}^L$
with respect to $(Q^{\prime},Q_{K},R)$, that is,
 \beq\label{eq: deriv1}
D_{Q^{\prime}, Q_{K}, R} {\mathcal S}_\sigma(\hat \phi, \hat \nabla
 \hat \phi,Q^{\prime},Q_{K},R,\hat g):\R^{K+4}\to \R^L
 \eeq
 is surjective.




Let us next combine the above constructions that were done for a single $\sigma$. 
Let $\psi_\sigma\in C^\infty(\hbox{cl}(\hat U))$ be the partition of unity
such that $\supp(\psi_\sigma)\subset  V_\sigma$ and $\sum_{\sigma\in\Sigma(L)}\psi_\sigma(x)=1$
in $\hbox{cl}(\hat U)$.

We define
\ba
& &{\mathcal S}_\ell(\phi, \nabla
 \phi,Q^{\prime},Q_{K},\nabla Q_{K},\nabla^g P,g)=\\
 & &\sum_{\sigma\in\Sigma(L)}\psi_\sigma(x)\,
 {\mathcal S}_{\sigma,\ell}(\phi, \nabla
 \phi,Q^{\prime},Q_{K},\nabla Q_{K},\nabla^g P,g),
\ea

Let us next denote 
\ba
 Q_K&=&Z,\\
 S_\ell&=&{\mathcal S}_\ell(\phi, \nabla
 \phi,Q^{\prime},Q_{K},\nabla Q_{K},\nabla^g P,g).
 \ea
Then we see,  using the partition of unity, that 
the functions $ S_\ell$ satisfy the equations
(\ref{meson number}) and (\ref{meson number2})
for all $x\in \hat U$.

Using the fact that 
 $Q_{\tilde \sigma,\ell}$, with $\tilde \sigma\not=\sigma$, do not appear 
in the formula (\ref{S sigma formulas}) and
that derivatives (\ref{eq: deriv1}) are surjective,
 we see that in
 the derivative of 
$ {\mathcal S}=( {\mathcal S}_\ell)_{\ell=1}^L$, c.f.\ Assumption S,
with respect to $(Q^{\prime},Q_{K},R)$, that is,
 \ba
D_{Q^{\prime}, Q_{K}, R} {\mathcal S}(\hat \phi, \hat \nabla
 \hat \phi,Q^{\prime},Q_{K},R,\hat g):\R^{K+4}\to \R^L
 \ea
 is surjective\footnote{ 
We make the following note related 
the case considered in the main part of the paper when $Q,P,R\in \I^m(Y)$, where $Y$ is 2-dimensional sufrace,
and we  need to use the principal symbol ${\bf r}$ of $R$
as an independent variable:
Let us consider also a point $x_0\in \hattuM _0$ and $\eta$
be a light-like covector choose
coordinates so that $g=\diag(-1,1,1,1)$ and $\eta=(1,1,0,0).$
When $c=(c^k)_{k=1}^4\in \R^4$ and $P^{jk}=C^{jk}(x\,\cdotp \eta)^a_+$, where $C^{jk}$ is such a symmetric
matrix that
$C^{11}=c_1$, $
C^{12}=C^{21}=\frac  12 c_2$,
 $C^{13}=C^{31}=c_3$, and  $C^{14}=C^{41}=c_4$ and other $C^{jk}$ are zeros. 
Then we have
\ba
\nabla^g_j (C^{jk}(x\,\cdotp \eta)_+^a)&=&\eta_jP^{jk}=(C^{1k}+ C^{2k})(x\,\cdotp \eta)_+^{a-1}=c_k(x\,\cdotp \eta)_+^{a-1}.
\ea
As we can always choose coordinates
that $\hat g$ and $\eta$ have at a given point the above forms,
we can obtain arbitrary vector ${\bf r}$, as
principal symbol of $R$, by considering
$Q^{\prime},Q_{L+1},P\in \I^m(Y)$, where $Y$ is 2-submanifold, with principal symbols of
$\tilde {\bf p}$ and $\tilde {\bf z}$ satisfing equations corresponding to
equations
$\hat g^{jk}\hat \nabla_j ({\bf p}_{kn}+ {\bf z}\hat g_{kn})\in \I^m(Y)$, 
and the sub-principal symbols of ${\bf p}_{kn}$ and ${\bf z}$
vary arbitrarily.} at all $x\in $cl$(U_{\hat g})$.
 Hence (iii) in the Assumption S is valid.
 \HOX{Remove the footnote in the final version of the paper.}

\subsection*{Appendix C: Stability and existence of the  direct problem}

 {

Let us  start by explaining  how we can choose a  $C^\infty$-smooth metric $\tilde g$ 
such that $\hat g<\tilde g$ and  $(M,\tilde g)$ is globally hyperbolic:
When $v(x)$ the eigenvector corresponding to the negative eigenvalue
of $\hat g(x)$, we can choose a smooth, strictly positive function
$\eta:\hattuM \to \R_+$ such that $\tilde g^\prime:=\hat g-\eta v\otimes v<\tilde g$. Then
$(\hattuM ,\tilde g^\prime)$ is globally hyperbolic, $\tilde g^\prime$ is smooth
and $\hat g<\tilde g^\prime$. Thus we can replace $\tilde g$ by the smooth metric
 $\tilde g^\prime$ having the same properties that are required for $\tilde g$.
 
 Let us now return to consider existence and stability of the solutions of the
 Einstein-scalar field equations.
Let $t={\bf t}(x)$ be local time so that there is a diffeomorphism
$\Psi:M\to \R\times N$, $\Psi(x)=({\bf t}(x),Y(x))$, and
 $S(T)=\{x\in \hattuM _0; \ {\bf t}(x)=T\}$, $T\in \R$ are Cauchy surfaces.
 Let $t_0>0$. Next we identify $M$ and $\R\times N$ via the map $\Psi$ and just
 denote $M=\R\times N$.
Let us denote $M(t)=(-\infty,t)\times N$,  and $M_0=M(t_0)$.
By \cite[Cor.\ A.5.4]{BGP} the set $\K=J^+_{\tilde g}(\hat p^-)\cap \hattuM _0$,
where $\hat p^0\in M_0$,
is compact.
Let $N_1,N_2\subset N$ be such
 open relatively compact sets with smooth boundary that
 $N_1\subset N_2$ and  $Y(J_{\tilde g}^+(\hat p^0)\cap \hattuM _0)\subset N_1$.  \HOX{We could improve explanation on $\tilde N$}
 
 To simplify citations to existing literature, 
let us define $\tilde N$ to be a compact manifold without boundary such
that $N_2$ can be considered as a subset of $\tilde N$. Using a construction based on
a suitable partition of unity, the Hopf double of the manifold $N_2$, 
and the Seeley extension of the metric tensor,  we can
endow $\tilde M=(-\infty,t_0)\times \tilde N$ with a smooth Lorentzian
metric $\hat g^e$ (index $e$ is for "extended")
so that $\{t\}\times \tilde N$ are Cauchy surfaces of $\tilde N$
and that $\hat g$ and $\hat g^e$ coincide in  the set $\Psi^{-1}((-\infty,t_0)\times N_1)$ that contains the set $J^+_{\hat g}(\hat p^0)\cap \hattuM _0$.} We extend the metric $\tilde g$ to a (possibly non-smooth) globally
hyperbolic metric $\tilde g^e$ on $\tilde M_0=(-\infty,t_0)\times \tilde N)$
such that $\hat g^e<\tilde g^e$.

 To simplify notations below we denote $\hat g^e=\hat g$ and
 $\tilde g^e=\tilde g$ 
 on the whole
 $\tilde M_0$.
 Our aim is to prove the estimate (\ref{eq: Lip estim}).

Let us denote by $t={\bf t}(x)$ the local time. 
Recall that when $(g,\phi)$ is a solution of 
the scalar field-Einstein equation, we denote $u=(g-\hat g,\phi-\hat \phi)$.
We will consider the equation for $u$, and to emphasize that the metric
depends on $u$, we denote $g=g(u)$ and assume below that  both the metric $g$ is  dominated
by $\tilde g$, that is, $ g<\tilde g$.
We use the pairs ${\bf u}(t)=(u(t,\cdotp),\p_t u(t,\cdotp))\in
H^1(\tilde N)\times L^2(\tilde N)$ and the notations
 ${\bf v}(t)=(v(t),\p_t v(t))$ etc. 
Let us consider a generalization of the system (\ref{eq: notation for hyperbolic system 1}) of the form 
\beq\label{eq: notation for hyperbolic system}
& &\square_{g(u)}u+V(x,D)u+H(u,\p u) =R(x,u,\p u)F+K,\ \ x\in \tilde M_0,\hspace{-1cm}\\
& & \nonumber \supp(u)\subset \K,
\eeq
where $\square_{g(u)}$ is the Lorentzian Laplace operator operating on the sections 
of the bundle $\B^L$ on $M_0$ and 
$\supp(F)\cup \supp(K)\subset \K$.
Note that above  $u=(g-\hat g,\phi-\hat \phi)$ and $g(u)=g$.
Also, $F\mapsto R(x,u,\p u)F$ is a 
 linear first order differential operator which coefficients at $x$ are depending smoothly on 
 $u(x)$, $\p_j u(x)$ and the  derivatives 
of $(\hat g,\hat \phi)$ at $x$ and
\ba V(x,D)=V^j(x)\p_j+V(x)\ea is a linear first order differential operator which coefficients 
at $x$ are depending smoothly on the  derivatives 
of $(\hat g,\hat \phi)$  at $x$, and finally, $H(u,\p u)$ is a polynomial
of $u(x)$ and $\p_j u(x)$ which coefficients 
at $x$ are depending smoothly on the  derivatives 
of $(\hat g,\hat \phi)$   such that
 $\p^\a_v\p^\beta_w H(v,w)|_{v=0,w=0}=0$ for
$|\a|+|\b|\leq 1$. 
By \cite[Lemma 9.7]{Ringstrom}, the equation  (\ref{eq: notation for hyperbolic system})
has at most one solution with given  $C^2$-smooth source functions $F$ and $K$.
Next we consider the existence of $u$ and its dependency on $F$ and $K$.


%

Below we use notations, c.f. (\ref{eq: notation for hyperbolic system 1}) and
(\ref{eq: notation for hyperbolic system 1b}) \ba
{\mathcal R}({\bf u},F)=R(x,u(x),\p u(x))F(x),\quad
{\mathcal H}({\bf u})=H(u(x),\p u(x)).\ea
Note that $u=0$, i.e., $g=\hat g$ and $\phi=\hat\phi$ satisfies (\ref{eq: notation for hyperbolic system}) with $F=0$
and $K=0$.
Let us use the same notations as in  \cite {HKM} cf. also \cite[section 16]{Kato1975}, 
 to consider quasilinear wave equation on $[0,t_0]\times \tilde N$.
Let $\H^{(s)}(\tilde N)=H^s(\tilde N)\times H^{s-1}(\tilde N)$ and 
\ba
Z=\H^{(1)}(\tilde N),\quad  Y=\H^{(k+1)}(\tilde N),\quad X=\H^{(k)}(\tilde N).
\ea
The norms on these space are defined invariantly using the smooth Riemannian 
metric $h=\hat g|_{\{0\}\times \tilde N}$ on $\tilde N$. \HOX{Maybe we should explain the 
bundles, the connection and the wave operator on section in a detailed way.}
Note that $\H^{(s)}(\tilde N)$ are in fact the Sobolev spaces of sections on the bundle $\pi:\B_K\to
\tilde N$, where $\B_K$ denotes also the  pull back bundle of $\B_K$ on $\tilde M$ in the map
$id:\{0\}\times \tilde N\to \tilde M_0$,
or on the bundle $\pi:\B_L\to \tilde N$. Below, $\nabla_h$ denotes the  standard connection of the bundle  $\B_K$
or  $\B_L$ associated to the metric $h$.

Let $k\geq 4$ be an even integer.
By definition of $H$ and $R$ we see that  there are  $0<r_0<1$ and $L_1,L_2>0$,
all depending on $\hat g$, $\hat \phi$, $\K$, and $t_0$, such that if 
$0<r\leq r_0$ and
\beq\label{new basic conditions}
& &\|{\bf v} \|_{C([0,t_0];\H^{(k+1)}(\tilde N))}\leq r,\quad \|{\bf v}^\prime\|_{C([0,t_0];\H^{(k+1)}(\tilde N))}\leq r,\\
& &\|F\|_{C([0,t_0];H^{(k+1)}(\tilde N))}\leq r^2,\quad \|K\|_{C([0,t_0];H^{(k+1)}(\tilde N))}\leq r^2 \nonumber\\
& &\|F^\prime\|_{C([0,t_0];H^{(k+1)}(\tilde N))}\leq r^2,\quad \|K^\prime\|_{C([0,t_0];H^{(k+1)}(\tilde N))}\leq r^2 \nonumber
\eeq
then
\beq\label{new f1f2 conditions}
& &\quad \quad \|g(\cdotp;v)^{-1}\|_{C([0,t_0];H^s(\tilde N))}\leq L_1 ,\\ \nonumber
& &\|{\mathcal H} ({\bf v})\|_{C([0,t_0];H^{s-1}(\tilde N))}\leq L_2r^2,
\quad \|{\mathcal H} ({\bf v}^\prime)\|_{C([0,t_0];H^{s-1}(\tilde N))}\leq L_2r^2,\\ \nonumber
& &\|{\mathcal H} ({\bf v})-{\mathcal H} ({\bf v}^\prime)\|_{C([0,t_0];H^{s-1}(\tilde N))}\leq  L_2r\, \|{\bf v}-{\bf v}^\prime\|_{C([0,t_0];\H^{(s)}(\tilde N))}, 
\\ \nonumber
& &\|\mathcal R ({\bf v}^\prime,F^\prime)\|_{C([0,t_0];H^{s-1}(\tilde N))}\leq L_2r^2,
\quad \|\mathcal R ({\bf v},F)\|_{C([0,t_0];H^{s-1}(\tilde N))}\leq L_2r^2,\hspace{-1cm} \\ \nonumber
& &\|\mathcal R ({\bf v},F)-\mathcal R ({\bf v}^\prime, F^\prime)\|_{C([0,t_0];H^{s-1}(\tilde N))}
\\ \nonumber & &\leq  
L_2r\, \|{\bf v}-{\bf v}^\prime\|_{C([0,t_0];\H^{(s)}(\tilde N))}+L_2\|F-{F}^\prime\|_{C([0,t_0];H^{s+1}(\tilde N))\cap C^1([0,t_0];H^{s}(\tilde N))},\hspace{-2cm}
%
  \eeq
 for all $s\in [1,k+1]$. 
%

%
%

Next we write (\ref{eq: notation for hyperbolic system}) as a first order
system.  To this end,
let $\mathcal A(t,{\bf v}):\H^{(s)}(\tilde N)\to \H^{(s-1)}(\tilde N)$ be the operator 
$\mathcal A(t,{\bf v})=\mathcal A_0(t,{\bf v})+\mathcal A_1(t,{\bf v})$
where in local coordinates and in the local trivialization of the bundle $\B^L$
\ba
\mathcal A_0(t,{\bf v})=-\left(\begin{array}{cc} 0 & I\\
 \frac 1{g^{00}(v)}\sum_{j,k=1}^3 g^{jk}(v)\frac \p{\p x^j}\frac\p{ \p x^k}\quad &
 \frac 1{g^{00}(v)}\sum_{m=1}^3 g^{0m}(v)\frac {\p}{\p x^m} \end{array}\right)
\ea 
with $g^{jk}(v)=g^{jk}(t,\cdotp;v)$ is a function on $\tilde N$ and 
\ba
\mathcal A_1(t,{\bf v})= \frac {-1}{g^{00}(v)}\left(\begin{array}{cc} 0 & 0\\
\sum_{j=1}^3 
 B^j(v)
   \frac \p{\p x^j}\quad &
B^{0}(v)\end{array}\right)
\ea 
where $B^j(v)$ depend on $v(t,x)$ and its first derivatives, and 
the connection coefficients (the Christoffel  symbols) corresponding  to $g(v)$. We denote $\S=(F,K)$ and 
 \ba
& &f_\S(t,{\bf v})=(f_\S^1(t,{\bf v}),f_\S^2(t,{\bf v}))\in \H^{(k)}(\tilde N),
\quad\hbox{where}\\
& & f_\S^1(t,{\bf v})=0,\quad f_\S^2(t,{\bf v})=\mathcal R({\bf v},F)(t,\cdotp)-\mathcal H({\bf v})(t,\cdotp)+K(t,\cdotp).
 \ea
 
 Note that when (\ref{new basic conditions}) are satisfied with $r<r_0$, inequalities (\ref{new f1f2 conditions}) imply that there exists $C_2>0$ so that
 \beq\label{fF bound r2}
 & &\|f_\S(t,{\bf v})\|_Y+\|f_{\S^\prime}(t,{\bf v}^\prime)\|_Y\leq C_2r^2,\\
 & &\|f_\S(t,{\bf v})-f_\S(t,{\bf v}^\prime)\|_Y\leq C_2
 r\,\|{\bf v}-{\bf v}^\prime\|_{C([0,t_0];Y)}. \nonumber
 \eeq

%
%

Let $U^{\bf v}(t,s)$ be the wave propagator corresponding to metric $g(v)$, that is,
$U^{\bf v}(t,s): {\bf h}\mapsto {\bf w}$, where ${\bf w}(t)=(w(t),\p_t w(t))$  solves
\ba
(\square_{g(v)}+V(x,D))w=0\quad\hbox{for } (t,y)\in [s,t_0]\times \tilde N,\quad\hbox{with }
 {\bf w}(s,y)={\bf h}.
 \ea 
 Let $S=(\nabla_{h}^*\nabla_{h}+1)^{k/2}:Y\to Z$ be an isomorphism.
As $k$ is an even integer, we see using multiplication estimates for Sobolev
spaces, see e.g.\
  \cite [Sec.\ 3.2, point (2)]{HKM}, that there exists $c_1>0$ (depending on
  $r_0,L_1,$ and $L^2$) so that  $\A(t,{\bf v})S-S\A(t,{\bf v})=C(t,{\bf v})$,
 where $\|C(t,{\bf v})\|_{Y\to Z}\leq c_1$ for all ${\bf v}$ satisfying
 (\ref{new basic conditions}). This yields that 
 the property (A2) in  \cite {HKM} holds, namely
 that $S\A(t,{\bf v})S^{-1}=\A(t,{\bf v})+B(t,{\bf v})$ where
$B(t,{\bf v})$ extends to a bounded operator in $Z$
 for which $\|B(t,{\bf v})\|_{Z\to Z}\leq c_1$ for all ${\bf v}$ satisfying
 (\ref{new basic conditions}).
Alternatively, to see the mapping properties of $B(t,{\bf v})$ we could use the fact that 
 $B(t,{\bf v})$ is a zeroth order pseudodifferential operator with
$H^{k}$-symbol.\hiddenfootnote{On mapping properties of such
pseudodifferential operators, see J. Marschall, Pseudodifferential operators with coefficients in Sobolev spaces. Trans. Amer. Math. Soc. 307 (1988), no. 1, 335-361.} 
 
 Thus the proof of  \cite [Lemma 2.6]{HKM} shows\hiddenfootnote{ Alternatively, as $U^{\bf v}(t,s)$ are propagators
  for linear wave equations having finite speed of wave propagation, one can prove 
the  estimates (\ref{U bounds}) by considering the wave equation
in local coordinate neighborhoods $W_j\subset W_j^\prime \subset \tilde N$,
$\Phi_j:  W_j^\prime\to \R^3$ and a partition of unity $\psi_j\in C^\infty_0(W_j)$ and
cut-off functions $\psi^\prime_j\in C^\infty_0(W_j^\prime)$.
Then  \cite [Lemma 2.6]{HKM}, used  in local coordinates, shows
that when $t_k-t_{k-1}$ is small enough then 
\ba 
\| \psi^\prime_j U^{\bf v}(t_k,t_{k-1})(\psi_j {\bf w}) \|_{Y}\leq C_{3,jk}^{\prime}e^{C_4(t-s)}\| \psi_j {\bf w} \|_{Y}.
\ea
Using  sufficiently small time steps $t_k-t_{k-1}$ and combining the estimates
for different $j$:s
together, one obtain estimates (\ref{U bounds}).}
 that 
%
there is a constant  $C_3>0$
   so that 
  \beq
  \label{U bounds}\|U^{\bf v}(t,s)\|_{Z\to Z}\leq C_3\quad\hbox{and}\quad
  \|U^{\bf v}(t,s) \|_{Y\to Y}\leq C_3\hspace{-1.5cm}
  \eeq 
  for $0\leq s<t\leq t_0$.
    By interpolation of estimates (\ref{U bounds}), we see also that 
    \beq
  \label{U bounds2}
  \|U^{\bf v}(t,s)\|_{X\to X}\leq C_3,
  \eeq  for $0\leq s<t\leq t_0$.

Let us next modify the reasoning given in \cite{Kato1975}: let  $r_1\in (0,r_0)$
be a parameter which value will be chosen later,  $C_1>0$
and 
$E$ be the space of functions ${\bf u}\in C([0,t_0];X)$ for which
\beq
\label{1 bounded}& &\|{\bf u}(t)\|_Y\leq r_1\quad\hbox{and}\\
\label{2 Lip}& &\|{\bf u}(t_1)-{\bf u}(t_2)\|_X\leq C_1|t_1-t_2|\eeq
for all $t,t_1,t_2\in [0,t_0]$.
The set $E$ is endowed by the metric of  $C([0,t_0];X)$. We 
note that  by \cite[ Lemma 7.3] 
{Kato1975},
a convex $Y$-bounded, $Y$-closed set is closed also in $X$.
Similarly, functions $G:[0,t_0]\to X$  satisfying (\ref{2 Lip}) 
form a closed subspace of $C([0,t_0];X)$. Thus $E\subset X$ is a closed set implying
that $E$ is a complete metric space.

Let
 \ba 
 & &\hspace{-1cm}W=\\ & &\hspace{-1cm}\{(F,K)\in C([0,t_0];H^{k+1}(\tilde N))^2;\ \sup_{t\in [0,t_0]}\|F(t)\|_{H^{k+1}(\tilde N)}+\|K(t)\|_{H^{k+1}(\tilde N)}< r_1\}.\hspace{-1cm}
 \ea

Following \cite[p. 44]{Kato1975}, we see that
the solution of
equation (\ref{eq: notation for hyperbolic system}) with the source $\S\in W$  is found as
a fixed point, if it exists, of the map $\Phi_\S:E\to C([0,t_0];Y)$ where 
$\Phi_\S({\bf v})={\bf u}$ is  given by
$$
{\bf u}(t)=\int_0^t U^{\bf v}(t,\tilde t)f_\S(\tilde t,{\bf v})\,d\tilde t,\quad 0\leq t\leq t_0.
$$ 
Below, we denote  ${\bf u}^{\bf v}=\Phi_\S({\bf v})$.

As $\Phi_{\S_0}(0)=0$ where
$\S_0=(0,0)$, we see using 
the above and the inequality $\|\,\cdotp\|_X\leq \|\,\cdotp\|_Y$ that the
 function ${\bf u}^{\bf v}$ satifies
  \ba
& &  \|{\bf u}^{\bf v}\|_{C([0,t_0];Y)}\leq C_3C_2t_0 r_1^2,\\
& &  \|{\bf u}^{\bf v}(t_2)-{\bf u}^{\bf v}(t_1)\|_{X}\leq C_3C_2r_1^2|t_2-t_1|,\quad t_1,t_2\in [0,t_0].
  \ea
  When $r_1>0$ is so small that $C_3C_2(1+t_0)<r_1^{-1}$ 
  and $C_3C_2r_1^2<C_1$
  we see that $  \|\Phi_\S({\bf v})\|_{C([0,t_0];Y)}<r_1$
  and $ \|\Phi_\S({\bf v})\|_{C^{0,1}([0,t_0];X)}<C_1$.
  Hence $ \Phi_\S(E)\subset E$
 and we can consider $\Phi_\S$ as a map
 $\Phi_\S:E\to E$.
  
  As $k>1+\frac 32$, it follows from Sobolev embedding theorem
  that $X=\H^{(k)}(\tilde N)\subset C^1(\tilde N)^2$. This yields that
  by \cite[Thm. 3]{Kato1975}, for the original reference, see Theorems III-IV
  in \cite{Kato1973},
  \ba
 & & \|(U^{\bf v}(t,s)- U^{\bf v^\prime}(t,s)){\bf h}\|_X \\
 & &\leq 
 C_3\bigg(\sup_{t^\prime \in [0,t]}\|\A(t^\prime,{\bf v})-
 \A(t^\prime,{\bf v}^\prime)\|_{Y\to X}  \|U^{\bf v}(t^\prime,0) {\bf h}\|_Y\bigg)\\
 & &\leq 
 C_3^2 \|{\bf v}-{\bf v}^\prime\|_{C([0,t_0];X)}  \|{\bf h}\|_Y.
  \ea
  Thus,
  \ba
&&\|  U^{\bf v}(t,s)f_\S(s,{\bf v})- U^{\bf v^\prime}(t,s)f_{\S}(s,{\bf v}^\prime)\|_X\\
&\leq&\|  (U^{\bf v}(t,s)- U^{\bf v^\prime}(t,s))f_{\S}(s,{\bf v})\|_X+
\|  U^{\bf v^\prime}(t,s)(f_\S(s,{\bf v})- f_{\S}(s,{\bf v}^\prime))\|_X\\
&\leq&
(1+C_3)^2C_2r_1^2 \|{\bf v}-{\bf v}^\prime\|_{C([0,t_0];X)} .
  \ea
This implies that
\ba
& & \| \Phi_\S({\bf v})-\Phi_{\S}({\bf v}^\prime)\|_{C([0,t_0];X)}
\leq t_0(1+C_3)^2C_2r_1^2 \|{\bf v}-{\bf v}^\prime\|_{C([0,t_0];X))}.
\ea

Assume next that  $r_1>0$ is so small that we have also
\ba
t_0(1+C_3)^2C_2r_1^2<\frac 12.
\ea
 cf.\
  Thm.\ I in \cite{HKM} (or
  (9.15)  and (10.3)-(10.5) in \cite{Kato1975}). For $\S\in W$ this 
 implies that $\Phi_\S:E\to E$ is a contraction with a contraction constant 
 $C_L\leq \frac 12$,
  and thus 
 $\Phi_\S$ has a unique fixed point ${\bf u}$  in the space $E\subset C^{0,1}([0,t_0];X)$.
 
 Moreover, elementary considerations related to fixed point
 of the map $\Phi_\S$ show that ${\bf u}$  in  $C([0,t_0];X)$ depends in $E\subset C([0,t_0];X)$
 Lipschitz-continuously on
 $\S\in  W\subset C([0,t_0];H^{k+1}(\tilde N))^2$. Indeed, if $\|\S-\S^\prime\|_{C([0,t_0];H^{k+1})^2}<\e$,
 we see that
  \beq\label{fF bound r2 B}
 \|f_\S(t,{\bf v})-f_{\S^\prime}(t,{\bf v})\|_Y\leq C_2\e,\quad t\in [0,t_0],
 \eeq
 and when  (\ref{fF bound r2}) and (\ref{U bounds}) are satisfied 
 with $r=r_1$,
 we have
  \ba
& & \| \Phi_\S({\bf v})-\Phi_{\S^\prime}({\bf v}^\prime)\|_{C([0,t_0];Y)}
\leq C_3C_2t_0 r_1^2.
\ea
 Hence 
 \ba
 \|\Phi_\S({\bf v})-\Phi_{\S^\prime}({\bf v})\|_{C([0,t_0];Y)}\leq  t_0C_3C_2\e.
 \ea
 This and standard estimates for fixed points, yield that when $\e$ is small enough 
 the fixed point ${\bf u}^\prime$ of the map  $\Phi_{\S^\prime}:E\to E$ 
 corresponding to the source ${\S^\prime}$ and
  the fixed point ${\bf u}$ of the map $\Phi_{\S}:E\to E$ 
  corresponding to the source ${\S}$ satisfy
  \beq\label{stability}
  \|{\bf u}-{\bf u}^\prime\|_{C([0,t_0];X)}\leq  \frac 1{1-C_L}t_0C_3C_2\e.
  \eeq
{
Thus the solution ${\bf u}$ 
  depends in $C([0,t_0];X)$ Lipschitz continuously on
 $\S\in  C([0,t_0];H^{k+1}(\tilde N))^2$ (see also \cite[Sect.\ 16]{Kato1975}, 
 and \cite{Ringstrom})}. 
In fact,  for analogous systems it is possible to show that $u$ is in $C([0,t_0];Y)$,
but one can not obtain  
 Lipschitz or H\"older  stability for $u$ in the $Y$-norm, see  \cite{Kato1975}, Remark 7.2. 
 
 Finally, we note that the fixed point ${\bf u}$ of $\Phi_\S$ can be found as a limit
 $ {\bf u}=\lim_{n\to \infty}{\bf u}_n$ in $C([0,t_0];X)$, where ${\bf u}_0=0$ and  ${\bf u}_n=\Phi_\S({\bf u}_{n-1})$.
Denote ${\bf u}_n=(g_n-\hat g,\phi_n-\hat \phi)$. We see that if
 $\supp({\bf u}_{n-1})\subset J_{\hat g}(\supp(\S))$ 
then also $\supp(g_{n-1}-\hat g)\subset J_{\hat g}(\supp(\S))$.
Hence  for all $x\in M_0\setminus J_{\hat g}(\supp(\S))$ we see that $J_{g_{n-1}}^-(x)
\cap J_{\hat g}(\supp(\S))=\emptyset.$ Then, using the definition of the map $\Phi_\S$ we see that
 $\supp({\bf u}_n)\subset J_{\hat g}(\supp(\S))$. Using induction we see that this holds
for all $n$ and hence we see that the solution ${\bf u}$ satisfies
 \beq\label{eq: support condition for u}
 \supp({\bf u})\subset J_{\hat g}(\supp(\S)).
 \eeq 
 \hiddenfootnote{REMOVE MATERIAL IN THIS FOOTNOTE OR MOVE IT ELSEWHERE:
  Let us next use the above considerations to study system 
   (\ref{eq: notation for hyperbolic system 1}) on $\R\times N$.
   For this end, we denote below the background metric of  $\tilde M=\R\times \tilde N$
   by $\tilde g$ and the background metric of  $\hattuM _0=\R\times  N$
   by $\hat g$.
 Consider a source 
\ba
\S\in {C^1([0,T];H^{s_0-1}(\tilde N))\cap
C^0([0,t_0];H^{s_0}(\tilde N))}
\ea  that is supported in a compact set $\K$ and
small enough in the norm of this space and let $\tilde u$ 
be the solution of the  
system (\ref{eq: notation for hyperbolic system 1}) on $\R\times \tilde N$.
Assume that  $u$ is some $C^2$-solution of  
(\ref{eq: notation for hyperbolic system 1}) on $\R\times N$
with the same source $\S$. Our aim is to show that it
coincides with $\tilde u$ in its support.
%
For this end, let $\tilde g^\prime$ be a metric tensor which
has the same eigenvectors as $\tilde g$ and the same
eigenvalues except that  the unique negative eigenvalue of $\tilde g$
 is multiplied by $(1-h_1)$, $h_1>0$. 
Assume that $h_1>0$ is so small
that $J^+_{\hattuM _0,\tilde g^\prime}(p^+)\subset  \hattuM _0\subset N_0\times (-\infty,t_0]$
and assume that $\e$ is so small that $g(\tilde u)<\tilde g^\prime$ for all 
sources $\S$ with $\|\S\|_{C^1([0,t_0];H^{s_0-1}(\tilde N))\cap
C^0([0,t_0];H^{s_0}(\tilde N))}\leq \e$. Then
 $J^+_{\tilde M,g(\tilde u)}(p^+)\subset  (-\infty,t_0]\times N_0$.
This  proves the estimate (\ref{eq: Lip estim}).

Next us we consider solution $u$ on $\hattuM _0$ corresponding to source $\S$.
Let  $T(h_1,\S)$
be the supremum of all $T^\prime\leq t_0$ for which
$J^+_{g(u)}(p^+)\cap \hattuM _0(T^\prime)\subset N_1\times (-\infty,T^\prime]$.
Assume that $ T(h_1,\S)<t_0$.
Then for $T^\prime= T(h_1,\S)$ we see that
$u$ can be continued by zero from $ (-\infty,T^\prime]\times N_1$
to a function on $(-\infty,T^\prime]\times \tilde N$ that 
satisfies the system (\ref{eq: notation for hyperbolic system 1}). Thus
$u$ coincides in $ (-\infty,T^\prime]\times N_1$
with $\tilde u$ and vanishes outside this set.
As  $J^+_{\tilde M,g(\tilde u)}(p^+)\subset  (-\infty,t_0]\times N_0$, this implies that
also $u$ satisfies 
$J^+_{g(u)}(p^+)\cap \hattuM _0(T^\prime)\subset (-\infty,T^\prime]\times N_0$.
As $u$ is $C^2$-smooth and $N_0\subset \subset N_1$
we see that there is $T^\prime_1> T(h_1,\S)$ so 
that $J^+_{g(u)}(p^+)\cap \hattuM _0(T^\prime_1)\subset  (-\infty,T^\prime_1]
\times N_1$ that is in contradiction with definition of $ T(h_1,\S)$. Thus
we have to have  $T(h_1,\S)=t_0$.


This shows that when the norm of $\S$ is smaller than the above chosen $\e$
then $J^+_{g(u)}(p^+)\cap \hattuM _0\subset N_0\times (-\infty,t_0]$.
Then $u$ vanishes outside $N_1\times (-\infty,t_0]$ by
the support condition in (\ref{eq: notation for hyperbolic system 1}).
Thus $u$ and $\tilde u$ are both supported on  $N_1\times (-\infty,t_0]$ 
and coincide there. As $\tilde u$ is unique, 
this shows that for a sufficiently small sources $ \S$
 supported on  $\K$   the solution $u$ of
(\ref{eq: notation for hyperbolic system 1})  in $M_0$
exists and is unique.}

\bigskip

\subsection*{Appendix D: An inverse problem for a non-linear wave equation}

In this appendix  we explain how a problem for a scalar wave equation can
be solved with the same techniques that we used for the Einstein equations.

Let $(M_j,g_j)$, $j=1,2$ be two globally hyperbolic  $(1+3)$  
dimensional Lorentzian
manifolds represented using
  global smooth time functions as $M_j=\R\times N_j$, $\mu_j=\mu_j([-1,1])
\subset M_j$ be
a time-like geodesic and $U_j\subset M_j$ be open, relatively compact  
neighborhood
of $\mu_j([s_-,s_+])$, $-1<s_-<s_+<1$. Let $M_j^0=(-\infty,T_0)\times  
N_j$ where $T_0>0$ is such
that $U_j\subset  M_j^0$.
{\mltext Consider the
non-linear wave equation
\beq\nonumber
& &\hspace{-.5cm}  \square_{g_j}u(x)+a_j(x)\,u(x)^2=f(x)
\quad\hbox{on }M_j^0,\hspace{-.5cm}\\
\label{eq: wave-eq general} & &\quad \supp(u)\subset J^+_{g_j}(\supp(f)),
\eeq
where $\supp(f)\subset U_j,$
\ba
\square_gu=
-\sum_{p,q=1}^4\det(-g(x))^{-1/2}\frac  \p{\p x^p}
\left ((-\det(g(x))^{1/2}
g^{pq}(x)\frac \p{\p x^q}u(x)\right),
\ea
$\det(g)=\det((g_{pq}(x))_{p,q=1}^4)$,
  $f\in C^6_0(U_j)$ is a controllable source, and $a_j$ is
a non-vanishing $C^\infty$-smooth  function.}
Our goal is to prove the following result:

\begin{theorem}\label{main thm3}
Let $(M_j,g_j)$, $j=1,2$ be two
open,  smooth, globally hyperbolic    Lorentzian manifolds of  
dimension $(1+3)$.
Let
$p^+_j=\mu_j(s_+), p^-_j=\mu_j(s_-)\in M_j$ the points of a  time-like  
geodesic  $\mu_j=
\mu_j([-1,1])\subset M_j$, $-1<s_-<s_+<1$,
and let  $U_j\subset M_j$ be an open  relatively compact neighborhood
of $\mu_j([s_-,s_+])$  given in (\ref{eq: Def Wg with hat}).  Let  
$a_j:M_j\to \R$, $j=1,2$ be $C^\infty$-smooth functions that are
non-zero on $M_j$.

Let   $L_{U_j}$, $j=1,2$ be measurement operators defined
in an open set $\mathcal W_j\subset C^6_0(U_j)$ containing the zero  
function by setting
\beq\label{measurement operator}
L_{U_j}: f\mapsto u|_{U_j},\quad f\in C^6_0(U_j),
\eeq
where $u$ satisfies the wave equation (\ref{eq: wave-eq general}) on  
$(M_j^0,g_j)$.

Assume that there is a diffeomorphic isometry
$\Phi:U_1\to U_2$ so that $\Phi(p^-_1)=p^-_2$ and
$\Phi(p^+_1)=p^+_2$ and the measurement maps
satisfy
\ba
((\Phi^{-1})^*\circ L_{U_1}\circ \Phi^*) f =L_{U_2}f
\ea
for all $f\in \W$ where $\W$ is  some  neighborhood of the zero function in $C^6_0(U_2)$.

Then there is a diffeomorphism $\Psi:I(p^-_1,p^+_1)\to I(p^-_2,p^+_2)$,
and the metric $\Psi^*g_2$ is conformal to $g_1$ in  
$I(p^-_1,p^+_1)\subset M_1$,
that is, there is $\beta(x)$ such that $g_1(x)=\beta(x)(\Psi^*g_2)(x)$ in  
$I(p^-_1,p^+_1)$.
\end{theorem}

We note that the  smoothness assumptions assumed above on the functions
$a$ and the source $f$ are not optimal.
The proof, presented below, is based on using the interaction
of singular waves.  The techniques used can be modified used to study
different non-linearities, such as the equations
$\square_{\hat g} u+a(x)u^3=f$, $\square_{\hat g} u+a(x)u_t^2=f$,
or $\square_{g(x,u(x))} u=f$, but these considerations are outside the scope
of this paper.

Theorem \ref{main thm3} can be applicable
  for example in the mathematical analysis of  non-destructive
testing or imaging in non-linear medium e.g,
in imaging  the non-linearity  of the acoustic material parameter  inside a
given body  when it
is under large, time-varying, possibly periodic, changes of the external
pressure and at the same time the body is probed with small-amplitude fields.
  Such acoustic measurements are analogous to
  the recently developed Ultrasound Elastography imaging technique where the interaction
  of the elastic shear and pressure
waves is used for medical imaging, see  
e.g.\ \cite{Hoskins,McLaughlin1,McLaughlin2,Ophir}. There, the slowly  
progressing
shear wave is imaged using a pressure wave and the
image of the shear wave inside the body is used to determine
approximately the material parameters. In other words,
the changes which the elastic wave causes in the medium
are imaged using the interaction of the s-wave and p-wave
components of the  elastic wave.

Let us also consider some implications of theorem \ref{main thm3}  for inverse
problems for a non-linear equation involving a time-independent metric
\beq\label{product metric} & &g(t,y)=-dt^2+
\sum_{\alpha,\beta=1}^3 h_{\alpha\beta}(y)dy^\alpha dy^\beta,\quad
(t,y)\in \R\times N.
\eeq
The metric (\ref{product metric}) corresponds to
the hyperbolic operator $\p_t^2- \Delta_{h}$,
with a time-independent Riemannian metric  
$h=(h_{\alpha\beta}(y))_{\alpha,\beta=1}^3,$
$y\in N$, where $N$ is a 3-dimensional manifold and
\ba
  \Delta_{h}u(y,t)=\sum_{\alpha,\beta=1}^3\frac \p{\p y^\alpha}
\left( h^{\alpha\beta}(y) \frac \p{\p y^\beta}u(t,y)\right).
\ea

\begin{corollary}\label{coro thm3}
Let $(M_j,g_j)$, $M_j=\R\times N_j$, $j=1,2$ be two
open,  smooth, globally hyperbolic    Lorentzian manifolds of  
dimension $(1+3)$.
Assume that $g_j$ is the product metric of the type (\ref{product metric}),
$g_j=-dt^2+h_j(y)$, $j=1,2$.
Let
$p^+_j=\mu_j(s_+),p^-_j=\mu_j(s_-)\in (0,T_0)\times N_j$ be points of a  time-like  
geodesic  $\mu_j=
\mu_j([-1,1])\subset M_j$, $-1<s_-<s_+<1$, and to fix
the time variable, assume that $\mu_j(s_-)\in \{1\}\times N_j$.

Let  $U_j\subset (0,T_0)\times N_j$ be an open  relatively compact neighborhood
of $\mu_j([s_-,s_+])$ given in (\ref{eq: Def Wg with hat}). 
  Let $a_j:M_j\to \R$, $j=1,2$ be  
$C^\infty$-smooth functions that are
non-zero on $M_j$ and $x=(t,y)\in \R\times N$.

For $j=1,2$, consider the
non-linear wave equations
\beq\nonumber
& &\hspace{-.5cm}  (\frac {\p^2}{\p  
t^2}-\Delta_{{h_j}})u(t,y)+a_j(y,t)(u(t,y))^2=f(t,y)
\quad\hbox{on }(0,T_0)\times N_j,\hspace{-.5cm}\\
\label{eq: wave-eq} & &\quad \supp(u)\subset J^+_{g_j}(\supp(f)),
\eeq
where $f\in C^6_0(U_j)$, $j=1,2$.
Let $L_{U_j}:f\mapsto u|_{U_j}$ be the measurement operator (\ref{measurement operator})
for the wave equation (\ref{eq: wave-eq}) with the  Riemannian metric $h_j(x)$ and the coefficient
$a_j(x,t)$ for $j=1,2$,  defined
in some $C^6_0(U_j)$ neighborhood of the zero function.

Assume that there is a diffeomorphism
$\Phi:U_1\to U_2$
of the form $\Phi(t,y)=(t,\phi(y))$  so that
\ba
((\Phi^{-1})^*\circ L_{U_1}\circ \Phi^*)f=L_{U_2} f
\ea
for all $f\in \W$ where $\W$ is  some  neighborhood of the zero function in $C^6_0(U_2)$.

Then there is a diffeomorphism $\Psi:I^+(p^-_1)\cap I^-(p^+_1)\to  
I^+(p^-_2)\cap I^-(p^+_2)$ of the form
$\Psi(t,y)=(\psi(y),t)$,
the metric $\Psi^*g_2$ is isometric to $g_1$ in $I^+(p^-_1)\cap I^-(p^+_1)$,
and $a_1(t,y)=a_2(t,\psi(y))$ in $I^+(p^-_1)\cap I^-(p^+_1)$.
\end{corollary}

\medskip

Next we consider the proofs.

\medskip

 {\bf Proof.} (of Theorem \ref{main thm3}).
We will explain how the proof of
Theorem \ref{main thm Einstein} for the Einstein equation needs to
be modified to obtain the similar result for the non-linear wave equation.

Let  $(M,\hat g)$ be a  smooth
globally hyperbolic Lorentzian manifold that we represent
using a global smooth time function as $M=(-\infty,\infty)\times N$,
and consider   $M^0=(-\infty,T)\times N\subset M$.
Assume that
the set $U$, where the sources are supported and where
we observe the waves, satisfies
$U\subset [0,T]\times N$.

The results of  section \ref{subsec: Direct problem} concerning the  
direct problem
for Einstein equations can be modified for the wave equation
\beq\label{PABC eq}
& &\square_{\hat g}u+au^2=f,\quad\hbox{in }M^0=(-\infty,T)\times N,\\
& &u|_{(-\infty,0)\times N}=0,\nonumber
\eeq
where
$a=a(x)$ is a smooth, non-vanishing function. Here we denote
the metric by $\hat g$ to emphasize the fact that it is independent
on the solution $u$.  Below,
let  $Q$ be the causal inverse  operator of $\square_{\hat g}$.

When $f$ in $C_0([0,t_0];H^6_0(B))\cap
C_0^1([0,t_0];H^5_0(B))$
is small enough,
we see by using  \cite[Prop.\ 9.17]{Ringstrom} and   \cite [Thm.\ III]{HKM},
see also (\ref{stability}) in Appendix C, that  the
equation (\ref{PABC eq}) has a unique solution
$u\in  C_0([0,t_0];H^{5}(N))\cap C_0^1([0,t_0];H^{4}(N))$.
Moreover, 
we  can consider the case when $f=\e f_0$ where $\e>0$
is small.
Then, we can write
\ba
u=\e w_1+\e^2 w_2+\e^3 w_3+\e^4 w_4+E_\e
\ea
where $w_j$ and the reminder term $E_\e$ satisfy
\ba
w_1&=&Qf,\\
w_2&=&-Q(a\,w_1\,w_1),\\
w_3&=&-2Q(a\, w_1\,w_2)\\
&=&2Q(a\, w_1\,Q(a\,w_1\,w_1))
,\\
w_4
&=&-Q(a\, w_2\,w_2)-2Q(a\, w_1\,w_3)\\
&=&-Q(a\, Q(a\,w_1\,w_1)\,Q(a\,w_1\,w_1))\\
& &+4Q(a\, w_1\,Q(a\, w_1\,w_2))
\\
&=&-Q(a\, Q(a\,w_1\,w_1)\,Q(a\,w_1\,w_1))\\
& &-4Q(a\, w_1\,Q(a\, w_1\,Q(a\,w_1\,w_1))),\\
& &\hspace{-1.5cm}\|E_\e\|_{C([0,t_0];H^{4}_0(N))\cap  
C^1([0,t_0];H^{3}_0(N))}\leq C\e^5.
\ea

If we consider sources  $f_{\vec\e}(x)=\sum_{j=1}^4\e_j f_{(j)}(x)$,
$\vec\e=(\e_1,\e_2,\e_3,\e_4),$
and the corresponding solution $u_{\vec \e}$ of (\ref{PABC eq}), we see that
\beq \nonumber
\M^{(4)}&=&\p_{\vec \e}^4u_{\vec \e}|_{\vec\e=0}\\
&=& \nonumber
\p_{\e_1}\p_{\e_2}\p_{\e_3}\p_{\e_4}u_{\vec \e}|_{\vec\e=0}\\
\label{4th interaction for wave eq B}
&=&-\sum_{\sigma\in \Sigma(4)}\bigg(Q(a\,  
Q(a\,u_{(\sigma(1))}\,u_{(\sigma(2))})\,Q(a\,u_{(\sigma(3))}\,u_{(\sigma(4))}))\\  
\nonumber
& &\quad+ 4Q(a\, u_{(\sigma(1))}\,Q(a\,  
u_{(\sigma(2))}\,Q(a\,u_{(\sigma(3))}\,u_{(\sigma(4))})))\bigg),
\eeq
where $u_{(j)}=Qf_{(j)}$ and $\Sigma(\ell)$ is the set
of permutations of the set $\{1,2,3,\dots,\ell\}$.

The results of Lemma \ref{lem: lagrangian 1} can
be replaced by the results of \cite[Prop.\ 2.1]{GU1} as follows.
Using the same notations as in  Lemma \ref{lem: lagrangian 1}, let
$Y=Y(x_0,\zeta_0;t_0,s_0)$, $K=K(x_0,\zeta_0;t_0,s_0)$, and  
$\Lambda_1=\Lambda(x_0,\zeta_0;t_0,s_0)$, and consider a source $f\in  
\I^{n-3/2}(Y)$.  Then $u=Qf$ satisfies
$u|_{M_0\setminus Y}\in
  \I^{n-1/2} ( M_0\setminus Y;\Lambda_1)$. Assume that  
$(x,\xi),(y,\eta)\in L^+M$  are on the same bicharacteristics of  
$\square_{\hat g}$,
  and $x<y$, that is, $((x,\xi),(y,\eta))\in \Lambda_{\hat g}^\prime$.  
Moreover, assume
  that $(x,\xi)\in N^*Y$.
Let   $\tilde b(x,\xi)$ be the principal
   symbol of $f$ at $(x,\xi)$ and
    $\tilde a(y,\eta)$ be the principal
   symbol of $u$ at $(y,\eta)$. Then $\tilde a(y,\eta)$
   depends linearly on $ \tilde f(x,\xi)$ and
    $\tilde a(y,\eta)$ vanishes if and only if
   $ \tilde f(x,\xi)$ vanishes.

Analogously to the Einstein equations,
we consider the indicator function
\beq\label{test sing 2}
\Theta_\tau^{(4)}=\bra F_{\tau},\M^{(4)}\cet_{L^2(U)},
\eeq
where
$\M^{(4)}$ is given by (\ref{4th interaction for wave eq B})
with  $u_{(j)}=Qf_{(j)}$, $j=1,2,3,4$, where  $f_{(j)}\in  
\I^{n-3/2}(Y(x_j,\xi_j;t_0,s_0))$, $n\leq -n_1$,
  and $F_\tau$ is the source producing a gaussian beam $Q^*F_\tau$
that propagates to the past along the geodesic $\gamma_{x_5,\xi_5}(\R_-)$,
see (\ref{Ftau source}).

Similar results to the ones given in Proposition \ref{lem:analytic limits A}
are
valid. Let us consider next the case when $(x_5,\xi_5)$ comes from the  
4-intersection
of  rays corresponding to $(\vec x,\vec \xi)=((x_j,\xi_j))_{j=1}^4$ and
$q$ is the corresponding intersection point, that is, $q=\gamma_{x_j,\xi_j}(t_j)$ for
all $j=1,2,3,4,5$.
Then 
\beq\label{indicator 2}
\Theta^{(4)}_\tau\sim 
\sum_{k=m}^\infty s_{k}\tau^{-k}
  \eeq
as $\tau\to \infty$ where   $m=-4n+4$.
Moreover,
let $b_j=(\dot\gamma_{x_j,\xi_j}(t_j))^\flat$ and
  $\bsequence=(b_{j})_{j=1}^5\in (T^*_q\hattuM _0)^5$,
  $w_j$ be the principal symbols of the waves $u_{(j)}$
  at $(q,b_j)$, and ${\bf w}=(w_j)_{j=1}^5$.
Then we see as in Proposition \ref{lem:analytic limits A}  that there is
  a real-analytic function $\mathcal G(\bsequence,{\bf w})$ such that
  the leading order term in (\ref{indicator})  satisfies
  \beq\label{definition of G 2}
s_{m}=
\mathcal  G(\bsequence,{\bf w}).
\eeq

The proof of Prop.\ \ref{singularities in Minkowski space} dealing
with Einstein equations needs significant changes and we need to prove  
the following:

\begin{proposition}\label{singularities in Minkowski space for wave equation}
The  function $\ \mathcal  G(\bsequence,{\bf w})
$ given in (\ref{definition of G 2}) for the non-linear wave equation
is a non-identically vanishing real-analytic function.
\end{proposition}

\noindent{\bf Proof.}
Let us use the notations introduced in Prop.\ \ref{singularities in  
Minkowski space}.

As for the Einstein equations, we consider light-like vectors
\ba
b_5=(1,1,0,0),\quad b_j=(1,1-\frac  
12\rhoepsilon_j^2,\rhoepsilon_j+O(\rhoepsilon_j^3),\rhoepsilon_j^3),\quad  
j=1,2,3,4,
\ea
  in the Minkowski space $\R^{1+3}$, endowed with the standard metric $g=\diag(-1,1,1,1)$, where the terms
  $O(\rhoepsilon_k^3)$ are such that the vectors $b_j$, $j\leq 5$,
are  light-like. Then
\ba
g(b_5,b_j)= -\frac 12 \rhoepsilon_j^2,\quad
g(b_k,b_j)=-\frac 12 \rhoepsilon_k^2-\frac 12  
\rhoepsilon_j^2+O(\rhoepsilon_k\rhoepsilon_j).
\ea
Below, we denote $\omega_{kj}=g(b_k,b_j)$.
Note that if $\rhoepsilon_j<\rhoepsilon_k^4$, we have
$\omega_{kj}=-\frac 12 \rhoepsilon_k^2+O(\rhoepsilon_k^3).$

{For the wave equation,
we use different parameters $\rhoepsilon_j$ than for the Einstein  
equations, and
define (so, we use here the "unordered" numbering 4-2-1-3)
\beq\label{eq: ordering of epsilons wave equation}
\rhoepsilon_4=\rhoepsilon_2^{100},\  
\rhoepsilon_2=\rhoepsilon_1^{100},\hbox{ and  
}\rhoepsilon_1=\rhoepsilon_3^{100}.
\eeq
Below in this proof, we denote $\vec\rhoepsilon\to 0$ when
$\rhoepsilon_3\to 0$ and
$\rhoepsilon_4,$ $\rhoepsilon_2,$ and $\rhoepsilon_1$ are defined using
$\rhoepsilon_3$ as in (\ref{eq: ordering of epsilons wave equation}).

Let us next consider in Minkowski space
the coordinates $(x^j)_{j=1}^4$ such that
$K_j=\{x^j=0\}$ are light-like hyperplanes and  the waves $u_j=u_{(j)}$ that
satisfy in the Minkowski space $\square u_j=0$ and can be written
as
\ba
u_j(x)=\int_{\R}e^{i x^j \theta }a_j(x,\theta)\,d\theta,
\quad
a_j(x,\theta^{\prime})\in S^{n}(\R^4;\R\setminus 0),\quad j\leq 4,
\ea
and
\ba
u^\tau(x) =\chi(x^0)w_{(5)} \exp(i\tau b^{(5)}\,\cdotp x).
\ea
Note that the singular supports of the waves $u_j$, $j=1,2,3,4,$  
intersect then at the point
$\cap_{j=1}^4 K_j=\{0\}$.
Analogously to the definition  (\ref{definition of G}) we considered
  for  Einstein equations, we
   define  the (Minkowski) indicator function
\ba
\mathcal G^{({\bf m})}(v,{\bf b})=
\lim_{\tau\to\infty} \tau^{m}(\sum_{\b\leq n_1}
\sum_{\sigma\in \Sigma(4)}
T^{({\bf m}),\b}_{\tau,\sigma}+\tilde T^{({\bf m}),\b}_{\tau,\sigma}),
\ea
where
\ba
T^{({\bf m}),\beta}_{\tau,\sigma}
&=&\bra Q_0(u^\tau\, \cdotp \a u_{\sigma(4)}), h\,\cdotp  \a  
u_{\sigma(3)}\,\cdotp Q_0(\a u_{\sigma(2)}\,\cdotp u_{\sigma(1)})\cet,\\
\tilde T^{({\bf m}),\beta}_{\tau,\sigma}
&=&\bra u^\tau,h\a\,Q_0(\a u_{\sigma(4)}\,\cdotp  \,u_{\sigma(3)})\,\cdotp
  Q_0(\a u_{\sigma(2)}\,\cdotp \, u_{\sigma(1)})\cet.
\ea

As for Einstein equations, we
see that when $\alpha$ is equal to the value of the function $a(t,y)$
at the intersection point $q=0$ of the waves,
we have
  $\mathcal G^{({\bf m})}(v,{\bf b})=\mathcal G(v,{\bf b})$.

Similarly to the Lemma \ref{lem:analytic limits A}  we analyze next  
the functions
\ba
\Theta_\tau^{({\bf m})}=\sum_{\beta\in J_\ell}\sum_{\sigma\in  
\Sigma(4)}(T_{\tau,\sigma}^{({\bf m}),,\beta}+\tilde  
T_{\tau,\sigma}^{({\bf m}),\beta}).
\ea
Here $({\bf m})$ refers to "Minkowski".
We denote $T_{\tau}^{({\bf m}),\beta}=T_{\tau,id}^{({\bf m}),\beta}$
and  $\tilde T_{\tau}^{({\bf m}),\beta}=\tilde T_{\tau,id}^{({\bf m}),\beta}$.

Let us first consider the case when the permutation
$\sigma=id$. Then, as in the proof of Prop.\ \ref{singularities in  
Minkowski space},
in the case  when $\vec S^\beta =(Q_0,Q_0)$, we have
\ba
T^{({\bf m}),\beta}_\tau\\
&& \hspace{-2cm}=C_1 \det(A) \cdotp
(i\tau)^{m}(1+O(\frac 1\tau))
\vec\rhoepsilon^{\,2\vec n}
(\omega_{45}\omega_{12})^{-1}\rhoepsilon_4^{-4}\rhoepsilon_2^{-4}\rhoepsilon_1^{-4}\rhoepsilon_3^{2}\cdotp\P\\
&& \hspace{-2cm}=C_2\det(A) \cdotp
(i\tau)^{m}(1+O(\frac 1\tau))
\vec\rhoepsilon^{\,2\vec n}
\rhoepsilon_4^{-4-2}\rhoepsilon_2^{-4}\rhoepsilon_1^{-4-2}\rhoepsilon_3^{-2}\cdotp\P
\ea
where $\P$ is the product of the principal symbols of the waves $u_j$
at zero, $\vec\rhoepsilon^{\,2\vec n}=
\rhoepsilon_1^{2n}\rhoepsilon_2^{2n}\rhoepsilon_3^{2n}\rhoepsilon_4^{2n}$, and  
$C_1$ and $C_2$ are non-vanishing.
Similarly, a direct computation yields
\ba
\tilde T^{({\bf m}),\beta}_\tau\\
&& \hspace{-2cm}=C_1 \det(A) \cdotp
(i\tau)^{n}(1+O(\frac 1\tau))
\vec\rhoepsilon^{\,2\vec n}
(\omega_{43}\omega_{21})^{-1}\rhoepsilon_4^{-4}\rhoepsilon_2^{-4}\rhoepsilon_1^{-4}\rhoepsilon_3^{-4}\cdotp\P\\
&& \hspace{-2cm}=C_2\det(A) \cdotp
(i\tau)^{m}(1+O(\frac 1\tau))
\vec\rhoepsilon^{\,2\vec n}
\rhoepsilon_4^{-4}\rhoepsilon_2^{-4}\rhoepsilon_1^{-4-2}\rhoepsilon_3^{-4-2}\cdotp\P,
\ea
where again, $\P$ is the product of the principal symbols of the waves $u_j$
at zero and $C_1$ and $C_2$ are non-vanishing.

Considering formula (\ref{4th interaction for wave eq B}), we
see that for the wave equation we do not need to consider the
terms  that for the Einstein equations correspond to the
cases when $\vec S^\beta =(Q_0,I)$, $\vec S^\beta =(I,Q_0)$,
or $\vec S^\beta =(I,I)$ as the corresponding terms do not appear
in formula (\ref{4th interaction for wave eq B}).

Let us now consider permutations $\sigma$ of the indexes $(1,2,3,4)$
and compare the terms
\ba
& &L^{({\bf m}),\beta}_{\sigma}=\lim_{\tau\to \infty} \tau^{m}  
T^{({\bf m}),\beta}_{\tau,\sigma},\\
&& \tilde L^{({\bf m}),\beta}_{\sigma}=\lim_{\tau\to \infty}  \tau^{m} \tilde  
T^{({\bf m}),\beta}_{\tau,\sigma}.
\ea
Due to the presence of $\omega_{45}\omega_{12}$ in the above computations,
we observe that  all the terms
  $ \tilde  L^{({\bf m}),\beta}_{\tau,\sigma}/L^{({\bf m}),\beta}_{\tau,id}\to 0$
  as $\vec \rhoepsilon\to 0$, see (\ref{eq: ordering of epsilons wave  
equation}).
Also, if $\sigma\not=(1,2,3,4)$ and  $\sigma\not=\sigma_01=(2,1,3,4)$,
we see that  $ \tilde   
L^{({\bf m}),\beta}_{\tau,\sigma}/L^{({\bf m}),\beta}_{\tau,id}\to 0$
  as $\vec \rhoepsilon\to 0$.
  Also, we observe that $ L^{({\bf m}),\beta}_{\tau,\sigma_1}=
L^{({\bf m}),\beta}_{\tau,id}$.
Thus we see that the equal terms
$ L^{({\bf m}),\beta}_{\tau,\sigma_1}=
  L^{({\bf m}),\beta}_{\tau,id}$
that give the largest contributions as $\vec \rhoepsilon\to 0$
and that when $\P\not =0$ the sum
\ba
S( \vec \rhoepsilon,\P)=\sum_{\sigma\in \Sigma(4)}(
  L^{({\bf m}),\beta}_{\tau,\sigma}+\tilde L^{({\bf m}),\beta}_{\tau,\sigma})
\ea
is non-zero when
$\rhoepsilon_3>0$ is small enough   and
$\rhoepsilon_4,$ $\rhoepsilon_2,$ and $\rhoepsilon_1$ are defined using
$\rhoepsilon_3$ as in (\ref{eq: ordering of epsilons wave equation}).
As the indicator
function is real-analytic, this shows that the indicator function
in non-vanishing in a generic set.   \hfill \Box \medskip

We need also to change the
  the singularity {\it detection condition} ({D}) with light-like  
directions $(\vec x,\vec \xi)$
  as follows:  We define that point $y\in \hat U$,
   satisfies the singularity {\it detection condition} (${D}^\prime$)  
with light-like directions $(\vec x,\vec \xi)$
   and $t_0,\hat s>0$
  if
  \medskip  
  
($D^\prime$) For any $s,s_0\in (0,\hat s)$  there are  
$(x_j^{\prime},\xi_j^{\prime})\in \W_j(s;x_j,\xi_j)$, $j=1,2,3,4$,
and ${f}_{(j)}\in {\mathcal  
I_{S}}^{n-3/2}(Y((x_j^{\prime},\xi_j^{\prime});t_0,s_0))$, and such  
that if  $u_{\vec \e}$ of is the solution  of (\ref{PABC eq})
with the source ${f}_{\vec \e}=\sum_{j=1}^4
\e_j{f}_{(j)}$, then
the function
$\p_{\vec \e}^4u_{\vec \e}|_{\vec \e=0}$ is not
$C^\infty$-smooth in any neighborhood of $y$.
  \medskip

When condition (D) is replace by ($D^\prime$), the considerations in  
the Sections
  \ref{sec: normal coordinates} and \ref{subsection combining}  show  
that we can recover
the conformal class of the metric. This proves Theorem \ref{main  
thm3}. \hfill \Box \medskip

\noindent
{\bf Proof.} (of Corollary \ref{coro thm3}).
Let us denote $W_j=I^+(p^-_j)\cap I^-(p^+_j)\subset M_j$.
By Theorem \ref{main thm3},  there is a map $\Psi:W_1\to W_2$
such that  the product type metrics $g_1=-dt^2+h_1(y)$
and $g_2=-dt^2+h_2(y)$ are conformal. At the vector field
$V=\p/\p x^0$ satisfies $\nabla^{g_j}V=0$,
for given $\mu_1(s)=x_0=(y_0,t_0),$ $s\in [s_-,s_+]$ we can consider  
all smooth paths $a:[0,1]\to W_1$
that satisfy $a(0)=x_0=(y_0,t_0)$ and $g_1(\dot a(s),V)=0$. The set of  
the end points
of such paths is equal to the set $W_1\cap (\{t_0\}\times N_1)$.
Considering all such paths on $W_1$ and the analogous paths
on $W_2$, we see that $\Psi:W_1\cap (\{t_0\}\times N_1)\to
W_2\cap (\{t_0\}\times N_2)$ is a diffeomorphism. Hence
$\Psi:W_1\to W_2$ has the form $\Psi(t,y)=(t,\tilde \psi(y,t))$.
This means that we can determine uniquely the foliation given by
the $t$-coordinate. As the metric tensors $g_1$ and $g_2$
are conformal and their  $(0,0)$-components in the
$(t,y)$ coordinates satisfy $(g_1)_{00}=-1$ and $(g_2)_{00}=-1$,
we conclude that $g_1$ and $g_2$ are isometric. Moreover, as
$g_1$ and $g_2$ are independent
of $t$, we see that there is a diffeomorphism $\Psi:W_1\to W_2$
of the form $\Psi(t,y)=(t,\psi(y))$ such that $g_1=\Psi^*g_2$.
Note also that if $\pi_2:(t,y)\mapsto y$,
then $h_1=\psi^*h_2$ on $\pi_2(W_1)$.
Thus, we can assume next that the metric tensors $g_1$ and $g_2$
are isometric and identify the sets $W_1$ and $W_2$
denoting $W=W_1=W_2$. 

As the linearized waves $u_{(j)}=Qf_{(j)}$ depend only on the
metric $g$, using the proof of Theorem \ref{main thm3}
we see that the indicator functions
$\mathcal  G(\bsequence,{\bf w})$ for $(U_1,g_1,a_1)$
and $(U_1,g_1,a_1)$ coincide for all $\bsequence$ and ${\bf w}$.
This implies that $a_1(t,y)^3=a_2(t,y)^3$ for all $(t,y)\in W$.
Hence $a_1(t,y)=a_2(t,y)$ for all $(t,y)\in W$.
  \hfill \Box \medskip

\subsection*{Acknowledgements}
The authors did part of this work at MSRI in Fall 2010, The Newton  
Institute in Fall 2011, the Fields Institute in the Spring 2012 and  
the Mittag-Leffler Institute in the Spring of 2013. They would like to  
express their gratitude to these institutions for their hospitality.

YK was partly supported by EPSRC.

ML was partly supported by the Finnish Centre of
Excellence in Inverse Problems Research 2006-2011 and 2012-2017.

GU was partly supported  by NSF. He was also partly supported by a  
Clay Senior Award at MSRI,  a Chancellor Professorship at UC Berkeley,  
a Rothschild Distinguished Visiting Fellowship at the Newton Institute  
and by the Fondation de Sciences Math\'ematiques de Paris.

\HOX{The style of the bibliography needs some more work. }

 \end{document}